\documentclass[10pt,leqno]{article}

\usepackage[francais]{babel}
\usepackage{amsmath,amssymb,amscd}
\usepackage{a4wide}
\usepackage[latin1]{inputenc}
\usepackage[all]{xy} 

\newenvironment{paragr}[1][]{\refstepcounter{subsection} \noindent \textbf{\thesubsection . \ #1}}{\medskip}

\newenvironment{theoreme}{ \medskip\refstepcounter{theo}  \noindent\textbf{Th\'eor\`eme \thetheo}. ---\em}{\em \medskip}
\newenvironment{proposition}{\medskip\refstepcounter{theo}   \noindent\textbf{Proposition \thetheo}. ---\em}{\em\medskip}
\newenvironment{corollaire}{\medskip\refstepcounter{theo}  \noindent\textbf{Corollaire \thetheo}. ---\em}{\em\medskip}

\newenvironment{lemme}{\medskip\refstepcounter{theo}   \noindent\textbf{Lemme \thetheo}. ---\em}{\em\medskip}

\newenvironment{definition}{\medskip\refstepcounter{theo}  \noindent\textbf{D\'efinition \thetheo}. ---}{\medskip}

\newenvironment{exemple}{\medskip\refstepcounter{theo}  \noindent\textbf{Exemple \thetheo}. ---}{\medskip}
\newenvironment{preuve}[1][]{\noindent \textbf{Démonstration.} #1 --- }{\hfill
  \ensuremath{\square} \medskip}

\newenvironment{remarque}{\medskip\refstepcounter{theo}  \noindent\textbf{Remarque \thetheo}. ---}{\medskip}
\newenvironment{remarques}{\medskip\refstepcounter{theo}  \noindent\textbf{Remarques \thetheo}. ---}{\medskip}

\DeclareMathOperator{\Lie}{Lie}

\DeclareMathOperator{\reg}{reg}
\DeclareMathOperator{\-reg}{-reg}

\DeclareMathOperator{\vol}{vol}
\DeclareMathOperator{\rang}{rang}
\DeclareMathOperator{\Ad}{Ad}

\DeclareMathOperator{\ad}{ad}
\DeclareMathOperator{\ab}{ab}

\DeclareMathOperator{\der}{der}
\DeclareMathOperator{\scnx}{sc}

\DeclareMathOperator{\Norm}{Norm}

\DeclareMathOperator{\el}{ell}

\DeclareMathOperator{\Aut}{Aut}
\DeclareMathOperator{\Gal}{Gal}
\DeclareMathOperator{\Hom}{Hom}
\DeclareMathOperator{\homc}{\mathcal{H}\mathit{om}}
\DeclareMathOperator{\inv}{inv}
\DeclareMathOperator{\Int}{Int}

\DeclareMathOperator{\aff}{aff}

\DeclareMathOperator{\Ker}{Ker}

\DeclareMathOperator{\Coker}{Coker}

\DeclareMathOperator{\val}{val}

\DeclareMathOperator{\Spec}{Spec}

\DeclareMathOperator{\Ext}{Ext}
\DeclareMathOperator{\trace}{trace}

\DeclareMathOperator{\car}{\mathfrak{car}}
\DeclareMathOperator{\lng}{long}

\DeclareMathOperator{\tors}{tors}
\DeclareMathOperator{\parab}{par}
\DeclareMathOperator{\sep}{sep}
\newcommand{\ZZ}{\mathbb{Z}}

\newcommand{\Gmk}{\mathbb{G}_{m,k}}

\newcommand{\NN}{\mathbb{N}}
\newcommand{\RR}{\mathbb{R}}
\newcommand{\AAA}{\mathbb{A}}
\newcommand{\CC}{\mathbb{C}}
\newcommand{\BBB}{\mathbb{B}}
\newcommand{\QQ}{\mathbb{Q}}
\newcommand{\FF}{\mathbb{F}}
\newcommand{\Fq}{\mathbb{F}_q}

\newcommand{\Ql}{\mathbb{Q}_{\ell}}
\newcommand{\Qlb}{\overline{\mathbb{Q}}_{\ell}}
\newcommand{\Qb}{\overline{\mathbb{Q}}}

\newcommand{\oc}{\mathcal{O}}

\newcommand{\Sc}{\mathcal{S}}
\newcommand{\tc}{\mathcal{T}}
\newcommand{\uc}{\mathcal{U}}

\newcommand{\wc}{\mathcal{W}}
\newcommand{\ec}{\mathcal{E}}

\newcommand{\mc}{\mathcal{M}}
\newcommand{\yc}{\mathcal{Y}}
\newcommand{\lc}{\mathcal{L}}

\newcommand{\fc}{\mathcal{F}}
\newcommand{\pc}{\mathcal{P}}
\newcommand{\hc}{\mathcal{H}}

\newcommand{\Ac}{\mathcal{A}}

\newcommand{\kc}{\mathcal{K}}

\newcommand{\Jc}{\mathcal{J}}

\newcommand{\xc}{\mathcal{X}}

\newcommand{\Gc}{\widehat{G}}
\newcommand{\Tc}{\widehat{T}}

\newcommand{\Mc}{\widehat{M}}
\newcommand{\Lc}{\widehat{L}}
\newcommand{\Rc}{\widehat{R}}

\newcommand{\Hc}{\widehat{H}}

\newcommand{\ggo}{\mathfrak{g}}

\newcommand{\mgo}{\mathfrak{m}}
\newcommand{\ngo}{\mathfrak{n}}
\newcommand{\ago}{\mathfrak{a}}

\newcommand{\pgo}{\mathfrak{p}}
\newcommand{\qgo}{\mathfrak{q}}

\newcommand{\tgo}{\mathfrak{t}}
\newcommand{\lgo}{\mathfrak{l}}

\newcommand{\zgo}{\mathfrak{z}}

\newcommand{\bgo}{\mathfrak{b}}
\newcommand{\Xgo}{\mathfrak{X}}

\newcommand{\al}{\alpha}

\newcommand{\la}{\lambda}

\newcommand{\back}{\backslash}

\newcommand{\bg}{\langle}
\newcommand{\bd}{\rangle}

\newcommand{\eps}{\varepsilon}

\renewcommand{\leq}{\leqslant}
\renewcommand{\geq}{\geqslant}

\title{Le lemme fondamental pondéré. II.\\ Énoncés cohomologiques}
\author{Pierre-Henri Chaudouard et Gérard Laumon}
\date{}

\begin{document}

\makeatletter
\@addtoreset{equation}{subsection}         
\def\theequation{\thesubsection.\arabic{equation}}
\makeatother

\maketitle

\begin{abstract}
Dans cet article, on étudie la cohomologie de la fibration de Hitchin tronquée introduite dans un article précédent. On étend les principaux théorèmes de Ngô sur la cohomologie de la partie elliptique de la fibration de Hitchin. Comme conséquence, on obtient une démonstration du lemme fondamental pondéré d'Arthur.
\end{abstract}

\renewcommand{\abstractname}{Abstract}
\begin{abstract}
In this paper, we study the cohomology of the  truncated Hitchin fibration, which was introduced in a previous paper. We extend Ngô's main theorems on the cohomology of the elliptic part of the Hitchin fibration. As a consequence, we get a proof of Arthur's weighted fundamental lemma.
\end{abstract}
\tableofcontents

\section{Introduction}

\begin{paragr} L'objet de ce travail, qui fait suite à \cite{LFPI}, est de démontrer \emph{le lemme fondamental pondéré}. Cet énoncé a été conjecturé par Arthur dans ses travaux sur la stabilisation des traces d'Arthur-Selberg  (cf. conjecture 5.1 de \cite{StableI}). Il généralise le lemme fondamental de Langlands-Shelstad et se compose d'une famille d'identités combinatoires entre intégrales orbitales pondérées $p$-adiques. Après les travaux de Ngô sur le lemme fondamental, il est le seul ingrédient manquant dans la stabilisation de la formule des traces. Selon une stratégie due à Langlands, on devrait obtenir, au moyen de la formule des traces stable, toutes les fonctorialités de type «endoscopique». Par exemple, des travaux en cours d'Arthur devraient donner une classification du spectre automorphe des groupes classiques en fonction de celui des groupes linéaires généraux. Hormis des travaux de Kottwitz (cf. \cite{KoBC}) et de Whitehouse (cf. \cite{Whitehouse}), le lemme fondamental pondéré était jusqu'à présent un problème complètement ouvert.
\end{paragr}

\begin{paragr} Plus précisément,  nous démontrons  une variante en caractéristique positive et pour les algèbres de Lie de l'énoncé d'Arthur. De plus, pour alléger autant que faire se peut l'exposition, nous nous sommes  limités dans cet article au cas des groupes déployés. Le cas général, qui comprend aussi les formes «non-standard» dues à Waldspurger du lemme fondamental pondéré,  s'obtient par des méthodes similaires et sera traité ultérieurement. Des travaux de Waldspurger (cf. \cite{Wlfp1} et \cite{Wlfp2}) montre que toutes ces variantes du lemme fondamental pondéré  impliquent l'énoncé original d'Arthur pour les groupes sur les corps $p$-adiques.
\end{paragr}

\begin{paragr} Désormais, nous réservons le terme de lemme fondamental ou lemme fondamental pondéré aux variantes pour les algèbres de Lie en  caractéristique positive des énoncés originaux. Dans \cite{ngo2}, Ngô a démontré le lemme fondamental ordinaire. Dans son approche, le lemme fondamental  résulte d'un énoncé  cohomologique sur la partie elliptique de la fibration de Hitchin. De même ici, le  lemme fondamental pondéré est la conséquence d'un théorème cohomologique (cf. théorème \ref{thm:cohomologique}) sur la fibration de Hitchin tronquée introduite dans \cite{LFPI}. Le résultat clef dans l'approche de Ngô est son théorème  qui décrit les supports des différents constituants de la cohomologie de la partie elliptique de la fibration de Hitchin. De manière imagée, ce théorème  montre que les intégrales orbitales globales sont les «limites» des intégrales orbitales associées aux éléments les plus réguliers possibles. En particulier, pour démontrer la variante globale du lemme fondamental, il suffit de le faire pour ces éléments très réguliers  pour lesquels une vérification  «à la main» est possible. L'énoncé local résulte ensuite de l'énoncé global. Notre principal résultat cohomologique est que, pour la cohomologie de la fibration de Hitchin tronquée, il n'y a aucun nouveau support qui apparaît. La conséquence la plus frappante est que   les intégrales orbitales pondérées globales sont aussi des «limites» d'intégrales orbitales ordinaires. 
\end{paragr}

\begin{paragr}[Le lemme fondamental pondéré.] --- Nous allons maintenant énoncer le théorème local que nous démontrons. Soit $\Fq$ un corps fini à $q$ éléments et $G$ un groupe réductif connexe, défini et déployé sur $\Fq$. Soit $T$ un sous-tore maximal de $G$ défini et déployé sur $\Fq$. On suppose que l'ordre du groupe de Weyl $W^G$ de $(G,T)$ est inversible dans $\Fq$. Soit $M$ un sous-groupe de Levi de $G$ défini sur $\Fq$ qui contient le tore $T$. On note par des lettres gothiques les algèbres de Lie correspondantes. Considérons le corps local $K=\Fq((\eps))$ et $\oc=\Fq[[\eps]]$ son anneau d'entiers. Soit $X\in \mgo(K)$ un élément semi-simple et $G$-régulier. Son centralisateur dans $G\times_{\Fq}K$ est un sous-$K$-tore maximal noté $J_X$ qui est inclus dans $M\times_{\Fq}K$. Dans ce contexte, une intégrale orbitale pondérée d'Arthur s'écrit
$$J_M^G(X)=|D^G(X)|^{1/2} \int_{J_X(K)\back G(K)} \mathbf{1}_{\ggo(\oc)}(\Ad(g^{-1})X) v_M^G(g) dg,$$
où
\begin{itemize}
\item $|D^G(X)|$ est la valeur absolue usuelle du discrimant de Weyl ;
\item $\mathbf{1}_{\ggo(\oc)}$ est la fonction caractéristique de $\ggo(\oc)$ ;
\item $\Ad$ désigne l'action adjointe de $G$ sur $\ggo$ ;
\item $v_M^G(g)$ est la fonction «poids» d'Arthur qui est $M(F)$-invariante à gauche ;
\item $dg$ est la mesure quotient déduite de mesures de Haar sur $G(K)$ et $J_X(K)$.
\end{itemize}
  Le poids et la mesure $dg$ dépendent de choix qu'il est possible de normaliser. Avec de bonnes conventions de mesures, l'intégrale $J_M^G(X)$ ne dépend de $X$ qu'à $M(K)$-conjugaison près.

Soit $\Gc$ et $\Tc$ les groupes complexes duaux au sens de Langlands de $G$ et $T$. Le groupe $\Tc$ s'identifie à un sous-tore maximal de $\Gc$. Le dual $\Mc$ de $M$ s'identifie à un sous-groupe de Levi de $\Gc$ qui contient $\Tc$. Soit $s\in \Tc$ et $\Mc'$ la composante neutre du centralisateur de $s$ dans $\Mc$. C'est un sous-groupe réductif connexe de $\Mc$ qui contient $\Tc$. On suppose  que les composantes neutres des centres de $\Mc$ et $\Mc'$ sont égales. Soit $M'$ le groupe sur $\Fq$ dual de $\Mc'$. Le tore $T$ s'identifie à un sous-tore maximal de $M'$ et le groupe de Weyl $W^{M'}$ s'identifie à un sous-groupe du groupe de Weyl  $W^M$. On a donc un morphisme canonique
\begin{equation}
  \label{eq:intro-morphcano}
  \mgo'//M' \simeq \tgo //W^{M'} \to \tgo //W^M \simeq \mgo//M,
\end{equation}
où les isomorphismes sont ceux de Chevalley et les quotients sont les quotients catégoriques. L'ouvert formé des éléments semi-simples $G$-réguliers définit par l'application canonique $\mgo\to \mgo//M$ un ouvert dans le quotient $\mgo//M$ dont l'image réciproque par le composé $\mgo'\to \mgo'//M'\to\mgo//M$ est l'ouvert des éléments semi-simples réguliers de $G$. À la suite de Waldspurger (cf. \cite{lem_fond}), on peut adapter la définition des facteurs de transfert de Langlands-Shelstad (cf. \cite{LS}) aux algèbres de Lie. En considérant $M'$ comme un «groupe endoscopique» de $M$, on obtient une fonction
$$\Delta_{M',M} : \mgo'(K) \times \mgo(K) \to \CC$$
définie sur l'ouvert formé des couples $(Y,X)$ dont chaque composante est un élément semi-simple et $G$-régulier. Le facteur de transfert est en fait invariant par l'action adjointe de $M'(K)\times M(K)$. Pour un couple $(Y,X)$ dont chaque composante est un élément semi-simple et $G$-régulier, le facteur 
 $\Delta_{M',M}(Y,X)$ est non nul si et seulement si $Y$ et $X$ ont même image dans $(\mgo//M)(K)$ par les applications canoniques $\mgo'(K)\to (\mgo//M)(K)$ déduite de (\ref{eq:intro-morphcano}) et $\mgo(K)\to (\mgo//M)(K)$. Dans ce dernier cas, lorsque de plus on a $M'=M$, le facteur de transfert vaut $1$.

Pour tout élement $Y\in \mgo'(K)$ semi-simple régulier, on peut former la combinaison linéaire
$$J^G_{M',M}(Y)=\sum_{X} \Delta_{M',M}(Y,X)J_M^G(X),$$
où la somme est prise sur l'ensemble des classes de $M(K)$-conjugaison des éléments $X \in \mgo(K)$ semi-simples et $G$-réguliers. La somme est à support fini.

On suppose désormais que $s$ \emph{n'est pas central} dans $\Mc$ c'est-à-dire $M'\not=M$. Soit $\ec(s)$ l'ensemble fini des sous-groupes réductifs connexes $\Hc$ de $\Gc$ qui vérifient
\begin{itemize}
\item il existe un entier $i\geq 2$ et des éléments $s_1\in s Z_{\Mc}$ (où $Z_{\Mc}$ est le centre de $\Mc$) et $s_2, \dots, s_n$ dans le centre de $\Mc'$ tels que $\Hc$ est la composante neutre de l'intersection des centralisateurs des $s_i$ dans $\Gc$ ;
\item les composantes neutres des centres de $\Gc$ et $\Hc$ sont égales.
\end{itemize}
Tout $\Hc\in\ec(s)$ contient $\Mc'$ comme sous-groupe de Levi. On remarque que l'ensemble $\ec(s)$ ne contient pas $\Gc$. Pour tout $s'\in \Tc$ et tout sous-groupe $\Hc$ de $\Gc$, on note $\Hc_{s'}$ la composante neutre du centralisateur de $s'$ dans $\Hc$. On dit que $s'$ est elliptique dans $\Hc$ si les composantes neutres des cenbtres de $\Hc_{s'}$ et $\Hc$ sont égales. Pour tous $\Hc\in \ec(s)$ et $s'\in Z_{\Mc'}$ elliptique dans $\Hc$ le groupe  $\Hc_s$ appartient à $\ec(s)$.  Par dualité de Langlands, on identifie les éléments de $\ec(s)$ à une famille de groupes réductifs connexes définis sur $\Fq$ qui contiennent $M'$ comme sous-groupe de Levi. Pour tout $Y\in \mgo'(K)$ semi-simple régulier, on définit pour tout $H\in \ec(s)$ la variante «stable» $S_{M'}^H(Y)$ des intégrales orbitales pondérées par la formule de récurrence suivante
$$S_{M'}^H(Y)=J_{M',M'}^H(Y)-\sum_{s'\in Z_{\Mc'}/Z_{\Hc}, \ s'\not=1 } |Z_{\Hc_{s'}}/Z_{\Hc}|^{-1} S_{M'}^{H_{s'}}(Y),$$ 
où par convention  $|Z_{\Hc_{s'}}/Z_{\Hc}|^{-1}$ vaut $0$ lorsque $s'$ n'est pas elliptique dans $\Hc$. On peut montrer que la somme est à support fini.

Pour tout $s'\in sZ_{\Mc}$ elliptique dans $\Gc$, le groupe $\Gc_{s'}$ appartient à $\ec(s)$. En particulier, l'expression $S_{M'}^{G_{s'}}(Y)$ est bien définie. On peut alors énoncer le lemme fondamental pondéré que nous démontrons.

\begin{theoreme}
Pour tout $Y\in \mgo'(K)$ semi-simple et $G$-régulier, on a l'identité suivante 
$$  J^G_{M',M}(Y)=\sum_{s'\in Z_{\Mc}/Z_{\Gc}} |Z_{\Gc_{s'}}/Z_{\Gc}|^{-1} S_{M'}^{G_{s'}}(Y),$$ 
où comme précédemment la somme est à support (fini) sur les éléments $s'$ elliptiques dans $\Gc$.
\end{theoreme}

\end{paragr}

\begin{paragr}[Le théorème de support.] --- Comme l'on dit plus haut, le théorème précédent repose sur une étude cohomologique de la fibration de Hitchin convenablement tronquée. Le pendant cohomologique du théorème précédent est le théorème \ref{thm:cohomologique} ci-dessous dont nous ne reproduirons pas dans cette introduction l'énoncé. Nous allons plutôt énoncer notre théorème de support qui est la clef du théorème \ref{thm:cohomologique} en ce sens qu'il ramène l'étude cohomologique de la fibration de Hitchin tronquée à celle de sa partie elliptique. Sur cette dernière le  théorème \ref{thm:cohomologique} est connu par un théorème fondamental de Ngô (cf. \cite{ngo2}).

La situation est la suivante. Soit $k$ une clôture algébrique d'un corps fini. Soit $G$ un groupe semi-simple défini sur $k$ et $T\subset G$ un sous-tore maximal. On suppose que la caractéristique de $k$ ne divise pas l'ordre du  groupe de Weyl $W$ de $(G,T)$. Soit $C$ une courbe projective, lisse et connexe muni d'un diviseur $D$ effectif et pair de degré supérieur au genre de $C$. Soit $\infty$ un point fermé de $C$ hors du support de $D$. Le $k$-champ de Hitchin $\mc$ classifie les triplets de Hitchin $(\ec,\theta,t)$ formés d'un $G$-torseur $\ec$, d'une section $\theta$ du fibré vectoriel $\Ad(\ec)(D)$ qu'on obtient en poussant le torseur $\ec$ par la représentation adjointe tordue par $D$ et d'un élément $G$-régulier $t\in \tgo$ tel que $t$ et $\theta(\infty)$ ont même image par le morphisme canonique
  \begin{equation}
    \label{eq:intro-cano2}
    \ggo \to \ggo//G\simeq \tgo//W.
  \end{equation}
 Soit $\Ac$ le $k$-schéma qui classifie les couples $(h_a,t)$ formés d'une section $h_a$ du fibré $\tgo//W \otimes_k \oc(D)$ et d'un  élément $G$-régulier $t\in \tgo$ tel que l'image de $t$ par $\tgo \to \tgo//W$ soit égale à $h_a(\infty)$.  La fibration de Hitchin $f : \mc \to \Ac$ est le morphisme qui à un triplet $(\ec,\theta,t)$ associe le couple $(\chi(\theta),t)$ où $\chi_D$ est une version tordue  du morphisme (\ref{eq:intro-cano2}).  On a introduit dans \cite{LFPI} un ouvert tronqué $\mc^\xi$ de $\mc$ formé de triplets «$\xi$-semi-stables» (pour des rappels, cf. section \ref{sec:Hitchin-tronque} ci-dessous) qui dépend d'un paramètre $\xi$ qui vit dans un espace réel. Lorsque $\xi$ est «en position générale» (au sens de la  remarque \ref{rq:position}), le morphisme $f^\xi$ restreint à l'ouvert $\mc^\xi$ est propre (cf. théorème 6.2.2 de \cite{LFPI} rappelé au \ref{thm:rappel-lfpI}). On définit à la section \ref{sec:Abon} un ouvert $\Ac^{\mathrm{bon}}$ de $\Ac$. Soit $\Ac^{\el}$ l'ouvert elliptique de $\Ac$ et 
$$j:\mathcal{A}^{\mathrm{ell}}\cap \mathcal{A}^{\mathrm{bon}} \hookrightarrow \mathcal{A}^{\mathrm{bon}}$$
l'immersion ouverte canonique. Soit $\ell$ un nombre premier inversible dans $k$ et $\Qlb$ une clôture algébrique de $\Ql$. Pour tout entier $i$, soit ${}^{\mathrm{p}}\mathcal{H}^{i}(Rf_{\ast}^\xi \overline{\mathbb{Q}}_{\ell})$ le $i$-ème faisceau de cohomologie perverse de $Rf_{\ast}^\xi \overline{\mathbb{Q}}_{\ell}$. Voici notre principal théorème cohomologique (cf. théorème \ref{thm:support} ci-dessous)

\begin{theoreme}\label{thm:intro-support}
   Pour tout entier $i$, on a un isomorphisme canonique
$${}^{\mathrm{p}}\mathcal{H}^{i}(Rf_{\ast}^\xi \overline{\mathbb{Q}}_{\ell})\cong j_{!\ast}j^{\ast}{}^{\mathrm{p}}\mathcal{H}^{i}(Rf_{\ast}^\xi \overline{\mathbb{Q}}_{\ell})$$
sur $\mathcal{A}^{\mathrm{bon}}$.
\end{theoreme}

On sait que le champ $\mc^\xi$ est lisse sur $k$. Par le théorème de Deligne, la propreté du morphisme $f^\xi$ implique que  le complexe $(Rf_{\ast}^\xi \overline{\mathbb{Q}}_{\ell}$ est pur. Il s'ensuit que les faisceaux pervers ${}^{\mathrm{p}}\mathcal{H}^{i}(Rf_{\ast}^\xi \overline{\mathbb{Q}}_{\ell})$ sont semi-simples. Le théorème précédent affirme que sur $\Ac^{\mathrm{bon}}$ les supports des constituants irréductibles rencontrent tous l'ouvert elliptique.

\end{paragr}

\begin{paragr}[Description de l'article.] --- Venons-en à une brève description du contenu de l'article. Les premières sections (sections  \ref{sec:notations} et \ref{sec:centralisateurs}) sont consacrées à quelques notations et compléments. On introduit la fibration de Hitchin tronquée à la section  \ref{sec:Hitchin-tronque}. À la section \ref{sec:Picard}, on rappelle l'action introduite par Ngô d'un champ de Picard sur les fibres de Hitchin. Comme cette action ne respecte pas les fibres tronquées, on introduit à la section \ref{sec:J1} un sous-champ qui, lui, respecte la troncature. La section \ref{sec:cameral} donne quelques rappels sur les courbes camérales et le faisceau des composantes connexes du champ de Picard introduit par Ngô. On en déduit à la section \ref{sec:sdecomposition} une décomposition en $s$-parties de la cohomologie de la fibration de Hitchin tronquée qui généralise celle de Ngô sur la partie elliptique (cf. théorème \ref{thm:sdecomposition}). À la section \ref{sec:Abon}, on introduit un certain ouvert «bon» de la base de la fibration de Hitchin et on donne une condition suffisante pour qu'un point de la base soit dans cet ouvert «bon» (théorème \ref{thm:Abon}). La section \ref{sec:cohomologie} contient les principaux résultats cohomologiques : théorème \ref{thm:support} de «support» et théorème \ref{thm:cohomologique} qui est l'analogue cohomologique du lemme fondamental pondéré. À la section \ref{sec:premiercalcul}, on exprime les $s$-intégrales orbitales pondérées globales  introduites à la section  \ref{sec:sIOP} en termes de trace de Frobenius sur la $s$-partie de la cohomologie (cf. théorème \ref{thm:scomptage}). À la section \ref{sec:LFPglobal}, on déduit des considérations précédentes une variante globale du lemme fondamental pondéré (cf. théorème \ref{thm:lfpglobal}). Le  reste de l'article est consacré au passage de cette forme  globale à l'énoncé local.

\end{paragr}

\begin{paragr}[Remerciements.] --- Nous remercions Luc Illusie, Ngô Bao Châu, Michel Raynaud et Jean-Loup Waldspurger pour des discussions utiles. Une partie de cet article a été écrit lors d'un séjour du premier auteur nommé à l'\emph{Institute for Advanced Study} de Princeton à l'automne 2008. Il souhaite remercier cet institut pour son hospitalité ainsi que la \emph{National Science Foundation}  pour le soutien (agreement No. DMS-0635607) qui a rendu ce séjour possible.
  
\end{paragr}

\section{Notations} \label{sec:notations}
\begin{paragr}\label{S:G}
Soit $k$ une clôture algébrique d'un corps fini $\FF_q$ à $q$ éléments. Soit $G$ un groupe algébrique réductif et connexe sur $k$.  Soit $G_{\der}$ le groupe dérivé de $G$. Soit $T$ un sous-tore maximal de $G$.  Soit 
$$W=W^G=W^G_T,$$
$$\Phi=\Phi^G=\Phi^G_T \ \ \text{    et    } \ \ \Phi^\vee=\Phi^{G,\vee}_T$$ 
respectivement le groupe de Weyl, l'ensemble de racines, resp. des coracines, de $T$ dans $G$. L'application qui à une racine associe sa coracine induit une bijection entre $\Phi$ et $\Phi^\vee$. Pour tout ensemble $X$, on note $|X|$ son cardinal (fini ou infini). On suppose dans toute la suite que l'ordre du groupe de Weyl $|W^G_T|$ est inversible dans $k$.

L'algèbre de Lie d'un groupe désigné par une lettre majuscule est notée par la lettre gothique correspondante. Ainsi  $\ggo$ et $\tgo$ sont les  algèbres de Lie respectives de $G$ et $T$. Soit
$$X^*(G)=\Hom_k(G,\Gmk)$$
le groupe des caractères de $G$. Soit
$$X_*(G)=\Hom_{\ZZ}(X^*(G),\ZZ)$$
le groupe dual. Dans le cas d'un tore, ce groupe s'identifie canoniquement au groupe des cocaractères.
On utilisera aussi par la suite les espaces vectoriels réels 
$$\ago_G=X_*(G)\otimes_{\ZZ} \RR$$
et
$$\ago_G^*=X^*(G)\otimes_{\ZZ} \RR.$$
\end{paragr}

\begin{paragr} Soit $\Int$ l'action de $G$ sur lui-même par automorphisme intérieur et $\Ad$ la représentation adjointe. Soit 
$$\car=\car_G=\ggo//G$$
le quotient catégorique de $\ggo$ par l'action adjointe de $G$. Soit 
$$ \chi=\chi_G \ :\ \ggo \to \car$$
le morphisme canonique, qu'on appelle \emph{morphisme caractéristique}. Soit
$$\chi^T=\chi_G^T \ :\ \tgo \to \car_G$$ 
la restriction de $\chi_G$ à $\tgo$. C'est un revêtement fini et galoisien de groupe $W$ de $\car_G$. Le morphisme  $\chi_G^T$ est $W$-invariant et induit un isomorphisme dû à Chevalley
$$\tgo//W\simeq \car.$$
On en déduit que le schéma $\car$ est un espace affine de dimension $n$, où $n$ est le rang de $G$, de coordonnées $n$ polynômes homogènes $G$-invariants sur $\mathfrak{g}$ dont les degrés, notés $e_{1}\leq \ldots \leq e_{n}$, sont uniquement déterminés.
\end{paragr}

\begin{paragr}[Section de Kostant.] --- \label{S:Kostant} On dit qu'un élément de $\ggo$ est $G$-régulier si la dimension de son centralisateur dans $G$ est égal au rang de $G$. Dans la suite, lorsque le contexte est clair, on parle simplement d'éléments réguliers. Soit $\mathfrak{g}^{\mathrm{reg}}\subset\mathfrak{g}$ l'ouvert des éléments réguliers de $\mathfrak{g}$.

Soit $B$ un sous-groupe de Borel de $G$ qui contient $T$ et $\Delta\subset \Phi$ l'ensemble des racines simples de $T$ dans $B$. Complétons ces données en un épinglage $(B,T,\{X_\al\}_{\al\in \Delta})$ de $G$. Soit $X_{\pm}=\sum_{\al \in \pm \Delta} X_\al$ où l'on a défini $X_{-\al}$, pour $\al\in \Delta$, comme l'unique vecteur radiciel associé à la racine $-\al$ qui vérifie $[X_{\al},X_{-\al}]=\al^\vee$. Soit $\ggo_{X_+}$ le centralisateur de $X_+$ dans $\ggo$. Kostant a montré que l'espace affine $X_- +\ggo_{X_+}$ est inclus dans l'ouvert régulier $\ggo^{\reg}$ et  que la restriction du morphisme  $\chi : \ggo \to \car$ induit un isomorphisme 
$$X_- +\ggo_{X_+} \to \car.$$
Soit
$$\eps=\eps_G  \ : \ \car \to  X_- +\ggo_{X_+}$$
le morphisme inverse.

\end{paragr}

\begin{paragr}[Sous-groupes paraboliques et sous-groupes  de Levi.] --- \label{S:parab}Pour tous sous-groupes $M$ et $Q$ de $G$, on note $\fc^Q(M)$ l'ensemble des  sous-groupes paraboliques $P$ de $G$ qui vérifient $M\subset P\subset Q$. Les éléments de $\fc=\fc^G(T)$ sont appelés les sous-groupes paraboliques \emph{semi-standard} de $G$.

On appelle \emph{sous-groupe de Levi} de $G$ un facteur de Levi d'un sous-groupe parabolique de $G$. Soit $\lc=\lc^G(T)$ l'ensemble des sous-groupes de Levi \emph{semi-standard} de $G$ (au sens où ils contiennent $T$).

  Pour tout $P\in\fc$, soit $N_P$ le radical unipotent de $P$ et $M_P\in \lc$ l'unique sous-groupe de Levi de $P$ qui contient $T$.  Pour tout sous-groupe de Levi $M\in \lc$, soit $\pc(M)$, resp. $\lc(M)$ le sous-ensemble des $P\in  \fc$ tels que $M_P=M$,  resp. des $L\in \lc$  tels que $M\subset L$. Si l'on veut rappeler le groupe ambiant $G$, on ajoute un exposant $G$ à ces notations.

Soit $M\in \lc$. Soit $\ago_T^M$ le sous-espace de $\ago_T$ engendré par les coracines de $T$ dans $M$. En considérant les racines, on définit de même le sous-espace $(\ago_T^M)^*\in \ago_T^*$. Le morphisme canonique $\ago_T \to \ago_M$ a pour noyau $\ago_T^M$. Ce dernier a pour supplémentaire dans $\ago_T$ l'orthogonal de $(a_T^M)^*$ ce qui identifie à $\ago_M$ à un sous-espace de $\ago_T$. De la sorte, on a une décomposition
$$\ago_T=\ago_T^M \oplus \ago_M.$$   
On  pose $\ago_M^G=\ago_M \cap \ago_T^G$. On a alors $\ago_M=\ago_M^G \oplus \ago_G$.

\medskip

On munit l'espace vectoriel réel $\ago_T$ d'un produit scalaire invariant sous l'action naturelle du groupe de Weyl $W^G$.  Les décompositions en sommes directes comme ci-dessus sont alors orthogonales. Tous les sous-espaces de $\ago_T$ sont munis de la mesure euclidienne (c'est-à-dire la mesure de Lebesgue qui donne le covolume $1$ aux réseaux de base une base orthonormale).
\end{paragr}

\begin{paragr}[Dual de Langlands.] --- \label{S:dual-Langlands} Soit $\Tc$ le tore sur $\CC$ défini par 
$$\Tc(\CC)=\Hom_{\ZZ}(X_*(T),\CC^\times).$$
Soit $\Gc$ un groupe réductif connexe sur $\CC$ muni d'un plongement de $\Tc$ comme sous-tore maximal de sorte qu'on ait l'identification suivante entre données radicielles
$$(X^*(T),\Phi^{G}_T,X_*(T),\Phi^{G,\vee}_T)=(X_*(\Tc),\Phi^{\Gc,\vee}_{\Tc},X\*(\Tc),\Phi^{\Gc}_{\Tc}).$$
On dit que $\Gc$ est «le» dual de Langlands de $G$. Réciproquement, tout sous-groupe réductif connexe $\Mc$ de $\Gc$ qui contient $\Tc$ détermine, à isomorphisme près, un groupe réductif connexe $M$ sur $k$ muni d'un plongement de $T$ dans $M$ pour lequel on ait une identification comme ci-dessus entre les données radicielles de $(M,T)$ et celles de $(\Mc,\Tc)$.  Les groupes $\Mc$ et $M$ sont dits duaux. On remarquera que le groupe de Weyl $W^M$ est alors naturellement un sous-groupe de $W^G$. Lorsque, de plus, $\Mc$ est un sous-groupe de Levi de $\Gc$, le groupe $M$ muni du plongement de $T$ s'identifie au sous-groupe de Levi $M\in \lc^G(T)$  défini par la condition 
$$\Phi^M_T=\Phi^{\Mc,\vee}_{\Tc}.$$

Dans la suite, par abus, on ne distinguera pas dans les notations  un groupe algébrique sur $\CC$ de son groupe de points à valeurs dans $\CC$.

\end{paragr}

\begin{paragr}
Soit $M\in \lc$ et $P\in \pc(M)$. Soit $\pgo$ l'algèbre de Lie de $P$. Soit
   $$\chi_P \ :\ \pgo \to \car_P=\pgo//P$$
le morphisme canonique de $\pgo$ sur son quotient adjoint. La restriction de ce morphisme à $\mgo$ est $M$-invariant et induit un isomorphisme de $\car_M=\mgo //M$ sur $\car_P$ (cf. \cite{LFPI} lemme 2.7.1). On a, de plus, un diagramme commutatif dont la second flèche horizontale est la projection donnée par la décomposition canonique $\pgo=\mgo\oplus \ngo_P$.
   
 \begin{equation}
    \label{eq:diagcar_P}
    \xymatrix{ \mgo \ar[r]  \ar[d]^{\chi_M} & \pgo \ar[r]  \ar[d]^{\chi_P} & \mgo \ar[d]^{\chi_M}  \\
 \car_M \ar[r]^{\sim} &\car_P\ar[r]^{\sim} & \car_M }
\end{equation}

Soit $\pgo^{\reg}=\pgo\cap \ggo^{\reg}$ l'ouvert de $\pgo$ formé des éléments $G$-réguliers. 

\begin{lemme}\label{lem:Pcj} Soit $X_1$ et $X_2$ dans $\pgo$ tel que $\chi_P(X_1)=\chi_P(X_2)$. Si $X_1$ et $X_2$ sont tous deux semi-simples ou tous deux $G$-réguliers, il existe $p\in P$ qui conjugue $X_1$ et $X_2$.   
\end{lemme}

\begin{preuve}
Tout élément semi-simple de $\pgo$ est évidemment conjugué sous $P$ à un élément de $\tgo$. Le résultat est donc évident si $X_1$ et $X_2$ sont tous deux semi-simples. 

Soit $X\in \pgo$ et $X=X_s+X_n$ sa décomposition de Jordan. Montrons qu'on a $\chi_P(X)=\chi_P(X_s)$. Quitte à conjuguer $X$ par un élément de $P$, on peut déjà supposer qu'on a  $X_s\in \mgo$. Soit $X'_n$ la projection de $X_n$ sur $\mgo$ (selon la décomposition $\pgo=\mgo\oplus\ngo_P$). Alors la projection de $X$ sur $\mgo$ admet $X_s+X_n'$ comme décomposition de Jordan. On sait bien qu'on a $\chi_M(X_s+X'_n)=\chi_M(X_s)$. Il s'ensuit qu'on a 
$$\chi_P(X_s)=\chi_M(X_s)=\chi_M(X_s+X'_n)=\chi_P(X)$$
par le diagramme commutatif (\ref{eq:diagcar_P}). 

Soit $X_1$ et $X_2$ deux éléments  $G$-réguliers de $\pgo$ tels que $\chi_P(X_1)=\chi_P(X_2)$. Il s'agit de montrer qu'ils sont conjugués sous $P$. D'après ce qui précède, les parties semi-simples dans les décompositions de Jordan de $X_1$ et $X_2$ sont conjuguées sous $P$. Quitte à conjuguer $X_2$ par un élément de $P$, on peut supposer qu'elles sont égales à un élément $Y\in \pgo$. 
 
Quitte maintenant à considérer $G_Y$ le centralisateur de $Y$ dans $G$, qu'on sait être un groupe réductif connexe, et son sous-groupe parabolique $G_Y\cap P$, on est ramené au cas où $X_1$ et $X_2$ sont nilpotents et $G$-réguliers. On sait bien que de tels éléments sont conjugués par un élément  $g\in G$. Montrons qu'on a nécessairement $g\in P$ ce qui conclura la démonstration. Soit $B_1$ et $B_2$ deux sous-groupes de Borel tels que, pour $i=1,2$, on ait $X_i\in \bgo_i$ et $B_i\subset P$. Ces groupes sont uniquement déterminés. En particulier, puisque $\Ad(g)X_1=X_2$, on doit avoir $gB_1g^{-1}=B_2$. Donc les sous-groupes paraboliques $gPg^{-1}$ et $P$ qui sont conjugués et contiennent le même sous-groupe de Borel à savoir $B_2$, sont égaux. Il s'ensuit qu'on a $g\in P$ ce qu'il fallait démontrer.   
 
\end{preuve}

\end{paragr}

\begin{paragr}[Factorisation du discriminant.] --- \label{S:discr} Soit $\Mc$ un sous-groupe réductif connexe de $\Gc$ qui contient $\Tc$. Soit $M$ un groupe réductif connexe sur $k$ muni d'un plongement de $T$ et dual de $\Mc$. On a des inclusions $\Phi^M\subset \Phi^G$ et $W^M\subset W^G$. On a donc un morphisme  
  $$\chi_G^M \ : \ \car_M \to \car_G$$
qui vérifie $\chi_G^M\circ \chi_M^T=\chi_G^T$. Lorsque, de plus, on a $M\in \lc^G(T)$, le morphisme  $\chi_G^M\circ \chi_M$ est la restriction de $\chi_G$ à $\mgo$.

 Soit $B\in \pc(T)$ un sous-groupe de Borel et $\Phi^B$ les racines de $T$ dans $B$. Pour toute racine $\al\in \Phi$, soit $d\al\in k[\tgo]$ sa dérivée. Les éléments de $k[\tgo]$ définis par
$$D^M=\prod_{\al\in \Phi^M} d\al$$
et 
$$R_M^B=\prod_{\al\in \Phi^B-\Phi^M} d\al$$
sont $W^M$-invariants : ils se descendent donc en des fonctions régulières sur $\car_M$.  C'est clair pour le premier qui n'est autre que le discriminant de Weyl. Pour le second, cela résulte de la formule
$$w(R_T^B)=(-1)^{\ell(w)} R_T^B$$
pour tout $w\in W$ où $\ell(w)$ est la longueur de $w$ (cf. \cite{Bki} chap. V, §5.4 proposition 5) et de la formule 
$$R_T^B = R_M^B R_{T}^{M\cap B}.$$
Notons que $R_M^B$ ne dépend du choix de $B$ qu'à un signe près. Soit $\car^{M\textrm{-}\reg}_M$ l'ouvert où $D^M$ ne s'annule pas. C'est précisément le lieu où le morphisme $\chi^T_M$ est étale. On ne le confondra pas avec l'ouvert $\car^{G\textrm{-}\reg}_M$ où le discriminant pour $G$
$$D^G=D^M (R_M^B)^2$$
ne s'annule pas. Seul ce dernier interviendra dans la suite et on pose
$$\car_M^{\reg}=\car_M^{G\textrm{-}\reg}.$$
On notera que la fibre de $\chi_M$ au-dessus de $\car_M^{\reg}$ est l'ouvert de $\mgo$ formé des éléments semi-simples $G$-réguliers. L'ouvert $\tgo^{\reg}$ est l'ouvert formé des éléments  $G$-réguliers.

\end{paragr}

\section{Compléments sur les centralisateurs}\label{sec:centralisateurs}

\begin{paragr}[Centralisateurs.] --- \label{S:centralisateurs} Pour tout sous-groupe parabolique $P\in \pc^G(T)$, soit $I^P$ le schéma en groupes des centralisateurs sur $\pgo$ défini par
$$I^P=\{(p,X)\in P \times_k \pgo \ | \ \Ad(p)X=X\}.$$
On a le lemme suivant.

\begin{lemme}\label{lem:IP}
 Le morphisme canonique $I^P \to I^G$ induit un isomorphisme 
\begin{equation}
  \label{eq:morphpreg}
  I^P_{|\pgo^{\reg}} \buildrel\sim\over\longrightarrow  I^G_{|\pgo^{\reg}}.
\end{equation}
sur l'ouvert  
$$\pgo^{\reg}=\pgo\cap \ggo^{\reg}.$$ 
En outre, les schémas en groupes  $I^P_{|\pgo^{\reg}}$ et $I^G_{|\pgo^{\reg}}$ sont commutatifs et lisses sur $\pgo^{\reg}$.
\end{lemme}

\begin{preuve}
 Le morphisme (\ref{eq:morphpreg}) est un isomorphisme si et seulement si, pour tout  $X\in \pgo^{\reg}$, le centralisateur $G_X$ de $X$ dans $G$ est inclus dans $P$. En utilisant la décomposition de Jordan $X=X_s+X_n$, on a $G_X=G_{X_s}\cap G_{X_n}$ et $G_{X_s}$ est un groupe réductif connexe puisque $X_s$ est semi-simple. Comme $P\cap G_{X_s}$ est un sous-groupe parabolique de $G_{X_s}$ et que $X_n$ est $G_{X_s}$-régulier, on est ramené au cas où $X$ est nilpotent, ce qu'on suppose désormais. Soit $B\subset P$ un sous-groupe de Borel dont l'algèbre de Lie contient $X$. Par $G$-régularité de $X$, un tel sous-groupe de Borel est unique donc $G_{X}$ est inclus dans le normalisateur de $B$ dans $G$ qui est $B$, donc qui est inclus dans $P$. 

On sait bien que  $I^G_{|\ggo^{\reg}}$ est un schéma en groupes commutatifs lisse sur $\ggo^{\reg}$. Par changement de base, il en est de même de $I^G_{|\pgo^{\reg}}$ et donc de $I^P_{|\pgo^{\reg}}$.
\end{preuve}
\end{paragr}

\begin{paragr}[Centralisateurs réguliers.] --- \label{S:centralisateurs-reg}Soit $M\in \lc^G(T)$ un sous-groupe de Levi semi-standard. Soit $\eps_M$ une section de Kostant du morphisme caractéristique $\chi_M  : \mgo \to\car_M$ (cf. § \ref{S:Kostant}). Cette section est à valeurs dans l'ouvert des éléments $M$-réguliers. Soit   
$$J^M=\eps_M^*I^M$$
qu'on appelle le schéma en groupes des centralisateurs réguliers sur $\car_M$. C'est un schéma en groupes commutatifs, lisse sur $\car_M$ (cf. lemme  \ref{lem:IP}). 

À la suite de Ngô, on introduit l'isomorphisme
\begin{equation}
  \label{eq:morphNgo} (\chi_G^*J^G)_{|\ggo^{\reg}} \to I^G_{\ggo^{\reg}}
\end{equation}
qui, à un couple $(g,X)$ dans $G\times_k \ggo^{\reg}$ tel que $\Ad(g)\eps_G(\chi_G(X))=\eps_G(\chi_G(X))$, associe $(hgh^{-1},X)$ où $h$ est un élément de $G$ tel que $\Ad(h)\eps_G(\chi_G(X))=X$. Comme $X$ et $\eps_G(\chi_G(X))$ sont $G$-régulier et ont même image dans $\car_G$, un tel $h$ existe et est bien défini à translation près à gauche par un élément du centralisateur $G_X$ de $X$. Par $G$-régularité de $X$, on sait que $G_X$ est commutatif. Le morphisme  (\ref{eq:morphNgo}) est donc bien défini et c'est clairement un isomorphisme. Voici une proposition qui généralise légèrement la proposition  3.2 de \cite{Ngo1}.

\begin{proposition} \label{prop:IP} Pour tout $P\in \fc^G(T)$, il existe un unique morphisme de schémas en groupes 
$$(\chi_G^* J^G)_{|\pgo} \to I^P$$
qui prolonge l'isomorphisme
$$(\chi_G^*J^G)_{|\pgo^{\reg}} \to I^G_{|\pgo^{\reg}}=I^P_{|\pgo^{\reg}}$$
obtenu par composition de la restriction à $\pgo^{\reg}$ du morphisme (\ref{eq:morphNgo}) avec  l'isomorphisme du lemme \ref{lem:IP}.

En outre, le diagramme ci-dessous est commutatif 

 \begin{equation}
   \label{eq:triangle}
 \xymatrix{    (\chi_G^* J^G)_{|\pgo}   \ar[r]     \ar[rd] &   I^P  \ar[d] \\ & I^G_{|\pgo}}
\end{equation}
où la flèche verticale est l'injection canonique et la flèche diagonale est la restriction à $\pgo$ du morphisme $\chi_G^* J^G\to I^G$.
\end{proposition}

\begin{preuve} Comme le schéma $(\chi_G^* J^G)_{|\pgo}$ est normal, que le complémentaire $\pgo-\pgo^{\reg}$ est de codimension $\geq 2$ dans $\pgo$ et que  $I^P$ est affine, le morphisme 
$$(\chi_G^*J^G)_{|\pgo^{\reg}} \to I^G_{|\pgo^{\reg}}=I^P_{|\pgo^{\reg}}$$
se prolonge de manière unique en un morphisme de schémas en groupes sur $\pgo$. Le même argument montre que, puisque le triangle (\ref{eq:triangle}) est commutatif sur $\pgo^{\reg}$, il l'est aussi sur $\pgo$.
 \end{preuve}
\end{paragr}

\begin{paragr} Dans ce paragraphe, on démontre la proposition suivante.

  \begin{proposition}\label{prop:JGXM}
    Soit $M\in \lc^{G}(T)$ un sous-groupe de Levi. Il existe un unique morphisme de schémas en groupes sur $\car_M$
    \begin{equation}
      \label{eq:morphcocar}
      (\chi_G^M)^* J^G \to X_*(M)\otimes_{\ZZ} \mathbb{G}_{m,\car_M}
    \end{equation}
    tel que pour tout $P\in \pc^G(M)$ on ait un diagramme commutatif
 \begin{equation}
   \label{eq:triangle2}
 \xymatrix{   \chi_P^*(\chi_G^M)^* J^G=(\chi_G^*J^G)_{|\pgo}   \ar[r]     \ar[rd] &   I^P  \ar[d] \\  & X_*(M)\otimes_{\ZZ} \mathbb{G}_{m,\pgo}}    
\end{equation}
où 
\begin{itemize}
\item la flèche diagonale est déduite du morphisme (\ref{eq:morphcocar}) par le changement de base 
$$\chi_P : \pgo\to \car_M \ ;$$
\item la flèche horizontale est celle définie dans la proposition \ref{prop:IP} ;
\item la flèche verticale est induite par le morphisme canonique 
$$P \to X_*(P)\otimes_{\ZZ}\Gmk =X_*(M)\otimes_{\ZZ}\Gmk.$$
\end{itemize}
  \end{proposition}

  \begin{remarque}
    Le $\ZZ$-module $X_*(M)$ est libre de rang fini égal à la dimension du centre de $M$. Par conséquent, $X_*(M)\otimes_{\ZZ}\Gmk$ est un tore sur $k$ du même rang.
  \end{remarque}

  \begin{preuve}Considérons les deux projections $p_1$ et $p_2$
$$ \pgo\times_{\car_M} \pgo \stackrel{\displaystyle \longrightarrow}{\longrightarrow} \pgo.$$
Soit $p=\chi_P\circ p_1=\chi_P\circ p_2$. Par descente fidèlement plate par le morphisme $\chi_P$, il s'agit de voir que le morphisme
\begin{equation}
  \label{eq:adescendre}
  \chi_P^*(\chi_G^M)^* J^G\to X_*(M)\otimes_{\ZZ} \mathbb{G}_{m,\pgo},
\end{equation}
qui fait commuter le diagramme (\ref{eq:triangle2}), appartient au noyau de 
$$\Hom_{\pgo}( (\chi_G^* J^G)_{|\pgo},  X_*(M)\otimes_{\ZZ} \mathbb{G}_{m,\pgo})  \stackrel{\displaystyle \longrightarrow}{\longrightarrow}  \Hom_{\pgo\times_{\car_M} \pgo} ( p^* (\chi_G^M)^* J^G,  X_*(M)\otimes_{\ZZ} \mathbb{G}_{m,\pgo\times_{\car_M} \pgo}).$$
Il suffit de le vérifier au-dessus de l'ouvert $\pgo^{\reg}\times_{\car_M}\pgo^{\reg}$. Or, pour $i=1,2$,  l'image du morphisme (\ref{eq:adescendre}) par $p_i^*$ s'explicite sur l'ouvert  $\pgo^{\reg}\times_{\car_M}\pgo^{\reg}$ comme l'application qui, à $(g,X_1,X_2) \in G\times_k \pgo^{\reg}\times_{\car_M}\pgo^{\reg}$ tel que $\Ad(g)\eps_G(\chi_G(X_i))=\eps_G(\chi_G(X_i))$, associe 
$$\la\in X^*(P) \mapsto  \la(h_igh_i^{-1}),$$
où $h_i$ conjugue $\eps_G(\chi_G(X_i))$ et $X_i$. Notons que $h_igh_i^{-1}$ appartient au centralisateur  $G_{X_i}$ de $X_i$ dans $G$ donc appartient à $P$ (rappelons l'inclusion $G_{X_i}\subset P$, cf. lemme \ref{lem:IP}). Or on sait que $X_1$ et $X_2$ sont conjugués dans $P$ (cf. lemme \ref{lem:Pcj}). Il s'ensuit qu'on a $h_2\in P G_{X_1}h_1 \subset Ph_1$. On a donc, pour tout $\la\in X^*(P)$,
$$ \la(h_1gh_1^{-1})=\la(h_2gh_2^{-1}),$$
ce qui était l'égalité à vérifier.

On a donc montré l'existence du morphisme (\ref{eq:morphcocar}) qui rend le diagramme (\ref{eq:triangle2}) commutatif pour un sous-groupe  parabolique $P\in\pc(M)$. On va maintenant vérifier que le morphisme (\ref{eq:morphcocar}) qu'on vient de  définir ne dépend pas du choix de $P$. Sur l'ouvert  dense $\car_M^{\reg}$ où $D^G$ ne s'annule pas, le morphisme (\ref{eq:morphcocar}) est l'application qui, à un couple $(g,a)\in G\times_k \car_M$ tel que $\Ad(g)\eps_G(\chi_G^M(a))=\eps_G(\chi_G^M(a))$, associe $\la\in X^*(P) \mapsto \la(hgh^{-1})$ où $h\in G$ est tel que $\Ad(h) \eps_G(\chi_G^M(a))$ appartienne à la fibre de $\chi_P$ au-dessus de $a$. Un tel élément $h$ existe : soit $Y\in \mgo$ tel que $\chi_P(Y)=a$. Alors $Y$ est semi-simple et $G$-régulier donc $G$-conjugué à $\eps_G(\chi_G^M(a))$. Un tel $h$ est bien défini à un élément de $P$ à droite près : la fibre de $\chi_P$ au-dessus de $a$ est une $P$-classe de conjugaison (cf. lemme \ref{lem:Pcj}).  En particulier,  $\la(hgh^{-1})$ est parfaitement défini pour tout $\la\in X^*(P)$ (par un argument déjà évoqué, on a  $hgh^{-1}\in P$). Il est clair qu'on aurait pu tout aussi bien définir $h$ comme un élément de $G$ tel $\Ad(h) \eps_G(\chi_G^M(a))$ appartienne à la fibre de $\chi_M$ au-dessus de $a$, d'où l'indépendance en $P$.
 \end{preuve}
\end{paragr}

 \begin{paragr}[Description galoisienne de $J^G$.] --- On laisse le soin au lecteur de vérifier le lemme suivant.

 \begin{lemme}\label{lem:IT} 
Le morphisme 
$$(\chi_G^T)^*J^G \to X_*(T)\otimes_\ZZ \mathbb{G}_{m,\tgo}$$
de la proposition \ref{prop:JGXM} est $W$-équivariant où l'action sur le but est l'action diagonale de $W$.
\end{lemme}

Pour tout $\tgo$-schéma $X$, soit $(\chi^T_G)_* X$ le $\car_G$-schéma obtenu par restriction des scalaires à la Weil de $X$ (rappelons que le $k[\car_G]$-module $k[\tgo]$ est libre de rang fini, égal à l'ordre de $W$). Le morphisme d'adjonction $J^G\to (\chi^T_G)_* (\chi^T_G)^* J^G$ se factorise par le morphisme
$$J^G \longrightarrow ((\chi^T_G)_* (\chi^T_G)^* J^G)^W$$
où l'exposant $W$ désigne le sous-schéma des points fixes sous $W$. En utilisant la $W$-équivariance du morphisme $(\chi_G^T)^*J^G \to X_*(T)\otimes_\ZZ \mathbb{G}_{m,\tgo}$ (cf. lemme \ref{lem:IT}),  on obtient un morphisme
\begin{equation}
  \label{eq:JGversT}
  J^G \longrightarrow (X_*(T)\otimes_\ZZ (\chi^T_G)_*\mathbb{G}_{m,\tgo})^W
\end{equation}
où le but est le sous-schéma fermé des points fixes sous $W$ dans  $X_*(T)\otimes_\ZZ \chi^T_*\mathbb{G}_{m,\tgo}$ ; c'est d'ailleurs un schéma en groupes commutatifs lisse sur $\car_G$ (cf. \cite{ngo2} lemme 2.4.1) dont l'algèbre de Lie, vue comme  un $\oc_{\car_G}$-module, est
$$\Lie((X_*(T)\otimes_\ZZ (\chi^T_G)_*\mathbb{G}_{m,\tgo})^W)=(X_*(T)\otimes_\ZZ (\chi^T_G)_*\oc_{\tgo})^W.$$
Rappelons la proposition suivante, due à Ngô, qui fournit une description «galoisienne» de l'algèbre de Lie de $J^G$.

\begin{proposition}\label{prop:J-galois}
Le morphisme 
$$J^G \longrightarrow (X_*(T)\otimes_\ZZ (\chi^T_G)_*\mathbb{G}_{m,\tgo})^W$$
est un isomorphisme sur l'ouvert $\car_G^{\reg}$. Il induit par ailleurs un isomorphisme entre les sous-schémas ouverts des composantes neutres des fibres. En particulier, ce morphisme induit un isomorphisme entre les algèbres de Lie vues comme $\oc_{\car_G}$-modules 
$$\Lie(J^G)=(X_*(T)\otimes_\ZZ (\chi^T_G)_*\oc_{\tgo})^W.$$
\end{proposition}

\begin{preuve}
  On renvoie le lecteur à \cite{ngo2} proposition 2.4.2 et corollaire 2.4.8.
\end{preuve}

\end{paragr}

\section{La fibration de Hitchin tronquée}\label{sec:Hitchin-tronque}

\begin{paragr}\label{S:cbe} Soit $C$ une courbe projective, lisse et connexe sur $k$ de genre $g$. Soit $D=2D'$ un diviseur  effectif de degré $>2g$ et $\infty$ un point de  $C(k)$ qui n'est pas dans le support de $D$.
 Soit $\lc_D$ le $\Gmk$-torseur sur $C$ associé à $D$, muni de sa trivialisation canonique sur l'ouvert complémentaire du support de $D$. Pour tout $k$-schéma $V$ muni d'une action de $\Gmk$, soit
$$V_D=\lc_D\times_k^{\Gmk} V$$
le produit contracté : c'est un schéma sur $C$ qui est une fibration localement triviale de fibre type $V$.  Soit $M\in \lc(T)$. L'action par homothétie du groupe multiplicatif $\Gmk$ sur $\tgo$ induit une action de $\Gmk$ sur $\car_M$ pour laquelle le morphisme caractéristique $\chi_M^T$ est  $\Gmk$-équivariant. Par produit contracté avec $\lc_D$, on obtient des $C$-morphismes 
$$\chi_{M,D}^T : \tgo_D \to \car_{M,D}$$
et, pour tout $P\in \pc(M)$,
$$\chi_{P,D} : \pgo_D  \to \car_{M,D}.$$
\end{paragr}

\begin{paragr}[Champ de Hitchin.] --- C'est le champ algébrique noté $\mc=\mc_G$ qui classifie les triplets de Hitchin $m=(\mathcal{E},\theta ,t)$ o\`{u}: 
\begin{itemize}
\item $\mathcal{E}$ est un $G$-torseur sur $C$ ;
\item $\theta \in H^{0}(C,\Ad_D(\ec))$ ;
\item $t\in \mathfrak{t}^{G-\textrm{reg}}$ avec $\chi(t)=\chi(\theta_{\infty})$ dans $\mathfrak{car}_G$.
\end{itemize}
 On a noté 
 $$\Ad_D(\ec)=\ec\times_k^{G,\Ad_D}\ggo_D$$
le fibré vectoriel sur $C$ obtenu lorsqu'on pousse le $G$-torseur $\ec$ par le morphisme 
$$\Ad_D \ : \ G\to \Aut(\ggo_D)$$ 
induit par l'action adjointe.

Soit $[\ggo_D/G]$ le champ sur $C$ défini comme le quotient de $\ggo_D$ par $G\times_k C$. Par définition, pour tout $C$-schéma test $S$, le groupoïde $[\ggo_D/G](S)$ est le groupoïde des couples $(\ec,\theta)$ où $\ec$ est un $G$-torseur sur $S$ et $\theta$ est une section sur $S$ de $\Ad_D(\ec)=\ec\times_k^{G,\Ad_D}\ggo_D$. 

 Le morphisme $\chi_G$ est à la fois $\Gmk$-équivariant et $G$-invariant. Il induit donc un morphisme de $C$-champs
 \begin{equation}
   \label{eq:[chiG]}
   [\chi_{G,D}] \ :\ [\ggo_D/G] \to \car_{G,D}.
 \end{equation}
 
On peut alors donner une définition équivalente de $\mc$. Soit $S$ un $k$-schéma test. La donnée d'un triplet  $m=(\ec,\theta,t)\in \mc(S)$ est équivalente à la donnée d'un couple $m=(h_m,t)$ formé d'un $C$-morphisme 
$$h_m  \ : \ C\times_k S \longrightarrow [\ggo_D/G]$$
et d'un point $t\in \tgo^{\reg}(S)$ qui vérifie 
$$ ([\chi_{G,D}]\circ h_m)(\infty_S)=\chi^T_G(t)$$
où $\infty_S$ se déduit de $\infty$ par le changement de base $S\to k$.
 
\end{paragr}

\begin{paragr}[Lissité de $\mc$ sur $k$.] --- On a la proposition suivante due à Biswas et Ramamanan.

\begin{proposition} (cf. \cite{Biswas-Ramanan} et \cite{ngo2} théorème 4.14.1)
Le champ algébrique $\mathcal{M}$ est lisse sur $k$.
\end{proposition}
  
\end{paragr}

\begin{paragr}[Espace caractéristique.] ---  \label{S:esp-car}Dans la (seconde) définition de $\mc$, si l'on remplace le champ quotient  $[\ggo_D/G]$ par le quotient catégorique de $\ggo_D$ par $G$, qui n'est autre que $\car_{G,D}$, on obtient un schéma quasi-projectif, lisse sur $k$, noté $\Ac=\Ac_G$ et appelé \emph{espace caractéristique}, qui classifie les couples $a=(h_a,t)$ où
\begin{itemize}
\item $h_a\in H^0(C,\car_{G,D})$ ;
\item $t\in \mathfrak{t}^{G-\textrm{reg}}$ avec $\chi(t)=h_a(\infty)$.
\end{itemize}

\medskip

Pour chaque $M\in \mathcal{L}(T)$,  on dispose de la base $\mathcal{A}_{M}$ de la fibration de Hitchin pour $M$ et d'un ouvert
$$
\mathcal{A}_{M}^{G\-reg}\subset \mathcal{A}_{M}
$$
o\`{u} l'on demande que $t$ soit dans $\mathfrak{t}^{G\-reg}\subset \mathfrak{t}^{M\-reg}$. Le morphisme
\begin{equation}
  \label{eq:AM-AG}
  \chi^M_{G} : \mathcal{A}_{M}^{G\-reg}\hookrightarrow \mathcal{A}
\end{equation}
qui envoie $(h,t)$ sur $(\chi_{G,D}^M \circ h,t)$ est une immersion fermée (cf. \cite{LFPI} proposition 3.5.1).
\smallskip

\textit{Dans la suite, seul interviendra l'ouvert $\mathcal{A}_{M}^{G\-reg}$ que l'on notera dorénavant simplement $\mathcal{A}_{M}$.}
\smallskip

\begin{definition}\label{def:partieell}
La partie elliptique $\mathcal{A}^{\mathrm{ell}}=\mathcal{A}_{G}^{\mathrm{ell}}$ de $\mathcal{A}$ est l'ouvert complémentaire de la réunion des $\mathcal{A}_{M}$ pour $M\in \mathcal{L}(T)$, $M\not=G$. 
\end{definition}

On notera que l'ouvert elliptique est non vide pour des raisons de dimension.
On définit de m\^{e}me $\mathcal{A}_{M}^{\mathrm{ell}}$. On a la proposition suivante (cf. \cite{LFPI} proposition 3.6.2). 

\begin{proposition}  \label{prop:reunionAM}
 Le schéma $\Ac_G$ est la réunion disjointe des sous-schémas localement fermés $\Ac_{M}^{\el}$
pour $M\in \lc(T)$  
$$\Ac_G=\bigcup_{M\in \lc^G(T)} \Ac_{M}^{\el}.$$
\end{proposition}

\end{paragr}

\begin{paragr}[Fibration de Hitchin.] --- On la définit comme suit à l'aide du morphisme $[\chi_{G,D}]$ de la ligne (\ref{eq:[chiG]}).

\begin{definition}
La fibration de Hitchin est le morphisme 
$$
f=f_{G}:\mathcal{M}\rightarrow \mathcal{A}
$$
qui envoie $m=(h_m,t)$ sur $a=([\chi_{G,D}]\circ h_m,t)$.
\end{definition}

Soit $a=(h_a,t)\in \Ac$. La fibre de Hitchin $\mc_a=f^{-1}(a)$ classifie les sections $h_m$ du champ $[\ggo_D/G]$ qui rendent le diagramme ci-dessous commutatif
\begin{equation}
  \label{eq:diag-Hitchin}
  \xymatrix{ C \ar[r]^{h_m}  \ar[rd]_{h_a} & [\ggo_D/G] \ar[d]^{[\chi_{G,D}]}  \\ &\car_{G,D} }
\end{equation}
\end{paragr}

\begin{paragr}[Réduction d'un triplet de Hitchin.] --- Soit $m=(h_m,t)\in \mc_G$ et $P\in \pc(M)$ un sous-groupe parabolique de $G$ de Levi $M\in \lc(T)$. Une réduction de $m$ à $P$ est, par définition, la donnée d'une section 
$$h_{P} :   C \to [\pgo_D/P]$$
du $C$-champ quotient, noté   $[\pgo_D/P]$ de $\pgo_D$ par $P\times_k C$ qui  vérifie les deux conditions suivantes :
\begin{enumerate}
\item on a $h_m=i_P\circ h_P$ où
$$i_P : [\pgo_D/P] \to [\ggo_D/G]$$
 est le morphisme canonique ;
\item le couple $(h_{a_P},t)$ défini par
$$h_{a_P}=[\chi_{P,D}]\circ h_P,$$
appartient à $\Ac_M$, autrement dit on a 
$$h_{a_P}(\infty)=\chi_P^T(t).$$
\end{enumerate}

Avec les notations ci-dessus, si $h_P$ est une réduction de $m$ à $P$, on a le diagramme commutatif suivant :

\begin{equation}
  \label{eq:diag-reduction}
\xymatrix{      C \ar[r]^{h_P}   \ar@/^3pc/[rr]_{h_m}  \ar@/_4pc/[rrd]_{h_a} \ar[rd]_{h_{a_P}} & [\pgo_D/P] \ar[d]^{[\chi_{P,D}]} \ar[r]^{i_P} &[\ggo_D/G] \ar[d]^{[\chi_{G,D}]}  \\  &\car_{P,D} \ar[r]_{\chi^M_{G,D}}     &\car_{G,D} }.
\end{equation}

Il est clair sur le diagramme ci-dessus que pour que $m$ admette une réduction à $P$ il faut que $f(m)$ appartienne à $\Ac_M$. En fait, cela suffit. Plus précisément, on a la proposition suivante (cf. proposition 4.10.1 de \cite{LFPI}).

\begin{proposition} \label{prop:reduc}Soit $m\in \mc_G$ et $M\in \lc(M)$ tel que $f(m)\in \Ac_{M}^{\el}$ (cf.  proposition \ref{prop:reunionAM}). Pour tout $P\in \fc(M)$, il existe une et une seule réduction $h_P$ de $m$ à $P$.  
\end{proposition}
\end{paragr}

\begin{paragr}[Degré $\deg_P$.] --- \label{S:degP} Pour tout schéma en groupes $H$ sur $C$, soit $\mathrm{B}(H)$ le $C$-champ classifiant de $H$ dont le groupoïde $\mathrm{B}(H)(S)$, pour tout $C$-schéma test $S$, est le groupoïde des $H\times_C S$-torseurs sur $S$. Ce champ est encore le quotient de $C$ par $H$ agissant trivialement. 

Le morphisme structural $\pgo_D \to C$ induit un morphisme de $C$-champs
\begin{equation}
  \label{eq:B1}
  [\pgo_D/P] \longrightarrow \mathrm{B}(P\times_k C).
\end{equation}
Par ailleurs, le morphisme évident $P\to X_*(P)\otimes_\ZZ \Gmk$ induit un morphisme de $C$-champs 

\begin{equation}
  \label{eq:B2}
\mathrm{B}(P\times_k C) \longrightarrow \mathrm{B}( X_*(P)\otimes_\ZZ \mathbb{G}_{m,C} ).
\end{equation}
On a un morphisme 
$$\mathrm{B}(\mathbb{G}_{m,C})(C)  \longrightarrow \ZZ$$
donné par le degré d'un $\mathbb{G}_{m,C}$-torseur. On en déduit de manière évidente un morphisme 
\begin{equation}
  \label{eq:BdegP}
   \mathrm{B}( X_*(P)\otimes_\ZZ \mathbb{G}_{m,C})(C)  \longrightarrow X_*(P).
\end{equation}
En composant le morphisme  (\ref{eq:BdegP}) avec les morphismes (\ref{eq:B1}) et (\ref{eq:B2}), on obtient un morphisme degré
\begin{equation}
  \label{eq:degP}
  \deg_P : [\pgo_D/P](C)  \longrightarrow X_*(P).
\end{equation}
\end{paragr}

\begin{paragr}[Convexe associé à $m\in\mc$.] --- \label{S:cvx} Soit $M\in \lc(T)$ et $m\in \mc$ tel que $f(m)\in \Ac_M^{\el}$. Pour tout $P\in \pc(M)$, soit $h_P$ l'unique réduction de $m$ à $P$ (cf. proposition \ref{prop:reduc}). On a défini en (\ref{eq:degP}) le degré $\deg_P(h_P)$ de la section $h_P$ de $[\pgo_D/P]$. 

Soit $\mathcal{C}_m\subset \ago_T$, le
 convexe obtenu lorsqu'on translate l'enveloppe convexe des points 
$$-\deg_P(h_P)\in X_*(P)=X_*(M) \subset \ago_M$$
pour $P\in \pc(M)$  par les vecteurs du sous-espace $\ago_T^M\oplus \ago_G$.
\end{paragr}

\begin{paragr}[La $\xi$-stabilité.] --- \label{S:xistab}Soit $\xi\in \ago_T$. Soit $m\in \mc$. On dit que $m$ est $\xi$-semi-stable si le convexe $\mathcal{C}_m$ du paragraphe \ref{S:cvx} contient $\xi$. Soit $\mc^\xi$ le sous-champ de $\mc$ formé de points $\xi$-stables et 
$$\mc^{\el}=\mc\times_{\Ac_G} \Ac_G^{\el}$$
l'ouvert elliptique de $\mc$. Soit 
$$f^\xi : \mc^\xi \to \Ac_G$$ 
et 
$$f^{\el} : \mc^{\el}\to\Ac_G^{\el}$$
les restrictions du morphisme de Hitchin $f$ à $\mc^\xi$ et $\mc^{\el}$.

Voici le principal résultat de \cite{LFPI} (\emph{ibid.} proposition 6.1.4 et théorème 6.2.2).

  \begin{theoreme}\label{thm:rappel-lfpI}
Le champ $\mc^\xi$ est un sous-champ ouvert de $\mc$, donc lisse sur $k$, qui contient l'ouvert elliptique $\mc^{\el}$. 

Si $G$ est semi-simple et si $\xi$ est \emph{en position générale} , le champ $\mc^\xi$  est un champ de Deligne-Mumford et le morphisme $f^\xi$ est propre.
  \end{theoreme}

  \begin{remarque} \label{rq:position}
    L'expression $\xi$ \emph{en position générale} signifie que $\xi$ n'appartient pas   à certains hyperplans rationnels de $\mathfrak{a}_T$ (cf. \cite{LFPI} définition 6.1.3). Elle implique en particulier que pour un triplet de Hitchin les notions de $\xi$-semi-stabilité et $\xi$-stabilité sont équivalentes.
  \end{remarque}

\end{paragr}

\section{Action du champ de Picard  $\Jc_a$}\label{sec:Picard}

\begin{paragr}[Action de $\Gmk$ sur les centralisateurs.] --- \label{S:actionGm}On a défini  au § \ref{S:centralisateurs-reg} le schéma en groupes  $J^G$ des centralisateurs réguliers sur $\ggo$. Soit $(g,a)\in J^G$ et $x\in \Gmk$. Les éléments de $\ggo^{\reg}$ définis par $\eps_G(x\cdot a)$ et $x\cdot \eps_G(a)$ (où $\cdot$ désigne respectivement l'action de $\Gmk$ sur $\car_G$ définie au § \ref{S:cbe} et l'action de $\Gmk$ sur $\ggo$ par homothétie) sont $G$-conjugués puisqu'ils ont même image par $\chi_G$. Il existe donc $h_x\in G$ tel que 
$$\eps_G(x\cdot a)=\Ad(h_x)(x\cdot \eps_G(a)).$$
On peut vérifier que $h_x g h_x^{-1}$ centralise $\eps_G(x\cdot a)$ et que cet élément ne dépend pas du choix de $h_x$ (les centralisateurs réguliers étant commutatifs). On pose alors 
$$x\cdot (g,a)= (h_x g h_x^{-1},x\cdot a)$$
ce qui définit une action de $\Gmk$ sur $J^G$ pour laquelle le morphisme canonique $J^G \to \car_G$ est $\Gmk$-équivariant. Par la construction du § \ref{S:cbe}, on en déduit un schéma en groupes $J^G_D$ au-dessus de $\car_{G,D}$.

Pour tout sous-groupe parabolique $P\in \fc(T)$, le schéma en groupes $I^P$ sur $\pgo$ des centralisateurs est muni de l'action de $\Gmk$ définie pour tous $x\in \Gmk$ et $(p,X)\in I^P$ par
$$x\cdot(p,X)=(p,x\cdot X)$$
où $\cdot$ est l'action par homothétie de $\Gmk$ sur $\pgo$. De nouveau, la construction du § \ref{S:cbe} donne un schéma en groupes $I^P_D$ au-dessus de $\pgo_{D}$.

Il résulte de la proposition \ref{prop:IP} que le morphisme  $(\chi_G^* J^G)_{|\pgo} \to I^P$ qui y est défini  est $\Gmk$-équivariant. On en déduit un morphisme de schémas en groupes 
\begin{equation}
  \label{eq:morph_IPD}
  (\chi_{G,D}^* J^G_D)_{|\pgo_D} \to I^P_D
\end{equation}
au-dessus de $\pgo_D$.

Soit $M\in \lc(T)$ et $\chi^M_{G,D}: \car_{M,D}\to \car_{G,D}$ le morphisme déduit du morphisme $\Gmk$-équivariant $\chi^M_G : \car_M \to \car_G$. Il résulte de la proposition \ref{prop:JGXM} qu'il existe un unique morphisme de schémas en groupes 
\begin{equation}
  \label{eq:morph_XMcarD}
(\chi^M_{G,D})^*J_D^G\to X_*(M)\otimes_\ZZ \mathbb{G}_{m,\car_{M,D}}
\end{equation}
qui, pour tout $P\in \pc(M)$, est le composé de (\ref{eq:morph_IPD}) avec le morphisme évident
$$I^P_D \to X_*(M)\otimes_\ZZ \mathbb{G}_{m,\car_{M,D}}.$$

\end{paragr}

\begin{paragr}[Le champ de Picard $\Jc_a$.] --- \label{S:Jc_a} Soit $a=(h_a,t)\in\Ac$. Soit 
$$J_a=h_a^* J^G_D.$$
C'est un schéma en groupes commutatifs lisse sur $C$. Soit $\Jc_a$ le $k$-champ des $J_a$-torseurs au-dessus de $C$. Plus exactement, pour tout $k$-schéma test $S$, le groupoïde $\Jc_a(S)$ est le groupoïde des $C$-morphismes  
$$C\times_k S \to \mathrm{B}(J_a).$$

Comme le schéma en groupes $J_a$ est commutatif, le produit contracté de deux $J_a$-torseurs est encore un $J_a$-torseur. Ce produit fait de $\Jc_a$ un champ de Picard au-dessus de $\Spec(k)$.
\end{paragr}

\begin{paragr}[Morphisme degré.] --- \label{S:morphdeg}Soit $M\in \lc(T)$ et $a_M=(h_{a_M},t)\in \Ac_M$. Soit $a$ l'image de $a_M$ dans $\Ac_G$ par le morphisme $\chi^M_G$, cf. l. (\ref{eq:AM-AG}). On a donc $h_a=\chi^M_{G,D}\circ  h_{a_M}$ et $J_a=h_{a_M}^* (\chi^M_{G,D})^* J^G_D$. Par le changement de base $h_{a_M} : C \to \car_{M,D}$, on déduit du morphisme (\ref{eq:morph_XMcarD}) le morphisme de schémas en groupes
$$J_a \to X_*(M)\otimes_\ZZ \mathbb{G}_{m,C}$$
et le morphisme de $C$-champs de Picard entre les classifiants
\begin{equation}
  \label{eq:classifiant}
  \mathrm{B}(J_a) \to \mathrm{B}( X_*(M)\otimes_\ZZ \mathbb{G}_{m,C}).
\end{equation}
Comme on dispose d'un morphisme degré, cf. l.(\ref{eq:BdegP}),
$$\mathrm{B}( X_*(M)\otimes_\ZZ \mathbb{G}_{m,C})(C) \to X_*(M),$$
on en déduit, par composition, un morphisme degré
\begin{equation}
  \label{eq:degreJ_a}
  \deg_M : \Jc_a \to X_*(M)
\end{equation}
pour tout $a\in \Ac_M$, qui est un homomorphisme de champs de Picard (si l'on voit  $X_*(M)$ comme un $k$-champ de Picard constant).

\begin{proposition}\label{prop:compatibilite} Soit $M\subset L$ deux sous-groupes de Levi dans $\lc^G(T)$. Soit 
$$p^M_L : X_*(M)\to X_*(L)$$
la projection duale du morphisme de restriction
$$X^*(L)\to X^*(M).$$
Soit $a\in \Ac_M$. On a l'égalité suivante entre morphismes de $\Jc_a$ dans $X_*(L)$
$$p_L^M\circ \deg_M=\deg_L.$$
  \end{proposition}

  \begin{preuve} Soit $Q\in \pc(L)$ et $P\in \pc(M)$ tels que $P\subset Q$. Partons du diagramme commutatif
$$\xymatrix{\pgo \ar[r] \ar[d]_{\chi_P}  & \qgo \ar[d]^{\chi_Q }\\ \car_M \ar[rd]_{\chi^M_G} \ar[r]^{\chi^M_L} & \car_L \ar[d]^{\chi_G^L} \\ & \car_G        }.$$
Rappelons qu'on a défini au § \ref{S:centralisateurs-reg} un schéma en groupes $J^G$ sur $\car_G$. D'après le diagramme ci-dessus, on a l'égalité 
$$\chi_P^* (\chi_G^M)^* J^G=(\chi_Q^* (\chi_G^L)^* J^G)_{|\pgo}.$$
On a donc un diagramme
\begin{equation}
  \label{eq:diagM-L}
  \xymatrix{\chi_P^* (\chi_G^M)^* J^G  \ar[r]  \ar[rd] &  X_*(M)\otimes_{\ZZ} \mathbb{G}_{m,\car_M} \ar[d] \\ & X_*(L)\otimes_{\ZZ} \mathbb{G}_{m,\car_L}}
\end{equation}
où la flèche horizontale est celle définie en (\ref{eq:morphcocar}) proposition \ref{prop:JGXM}, la flèche diagonale est celle définie  (\ref{eq:morphcocar}) relativement à $L$ et restreinte à $\pgo$, enfin la flèche verticale est induite par $p_M^L$. Montrons que le diagramme (\ref{eq:diagM-L}) est commutatif. Il suffit de le vérifier au-dessus de l'ouvert dense $\pgo^{\reg}$. Le morphisme horizontal, resp. diagonal, de (\ref{eq:diagM-L}) se factorise par $I^P$, resp. $I^Q_{\pgo}$ (cf.  proposition \ref{prop:JGXM}). On a donc un diagramme
$$ \xymatrix{(\chi_P^* (\chi_G^M)^* J^G)_{|\pgo^{\reg}}  \ar[r]  \ar[rd] & I^P_{|\pgo^{\reg}} \ar[r] \ar[d] &  X_*(M)\otimes_{\ZZ} \mathbb{G}_{m,\car_M} \ar[d] \\ &   I^Q_{|\pgo^{\reg}} \ar[r]   &  X_*(L)\otimes_{\ZZ} \mathbb{G}_{m,\car_L}}$$
dont il s'agit de vérifier la commutativité. Or, le premier triangle est formé uniquement d'isomorphismes et est clairement commutatif (cf. lemme \ref{lem:IP} et (\ref{eq:morphNgo})). On vérifie sans peine la commutativité du  carré de droite.

Pour tout $a\in \Ac_M$, on déduit de la commutativité de (\ref{eq:diagM-L}) un diagramme commutatif
$$ \xymatrix{ \mathrm{B}(J_a) \ar[r] \ar[rd] &  \mathrm{B}(X_*(M)\otimes_{\ZZ} \mathbb{G}_{m,C}) \ar[d]\\ &  \mathrm{B}(X_*(L)\otimes_{\ZZ} \mathbb{G}_{m,C})}$$
où la flèche verticale est induite par $p^M_L$ et les autres flèches sont données par la ligne (\ref{eq:classifiant}). Pour conclure, il suffit d'invoquer la commutativité du diagramme
$$ \xymatrix{\mathrm{B}(X_*(M)\otimes_{\ZZ} \mathbb{G}_{m,C})(C) \ar[r] \ar[d] & X_*(M) \ar[d] \\  \mathrm{B}(X_*(L)\otimes_{\ZZ} \mathbb{G}_{m,C})(C)  \ar[r] & X_*(L)}$$
où les morphismes horizontaux sont les degrés et les morphismes verticaux sont induits par $p_L^M$.

  \end{preuve}

\end{paragr}

\begin{paragr}[Quotient champêtre des centralisateurs.] --- \label{S:centr-champetre}Soit $P\in \fc^G(T)$ un sous-groupe parabolique semi-standard. On définit une action du groupe $P$ sur le schéma en groupes $I^P$ ainsi : pour tout $h\in P$ et tout $(p,X)\in I^P$, on a
$$h\cdot(p,X)= (hph^{-1},\Ad(h)X).$$
Cette action commute à l'action de $\Gmk$ par homothétie sur le premier facteur. On en déduit donc une action de $P$ sur $I^P_D$ pour laquelle le morphisme $I^P_D \to \pgo_D$ est $P$-équivariant.

Le champ quotient $[I^P_D/P]$ définit un schéma en groupes au-dessus du champ $[\pgo_D/P]$. Soit $h_P$ une section au-dessus de $C$ du champ $[\pgo_D/P]$. En tirant en arrière le champ  $[I^P_D/P]$ par la section $h_P$, on obtient un schéma en groupes sur $C$ noté $h_P^*[I^P_D/P]$. Donnons une description concrète de ce groupe : la section $h_P$ s'identifie à la donnée d'un $P$-torseur sur $C$ ainsi que d'une section $\theta$ de $\Ad{P,D}(\ec)=\ec\times^{P,\Ad}\pgo_D$. Le schéma en groupes $h_P^*[I^P_D/P]$ est le schéma en groupes $\Aut_P(\ec,\theta)$ des automorphismes de $\ec$ qui centralisent $\theta$.

On définit une action du groupe $G$ sur le schéma en groupes $\chi_G^* J^G$ de la manière suivante : pour tout $h\in G$ et tout $(g,X)\in \chi_G^* J^G$, on a
$$h\cdot(g,X)= (g,\Ad(h)X).$$
On vérifie que cette action commute à celle de $\Gmk$ (déduite de celle définie sur $J^G$ au paragraphe \ref{S:actionGm}. On en déduit une action de $G$ sur $\chi_{G,D}^* J^G_D$. On vérifie que le morphisme  $$(\chi_{G,D}^* J^G_D)_{|\pgo_D} \to I^P_D$$
défini au paragraphe \ref{S:actionGm} est $P$-équivariant. 

Le champ quotient $[(\chi_{G,D}^* J^G_D)/G]$ définit un schéma en groupes commutatifs lisse  sur $[\ggo_D/G]$. Comme l'action de $G$ sur le premier facteur de $\chi_G^* J^G$ est triviale, on a l'identification naturelle
\begin{equation}
  \label{eq:JD-champetre}
  [(\chi_{G,D}^* J^G_D)/G]=[\chi_{G,D}]^*J_D^G,
\end{equation}
où l'on distingue les morphismes $\chi_{G,D} : \ggo_D \to \car_{G,D}$ (cf. §\ref{S:cbe}) et $[\chi_{G,D}] : [\ggo_D/G] \to \car_{G,D}$, cf. (\ref{eq:[chiG]}).
 
\end{paragr}

\begin{paragr}[Action de $\Jc_a$ sur $\mc_a$.] --- \label{S:ActionJa}Soit $a=(h_a,t)\in \Ac$ et $m=(h_m,t)\in \mc_a$. On a défini au § \ref{S:Jc_a} un schéma en groupes $J_a$ sur $C$. D'après le diagramme (\ref{eq:diag-Hitchin}) et la ligne (\ref{eq:JD-champetre}) du paragraphe  \ref{S:centr-champetre}, on a aussi
\begin{equation}
  \label{eq:Ja-champetre}
J_a=h_m^* [\chi_{G,D}]^*J_D^G= h_m^*[(\chi_{G,D}^* J^G_D)/G].
\end{equation}
Le morphisme  $\chi_{G,D}^* J^G_D \to I^G_D$ défini en (\ref{eq:morph_IPD}) du paragraphe \ref{S:actionGm} est $G$-équivariant (cf. § \ref{S:centr-champetre}). On en déduit un morphisme de schémas en groupes
\begin{equation}
  \label{eq:morph-Aut(m)}
  J_a=h_m^*[(\chi_{G,D}^* J^G_D)/G] \to h_m^*  [I^G_D/G]=\Aut(m).
\end{equation}
La dernière égalité identifie $h_m^*  [I^G_D/G]$ au schéma en groupes $\Aut(m)$ des automorphismes de $m$ (cf. § \ref{S:centr-champetre}).

On fait alors agir le champ de Picard $\Jc_a$ sur le fibre de Hitchin $\mc_a$ de la façon suivante : à tout $m\in \mc_a$ et tout $J_a$-torseur $j$ on associe le produit contracté
$$j\cdot m=j\times^{J_a} m$$
où $J_a$ agit sur $m$ via le morphisme $J_a \to \Aut(m)$ (cf. l.(\ref{eq:morph-Aut(m)})). On vérifie que $j\times^{J_a} m \in \mc_a$.

\begin{proposition} \label{prop:translation}Soit $a\in \Ac_M^{\el}$ et $m\in \mc_a$. Pour tout $j\in \Jc_a$, les convexes associés à $m$ et $ j\cdot m$ (cf. \ref{S:cvx}) sont translatés l'un de l'autre. Plus précisément, on a 
$$\mathcal{C}_m=\mathcal{C}_{ j\cdot m}+\deg_M(j)$$
où $\deg_M$ est le morphisme degré défini en (\ref{eq:degreJ_a}) du §\ref{S:morphdeg}.  
\end{proposition}

\begin{preuve} Soit $P\in \pc(M)$ et $h_P$ la réduction correspondante de $m$ à $P$ (cf. proposition \ref{prop:reduc}). En examinant le diagramme \ref{eq:diag-reduction}, on s'aperçoit qu'on a 
$$ J_a= h_P^* i_P^* [\chi_{G,D}]^* J_D^G.$$
On vérifie qu'on a l'identification suivante de schémas en groupes sur $[\pgo_D/P]$
$$ i_P^* [\chi_{G,D}]^* J_D^G = [(J^G_D)_{|\pgo_D}/P].$$
 
Le morphisme de schémas en groupes sur  $[\ggo_D/G]$ déduit du morphisme $G$-équivariant (\ref{eq:morph_IPD}) (pour $P=G$) donne un morphisme (cf. aussi (\ref{eq:JD-champetre}))
$$ [\chi_{G,D}]^* J_D^G \to [I^G_D/G].$$
Il résulte de la proposition \ref{prop:IP} que le morphisme qui s'en déduit
$$ i_P^* [\chi_{G,D}]^* J_D^G  = [(J^G_D)_{\pgo_D}/P]  \to  i_P^*[I^G_D/G]$$
se factorise par
$$ [(J^G_D)_{|\pgo_D}/P]  \to [I^P_D/P],$$
où cette dernière flèche se déduit du morphisme $P$-équivariant  (\ref{eq:morph_IPD}) (l'équivariance est discutée au § \ref{S:centr-champetre}). Par conséquent, le morphisme
$$J_a \to h_m^*  [I^G_D/G]=\Aut(m)$$
se factorise par
\begin{equation}
  \label{eq:Aut(h_P)}
  J_a \to h_P^* [I^P_D/P]=\Aut(h_P).
\end{equation}

En particulier, si $j$ est un $J_a$-torseur, l'action de $j$ sur $m$ envoie la $P$-réduction $h_P$ sur $j\times^{J_a} h_P$ (où $J_a$ agit sur $h_P$ via (\ref{eq:Aut(h_P)})) qui est nécessairement l'unique $P$-réduction de $j\cdot m$. La proposition résulte alors de l'égalité, évidente sur les définitions,
$$\deg_P(j\times^{J_a} h_P)=\deg_M(j)+\deg_P(h_P).$$

\end{preuve}

\end{paragr}

\section{Morphisme degré et champ de Picard  $\Jc^1$}\label{sec:J1}

\begin{paragr}[Espace topologique $\bg\lc^G\bd$.] ---Soit $\bg\lc^G\bd$ l'espace topologique sobre dont l'ensemble sous-jacent est l'ensemble $\lc^G=\lc^G(T)$ et dont la topologie est la moins fine pour laquelle les parties $\lc^M$ pour  $M\in \lc^G$ sont fermées.

Soit $M\subset L$ deux éléments de $\lc^G(T)$. Soit $X_*(L)_{\lc^M}$ le faisceau constant sur $\bg \lc^M \bd$ de valeur $X_*(L)$. Soit
$$i_M : \bg\lc^M\bd \hookrightarrow \bg\lc^G\bd$$
l'inclusion canonique. Soit $i_{M,*}(X_*(L)_{\lc^M})$ le faisceau sur  $\bg \lc^G \bd$ à support dans $\lc^M$ et dont la fibre en tout point de $\lc^M$ est $X_*(L)$. On a un morphisme
\begin{equation}
  \label{eq:morph1}
  i_{M,*}(X_*(M)_{\lc^M})\to i_{M,*}(X_*(L)_{\lc^M})
\end{equation}
induit par le morphisme $p^M_L :  X_*(M)\to X_*(L)$ (cf. \ref{prop:compatibilite}). D'autre part, on a un morphisme
\begin{equation}
  \label{eq:morph2}
  i_{L,*}(X_*(L)_{\lc^L})\to i_{M,*}(X_*(L)_{\lc^M})
\end{equation}
qui, fibre à fibre, est le morphisme identique de $X_*(L)$ en un point de $\lc^M$ et qui est nul en dehors de $\lc^M$.

On définit alors un double morphisme 

\begin{equation}
  \label{eq:morph-dble}
  \bigoplus_{M\in \lc^G(T)}  i_{M,*}(X_*(M)_{\lc^M}) \stackrel{\displaystyle \longrightarrow}{\longrightarrow}  \bigoplus_{\stackrel{ M,L\in \lc^G(T)}{ M\subset L}}  i_{M,*}(X_*(L)_{\lc^M})
\end{equation}
  de la manière suivante : sur chaque composante $i_{M,*}(X_*(M)_{\lc^M})$, pour $M\in \lc^G$, la première flèche est obtenue par composition de la somme des morphismes (\ref{eq:morph1}) pour $L\in \lc^G(M)$ avec l'inclusion canonique
$$ \bigoplus_{L\in \lc^G(M)} i_{M,*}(X_*(L)_{\lc^M})\to  \bigoplus_{\stackrel{ M,L\in \lc^G(T)}{ M\subset L}}  i_{M,*}(X_*(L)_{\lc^M}) \ ; $$
le seconde est obtenue,  sur la composante $i_{L,*}(X_*(L)_{\lc^L})$, pour $L\in \lc^G$, par 
composition de la somme des morphismes (\ref{eq:morph2}) pour $M\in \lc^L(T)$ avec l'inclusion canonique
$$ \bigoplus_{M \in \lc^L(T)} i_{M,*}(X_*(L)_{\lc^M})\to  \bigoplus_{\stackrel{ M,L\in \lc^G(T)}{ M\subset L}}  i_{M,*}(X_*(L)_{\lc^M}).$$
On notera que, dans la double flèche (\ref{eq:morph-dble}), chaque composante  $i_{M,*}(X_*(L)_{\lc^M})$ du but (avec $M\subset L$) ne reçoit par chaque flèche qu'une composante de la source, à savoir $  i_{M,*}(X_*(M)_{\lc^M})$ par la première flèche et $i_{L,*}(X_*(L)_{\lc^L})$ par la seconde.

Soit $\mathcal{K}$ le noyau du morphisme double (\ref{eq:morph-dble}). On laisse au lecteur la vérification du lemme suivant.

\begin{lemme} \label{lem:concretement} Le faisceau  $\mathcal{K}$ est le faisceau constructible sur $\bg \lc^G \bd$ dont la fibre $\mathcal{K}_M$ en $M\in \lc^G$ est le $\ZZ$-module $X_*(M)$ et dont la flèche de spécialisation $\mathcal{K}_M \to \mathcal{K}_L$, pour $M$ et $L $ dans $\lc^G$ tels que $M\subset L$, est le morphisme $p^M_L : X_*(M) \to X_*(L)$ (cf. proposition \ref{prop:compatibilite}).  
\end{lemme}

\end{paragr}

\begin{paragr}[Le faisceau $\xc$ sur $\Ac_G$.] --- Soit
  \begin{equation}
    \label{eq:AG-LG}
    \Ac_G \to  \bg \lc^G \bd
  \end{equation}
  l'application de valeur $M$ sur chaque partie localement fermée $\Ac^{\el}_M$ pour $M\in \lc^G(T)$ (cf. proposition \ref{prop:reunionAM}).

  \begin{lemme}\label{lem:continue}
    L'application (\ref{eq:AG-LG}) est continue.
  \end{lemme}

  \begin{preuve}
    Soit $M\in \lc^G(T)$. L'image réciproque de la partie fermée $\lc^M(T)$ est la réunion des $\Ac_L^{\el}$ pour $L\in \lc^M(T)$ qui s'identifie au fermé $\Ac_M$ (cf. (\ref{eq:AM-AG}) et proposition  \ref{prop:reunionAM}).
  \end{preuve}

Soit $\xc$ le faisceau constructible sur $\Ac_G$ définie comme l'image réciproque par l'application continue (\ref{eq:AG-LG}) (cf. lemme \ref{lem:continue}) du faisceau $\mathcal{K}$ sur $  \bg \lc^G \bd$. Ce faisceau est donc constant de valeur $X_*(M)$ sur $\Ac^{\el}_M$, pour $M\in \lc^G(T)$.

\end{paragr}

\begin{paragr}[Le champ de Picard $\Jc$.] --- Soit
  \begin{equation}
    \label{eq:morphCA-car}
h_\Ac :    C\times_k \Ac \longrightarrow \car_{G,D}
  \end{equation}
  la flèche canonique qui à un triplet $(c,h_a,t)\in C\times_k \Ac$ associe $h_a(c)\in \car_{G,D}$. Soit 
$$J_{\Ac}=h^*_\Ac J^G_D$$ 
le schéma en groupes commutatifs, lisse sur $C\times_k\Ac$, obtenu comme l'image réciproque de $J^G_D$ par le morphisme $h_\Ac$ de (\ref{eq:morphCA-car}). Tout $a=(h_a,t)\in \Ac(k)$ définit par changement de base un morphisme $a_C : C\to C\times_k \Ac$ qui, composé avec $h_\Ac$, donne $h_a$. Le schéma en groupes sur $C$ défini par $a_C^* J_{\Ac}$ n'est autre que le schéma en groupes noté $J_a$ défini au §\ref{S:Jc_a}.

Soit $\mathrm{B}(J_{\Ac})$  le $C\times_k \Ac$-champ classifiant de $J_{\Ac}$. Soit $\Jc$ le champ sur $\Ac$ dont le groupoïde associé $\Jc(S)$ pour tout $\Ac$-schéma test $S$ est le groupoïde des morphismes
$$C\times_k S \to \mathrm{B}(J_{\Ac}).$$
Comme $J_\Ac$ est commutatif, le produit contracté confère à $\Jc$ une structure de champ de Picard au-dessus de $\Ac$.

\end{paragr}

\begin{paragr}[Le degré de $\Jc$.] --- Il est défini par la proposition suivante.
 
  \begin{proposition}\label{prop:degJ}
    Il existe un unique morphisme 
$$\deg : \Jc \to \xc$$
qui, pour tout $M\in \lc^G(T)$ et tout $a\in \Ac_M^{\el}$, induit le morphisme
$$\deg_M : \Jc_a \to \xc_a=X_*(M)$$
défini en (\ref{eq:degreJ_a}) du §\ref{S:morphdeg}.
  \end{proposition}

  \begin{preuve}
    L'unicité étant évidente, nous allons nous concentrer sur l'existence. Soit $M\in \lc^G(T)$. L'homomorphisme $\deg_M$ de (\ref{eq:degreJ_a}) induit un homomorphisme
$$(\chi_G^M)^* \Jc  \longrightarrow X_*(M)_{\Ac_M}$$
où $\chi_G^M$ est l'immersion fermée de $\Ac_M$ dans $\Ac_G$ définie en (\ref{eq:AM-AG}) et $X_*(M)_{\Ac_M}$ est le faisceau constant sur $\Ac_M$ de valeur $X_*(M)$. Par adjonction, on a donc un homomorphisme
$$\Jc \longrightarrow (\chi_G^M)_* (X_*(M)_{\Ac_M}).$$

En prenant la somme sur $M\in \lc^G(T)$, on obtient un homomorphisme
$$\Jc \longrightarrow \bigoplus_{M\in \lc^G(T)} (\chi_G^M)_* (X_*(M)_{\Ac_M}).$$
Nous allons montrer que cet homomorphisme est à valeurs dans le sous-faisceau $\xc$ et vérifie les propriétés annoncées. Pour voir qu'il est à valeurs dans $\xc$, il suffit de voir qu'il tombe dans le noyau de la double flèche $$\bigoplus_{M\in \lc^G(T)} (\chi_G^M)_* (X_*(M)_{\Ac_M})  \stackrel{\displaystyle \longrightarrow}{\longrightarrow}  \bigoplus_{\stackrel{ M,L\in \lc^G(T)}{ M\subset L}}  (\chi_G^M)_* (X_*(L)_{\Ac_M})$$
image réciproque de la double flèche (\ref{eq:morph-dble}) par l'application continue  $\Ac_G\to \bg \lc^G \bd$ définie en (\ref{eq:AG-LG}). Cela découle de la proposition \ref{prop:compatibilite}. Les autres propriétés résultent du lemme \ref{lem:concretement}.

  \end{preuve}

  \begin{paragr}[Action de $\Jc$ sur $\mc$.] --- \label{S:ActionJ} Il s'agit simplement de mettre en famille la construction du §\ref{S:ActionJa}. Soit $S$ un $k$-schéma test. Soit $m=(h_m,t)\in \mc(S)$ et $f(m)=(h_a,t)\in \Ac(S)$. Soit $J_{\Ac,S}=J_{\Ac}\times_{\Ac}S$ et soit $j\in \Jc(S)$ un $J_{\Ac,S}$-torseur. On a, cf. l. (\ref{eq:Ja-champetre}) du §\ref{S:ActionJa},
$$ J_{\Ac,S} = h_a^* J^G_D= h_m^* [\chi_{G,D}]^*J_D^G=h_m^*[(\chi_{G,D}^* J^G_D)/G].$$
Comme en (\ref{eq:morph-Aut(m)}), on en déduit un morphisme
$$J_{\Ac,S}\longrightarrow   h_m^*  [I^G_D/G]=\Aut(m).$$
L'action de $\Jc$ est alors donnée par le produit contracté
$$j\cdot m= j \times^{J_{\Ac,S}}m$$
où $J_{\Ac,S}$ agit sur $m$ via le morphisme vers $\Aut(m)$ que l'on vient d'introduire.
    
  \end{paragr}

  \begin{paragr}[Action de $\Jc^1$ sur $\mc^\xi$.] --- \label{S:actionJ1}Soit $\Jc^1$ le noyau de l'homomorphisme degré défini à la proposition \ref{prop:degJ}. C'est un sous-champ de Picard ouvert du champ $\Jc$ qui, fibre à fibre,  contient la composante neutre de $\Jc$. Comme $\Jc$ agit sur $\mc$, il en est de même de $\Jc^1$.

    \begin{proposition}
      Pour tout $\xi\in \ago_T$, l'action de $\Jc^1$ préserve l'ouvert $\mc^\xi$.
    \end{proposition}
    
    \begin{preuve}
      Soit $m\in \mc$ et $a=f(m)\in \Ac$. Soit $j\in \Jc_a$. D'après la proposition \ref{prop:translation}, on a $\mathcal{C}_m=\mathcal{C}_{ j\cdot m}+\deg_M(j)$. Si, de plus, $j\in \Jc_a^1$ alors on a $\deg_M(j)=0$ (cf. proposition \ref{prop:degJ}). Il s'ensuit que, pour tout  $j\in \Jc_a^1$, on a 
 $$\mathcal{C}_m=\mathcal{C}_{ j\cdot m}$$
et \emph{a fortiori} $j$ préserve la $\xi$-stabilité, qui est définie par la condition $\xi\in \mathcal{C}_m$.
    \end{preuve}

  \end{paragr}

\end{paragr}

\section{Courbe camérale et faisceau $\pi_0(\Jc)$.}\label{sec:cameral}

\begin{paragr}[Courbe camérale.] --- \label{S:cameral}Soit
$$\infty_\Ac : \Ac \to C\times_k \Ac$$
le morphisme déduit de $\infty$ par le changement de base  $\Ac\to \Spec(k)$. À la suite de Donagi et Gaitsgory (cf. \cite{DonagiGaitsgory}), on introduit la courbe camérale universelle, notée $Y$, qui est le revêtement de $C\times_k \Ac$  défini par le carré cartésien
$$\xymatrix{Y\ar[r]\ar[d]_{\pi_{Y}}&\tgo_{D}
\ar[d]^{\chi^T_{G,D}}\cr
C\times_k \Ac\ar[r]_{h_{\Ac}}&\car_{G,D}}
$$
où $h_\Ac$ est le morphisme canonique, cf. (\ref{eq:morphCA-car}), et $\chi^T_{G,D}$ est le morphisme $\chi^T_G$ du \ref{S:discr} «tordu par $D$». Comme le morphisme $\chi^T_{G,D}$ est $W^G$-invariant, on a une action naturelle de $W^G$ sur $Y$. Le revêtement $\pi_Y$ est fini, galoisien de groupe $W^G$ et  étale au-dessus d'un voisinage de l'image de $\infty_\Ac$. Pour tout $a=(h_a,t)\in \Ac$, le point $t$ vérifie $\chi_G^T(t)=h_a(\infty)$. On donc un morphisme canonique
$$\infty_Y : \Ac \to Y$$
au-dessus de $\infty_\Ac$. Soit $U\subset C\times_k \Ac$ l'image réciproque par $h_\Ac$ de l'ouvert $\car_{G,D}^{\reg}$ où le discriminant $D^G$ ne s'annule pas (cf. §\ref{S:discr}). Soit $V$ l'ouvert de $Y$ défini par $V=\pi_Y^{-1}(U)$. On notera que $\infty_\Ac$ est à valeurs dans $U$ et $\infty_Y$ est à valeurs dans $V$. Le revêtement $\pi_Y$ est étale précisément au-dessus de $U$.

Soit $a=(h_a,t)\in \Ac$ et $k(a)$ le corps résiduel de $a$. En tirant la construction précédente par le morphisme $C\times_k k(a)\to C\times_k\Ac$ déduit de $a$ par le changement de base $C\to \Spec(k)$, on obtient la courbe camérale $Y_a$ sur $k(a)$ définie par le diagramme cartésien
$$\xymatrix{Y_{a}\ar@{^{(}->}[r]\ar[d]_{\pi_{Y_a}}&\mathfrak{t}_{D}\times_k k(a)
\ar[d]^{\chi^T_{G,D}}    \\
C\times_k k(a)\ar@{^{(}->}[r]_{h_{a}}&\mathfrak{car}_{G,D}\times_k k(a)}.$$
Le groupe $W^G$ préserve $Y_a$ et le revêtement $\pi_{Y_a}$ est fini, galoisien de groupe $W^G$ et  étale au-dessus du point $\infty$. Dans la fibre de ce dernier, on a le point marqué $\infty_{Y_a}=\infty_\Ac\circ a$. Soit $U_a$ et $V_a$ les ouverts de $C\times_k k(a)$ et $Y_a$ déduits de $U$ et $V$ par changement de base.
\end{paragr}

\begin{paragr}[La fonction $W_a$.] --- \label{S:Wa} Dans tout ce paragraphe, la lettre $W$ désigne un sous-groupe \emph{quelconque} du groupe de Weyl $W^G$. Pour tout $a\in \Ac$, soit $W_a$ le sous-groupe de $W^G$ défini comme le normalisateur de la composante irréductible de $Y_a$ passant par $\infty_{Y_a}$. De manière équivalente, c'est le normalisateur de la composante connexe de $V_a$ qui contient $\infty_{Y_a}$.

Soit $\bg W^G \bd$ l'espace topologique sobre dont l'ensemble sous-jacent est formé des sous-groupes $W$ de $W^G$ et dont la topologie est la moins fine pour laquelle les parties formées des sous-groupes $W'$ inclus dans un sous-groupe $W$ sont fermées.

\begin{proposition}
  \label{prop:Wacontinue}
L'application $\Ac \to \bg W^G \bd$ définie par $a\mapsto W_a$ est continue.
\end{proposition}

\begin{preuve}
Montrons d'abord la constructibilité de cette application.   Le morphisme composé $Y\to C\times_k \Ac \to \Ac$ induit un morphisme surjectif $p : V \to \Ac$ dont $\infty_Y$ est une section. Pour tout $a\in \Ac$ soit $V_a^0$ la composante connexe de $V_a$ qui contient $\infty_Y(a)$. Soit $V^0$ la réunion des $V_a^0$ quand $a$ parcourt $\Ac$. D'après Grothendieck (cf. \cite{EGAIV3} proposition 15.6.4), on sait que $V^0$ est un ouvert. Pour tout sous-groupe $W$ de $W^G$, considérons l'ensemble constructible
$$\Ac_W= p( \bigcap_{w\in W} wV^0) - \bigcup_{w\in W^G-W} p( V^0\cap  wV^0).$$
Nous allons démontrer l'égalité
\begin{equation}
  \label{eq:AcW}
  \Ac_W=\{a\in \Ac \mid W_a=W \}
\end{equation}
qui implique la constructibilité de l'application $W_a$. Soit $a\in \Ac_W$. Il existe $y\in \bigcap_{w\in W} wV^0$ tel que $p(y)=a$. Par conséquent, la réunion $\bigcup_{w\in W} wV^0_a$ dont chaque constituant contient $y$  est connexe. Il s'ensuit qu'on a   $\bigcup_{w\in W} wV^0_a=V^0_a$ et donc $W\subset W_a$. Soit $w\in W_a$. On a donc $wV_a^0=V_a^0$ et $y\in V^0_a\cap  wV^0_a$. Nécessairement on a $w\in W$ d'où  $W_a\subset W$ et finalement $W_a=W$. Réciproquement soit $a\in \Ac$ tel que $W_a=W$. Soit $y\in V_a^0$. Comme $V_a^0$ est stable par $W_a=W$, on a $y\in \bigcap_{w\in W} wV^0$. Par ailleurs si $w\notin W=W_a$, on a  $V^0_a\cap  wV^0_a=\emptyset$. D'où $a\in \Ac_W$. On a donc démontré (\ref{eq:AcW}).

Montrons ensuite la continuité. Il suffit de montrer que pour tout sous-groupe $W$ de $W^G$ la réunion $\bigcup \Ac_{W'}$ prise sur les sous-groupes $W'$ de $W$ est fermée. Comme cette réunion est constructible, il suffit de montrer qu'elle est stable par spécialisation. Soit $a$ un point de $\Ac$ et $a'\in \Ac$ une générisation de $a$ dans $\Ac$. Nous allons montrer qu'on a  $W_{a}\subset W_{a'}$ ce qui implique la stabilité cherchée. Par définition de $W_{a'}$, on a 
$$V_{a'}^0 \subset V- \bigcup_{w\in W^G-W_{a'}} wV^0.$$
D'après le résultat déjà cité de Grothendieck, on sait que l'ensemble  $V- \bigcup_{w\in W^G-W_{a'}} wV^0$ est fermé dans $V$. On a donc
$$\overline{V_{a'}^0} \subset V- \bigcup_{w\in W^G-W_{a'}} wV^0.$$
où $\overline{V_{a'}^0}$ désigne l'adhérence de $V_a^0$ dans $V$. Toujours d'après Grothendieck (cf. \cite{EGAIV3} proposition 15.6.1), on a $V^0_{a}\subset \overline{V_{a'}^0}$ et donc
 $$V_{a}^0 \subset V- \bigcup_{w\in W^G-W_{a'}} wV^0$$
ce qui implique immédiatement qu'on a $W_a\subset W_{a'}$ comme voulu.
\end{preuve}
  
\end{paragr}

\begin{paragr}[Le morphisme de $X_*(T)_{\Ac}$ dans $\Jc$.] --- \label{S:X-J}Soit $S$ un schéma affine sur $k$ et $a=(h_a,t)\in \Ac(S)$. Soit $C_S=C\times_k S$ et $\mathcal{O}_{C_S,\infty}$ le complété de $\mathcal{O}_{C_S}$ le long de $\{\infty\}\times_k S \subset C_S$ que l'on peut identifier à $\mathcal{O}_{S}[[\varpi_{\infty}]]$ par le choix d'une uniformisante $\varpi_{\infty}$ sur $C$ au point $\infty$. Soit $F_{C_S,\infty}=\mathcal{O}_{C_S,\infty}[1/\varpi_{\infty}]$ (qui ne dépend pas du choix de $\varpi_{\infty}$). Soit $C^\infty= C-\{\infty\}$ et $C^\infty_S =C^\infty\times_k S$.

  \begin{proposition}\label{prop:Ja=T-isocan}
    Il existe un isomorphisme canonique 
\begin{equation}
  \label{eq:Ja=T-isocan}
   J_a \times_C  \Spec(\mathcal{O}_{C_S,\infty}) \simeq  T\times_k  \Spec(\mathcal{O}_{C_S,\infty}).
\end{equation}

  \end{proposition}

\begin{preuve}
    Dans un voisinage de $\infty$, le fibré $\car_{G,D}$ s'identifie canoniquement au fibré trivial de fibre  $\car_G$. Il s'ensuit que la section $h_a : C_S \to \car_{G,D} \times_k S$ se restreint en une section
$$h_{a,\infty} : \Spec(\mathcal{O}_{C_S,\infty}) \to \car^{\reg}_G\times_k \mathcal{O}_{C_S,\infty}.$$
Comme le morphisme  $\chi_G^T : \tgo^{\reg} \to \car^{\reg}$ est étale, il existe un unique morphisme $
h_{a_T,\infty}$ qui relève $t$ et qui rend le diagramme ci-dessous commutatif
$$\xymatrix{\Spec(\mathcal{O}_{C_S,\infty}) \ar[rr]^{h_{a_T,\infty}} \ar[drr]_{h_{a,\infty}} & &  \tgo^{\reg}\times_k \mathcal{O}_{C_S,\infty} \ar[d]^{\chi_T^G} \\ & & \car^{\reg}_G\times_k \mathcal{O}_{C_S,\infty}}.$$
L'isomorphisme (\ref{eq:morphNgo}) du \S \ref{S:centralisateurs-reg} induit l'isomorphisme suivant sur $\tgo^{\reg}$
$$\big((\chi_G^T)^*J^G\big)_{|\tgo^{\reg}}\simeq I^G_{\tgo^{\reg}}=T\times_k \tgo^{\reg}.$$
De cet isomorphisme et de l'égalité
$$ J_a \times_C  \Spec(\mathcal{O}_{C_S,\infty})=h_{a,\infty}^*J^G= h_{a_T,\infty}\* (\chi_G^T)^*J^G,$$
on déduit l'isomorphisme cherché.

\end{preuve}

D'après la descente formelle de Beauville-Laszlo (cf. \cite{BL}), pour tout $\la \in X_*(T)$, il existe un triplet $(j,\al,\beta)$, unique à un unique isomorphisme près, formé d'un $J_a$-torseur $j$ sur $C_S$ et d'isomorphismes $J_a$-équivariants 
$$\al : j_{|\Spec(\mathcal{O}_{C_S,\infty}) } \to J_a \times_C  \Spec(\mathcal{O}_{C_S,\infty})  $$
et
$$\beta : j_{|C^\infty_S} \to J_a \times_C   C_S^\infty$$
de sorte que l'isomorphisme sur  $\Spec(F_{C_S,\infty})$ induit par $\beta\circ \al^{-1}$ sur  
$$J_a \times_C  \Spec(F_{C_S,\infty}) \simeq T\times_k  \Spec(F_{C_S,\infty})$$
(isomorphisme déduit de  (\ref{eq:Ja=T-isocan})) soit donné par la translation à gauche par $\varpi_\infty^\la$.

On déduit de cette construction  un homomorphisme
$$X_*(T)_{\Ac} \to \Jc$$
où $X_*(T)_{\Ac}$ le faisceau constant sur $\Ac$ de valeurs $X_*(T)$. 

\end{paragr}

\begin{paragr}[Retour sur des résultats de Ngô.] --- \label{S:retour} Soit $K$ une extension de $k$ séparablement close. Soit $a=(h_a,t)\in \Ac(K)$ un point «géométrique». Dans tout ce paragraphe, $\tc$ désigne un schéma en tores sur l'ouvert $U_a$ de $C_K=C\times_k K$ (cf. \ref{S:cameral}). Un tel schéma en tores se prolonge en un schéma en groupes lisse sur $C_K$ \emph{à fibres connexes} et tout tel prolongement s'envoie dans l'unique prolongement maximal noté $\tilde{\tc}$ et appelé le modèle de Néron connexe de $\tc$.  On introduit les notations suivantes : soit $\tc_{\infty_K}$ la fibre au point $\infty_K=\infty\times_k K$ de $\tc$. Soit $X_*(\tc_{\infty_K})=\Hom_{K\textrm{-gp}}(\mathbb{G}_{m,K},\tc_{\infty_K})$. Le groupe fondamental de $U_a$ pointé par $\infty_K$ est noté  $\pi_1(U_a,\infty_K)$ et agit sur $X_*(\tc_{\infty_K})$. On pose
$$A(\tc)=X_*(\tc_{\infty_K})_{\pi_1(U_a,\infty_K)}$$
 où l'indice $\pi_1(U_a,\infty_K)$ désigne le groupe des co-invariants sous $\pi_1(U_a,\infty_K)$. Soit $B(\tc)$ le groupe des composantes connexes du champ de Picard des $\tilde{\tc}$-torseurs. Soit $A$ et $B$ les foncteurs de la catégorie des schémas en tores sur $U_a$ vers la catégorie des groupes abéliens définis par $\tc \mapsto A(\tc)$ et  $\tc \mapsto B(\tc)$. À l'aide d'une  variante d'un résultat de Kottwitz (cf. \cite{Kot-isocrystal} lemme 2.2), Ngô a montré que les foncteurs $A$ et $B$ sont isomorphes. Il est important pour nous de préciser un tel isomorphisme. Pour cela, on utilise la variante suivante «à la Ngô» du lemme de Kottwitz déjà cité, où l'on note $X_*$ le foncteur $\tc \mapsto X_*(\tc_{\infty_K})$.

 \begin{lemme}\label{lem:kottwitz-ngo}
 Les morphismes canoniques
$$\Hom(A,B) \to \Hom(X_*,B)$$
et 
$$ \Hom(X_*,B) \to \Hom(X_*(\mathbb{G}_{m,U_a}),B(\mathbb{G}_{m,U_a}))= B(\mathbb{G}_{m,U_a})$$
sont des isomorphismes. Le groupe  $B(\mathbb{G}_{m,U_a})$ est libre de rang $1$ et tout élément de $\Hom(A,B)$ dont l'image par le composé des deux isomorphismes précédents est un générateur  de $B(\mathbb{G}_{m,U_a})$ est un isomorphisme.
 \end{lemme}

 \begin{preuve}
   Les propriétés fondamentales du foncteur $B$ sont démontrées par Ngô (cf. \cite{Ngo1}). Le reste n'est qu'une transposition de résultats de Kottwitz (cf. \cite{Kot-isocrystal} pp.206-207).
 \end{preuve}

Définissons alors un morphisme de foncteurs de $X_*$ vers $B$. Pour tout schéma en tores $\tc$ sur $U_a$, on définit ainsi un morphisme de $X_*(\tc_{\infty_K})$ dans le champ des $\tilde{\tc}$-torseurs. Les notations étant celles du §\ref{S:X-J}, pour tout $\la \in X_*(\tc_{\infty_K})$, on recolle le $\tilde{\tc}$-torseur  trivial sur $C^\infty_K$  avec le $\tc$-torseur  trivial sur $\Spec(\oc_{C_K,\infty})$ à la Beauville-Laszlo avec la donnée de recollement $\varpi_\infty^\la$ (le tore $\tc$ est en fait constant sur $\oc_{C_K,\infty}$, isomorphe à $X_*(\tc_{\infty_K})\otimes_\ZZ \mathbb{G}_{m,\oc_{C_K,\infty}}$). On en déduit un morphisme 
$$X_* \to B$$
en composant le morphisme précédent avec le morphisme canonique du champ des $\tilde{\tc}$-torseurs dans son groupe de composantes connexes.

\begin{lemme}\label{lem:factorisation-foncteur}
  Le morphisme de foncteur $X_* \to B$ se factorise canoniquement selon le diagramme
$$\xymatrix{X_* \ar[rd] \ar[r] & A \ar[d] \\ &  B }$$
où le morphisme horizontal est le morphisme canonique et le morphisme vertical est un isomorphisme.
\end{lemme}

\begin{preuve}
D'après le lemme \ref{lem:kottwitz-ngo}, il s'agit de voir que le foncteur $X_* \to B$ envoie le générateur canonique de  $X_*(\tc_{\infty_K})$ (donné par l'homomorphisme identique) sur un générateur de 
$B(\mathbb{G}_{m,U_a})$. Or lorsque $\tc= \mathbb{G}_{m,U_a}$, on a $\tilde{\tc}=\mathbb{G}_{m,C_K}$ et le morphisme degré du champ de Picard  des $\tilde{\tc}$-torseurs, qui n'est autre que le champ de Picard des $\oc_{C_K}$-modules localement libres de rang $1$, identifie son groupe des composantes connexes $\ZZ$. Le générateur canonique de  $X_*(\tc_{\infty_K})$ définit par le recollement expliqué ci-dessus, un $\mathbb{G}_{m,C_K}$-torseur de degré $1$ d'où le lemme.
\end{preuve}

\end{paragr}

\begin{paragr}[Factorisation de Ngô.] --- Soit  $\mathcal{Y}$ le faisceau sur $\mathcal{A}$ image réciproque par l'application continue $\mathcal{A}\rightarrow \bg W^G \bd$ du faisceau dont la fibre en $W\subset W^G$ est $X_{\ast}(T)_{W}$, avec les morphismes de transition naturels. Ce faisceau est par définition un quotient du faisceau constant $X_{\ast}(T)_{\Ac}$. On a défini au §\ref{S:X-J} un morphisme $X_{\ast}(T)_{\Ac} \to \Jc$. D'après Ngô, il existe un faisceau constructible sur $\Ac$  qu'on  note $\pi_0(\Jc)$ et dont la fibre en un point géométrique $a$ de $\Ac$ est le groupe $\pi_0(\Jc_a)$ des composantes connexes de $\Jc_a$ (cf. proposition 6.2 de \cite{Ngo1}). Par composition avec le morphisme canonique $\Jc\to \pi_0(\Jc)$, on a obtient un morphisme  $X_{\ast}(T)_{\Ac}\to \pi_0(\Jc)$. Ngô a montré la factorisation suivante (cf. proposition 5.5.1 de \cite{ngo2}). 

  \begin{proposition} \label{prop:factorisation-faisceau}On a la factorisation suivante
$$\xymatrix{   X_{\ast}(T)_{\Ac}  \ar[r] \ar[rd] & \yc \ar[d]\\ & \pi_0(\Jc)}$$
où la flèche horizontale est la flèche canonique et les deux autres flèches sont des épimorphismes.
    
  \end{proposition}
  
  \begin{preuve}
    Pour la commodité du lecteur, on rappelle quelques points de la démonstration de Ngô. Il suffit de prouver la factorisation fibre à fibre. Soit un point géométrique $a=(h_a,t)\in \Ac(K)$, où $K$ est une extension de $k$ séparablement close. Soit $\Jc^0_a$ le champ des $J^0$-torseurs  sur $C\times_k K$ où $J^0_a$ est le schéma en groupes connexes qui, fibre à fibre, est la composante connexe de la section neutre de $J_a$. Comme la fibre de $J_a$ en $\infty_K$ est connexe, le morphisme $X_*(T)\to \Jc_a$ se factorise par un morphisme évident    $X_*(T)\to \Jc^0_a$. D'après Ngô  (cf. \cite{Ngo1} propositions 6.3), le morphisme canonique $\Jc^0_a \to \Jc_a$ induit un morphisme \emph{surjectif} $\pi_0(\Jc^0_a) \to \pi_0(\Jc_a)$. Par ailleurs, le morphisme canonique $J^0_a \to \tilde{J}_a$ de $J^0_a$ dans le modèle de Néron connexe  $\tilde{J}_a$ du schéma en tores $J_{|U_a}$ induit un morphisme de $\Jc^0_a$ dans le champ $\tilde{\Jc}_a$ des $\tilde{J}_a$-torseurs sur $C\times_k K$ et, d'après Ngô (cf. \cite{Ngo1} propositions 6.4), un isomorphisme $\pi_0(\tilde{\Jc}_a)=\pi_0(\Jc^0_a)$.

Il suffit donc de prouver la factorisation cherchée pour le morphisme $X_*(T)\to \pi_0(\tilde{\Jc}_a)$, c'est-à-dire la factorisation par l'épimorphisme $X_*(T)_{W_a}\to \pi_0(\tilde{\Jc}_a)$ (ce sera même un isomorphisme).

Avec les notations du § \ref{S:cameral}, on a un diagramme commutatif

$$\xymatrix{\Spec(K) \ar[rr]^{\infty_{Y_a}} \ar[rrd]_{\infty} &  & V_{a}\ar@{^{(}->}[rr]^{i_a}    \ar[d]^{\pi_{Y_a}}&& \mathfrak{t}^{\reg}_{D,K}
\ar[d]^{\chi^T_{G,D}}    \\
& &   U_a\ar@{^{(}->}[rr]_{h_{a}}& &\car^{\reg}_{G,D,K}}$$
où le carré est cartésien et où l'on a noté par un indice $K$ le changement de base de $k$ à $K$. 
 
Soit $(Z,z)$ un revêtement connexe, étale et galoisien de $U_a$ pointé par un point $z : \Spec(K)\to Z$ au-dessus de $\infty$. On suppose, ainsi qu'il est loisible, que pour chaque composante connexe de $V_a$, il existe un $U_a$-morphisme surjectif de $Z$ sur celle-ci. La composition à droite avec $z$ induit une bijection de $\Hom_{U_a}(Z,V_a)$ sur le fibre de $\pi_{Y_a}$ au-dessus  de $\infty$. On note $\al : Z\to V_a$ le morphisme correspondant à $\infty_{Y_a}$. Avec les notations de §\ref{S:Jc_a}, on a 
$$J_a=J_D^G\times_{\car^{\reg}_{G,D,K}} U_a$$
Or d'après l'isomorphisme (\ref{eq:morphNgo}) du §\ref{S:centralisateurs-reg}, on a l'isomorphisme suivant 
$$J_D^G \times_{\car^{\reg}_{G,D,K}}\mathfrak{t}^{\reg}_{D,K}\simeq  T \times_k \mathfrak{t}^{\reg}_{D,K}$$
d'où 
$$J_a\times_{U_a}V_a=J_D^G\times_{\car^{\reg}_{G,D,K}} V_a \simeq  T \times_k V_a.$$
En utilisant le morphisme $\al$, on obtient
\begin{equation}
  \label{eq:triv-de-T}
  J_a\times_{U_a} Z = (J_a\times_{U_a} V_a)\times_{V_a,\al} Z \simeq T \times_k Z.
\end{equation}
Pour tout $\Phi\in \Aut_{U_a}(Z)$, il existe un unique élément $w_{\phi}\in W^G$ tel $\al\circ \Phi=w_{\Phi}\cdot\al$ (où $\cdot$ désigne l'action de $W^G$ sur $V_a$). Soit $\pi_1(U_a,\infty)$ le groupe fondamental de $U_a$ pointé par $\infty$. L'application $\Phi \mapsto w_{\Phi}$ définit un morphisme de $\pi_1(U_a,\infty)$ dans $W^G$, dont l'image est précisément le sous-groupe $W_a$. L'action évidente de $\pi_1(U_a,\infty)$ sur  $J_a\times_{U_a} Z$ correspond, par l'isomorphisme (\ref{eq:triv-de-T}), à l'action diagonale de $\pi_1(U_a,\infty)$  sur $T \times_k Z$, où l'action sur le premier facteur est donné par le morphisme dans $W^G$. En particulier, on a l'identification 
$$\Hom_{Z\textrm{-gpe}}( \mathbb{G}_{m,Z},  J_a\times_{U_a} Z)=X_*(T)$$
qui est équivariante sous l'action de $\pi_1(U_a,\infty)$ lorsque ce groupe agit sur  $X_*(T)$ via son morphisme dans $W_a$. On en déduit les identifications
$$X_*(J_a)_{\infty}=X_*(T)$$
et
\begin{equation}
  \label{eq:identification}
(X_*(J_a)_{\infty})_{\pi_1(U_a,\infty) }=X_*(T)_{W_a}.
\end{equation}
Notons que cette identification vaut aussi pour $\tilde{J}_a$ puisque $J_a\simeq \tilde{J}_a$ sur $U_a$.

Ainsi, le morphisme  $X_*(T)\to \pi_0(\tilde{\Jc}_{\bar{a}})$ s'identifie au morphisme $X_*(\tilde{J_a})_{\infty} \to \pi_0(\tilde{\Jc}_{a})$ considéré au \S\ref{S:retour}. Le lemme \ref{lem:factorisation-foncteur} et la ligne (\ref{eq:identification}) impliquent que ce morphisme se factorise par $X_*(T)_{W_a}$ et que le morphisme $X_*(T)_{W_a}\to  \pi_0(\tilde{\Jc}_{a})$ est un isomorphisme. Cela conclut la démonstration.

  \end{preuve}

\end{paragr}

\section{La $s$-décomposition}\label{sec:sdecomposition}

\begin{paragr}
Dans \cite{Ngo1} section 8, Ngô a introduit une décomposition des faisceaux de cohomologie perverse de la fibration de Hitchin sur la partie elliptique. Le but de cette section est de généraliser à des ouverts non nécessairement elliptiques la décomposition de Ngô.
\end{paragr}

\begin{paragr} La proposition suivante sera utile pour la suite.

  \begin{proposition}\label{prop:AM-Wa}
    Soit $a\in \Ac_G$ et $M\in \lc^G$ un sous-groupe de Levi. Les deux assertions suivantes sont équivalentes :
    \begin{enumerate}
    \item $a\in \Ac_M$ ;
    \item $W_a\subset W^M$.
    \end{enumerate}
  \end{proposition}

  \begin{preuve} Soit $a=(h_a,t)\in \Ac_G$ et $k(a)$ le corps résiduel de $a$. Montrons que l'assertion 2 implique la première. Avec les notations du § \ref{S:cameral}, on a un diagramme commutatif
    \begin{equation}
      \label{eq:AM-Wa}
      \xymatrix{\Spec(k(a)) \ar[rr]^{\infty_{Y_a}} \ar[rrd]_{\infty} &  & V_{a}\ar@{^{(}->}[rr]^{i_a}    \ar[d]^{\pi_{Y_a}}&& \mathfrak{t}^{\reg}_{D,k(a)}
\ar[d]^{\chi^T_{G,D}}    \\
& &   U_a\ar@{^{(}->}[rr]_{h_{a}}& &\car^{\reg}_{G,D,k(a)}}
\end{equation}
où le carré est cartésien. Soit $V_a^0$ la composante  connexe de $V_a$ qui contient $\infty_{Y_a}$. Par définition de $W_a$, la restriction de $\pi_{Y_a}$ à $V_a^0$ fait de $V_a^0$ un revêtement connexe, étale et galoisien de groupe $W_a$. Comme l'inclusion $i_a$ est $W$-équivariante, on obtient par quotient par $W_a$ un morphisme
$$U_a \to \mathfrak{t}^{\reg}_{D,k(a)}//W_a$$ au-dessus de $\car^{\reg}_{G,D,k(a)}$. Comme $W_a \subset W^M$, on peut composer ce morphisme avec le morphisme canonique $\mathfrak{t}^{\reg}_{D,k(a)}//W_a \to \car^{\reg}_{M,D,k(a)}$ de sorte qu'on obtient un morphisme $h_{a_M}$ qui s'inscrit dans le diagramme commutatif

$$\xymatrix{ U_{a}\ar@{^{(}->}[rr]^{h_{a_M}}    \ar[d] & & \car_{M,D,k(a)} \ar[d]^{\chi^M_{G,D}}    \\
C\times_k k(a) \ar@{^{(}->}[rr]_{h_{a}} & & \car_{G,D,k(a)} }$$
    Par propreté du morphisme $\chi^M_{G,D}$, le morphisme $h_{a_M}$ s'étend à $C\times_k k(a)$. On vérifie que $(h_{a_M},t)$ appartient à $\Ac_M$ et a pour image $a$ dans $\Ac_G$.

\medskip
La réciproque résulte du lemme suivant.
\end{preuve}

\begin{lemme}\label{lem:WaM}
 Soit $a_M \in \Ac_M$ et $a\in \Ac_G$ son image par l'immersion fermée $\Ac_M \to \Ac_G$ de (\ref{eq:AM-AG}). Soit $W_{a_M}$ le sous-groupe de $W^M$ défini au §\ref{S:Wa} relativement au groupe $M$. On a 
$$W_{a_M}=W_a$$ 
\end{lemme}

\begin{preuve}Soit $a_M=(h_{a_M},t)\in \Ac_M$ et $a=(h_a,t)$ son image dans $\Ac_G$. On adopte les notations de la démonstration de la proposition \ref{prop:AM-Wa} mais relativement au groupe $M$. On a donc un diagramme analogue à (\ref{eq:AM-Wa})

\begin{equation*}
            \xymatrix{\Spec(k(a_M)) \ar[rr]^{\infty_{Y_{a_M}}} \ar[rrd]_{\infty} &  & V_{a_M}\ar@{^{(}->}[rr]^{}    \ar[d]^{}&& \mathfrak{t}^{\reg}_{D,k(a_M)}
\ar[d]^{\chi^T_{M,D}}    \\
& &   U_{a_M}\ar@{^{(}->}[rr]_{h_{a_M}}& &\car^{\reg}_{G,D,k(a_M)}}
\end{equation*}
où l'exposant $\reg$ signifie $G$-régulier (conformément aux conventions du §\ref{S:discr}). On a d'ailleurs $U_a=U_{a_M}$ et une immersion à la fois ouverte et fermée $V_{a_M} \hookrightarrow V_a$, qui est $W^M$-équivariante et qui envoie $\infty_{Y_{a_M}}$ sur   $\infty_{Y_{a}}$. En particulier, la composante connexe de  $\infty_{Y_{a}}$ dans $V_a$ est égale à l'image de celle de $\infty_{Y_{a_M}}$ dans $V_{a_M}$. Si $w\in W^G$ vérifie $w V_{a_M}\cap V_{a_M}\not=\emptyset$, on a nécessairement $w\in W^M$. Il s'ensuit qu'on a $W_a\subset W^M$. Le résultat est alors clair.

\end{preuve}

\end{paragr}

\begin{paragr}\label{S:X1}
 On a défini respectivement aux propositions \ref{prop:degJ} et \ref{prop:factorisation-faisceau} un morphisme $\deg : \Jc \to \xc$, qui se factorise évidemment en un morphisme 
$$\deg : \pi_0(\Jc) \to \xc,$$
et  un morphisme $\yc \to \pi_0(\Jc)$. En composant ces deux morphismes, on obtient $\yc \to\xc$ et en composant encore avec le morphisme  $X_*(T)_{\Ac}\to \yc$, on obtient  un morphisme $X_*(T)_{\Ac}\to\xc$.  Soit $\yc^1$ et $X_*(T)_{\Ac}^1$ les faisceaux constructibles sur $\Ac$ définis respectivement  par
$$\yc^1=\Ker (\yc \to \xc)$$
et 
$$X_*(T)_{\Ac}^1=\Ker (X_*(T)_{\Ac}\to \xc).$$
Soit $\pi_0(\Jc^1)$ le faisceau des composantes connexes de $\Jc^1$ : c'est un sous-faisceau de $\pi_0(\Jc^1)$ dont la fibre en un point géométrique $a$ est le sous-groupe de $\pi_0(\Jc_a)$ formé des composantes de degré $0$ (où le degré est le morphisme ci-dessus). Le morphisme $X_*(T)_{\Ac} \to \Jc$ se restreint alors en un morphisme  $X_*(T)_{\Ac}^1 \to \Jc^1$ qui, par la proposition \ref{prop:factorisation-faisceau}, s'inscrit dans le  diagramme commutatif
\begin{equation}
  \label{eq:diag-X1}
  \xymatrix{ X_*(T)_{\Ac}^1 \ar[r] \ar[rd] & \yc^1  \ar[d]  \\ & \pi_0(\Jc^1)}
\end{equation}
où les deux flèches de but  $\pi_0(\Jc^1)$ sont des épimorphismes.

\begin{proposition} \label{prop:description}  Soit $a\in \Ac_M^{\el}$ avec $M\in \lc^G$ et $a'\in \Ac_G$  une spécialisation de $a$. Soit $L\in \lc^G$ défini par la condition $a'\in \Ac_L^{\el}$. 

  \begin{enumerate}
  \item Les  faisceaux $X_*(T)_{\Ac}^1$ et $\yc^1$ ont respectivement pour fibre en $a$ les groupes $X_*(T\cap M_{\der})$ et $X_*(T\cap M_{\der})_{W_a}$ ;
  \item on a $L\subset M$ et $W_{a'}\subset W_a$ ;
  \item les flèches de  spécialisation $  (X_*(T)_{\Ac}^1)_{a'} \to (X_*(T)_{\Ac}^1)_{a}$ et 
$\yc^1_{a'} \to \yc^1_{a}$ sont respectivement données par les morphismes canoniques
$$X_*(T\cap L_{\der})\to X_*(T\cap M_{\der})$$
et
$$X_*(T\cap L_{\der})_{W_{a'}}\to X_*(T\cap M_{\der})_{W_{a}}\ ;$$
\item les fibres du faisceau   $\yc^1$ sont finies ;
\item les fibres du faisceau   $\pi_0(\Jc^1)$ sont finies ;
\end{enumerate}
\end{proposition}

\begin{preuve}
  Montrons l'assertion 1. D'après le lemme \ref{lem:explicite} qui suit, la fibre en $a\in \Ac_M^{\el}$ de $X_*(T)_{\Ac}^1$ est le noyau du morphisme canonique $X_*(T)\to X_*(M)$ : c'est donc bien $X_*(T\cap M_{\der})$. On sait par la proposition \ref{prop:factorisation-faisceau} que cette flèche se factorise par $X_*(T)_{W_a}$ qui est la fibre en $a$ de $\yc$. La fibre en $a$ est donc le noyau de la flèche  $X_*(T)_{W_a}\to X_*(M)$ (qui est bien définie puisque $W_a\subset W^M$ d'après la proposition \ref{prop:AM-Wa} : c'est donc bien  $X_*(T\cap M_{\der})_{W_a}$.

L'assertion 2 résulte de la continuité des applications $\Ac_G \to \bg \lc^G \bd$ et $\Ac_G \to \bg W^G \bd$ (cf. lemme \ref{lem:continue} et proposition \ref{prop:Wacontinue}.

On laisse l'assertion 3 au lecteur. L'assertion 4 résulte du lemme \ref{lem:finitude} qui suit. L'assertion 5 est une conséquence de l'assertion 4 vu qu'on dispose d'un épimorphisme $\yc^1 \to \pi_0(\Jc_1)$.
\end{preuve}

\begin{lemme} \label{lem:explicite}Soit $M\in \lc^G$.  Le morphisme $X_*(T)_{\Ac} \to \xc$ a pour fibre en $a\in \Ac_M^{\el}$ le morphisme canonique 
$$X_*(T)\to X_*(M).$$
\end{lemme}

\begin{preuve}
Soit $a\in \Ac_M^{\el}$ et  $K$ son corps résiduel. Soit $h_{a_M}$ le morphisme de $C_K$ dans $\car_{M,D}$ qui se déduit de $a$. On reprend les notations du §\ref{S:X-J}. Soit $\la\in X_*(T)$ et $j_\la$ le $J_a$-torseur sur $C_{K}=C\times_k K$ obtenu à l'aide du recollement par $\varpi_\infty^\la\in T(\oc_{C_{K},\infty})\simeq J_a(\oc_{C_{K},\infty})$ (cf. l'isomorphisme de la proposition \ref{prop:Ja=T-isocan}) qui définit le morphisme $X_*(T) \to \Jc_a$. Le degré du $J_a$-torseur $j_\la$ est l'application qui, à $\mu\in X^*(M)$, associe le degré du $\mathbb{G}_{m,C_{K}}$-torseur $\mu(j_\la)$ obtenu lorsqu'on pousse $j_\la$ par le morphisme composé
\begin{equation}
  \label{eq:morph_compose}
  J_a \longrightarrow X_*(M)\otimes_{\ZZ} \mathbb{G}_{m,C_{K}} \stackrel{\mu }{\longrightarrow} \mathbb{G}_{m,C_{K}},
\end{equation}
où l'on obtient la première flèche est obtenue en tirant en arrière par $h_{a_M}$ le morphisme 
$$(\chi^M_{G,D})^* J_D^G \to X_*(M)\otimes_{\ZZ}   \mathbb{G}_{m,\car_{M,D}}$$
défini en (\ref{eq:morph_XMcarD}) et déduit de la proposition \ref{prop:JGXM}. Le degré de $j_\la$ est donc la valuation de l'image de $\varpi_\infty^\la$ par le morphisme composé (\ref{eq:morph_compose}).

En fait, sur $\Spec(\oc_{C_{K},\infty})$, la section $h_{a_M}$ est à valeurs régulières et se factorise par une section canonique notée $h_{a_T,\infty}$ de $\tgo^{\reg}_D$ (cf. la démonstration de la proposition \ref{prop:Ja=T-isocan}). L'isomorphisme  $J_a \times_C  \Spec(\mathcal{O}_{C_S,\infty}) \simeq  T\times_k  \Spec(\mathcal{O}_{C_S,\infty})$ s'obtient alors par image réciproque par $h_{a_T,\infty}$ de l'isomorphisme $\big((\chi_G^T)^*J^G\big)_{|\tgo^{\reg}}\simeq X_*(T)\otimes_\ZZ \mathbb{G}_{m,\tgo^{\reg}}$. Comme on peut le voir sur les définitions, le morphisme
$$X_*(T)\otimes_\ZZ \mathbb{G}_{m,\tgo^{\reg}}=\big((\chi_G^T)^*  J^G\big)_{\tgo^{\reg}}  \to X_*(M)\otimes_{\ZZ}   \mathbb{G}_{m,\tgo^{\reg}}$$
qu'on obtient en tirant en arrière par $\chi_M^T$ le morphisme de la proposition \ref{prop:JGXM} se déduit du morphisme canonique $X_*(T)\to X_*(M)$. Par conséquent, on a 
$$\deg(\mu(j_\la))=\mu(\la)$$
(l'accouplement est celui entre caractères et cocaractères de $T$) et le degré de $\mu(j_\la)$ est bien l'image de $\la$ par la projection $X_*(T)\to X_*(M)$.

\end{preuve}

\begin{lemme}
  \label{lem:finitude}
Pour tout $M\in \lc^G$ et tout $a\in \Ac^{\el}_M$
le groupe $X_*(T\cap M_{\der})_{W_a}$ est fini.
\end{lemme}

\begin{preuve}
  Nous allons d'abord prouver que le groupe $X_*(T\cap M_{\der})^{W_a}$ des invariants sous $W_a$ est trivial. Soit $\la$ un élément de ce groupe qu'on voit comme un morphisme de $\Gmk$ dans $M$. Soit $L$ le sous-groupe de Levi dans $\lc^M$ qui est le centralisateur du tore image. On a donc $W_a\subset W^L$ et par la proposition \ref{prop:AM-Wa} on a $a\in \Ac_L$. Comme $a$ est elliptique dans $M$, il vient $M=L$. Donc $\la$ est central dans $M_{\der}$ ce qui n'est possible que si $\la=0$ ce qu'il fallait démontrer. 

Comme  $X_*(T\cap M_{\der})$ est de type fini, la finitude de l'énoncé est équivalente à la nullité de $X_*(T\cap M_{\der})_{W_a}\otimes_\ZZ \QQ$. Or la moyenne sur $W_a$ induit un isomorphisme de ce $\QQ$-espace vectoriel sur   $X_*(T\cap M_{\der})^{W_a} \otimes_\ZZ \QQ$ dont on a prouvé la nullité. Cela conclut la démonstration.
\end{preuve}

\end{paragr}

\begin{paragr}[L'ouvert $\uc_M$.] --- \label{S:OmegaM} Soit $M\in \lc^G$ et $\uc_M\subset \Ac_G$ l'ouvert défini par
$$\uc_M=\bigcup_{L\in \lc^G(M)} \Ac^{\el}_L.$$
C'est encore le complémentaire du fermé
$$\bigcup_{L\in \lc^G- \lc^G(M)} \Ac_L.$$
L'ouvert $\uc_G$ n'est autre que l'ouvert elliptique $\Ac_G^{\el}$.
On définit la partie $M$-parabolique de $\mc^\xi$ notée $\mc^{\xi, M-\mathrm{par}}$ par le changement de base
$$\mc^{\xi,M\text{-}\mathrm{par}}=\mc^\xi\times_{\Ac} \uc_M.$$
Soit 
$$f^{\xi,M\text{-}\mathrm{par}} : \mc^{\xi,M\text{-}\mathrm{par}} \to \uc_M$$
 la restriction du morphisme de Hitchin $f^\xi$ à $\mc^{\xi,M\text{-}\mathrm{par}}$. Lorsque $M=G$, on remplace «$G\text{-}\mathrm{par}$» par «$\el$». 

\medskip

Soit  $\Gamma(\uc_M, \cdot)$ le foncteur des sections sur $\uc_M$.

\begin{lemme} \label{lem:quotient} La restriction de $\yc^1$ à $\uc_M$ est un quotient du faisceau constant de valeur $X_*(T\cap M_{\der})$. En particulier, on dispose d'homomorphismes de groupes
$$X_*(T\cap M_{\der}) \longrightarrow \Gamma(\uc_M,\yc^1)$$
et
\begin{equation}
  \label{eq:J1=quotient}
  X_*(T\cap M_{\der}) \longrightarrow \Gamma(\uc_M,\pi_0(\Jc^1)).
\end{equation}
\end{lemme}

\begin{preuve}
  La première assertion résulte de la description de $\yc^1$ donnée à la proposition \ref{prop:description}. Le reste s'en déduit.
\end{preuve}
\end{paragr}

\begin{paragr}[la $s$-décomposition.] --- \label{S:sdecomposition}Avec les notations de §\ref{S:dual-Langlands}, on note $\Tc$ le tore complexe de réseau de caractère $X^*(\Tc)=X_*(T)$. Dans toute la suite, on confond toujours dans les notations les groupes algébriques sur $\CC$ avec leurs groupes de points à valeurs dans  $\CC$. Soit $Z_{\Mc}^0$ le centre connexe du dual de Langlands $\Mc$ de $M$. C'est encore la composante neutre du sous-groupe de $\Tc$ défini comme l'intersection des noyaux des caractères $\al^\vee$ où $\al^\vee$ parcourt les coracines de $T$ dans $M$. Notons qu'on a 
$$\Tc/Z_{\Mc}^0=\Hom_{\ZZ}(X_*(T\cap M_{\der}),\CC^\times).$$
Soit $(\Tc/Z_{\Mc}^0)_{\tors}$ le sous-groupe des points de  torsion de $\Tc/Z_{\Mc}^0$. Soit $\overline{\QQ}$ la clôture algébrique de $\QQ$ dans $\CC$. On fixe un plongement de $\overline{\QQ}$ dans $\Qlb$. On a alors les inclusions naturelles
$$(\Tc/Z_{\Mc}^0)_{\tors}\subset \Hom_\ZZ(X_*(T\cap M_{\der}),\Qb^\times)\subset \Hom_\ZZ(X_*(T\cap M_{\der}),\Qlb^\times).$$

 L'action du champ de Picard $\Jc^1$ sur $\mc^{\xi}$ au-dessus de $\Ac$ induit une action de $\Jc^1$ sur  $\mc^{\xi, M-\mathrm{par}}$ au-dessus de $\uc_M$ donc une action sur chaque faisceau de cohomologie perverse ${}^{\mathrm{p}}\mathcal{H}^{i}( Rf^{\xi,M\text{-}\mathrm{par}}_* \Qlb)$ de $Rf^{\xi,M\text{-}\mathrm{par}}_* \Qlb$. Par le lemme d'homotopie (cf. \cite{laumon-Ngo} lemme 3.2.3), cette action se factorise en une action du faisceau $\pi_0(\Jc^1)$ sur ${}^{\mathrm{p}}\mathcal{H}^{i}( Rf^{\xi,M\text{-}\mathrm{par}}_* \Qlb)$.

 \begin{theoreme} \label{thm:sdecomposition}Pour tous $s\in (\Tc/Z_{\Mc}^0)_{\tors}$ et $i\in \ZZ$, soit
$${}^{\mathrm{p}}\mathcal{H}^{i}( Rf^{\xi,M\text{-}\mathrm{par}}_* \Qlb)_{s}$$
la somme des composantes $s'$-isotypiques du groupe commutatif fini  $\Gamma(\uc_M,\pi_0(\Jc^1))$ agissant sur le faisceau de  cohomologie perverse ${}^{\mathrm{p}}\mathcal{H}^{i}( Rf^{\xi,M\text{-}\mathrm{par}}_* \Qlb)$ pour les  caractères $s'\in \Hom_{\ZZ}(\Gamma(\uc_M,\pi_0(\Jc^1)),\Qlb^\times)$ qui induisent le  caractère $s$ sur $X_*(T\cap M_{\der})$ via le morphisme (\ref{eq:J1=quotient}). 

On a alors les assertions suivantes :
\begin{enumerate}
\item si $L$ est un élément de $\lc^G(M)$ et si $a$ est un point géométrique  de $\Ac_L^{\el}$ tel que la fibre ${}^{\mathrm{p}}\mathcal{H}^{i}( Rf^{\xi,M\text{-}\mathrm{par}}_* \Qlb)_{a,s}$ soit non nulle alors il existe un relèvement $\tilde{s}$ de $s$ à $\Tc/Z_{\Lc}^0$ tel que
$$W_a \subset \{w\in W^L \mid w\cdot \tilde{s}=\tilde{s} \}\ ;$$
\item il n'y a qu'un nombre fini de $s$ pour lesquels 
$${}^{\mathrm{p}}\mathcal{H}^{i}( Rf^{\xi,M\text{-}\mathrm{par}}_* \Qlb)_{s}\not=0 \ ;$$
\item on a la décomposition suivante
  \begin{equation}
    \label{eq:s-decomposition}
    {}^{\mathrm{p}}\mathcal{H}^{i}( Rf^{\xi,M\text{-}\mathrm{par}}_* \Qlb)=\bigoplus_{s\in (\Tc/Z_{\Mc}^0)_{\tors}}{}^{\mathrm{p}}\mathcal{H}^{i}( Rf^{\xi,M\text{-}\mathrm{par}}_* \Qlb)_{s}.
  \end{equation}
  
\end{enumerate}
    \end{theoreme}

    \begin{remarque}
      Lorsque $M=T$, on a $Z_{\Mc}=\Tc$ et la décomposition (\ref{eq:s-decomposition}) est tautologique.  Lorsque $M=G$, le théorème ci-dessus est dû à Ngô.
    \end{remarque}

    \begin{preuve} Comme $\Gamma(\uc_M,\pi_0(\Jc^1))$ est un groupe commutatif fini (cf. proposition \ref{prop:description} assertion 5) qui agit sur ${}^{\mathrm{p}}\mathcal{H}^{i}( Rf^{\xi,M\text{-}\mathrm{par}}_* \Qlb)$, on a une décomposition de la cohomologie en somme directe de composantes isotypiques indexées par les caractères de $\Gamma(\uc_M,\pi_0(\Jc^1))$ à valeurs dans $\Qb$. Un tel caractère induit un caractère d'ordre fini de $X_*(T\cap M_{\der})$ par le morphisme (\ref{eq:J1=quotient}) du lemme \ref{lem:quotient} donc un élément de $(\Tc/Z_{\Mc}^0)_{\tors}$. La décomposition (\ref{eq:s-decomposition} s'en déduit ainsi que la finitude de l'assertion 2.

Soit $L\in \lc^G(M)$ et $a$ un point géométrique de $\Ac_L^{\el}$. On a donc $W_a\subset W^L$ en vertu de la proposition \ref{prop:AM-Wa}. Le groupe  $\Gamma(\uc_M,\pi_0(\Jc^1))$ agit sur la fibre ${}^{\mathrm{p}}\mathcal{H}^{i}( Rf^{\xi,M\text{-}\mathrm{par}}_* \Qlb)_{a,s}$  via le morphisme de restriction $\Gamma(\uc_M,\pi_0(\Jc^1))\to \pi_0(\Jc^1_a)$. Soit $s'$ un caractère du groupe fini  $\pi_0(\Jc^1_a)$ tel que ${}^{\mathrm{p}}\mathcal{H}^{i}( Rf^{\xi,M\text{-}\mathrm{par}}_* \Qlb)_{a,s}$ possède une composante $s'$-isotypique non nulle. Alors, par définition de la $s$-partie, le caractère $s'$ composé avec le morphisme $X_*(T\cap M_{\der})\to \pi_0(\Jc^1_a)$ redonne $s$. Or ce dernier morphisme se factorise par $\pi_0(\yc^1_a)=X_*(T\cap L_{\der})_{W_a}$ (cf. le diagramme commutatif (\ref{eq:diag-X1})). Soit $\tilde{s}$ le caractère de $\Tc/Z_{\Lc}^0$ induit par $s'$. Cet élément $\tilde{s}$ appartient à $(\Tc/Z_{\Lc}^0)^{W_a}$ et relève $s$. Cela démontre l'assertion 1.
    \end{preuve}

    \begin{corollaire} \label{cor:dec-ell}Soit $s\in (\Tc/Z_{\Mc}^0)_{\tors}$ et $j: \Ac_G^{\el}\hookrightarrow \Ac_G$. On a alors la décomposition suivante 

$$j^*({}^{\mathrm{p}}\mathcal{H}^{i}( Rf^{\xi,M\text{-}\mathrm{par}}_* \Qlb)_s)=\bigoplus_{s'}{}^{\mathrm{p}}\mathcal{H}^{i}( Rf^{\el}_* \Qlb)_{s'}$$
      où la somme est prise sur les $s'\in  (\Tc/Z_{\Gc}^0)_{\tors}$ d'image $s$ par l'application canonique 
$$\Tc/Z_{\Gc}^0\to \Tc/Z_{\Mc}^0.$$
    \end{corollaire}

    \begin{preuve}
Sur l'ouvert $\Ac^{\el}$, le groupe $\Gamma(\uc_M,\pi_0(\Jc^1))$ agit sur $j^*({}^{\mathrm{p}}\mathcal{H}^{i}( Rf^{\xi,M\text{-}\mathrm{par}}_* \Qlb))$ via le morphisme de restriction $\Gamma(\uc_M,\pi_0(\Jc^1)) \to   \Gamma(\Ac^{\el},\pi_0(\Jc^1))$. Le corollaire est alors une conséquence évidente du diagramme commutatif 
$$\xymatrix{  X_*(T\cap M_{\der}) \ar[r] \ar[d] &    X_*(T\cap G_{\der})\ar[d] \\ \Gamma(\uc_M,\pi_0(\Jc^1))   \ar[r] &   \Gamma(\Ac^{\el},\pi_0(\Jc^1)) }$$
qui résulte du lemme \ref{lem:quotient}.
   \end{preuve}
\end{paragr}

\section{L'ouvert $\Ac^{\mathrm{bon}}_G$} \label{sec:Abon}

\begin{paragr}[L'énoncé principal.] --- \label{S:enonceAbon}Dans cette section, on reprend les notations des sections précédentes. En particulier, $C$ est un courbe projective, lisse et connexe sur $k$ de genre $g$ munie  d'un diviseur effectif de degré $>2g$ et d'un point fermé $\infty$ hors du support de $D$. Le couple $(G,T)$ est formé d'un groupe $G$ réductif et connexe sur $k$ et d'un sous-tore $T$ maximal. Le groupe de Weyl $W^G$ est noté simplement $W$ dans cette section. On rappelle également que la caractéristique de $k$ ne divise pas l'ordre de $W^G$. Soit $\Ac_G$ le $k$-schéma lisse défini au §\ref{S:esp-car}. À la suite de Ngô, on va définir une fonction semi-continue  supérieurement sur $\Ac_G$ qui sera notée $\delta$ (cf. §\ref{S:demo-Abon}). En $a\in \Ac_G(k)$, cette fonction vaut 
$$\delta_a= \dim_k(H^0(C,(\mathfrak{t}\otimes_{k}\pi_{Y_a,\ast}(\rho_{a,\ast}\mathcal{O}_{X_a}/ \mathcal{O}_{Y_a}))^{W}))$$
où $\pi_Y: Y\to C$ est le revêtement caméral de $C$ associé à $a$ (cf. §\ref{S:cameral}) et $\rho_a: X_a \to Y_a$ est le morphisme de normalisation. On pose alors la définition suivante.

\begin{definition}
  \label{def:Abon}
Soit $\Ac_G^{\mathrm{bon}}\subset \Ac_G$ le plus grand ouvert de $\Ac_G$ tel que, pour tout point de Zariski $a\in \Ac_G^{\mathrm{bon}}$, on ait l'inégalité suivante
$$\mathrm{codim}_{\Ac_G}(a)\geq \delta_a.$$
\end{definition}

\begin{remarque}
  Comme $\delta$ est semi-continue, la condition $\delta_a=\delta$ définit une partie localement fermée $\Ac_{\delta}$ de $\Ac_G$. Par noethérianité de $\Ac_G$, il n'y a qu'un nombre fini de $\delta$ pour lesquels $\Ac_\delta$ n'est pas vide. L'ouvert $\Ac_G^{\mathrm{bon}}$ est le complémentaire des composantes irréductibles $\Ac_\delta'$ des strates $\Ac_\delta$ qui vérifient
$$\mathrm{codim}_{\Ac_G}(\Ac_\delta')<\delta.$$
\end{remarque}

\emph{A priori}, $\Ac_G^{\mathrm{bon}}$ peut être vide. Si le corps de base est de caractéristique  nulle, Ngô a montré qu'on a en fait $\Ac_G^{\mathrm{bon}}=\Ac_G$. Pour notre corps de base $k$,  le théorème \ref{thm:Abon} (et la remarque qui suit) montre que l'ouvert $\Ac_G^{\mathrm{bon}}$ est assez gros. Avant d'énoncer ce théorème on va introduire certains diviseurs. Soit $M$ un groupe réductif muni d'un plongement de $T$, qui est dual d'un sous-groupe de $\Gc$ qui contient $\Tc$. Dans ce contexte, on a défini au §\ref{S:discr} le discriminant $D^M$ sur $\car_{M}$. Celui-ci détermine un diviseur de $\car_M$ qui est $\Gmk$-invariant. On en déduit un diviseur $\mathfrak{D}^M$ de $\car_{M,D}$. Soit $B\in \pc(T)$ un sous-groupe de Borel de $G$. Toujours au  §\ref{S:discr}, on a défini  une fonction homogène $R_M^B$ sur $\car_M$. Soit $ \mathfrak{R}_M^B$ le diviseur sur $\car_{M,D}$ qui s'en déduit. Comme $R_M^B$ ne dépend du choix de $B$ qu'à un signe près, le diviseur  $ \mathfrak{R}_M^B$ ne dépend pas du choix de $B$. Pour cette raison, on le note  $\mathfrak{R}_M^G$. 
 
\begin{theoreme}\label{thm:Abon}
  Soit $M$ un groupe réductif sur $k$ dual d'un sous-groupe de $\Gc$ qui contient $\Tc$ et Soit $a=(h_a,t)\in \Ac_M(k)$ tel qu'il existe un ensemble fini $S$ de points fermés de $C$ tels que
  \begin{enumerate}
  \item  la courbe $h_a(C)$ tracée dans $\car_{M,D}$ vérifie pour tout $c\in C-S$ l'une des trois hypothèses suivantes
    \begin{enumerate}
    \item $h_a(c)$ n'appartient pas au diviseur  $\mathfrak{R}_{M}^{G}\cup \mathfrak{D}^{M}$ ;
\item au point $h_a(c)$, la courbe $h_a(C)$ ne coupe pas  $\mathfrak{R}_{M}^{G}$ et coupe  transversalement le diviseur $\mathfrak{D}^{M}$ qui est lisse au point  $h_a(c)$ ;
\item  au point $h_a(c)$, la courbe $h_a(C)$ ne coupe pas  $\mathfrak{D}^{M}$ et coupe  transversalement le diviseur $\mathfrak{R}_{M}^G$ qui est lisse au point $h_a(c)$ ;
    \end{enumerate}

  \item pour tout $c$ dans $S$, la multiplicité d'intersection de $h_a(C)$ avec le diviseur $\mathfrak{D}^G$ est majorée par
$$2|S|^{-1}|W|^{-1} (\deg(D)-2g+2).$$
\end{enumerate}
Alors l'image de $a$ dans $\Ac_G$ appartient à l'ouvert $\Ac_G^{\mathrm{bon}}$.
\end{theoreme}

  \begin{remarque}
 Soit $a=(h_a,t)\in \Ac_G(k)$ tel que $h_a(C)$ coupe transversalement  le diviseur $\mathfrak{D}^G$. Alors $a=(h_a,t)$  appartient à l'ouvert $\Ac_G^{\mathrm{bon}}$. En particulier, $\Ac_G^{\mathrm{bon}}$ n'est pas vide. 
    \end{remarque}

Le reste de la section est consacrée à la démonstration de ce théorème.

\end{paragr}

\begin{paragr}[Le problème de modules $\BBB_G$.] --- On considère le problème de modules $\BBB_{G}$ suivant. À un $k$-schéma $S$ on associe l'ensemble des classes d'isomorphie de diagrammes commutatifs
$$
\xymatrix{ X\ar[r]^{f}\ar[d]_{\pi}&\mathfrak{t}_{D}\times_{k} S\ar[dl]\\
C\times_{k}S&\\}
$$
o\`{u} $X$ est une courbe projective et lisse sur $S$, à fibres géométriques non nécessairement connexes, munie d'une action de $W$, $\pi$ est un morphisme fini,  $W$-invariant et plat,  et $f$ est un morphisme fini et $W$-équivariant. En outre, ces données satisfont la condition suivante : le morphisme  $\pi$ est étale et galoisien de groupe de Galois $W$ au-dessus de $\infty\times_{k}S $ et le morphisme $f$ restreint à $\pi^{-1}(\infty\times_{k}S)$ est injectif.

On laisse au lecteur le soin de vérifier l'énoncé de représentabilité suivant.

\begin{lemme}
Le problème de modules $\BBB_{G}$ est représentable par un $k$-espace algébrique de type fini.
\end{lemme}

\end{paragr}

\begin{paragr}[Le schéma $\AAA_G$.] --- Soit $\AAA_G$ le schéma lisse tel que pour tout $k$-schéma $S$ l'ensemble $\AAA_G(S)$ est l'ensemble des sections $a$ de $\car_D\times_k S$ au-dessus de $C\times_k S$ qui envoie le $S$-point $\infty\times_k S$ dans l'ouvert $G$-régulier $\car^{\reg}_D$. On a la formule suivante pour la dimension de $\AAA_G$ (cf. \cite{ngo2} lemme 4.13.1)
  \begin{equation}
    \label{eq:dimAAA_G}
    \dim_k(\AAA_G)=\deg(D)(\dim(G)+\rang(G))/2 + \rang(G)(1-g).
  \end{equation}

\begin{lemme}\label{lem:3assertions}Soit  $S$ un $k$-schéma et $(X,\pi,f) \in \mathbb{B}_{G}(S)$. 
  \begin{enumerate}
  \item La formation du schéma quotient $X//W$ commute à tout changement de base $S'\rightarrow S$.
  \item Le morphisme $X//W\rightarrow C\times_{k}S$ induit par $\pi$ est un isomorphisme.
  \item Le morphisme $f$ définit par passage au quotient par $W$ un élément $a_f\in \AAA_G(S)$ ; l'application $f\mapsto a_f$ définit un morphisme noté $\psi$ de $\BBB_G$ dans $\AAA_G$. 
  \end{enumerate}
\end{lemme}

\begin{preuve}
  Soit $S'\rightarrow S$. Comme l'ordre de $W$ est inversible dans $k$, pour tout $S'$-algèbre munie d'une action de $W$ au-dessus de $S$, le foncteur des $W$-invariants est donné par la moyenne sur $W$ qui est $S'$-linéaire. L'assertion $1$ s'en déduit.

Soit $K$ le corps des fonctions de $X$. Il résulte de nos hypothèses que $\pi$ est génériquement étale galoisien de groupe $W$. Il s'ensuit que le sous-corps des points $W$-fixes de $K$ est le corps $F$ des fonctions de $C\times_k S$. En particulier, les sous-anneaux des points fixes sous $W$ des anneaux locaux de $X$ sont inclus dans la clôture intégrale des anneaux locaux de $C\times_k S$ : ils sont donc égaux. L'assertion 2 s'en déduit.

L'assertion $2$ implique que $f$ définit par passage au quotient un élément $a_f\in \AAA_G(S)$. On obtient bien un morphisme par l'assertion 1 d'où l'assertion 3.

\end{preuve}

\end{paragr}

\begin{paragr}[L'invariant $\delta$ et les strates $\AAA_{G,\delta}$.] --- Soit  $S$ un $k$-schéma et $(X,\pi,f) \in \BBB_{G}(S)$. Par l'assertion 3 du lemme \ref{lem:3assertions}, on en déduit une section $a$ de $\car_D\times_k S$. La plupart du temps on sous-entendra $a$ et on notera simplement $\pi : Y \to C\times_k S$ le revêtement caméral  $\pi_a : Y_a \to C\times_k S$  du § \ref{S:cameral} (on ne confondra pas avec le revêtement caméral universel noté $Y$ du §\ref{S:cameral}). On a donc un diagramme commutatif

  \begin{equation}
    \label{eq:memento}
    \xymatrix{ X \ar[r]_{\rho}    \ar[rd]_{\pi}  \ar@/^1pc/[rr]^{f} & Y \ar[d]^{\pi_Y}     \ar[r]_{f_Y}     & \tgo_D\times_k S \ar[d]^{\chi}  \\
  & C\times_k S   \ar@{=}[rd]    \ar[r]^{a}     & \car_{D}\times_k S \ar[d]^{p}\\
 & & C\times_k S      }
\end{equation}
de carré cartésien, où 

\begin{itemize}
\item $\rho$ est le morphisme de normalisation de la courbe camérale $Y$;
\item $p$ est la projection canonique ;
\item $\car_D=\car_{G,D}$ ;
\item $\chi$ est le morphisme $\chi_{G,D}^T$ du §\ref{S:cameral}.
\end{itemize}

Ce diagramme est obtenu par le changement de base $S \to \BBB_G$ associé à $(X,\pi,f)$ à partir du digramme universel 
 \begin{equation}
    \label{eq:memento-univ}
    \xymatrix{ \xc \ar[r]_{\rho_{\xc}}    \ar[rd]_{\pi_{\xc}}  \ar@/^1pc/[rr]^{f_{\xc}} & \yc \ar[d]^{\pi_{\yc}}     \ar[r]_{f_{\yc}}     & \tgo_D\times_k \BBB_G \ar[d]^{\chi}  \\
  & C\times_k \BBB_G   \ar@{=}[rd]    \ar[r]     & \car_{D}\times_k \BBB_G \ar[d]^{p}\\
 & & C\times_k \BBB_G      }
\end{equation}
où  $\xc$ la courbe universelle au-dessus de $\BBB_G$ et $\yc$ est la courbe camérale universelle au-dessus de $c \times_k \BBB_G$.

\begin{lemme} \label{lem:bij-points}
  Le morphisme $\psi : \BBB_G \to \AAA_G$ (défini au lemme \ref{lem:3assertions} assertion 3) induit une application bijective sur les espaces topologiques sous-jacents.
\end{lemme}

\begin{preuve}
  Montrons d'abord que l'application induite sur les espaces topologiques sous-jacents est surjective sur les points fermés. Soit $a$ un tel point fermé. Son corps résiduel est le corps algébriquement clos $k$. En particulier, $a$ s'identifie à une section de $\car_D$ au-dessus de $C$. Soit $Y$ la courbe camérale associée. Il résulte du fait que $a$ est génériquement régulière que $Y$ est géométriquement réduite. Soit $X$ la normalisée de $Y$. On sait alors que $X$  est lisse (cf. \cite{EGAIV4} proposition 17.15.14). Selon le diagramme (\ref{eq:memento}), on obtient un triplet $(X,\pi,f)$ qui est un $k$-point de $\BBB_G$ d'image $a$.  Cela prouve que tous les points fermés de $a$ sont atteints par le morphisme  $\BBB_G \to \AAA_G$.

Soit $a$ un point de l'espace topologique sous-jacent à $\AAA_G$ non fermé. Soit $\bar{a}$ son adhérence. Il s'agit de voir que $a$ appartient à $\psi(\psi^{-1}(\bar{a}))$. Par le théorème de Chevalley, il s'agit d'un ensemble constructible. Par \cite{EGAIII1} chap. 0  proposition 9.2.3, on a l'alternative suivante sur $\bar{a}\cap \psi(\psi^{-1}(\bar{a}))$ :
\begin{itemize}
\item soit cette partie contient un ouvert non vide de $\bar{a}$ et donc $a$ ;
\item soit cette partie est rare dans $\bar{a}$.
\end{itemize}
Il s'agit d'exclure la seconde possibilité. Celle-ci signifie que l'adhérence de $\psi(\psi^{-1}(\bar{a}))$ dans $\bar{a}$ est d'intérieur vide. En particulier, l'adhérence de $\psi(\psi^{-1}(\bar{a}))$ dans $\bar{a}$ est strictement plus petite que $\bar{a}$. Il existe donc un point fermé de $\bar{a}$ qui n'est pas dans l'adhérence de $\psi(\psi^{-1}(\bar{a}))$. Mais cela n'est pas possible puisque, comme on l'a vu, un tel point est toujours dans $\psi(\psi^{-1}(\bar{a}))$.  On a donc écarté  cette seconde possibilité et ainsi terminé la démonstration de la surjectivité.

L'injectivité est évidente puisque tout triplet $(X,\pi,f)$ de $\BBB_G$ d'image $a\in \AAA_G$ est obtenu, après extension convenable des scalaires, par normalisation de la courbe camérale associée à $a$.
\end{preuve}

Soit $p_2 : C\times_k \BBB_G \to \BBB_G$ la projection sur le second facteur. Par définition de $\BBB_G$, le $\oc_{\BBB_G}$-module 
$$p_{2,*} \pi_{\yc,*}( \rho_{\xc,*}\oc_{\xc}/ \oc_{\yc})$$
est localement libre. Comme la caractéristique de $k$ ne divise pas l'ordre de $W$, il en est de même du   $\oc_{\BBB_G}$-module
\begin{equation}
  \label{eq:ocBmodule}
  (\tgo \otimes_k p_{2,*} \pi_{\yc,*}( \rho_{\xc,*}\oc_{\xc}/ \oc_{\yc}))^W.
\end{equation}

\begin{definition} \label{def:delta}
Soit $\delta : \BBB_G \to \NN$ la fonction localement constante sur $\BBB_G$ donnée par le rang du $\oc_{\BBB_G}$-module localement libre (\ref{eq:ocBmodule}).

Soit  $\delta : \AAA_G \to \NN$ la fonction constructible sur $\AAA_G$ définie pour tout $a\in\AAA_G$ par
$$\delta(a)=\delta(b)$$
où $b\in \BBB_G$ est l'unique point tel que $\psi(b)=a$.
\end{definition}

\begin{remarques} L'existence et l'unicité de $b$ a été démontrée au lemme \ref{lem:bij-points}. La constructibilité de $\delta$ sur $\AAA_G$ résulte du théorème de Chevalley.
  
\end{remarques}

Soit $K$ est une extension de $k$ et $(X,\pi,f) \in \BBB_G(K)$. Alors, avec les notations  du diagramme (\ref{eq:memento}) pour $S=\Spec(K)$, on a 
$$ \delta (X,\pi,f)=\sum_{c} \delta_{c}(X,\pi,f)$$
où la somme est prise sur  tous les points fermés $c$ de $C\times_k K$ et l'invariant local $\delta_c$ est la longueur suivante
$$
\delta_{c}(X,\pi,f)=\mathrm{long}((\mathfrak{t}\otimes_{k}\pi_{Y,\ast}(\rho_{\ast}\mathcal{O}_{X}/\mathcal{O}_{Y}))_{c}^{W}).
$$
Cet invariant  est nul en un point $c$ tel que $Y$ est lisse en tout point $y\in Y$ au-dessus de $c$. On a également la formule suivante
$$ \delta (X,\pi,f)=\dim_K(H^0(C\times_k K, (\tgo\otimes_l \pi_{Y,*}(\rho_*\oc_X/\oc_Y))^W)).$$

\begin{proposition}\label{prop:semi-continue}
La fonction $\delta$ sur $\AAA_{G}$ est semi-continue supérieurement.
\end{proposition}

\begin{preuve}
Comme la fonction $\delta$ est constructible, il  suffit de montrer que pour tout morphisme d'un trait $S$ dans $\AAA_G$, on a $\delta(\eta)\leq \delta(s)$ où $s$ et $\eta$ sont les points fermé et générique de $S$.  Un tel  morphisme de $S$ dans $\AAA_G$ s'interprète comme une section de $\car_{D,S}$ au-dessus de $C\times_k S$. Soit $Y$ la courbe camérale relative sur $S$ qui s'en déduit. Soit $X$ la normalisée de $Y$. On est donc dans la situation du diagramme (\ref{eq:memento}) dont on utilise les notations. La fibre générique de $Y$ est génériquement réduite. Quitte à remplacer $\eta$ par une extension radicielle et $S$ par son normalisé dans cette extension, on peut supposer que la fibre générique $X_\eta$ de $X$ est lisse (cf. \cite{EGAIV4} proposition 17.15.14). Le triplet $(X_\eta,\pi_\eta,f_\eta)$ déduit de $(X,\pi,f)$ par changement de base à $\eta$ définit un point de $\BBB_G$, à valeurs dans le corps résiduel $k(\eta)$ de $\eta$, au-dessus de $a(\eta)$. On a donc 
 $$\delta(\eta)=\dim_{k(\eta)}( (\mathfrak{t}\otimes_{k} H^0(C_\eta,\pi_{Y_\eta,\ast}(\rho_{\eta,\ast}\mathcal{O}_{X_\eta}/ \mathcal{O}_{Y_\eta})))^{W}).$$
 Comme le $\oc_S$-module $\tgo\otimes_k H^0(C,\pi_{Y,\ast}(\rho_{\ast}\mathcal{O}_{X}/ \mathcal{O}_{Y}))$ est plat, le  $\oc_S$-module 
 $$(\tgo\otimes_k H^0(C,\pi_{Y,\ast}(\rho_{\ast}\mathcal{O}_{X}/ \mathcal{O}_{Y})))^{W},$$ qui en est un facteur direct, est aussi plat. Il s'ensuit qu'on a aussi
  $$\delta(\eta)=\dim_{k(s)}( (\mathfrak{t}\otimes_{k} H^0(C_s,\pi_{Y_s,\ast}(\rho_{s,\ast}\mathcal{O}_{X_s}/ \mathcal{O}_{Y_s})))^{W}),$$
 où l'on note par un indice $s$ le changement de base à $s$ et $k(s)$ le corps résiduel de $s$. De nouveau, quitte à remplacer $k(s)$ par une extension radicielle, on peut supposer que la normalisée, notée $\tilde{X}_s$, de $Y_s$ est lisse. Soit $\tilde{\rho}_s$ le morphisme $\tilde{X}_s\to Y_s$ de normalisation. On a alors
 $$\delta(s)=\dim_{k(s)}( (\mathfrak{t}\otimes_{k} H^0(C_s,\pi_{Y_s,\ast}(\tilde{\rho}_{s,\ast}\mathcal{O}_{\tilde{X}_s}/ \mathcal{O}_{Y_s})))^{W}).$$
 Comme on a $\rho_{s,*}\oc_{X,s}\subset \tilde{\rho}_{s,\ast}\mathcal{O}_{\tilde{X}_s}$, on en déduit l'inégalité cherchée $\delta(\eta)\leq \delta(s)$.

\end{preuve}

Pour terminer ce paragraphe, introduisons la définition suivante.

\begin{definition}\label{def:Adelta}
  Pour tout entier naturel  $\delta$, soit  $\AAA_{G,\delta}$ la partie localement fermée formée des $a \in \AAA_G$ tels que $\delta(a)=\delta$.
\end{definition}

\begin{remarques} D'après la proposition \ref{prop:semi-continue}, la réunion $\bigcup_{\delta\geq\delta'}\Ac_{\delta}$ est fermée. L'ouvert $\AAA_{G,0}$ est formé des $a$ pour lesquels la courbe camérale $Y_a$ est lisse.
  \end{remarques}

Dans la suite, sous certaines conditions, nous allons majorer la dimension des strates $\AAA_{G,\delta}$. 

\end{paragr}

\begin{paragr}[Complexe cotangent pour un point de $\BBB_G(k)$.] --- Soit $(X,\pi,f)$ un point de $\BBB_{G}(k)$. La théorie des déformations en ce point pour $\BBB_{G}$ 
est contr\^{o}lée par le complexe
$$
R\mathrm{Hom}_{\mathcal{O}_{X}}(L_{X/\tgo_D},\mathcal{O}_{X})^{W}
$$
o\`{u} le complexe cotangent $L_{X/\tgo_D}$ est ici un  complexe parfait de longueur $1$ concentré en degrés $-1$ et $0$,
$$
L_{X/\tgo_D}=[f^{\ast}\Omega_{\mathfrak{t}_{D}/k}\rightarrow \Omega_{X/k}^{1}]\cong [f^{\ast}\Omega_{\mathfrak{t}_{D}/C}\rightarrow \Omega_{X/C}^{1}]
$$
dont la flèche est $df$. Ce complexe a deux faisceaux de cohomologie, à savoir 
$$\mathcal{H}^{-1}(L_{X/\tgo_D})=\mathrm{Ker}(f^{\ast}\Omega^{1}_{\mathfrak{t}_D/C}\rightarrow \Omega^{1}_{X/C})$$
qui est un fibré vectoriel sur $X$ de rang le rang du groupe  $G$ et  
$$\mathcal{H}^{0}(L_{X/\tgo_D})=\Coker(f^{\ast}\Omega^{1}_{\mathfrak{t}_D/C}\rightarrow \Omega^{1}_{X/C})=\Omega_{X/Y}^{1}$$
 qui est un  un $\mathcal{O}_{X}$-module de torsion sur $X$.
\end{paragr}

\begin{paragr}[Complexe cotangent pour un point de $\AAA_G(k)$.] --- De manière analogue, la théorie des déformations de $\AAA_{G}$ en un point $a\in \AAA_G(k)$ est contrôlée par le complexe
$$
R\mathrm{Hom}_{\mathcal{O}_{Y}}(L_{Y/\tgo_D},\mathcal{O}_{Y})^{W}
$$
o\`{u} $Y=Y_a$ est la courbe camérale et  
$$
L_{Y/\tgo_D}=[f_Y^{\ast}\Omega^{1}_{\mathfrak{t}_D/k}\rightarrow \Omega^{1}_{Y/k}]\cong [f_Y^{\ast}\Omega^{1}_{\mathfrak{t}_D/C}\rightarrow \Omega^{1}_{Y/C}]
$$
est un complexe parfait concentré en degrés $-1$ et $0$, la flèche étant donnée par $df_{Y}$.  Ce complexe est quasi-isomorphe à $\mathcal{H}^{-1}(L_{Y/\tgo_D})[1]$.

\begin{lemme}\label{lem:HY}
On a 
$$
\mathcal{H}^{-1}(L_{Y/\tgo_D})=
f_{Y}^{\ast}\chi^{\ast}\Omega^{1}_{\mathfrak{car}_{D}/C}=
\pi_{Y}^{\ast}a^{\ast}\Omega^{1}_{\mathfrak{car}_{D}/C}.
$$
\end{lemme}

\begin{preuve}
Le morphisme $\chi$ induit la suite exacte 
\begin{equation}\label{eq:OmegatD}
0 \longrightarrow \chi^{\ast}\Omega^{1}_{\mathfrak{car}_{D}/C}\longrightarrow
\Omega^{1}_{\mathfrak{t}_{D}/C}\longrightarrow
\Omega^{1}_{\mathfrak{t}_{D}/\mathfrak{car}_{D}} \longrightarrow 0.
\end{equation}
Montrons l'exactitude à gauche, la seule non évidente. Comme $\chi$ est génériquement étale, la suite est génériquement exacte à gauche. Puisque $\mathfrak{car}_{D}$ est lisse sur $C$, le $\mathcal{O}_{\mathfrak{t}_{D}}$-module $\chi^{\ast}\Omega^{1}_{\mathfrak{car}_{D}/C}$ est localement libre et ne peut contenir un module de torsion non trivial.

Comme $\pi_{Y}$ est obtenu à partir de $\chi$ par le changement de
base donné par $a:C\rightarrow \mathfrak{car}_{D}$, on a 
$a\pi_{Y}=\chi f_{Y}$ et
\begin{equation}
\label{eq:omegaY}
\Omega^{1}_{Y/C}=f_{Y}^{\ast}\Omega^{1}_{\mathfrak{t}_{D}/\mathfrak{car}_{D}}.
\end{equation}
Par suite, en appliquant le foncteur $f^{\ast}_{Y}$ à
la suite (\ref{eq:OmegatD}), on obtient la suite exacte
\begin{equation}
\label{eq:OmegatD2}
0 \longrightarrow f_{Y}^{\ast}\chi^{\ast}\Omega^{1}_{\mathfrak{car}_{D}/C}\longrightarrow
f_{Y}^{\ast}\Omega^{1}_{\mathfrak{t}_{D}/C}\longrightarrow \Omega^{1}_{Y/C}
\longrightarrow 0
\end{equation}
o\`{u} l'exactitude à gauche se voit comme avant, d'o\`{u} le 
lemme.
\end{preuve}

\begin{lemme} \label{lem:pas-d-obstruction}Soit $a\in \AAA_G(k)$ et $Y=Y_a$.
  \begin{enumerate}
  \item Il n'y a pas d'obstruction à déformer  $Y\hookrightarrow \mathfrak{t}_{D}$ et $\AAA_{G}$ est lisse au point $a$ autrement dit on a 
 $$\mathrm{Ext}_{\mathcal{O}_{Y}}^{2}(L_{Y/\tgo_D},\mathcal{O}_{Y})^{W}=0 \ ;$$
\item L'espace tangent à $\AAA_{G}$ en le point $a$ est l'espace 
  $$
\mathrm{Ext}^{1}(L_{Y/\tgo_D},\mathcal{O}_{Y})^{W}$$
qui est de dimension $\dim(\AAA_G)$.
\end{enumerate}

\end{lemme}

\begin{preuve}
Soit $i\in \{0,1\}$.  Comme $\mathcal{H}^{-1}(L_{Y/\tgo_D})$ est localement
libre, on a
\begin{eqnarray*}
  \mathrm{Ext}^{i+1}(L_{Y/\tgo_D},\mathcal{O}_{Y})&=&\mathrm{Ext}^{i}(\mathcal{H}^{-1}(L_{Y/\tgo_D}),\mathcal{O}_{Y})\\
&=& H^{i}(Y,\mathcal{H}\mathit{om}_{\mathcal{O}_{Y}}(\mathcal{H}^{-1}(L_{Y/\tgo_D}),O_{Y})).
\end{eqnarray*}
D'après le lemme \ref{lem:HY}, on a 
$$\homc_{\oc_Y}(\hc^{-1}(L_{Y/\tgo_D}),\oc_Y)=\pi_Y^* \homc_{\oc_C}(a^*\Omega^1_{\car_D/C},\oc_C).$$
Soit $n$ le rang de $G$. Il existe $n$ éléments homogènes $(t_1,\ldots,t_n)$ de $k[\tgo]$ de degrés $e_1,\ldots,e_n$ tels que $k[\tgo]^W=k[t_1,\ldots,t_n]$ et $2\sum_{i=1}^n e_i =\dim(G)+n$ (cf. propriété 33 du résumé de \cite{Bki}). En utilisant ces éléments, on obtient une identification
$$a^*\Omega^1_{\car_D/C}\cong \bigoplus_{i=1}^n \oc_C(-e_iD)$$
donc un isomorphisme
$$\homc(L_{Y/\tgo_D},O_Y)=\pi_Y^*  \bigoplus_{i=1}^n \oc_C(e_iD).$$
Ainsi, on a
$$
\mathrm{Ext}^{i+1}(L_{Y/\tgo_D},\mathcal{O}_{Y})= H^{i}(Y,\pi_Y^*  \bigoplus_{i=1}^n \oc_C(e_iD))= H^{i}(C,\pi_{Y,*}\pi_Y^*  \bigoplus_{i=1}^n \oc_C(e_iD)).$$
Par la formule des projections, on a 
$$ \pi_{Y,*}\pi_Y^* \bigoplus_{i=1}^n \oc_C(e_iD)=\pi_{Y,*}\oc_Y  \otimes_{\oc_C}  \bigoplus_{i=1}^n \oc_C(e_iD) .$$
En passant aux invariants sous $W$, on trouve
$$\Ext^{i+1}(L_{Y/\tgo_D},\oc_Y)^W=H^i(C,(\pi_{Y,*}\oc_Y)^W  \otimes_{\oc_C}  \bigoplus_{i=1}^n \oc_C(e_iD))=H^i(C,\bigoplus_{i=1}^n \oc_C(e_iD))$$
puisque $(\pi_{Y,*}\oc_Y)^W=\oc_C$. Le résultat se déduit du théorème de Riemann-Roch.
\end{preuve}

\end{paragr}

\begin{paragr}[Invariant $\kappa$.] --- Considérons la situation universelle du diagramme commutatif (\ref{eq:memento-univ}). Soit $L_{\xc/\BBB_G}$ le complexe cotangent de $\xc$ relatif à $\BBB_G$. Soit 
$$\kc=p_{2,*} (\pi_{\xc,*}(\hc^{-1}(L_{\xc/\BBB_G})\otimes_{\oc_{\xc}} \Omega^1_{\xc/ C\times\BBB_G}))^W $$
qu'on voit comme un $\oc_{\BBB_G}$-module cohérent sur $\BBB_G$ (cf. \cite{Illusie} corollaire 2.3.11).

\begin{definition}\label{def:kappa} Soit $\kappa : \BBB_G \to \NN$ la fonction semi-continue supérieurement définie pour tout $b\in \BBB_G$ par
$$\kappa(b)=\dim_{k(b)} (\kc\otimes_{\oc_{\BBB_G}} k(b))$$
où $k(b)$ est le corps résiduel de $b$.

Soit $\kappa : \AAA_G \to \NN$ la fonction constructible définie par 
$$\kappa(a)=\kappa(b)$$
où $b$ est l'unique point de $\BBB_B$ d'image $a$ par la bijection $\psi$ du lemme \ref{lem:bij-points}.
\end{definition}

\begin{remarques}
  La semi-continuité de $\kappa$ sur $\BBB_G$ est une conséquence évidente du lemme de Nakayama. La constructibilité de $\kappa$ sur $\AAA_G$ résulte du théorème de Chevalley. En général, la fonction $\kappa$ n'est pas semi-continue sur $\AAA_G$. On peut montrer cependant qu'elle est semi-continue sur les strates $\AAA_{G,\delta}$ de la définition \ref{def:Adelta}.
 
\end{remarques}

Soit $(X,\pi,f)\in \BBB_G(k)$. Pour tout point fermé $c$ de $C$, soit
\begin{equation}
  \label{eq:kappalocal}
  \kappa_c(X,\pi,f)=\lng( \pi_*(\hc^{-1}(L_{X/\tgo_D})\otimes_{\oc_X}\Omega^1_{X/C})^W_c).
\end{equation}
 
\begin{remarque}
  Pour tout $\oc_X$-module $\lc$, on note, par abus,  $\pi_*(\lc)^W$ ce qui devrait être plus correctement noté $(\pi_*\lc)^W$.
  \end{remarque}

On a alors
\begin{equation}
  \label{eq:sumkappaloc}
  \kappa(X,\pi,f)=\sum_{c\in C}\kappa_c(X,\pi,f).
\end{equation}
Cette somme est bien entendu à support fini. Plus précisément, on a la proposition suivante.

\begin{proposition}\label{prop:kappanul}
Pour tout point fermé $c$ de $C$ où l'une des deux conditions suivantes est réalisée :
\begin{itemize}
\item $X$ est étale au-dessus de  $c$ ;
\item $Y$ est lisse au-dessus de  $c$,
\end{itemize}
on a 
$$\kappa_c(X,\pi,f)=0.$$
\end{proposition}

\begin{preuve}
  Soit $x$ un point fermé de $X$ tel que $\pi$ est étale en $x$. On a donc 
$$\Omega^1_{X/C,x}=0$$
et $\kappa_c(X,\pi,f)=0$ pour tout point fermé $c$ de $C$ tel que $\pi$ est étale au-dessus de $c$.

Plaçons-nous ensuite en un point $c$ au-dessus duquel  $Y$ est lisse. Dans ce cas, localement au-dessus de $c$ on a $X=Y$, $\pi=\pi_Y$ et $\hc^{-1}(L_Y)=\hc^{-1}(L_X)$. En utilisant le lemme \ref{lem:HY}, on a donc, localement au-dessus de $c$ 
$$\hc^{-1}(L_{X/\tgo_D})=\pi^*a^*\Omega^1_{\car_D/C}.$$ 
La formule des projections donne alors
$$(\pi_*(\hc^{-1}(L_{X/\tgo_D})\otimes_{\oc_X}\Omega^1_{X/C}))^W_c= (a^*\Omega^1_{\car_D/C}\otimes_{\oc_C}(\pi_*(\Omega^1_{X/C}))^W)_c$$
qui est nul puisque
\begin{equation}
  \label{eq:piOmegaW=0}
  \pi_*(\Omega^1_{X/C}))^W=0.
\end{equation}
En effet, ce dernier module est le quotient de  $\pi_*(\Omega^1_{X/k}))^W=\Omega^1_{C/k}$ (cf. lemme \ref{lem:annexes} assertion 2) par 
$$\pi_*( \pi^*\Omega^1_{C/k}))^W=  \Omega^1_{C/k} \otimes_{\oc_C} \pi_*(\oc_X)^W= \Omega^1_{C/k}$$
puisque $  \pi_*(\oc_X)^W=\oc_C$.

\end{preuve}

\end{paragr}

\begin{paragr}[Lissité de $\BBB_G$.]  ---Voici une condition, en termes de la fonction $\kappa$, suffisante  pour qu'un point de $\BBB_G(k)$ soit lisse. 

  \begin{proposition}\label{prop:lissite}
    Soit $(X,\pi,f)\in \BBB_G(k)$ tel que 
$$\deg(D)> 2g-2 +\kappa(X,\pi,f).$$
Alors on a 
$$\mathrm{Ext}_{\mathcal{O}_{X}}^{2}(L_{X/\tgo_D},\mathcal{O}_{X})^{W}=(0)$$
et $\BBB_G$ est donc lisse en un tel point.
  \end{proposition}

  \begin{preuve} On a 
$$\mathrm{Ext}_{\mathcal{O}_{X}}^{2}(L_{X/\tgo_D},\mathcal{O}_{X})=\mathrm{Ext}_{\mathcal{O}_{X}}^{1}( \hc^{-1}(L_{X/\tgo_D}),\mathcal{O}_{X}).$$
Par dualité de Serre, ce $k$-espace vectoriel est dual de $H^0(X,  \hc^{-1}(L_{X/\tgo_D})  \otimes_{\oc_X}  \Omega^{1}_{X/k} )$.
En passant au module des invariants sous $W$ (qui est isomorphe au module des co-invariants vu que l'ordre de $W$ est inversible dans $k$), on voit que l'annulation cherchée est équivalente à
$$H^0(X,  \hc^{-1}(L_{X/\tgo_D})  \otimes_{\oc_X}  \Omega^{1}_{X/k} )^W=0$$
soit encore
\begin{equation}
  \label{eq:annulationaprouver}
  H^0(C,  F_L(\Omega^{1}_{X/k}))=0,
\end{equation}
où $F_L$ est le foncteur exact de la catégorie des  $\oc_X$-modules  vers celles des $\oc_C$-modules défini par 
$$F_L(\lc)= \pi_*(\hc^{-1}(L_{X/\tgo_D})\otimes_{\oc_X} \lc )^W.$$
Appliquons  le foncteur $F_L$ à la suite exacte de $\oc_X$-modules 
\begin{equation}
  \label{eq:hurwitz}
   0 \longrightarrow \pi^* \Omega^{1}_{C/k} \longrightarrow \Omega^{1}_{X/k} \longrightarrow \Omega^{1}_{X/C} \longrightarrow 0.
 \end{equation}
 On obtient la suite exacte de $\oc_C$-modules 
 $$ 0 \longrightarrow F_L(\pi^* \Omega^{1}_{C/k}) \longrightarrow F_L(\Omega^{1}_{X/k}) \longrightarrow F_L(\Omega^{1}_{X/C}) \longrightarrow 0.$$
Le module $F_L(\Omega^{1}_{X/C})$ est de torsion et sa longueur est $\kappa(X,\pi,f)$ par définition même de cet invariant. Il existe donc un diviseur effectif $K$ sur $C$, de degré majoré par  $\kappa(X,\pi,f)$ tel que 
$$ F_L(\pi^* \Omega^{1}_{C/k})  \subset F_L(\Omega^{1}_{X/k}) \subset  F_L(\pi^* \Omega^{1}_{C/k})(K).$$
Pour prouver (\ref{eq:annulationaprouver}), il suffit donc de prouver qu'on a
\begin{equation}
  \label{eq:annulationaprouver2}
H^0(C, F_L(\pi^* \Omega^{1}_{C/k})(K))=0.
\end{equation}
Or, par la formule des projections, on a
\begin{equation}
  \label{eq:FpiOmegaC}
   F_L(\pi^* \Omega^{1}_{C/k})=    \pi_*(\hc^{-1}(L_{X/\tgo_D}))^W   \otimes_{\oc_C}    \Omega^{1}_{C/k}.
 \end{equation}
De l'injection évidente $\hc^{-1}(L_{X/\tgo_D})\hookrightarrow f^*\Omega^1_{\tgo_D/C}$, on déduit les injections 
\begin{equation}
  \label{eq:injection}
   \pi_*(\hc^{-1}(L_{X/\tgo_D}))^W  \hookrightarrow   \pi_*(\hc^{-1}(L_{X/\tgo_D}))  \hookrightarrow   \pi_*(f^*\Omega^1_{\tgo_D/C}).
 \end{equation}
d'où, par tensorisation par $\Omega^1_{C/k}(K)$ et l'isomorphisme  (\ref{eq:FpiOmegaC}), une injection
\begin{equation}
  \label{eq:injection2}
 F_L(\pi^* \Omega^{1}_{C/k})(K) \hookrightarrow   \pi_*(f^*\Omega^1_{\tgo_D/C})\otimes_{\oc_C}  \Omega^{1}_{C/k}(K).
\end{equation}
On va prouver qu'on a 
\begin{equation}
  \label{eq:annulationaprouver3}
H^0(C, \pi_*(f^*\Omega^1_{\tgo_D/C})\otimes_{\oc_C}  \Omega^{1}_{C/k}(K) )=0,
\end{equation}
ce qui prouvera (\ref{eq:annulationaprouver2}) et donc  (\ref{eq:annulationaprouver}).
L'annulation ci-dessus est équivalente à
\begin{equation}
  \label{eq:annulationaprouver4}
H^0(X, f^*\Omega^1_{\tgo_D/C} \otimes_{\oc_X}  \pi^*\Omega^{1}_{C/k}(K) )=0.
\end{equation}

Soit $\tgo^\vee$ le dual du $k$-espace vectoriel $\tgo$. On a alors un isomorphisme 
$$\Omega^1_{\tgo_D/C}=\tgo^\vee  \otimes_k \chi^* p^* \oc_C(-D) $$
d'où
  $$f^*\Omega^1_{\tgo_D/C} =\tgo^\vee  \otimes_k \pi^* \oc_C(-D).$$
Donc l'annulation (\ref{eq:annulationaprouver4}) est équivalente à l'annulation 
\begin{equation}
  \label{eq:annulationaprouver5}
H^0(X, \pi^*(\Omega^{1}_{C/k}(K-D)) )=0.
\end{equation}
Mais, pour cette dernière, il suffit d'avoir 
$$\deg( \Omega^{1}_{C/k}(K-D))<0$$
soit 
$$\deg(K)< \deg(D)-2g+2.$$
Comme $\deg(K)\leq \kappa(X,\pi,f)$, la condition 
$$\kappa(X,\pi,f)< \deg(D)-2g+2$$
est suffisante pour avoir l'annulation (\ref{eq:annulationaprouver5}) donc l'annulation (\ref{eq:annulationaprouver}) et donc $\mathrm{Ext}_{\mathcal{O}_{X}}^{2}(L_{X/\tgo_D},\mathcal{O}_{X})^{W}=(0)$.
  \end{preuve}

\end{paragr}

\begin{paragr}[Estimation de la dimension de $\BBB_G$ aux points lisses.] --- Le but de cette section est de démontrer la proposition suivante.

  \begin{proposition}\label{prop:inegalitedelta}
     Soit $(X,\pi,f)\in \BBB_G(k)$ tel que 
 $$\Ext^2(L_{X/\tgo_D},\oc_X)^W=0.$$
On a l'inégalité 
$$\dim_k(\AAA_G)-\dim_k(\Ext^1(L_{X/\tgo_D},\oc_X)^W)\geq \delta(X,\pi,f).$$
  \end{proposition}

Dans la démonstration de la proposition \ref{prop:inegalitedelta}, il sera commode de disposer de la définition suivante.
    
\begin{definition}\label{def:F}
  Soit $F_{\Omega}$ le foncteur  de la catégorie des $\oc_X$-modules vers celles des $\oc_C$-modules défini pour tout  $\oc_X$-module $\lc$ par
$$F_{\Omega}(\lc)=\pi_*(\lc\otimes_{\oc_X}\Omega^1_{X/k})^W.$$
 \end{definition}

 \begin{remarque}
    Le foncteur $F_{\Omega}$ est exact.
 \end{remarque}

 \begin{preuve}(de la proposition \ref{prop:inegalitedelta})  Sous la condition $\mathrm{Ext}_{\mathcal{O}_{X}}^2(L_{X/\tgo_D},\mathcal{O}_{X})^W=0$, l'inégalité cherchée est équivalente à 
$$\mathrm{deg}(F_{\Omega}(\rho^{\ast}\Omega^{1}_{Y/C}))-\mathrm{deg}(F_{\Omega}(\Omega^{1}_{X/C}))\geq \delta(X,\pi,f)$$
d'après la proposition \ref{prop:lacodim} du paragraphe \ref{S:auxiliaire} ci-dessous. En combinant les lemmes \ref{lem:FOmegaY}, \ref{lem:defaut} et \ref{lem:Delta} du  paragraphe \ref{S:auxiliaire} ci-dessous, on obtient 
$$\mathrm{deg}(F_{\Omega}(\rho^{\ast}\Omega^{1}_{Y/C}))=\delta(X,\pi,f) +\deg((\tgo^\vee\otimes_k\pi_*\Omega^{1}_{X/C})^W).$$
L'inégalité cherchée est donc équivalente à 
$$\deg((\tgo^\vee\otimes_k\pi_*\Omega^{1}_{X/C})^W)-\mathrm{deg}(F_{\Omega}(\Omega^{1}_{X/C}))\geq 0$$
Il suffit pour cela de vérifier qu'on a 
$$\dim_k((\tgo^\vee\otimes_k \pi_*\Omega^1_{X/C})^W_c)\geq \dim_k(F_{\Omega}(\Omega^1_{X/C})_c)$$
pour tout $c\in C(k)$. Or cette dernière inégalité résulte des assertions 3 et 4 du lemme \ref{lem:annexes} du paragraphe \ref{S:annexes}.
\end{preuve}
\end{paragr}

\begin{paragr} \label{S:auxiliaire} Dans ce paragraphe, on énonce et on démontre quelques résultats qui interviennent dans la démonstration de la proposition \ref{prop:inegalitedelta}.

\begin{proposition}\label{prop:lacodim}
Sous la condition 
$$
\mathrm{Ext}_{\mathcal{O}_{X}}^2(L_{X/\tgo_D},\mathcal{O}_{X})^W=0,
$$
on a 
$$
\mathrm{dim}_k(\mathrm{Ext}_{\mathcal{O}_{X}}^{1}(L_{X/\tgo_D},\mathcal{O}_{X})^W)=\dim(\AAA_G)+\mathrm{deg}(F_{\Omega}(\Omega^{1}_{X/C}))-\mathrm{deg}(F_{\Omega}(\rho^{\ast}\Omega^{1}_{Y/C}))
$$
\end{proposition}

\begin{preuve}
Pour tout $k$-espace vectoriel $V$, on note  $V^{\vee}$ son dual.  Par
dualité de Serre, on a
$$
\mathrm{Ext}_{\mathcal{O}_{X}}^{i}(L_{X/\tgo_D},\mathcal{O}_{X})^{\vee}
=H^{1-i}(X,L_{X/\tgo_D}\otimes_{\mathcal{O}_{X}}\Omega^{1}_{X/k})
$$
et donc
$$
(\mathrm{Ext}_{\mathcal{O}_{X}}^{i}(L_{X/\tgo_D},\mathcal{O}_{X})^{\vee})^W
=H^{1-i}(C, F_{\Omega}(L_{X/\tgo_D})).
$$
L'hypothèse de l'énoncé équivaut à
$$
H^{-1}(C,F_{\Omega}(L_{X/\tgo_D}))=0
$$
et la dimension cherchée, à savoir $\mathrm{dim}_k(\mathrm{Ext}_{\mathcal{O}_{X}}^{1}(L_{X/\tgo_D},\mathcal{O}_{X})^W)$, est égale à la caractéristique d'Euler-Poincaré notée $\chi(C,F_{\Omega}(L_{X/\tgo_D}))$ de $R\Gamma (C,F_{\Omega}(L_{X/\tgo_D}))$.
On a un triangle distingué
$$
\rho^{\ast}L_{Y/\tgo_D}\rightarrow L_{X/\tgo_D}\rightarrow [\rho^{\ast}\Omega_{Y/C}^{1}\rightarrow
\Omega_{X/C}^{1}]\rightarrow.
$$
et donc un triangle distingué
$$
F_{\Omega}(\rho^{\ast}L_{Y/\tgo_D})\rightarrow F_{\Omega}(L_{X/\tgo_D})\rightarrow [F_{\Omega}(\rho^{\ast}\Omega_{Y/C}^{1})\rightarrow F_{\Omega}(\Omega_{X/C}^{1})]\rightarrow.
$$
o\`{u} le dernier complexe est concentré en degrés $-1$ et $0$. Par suite, on a l'égalité entre caractéristiques d'Euler-Poincaré

$$\chi(C,F_{\Omega}(L_{X/\tgo_D}))=\chi(C,F_{\Omega}(\rho^*L_{Y/\tgo_D}))+ \chi(C,F_{\Omega}(\Omega_{X/C}^{1}))- \chi(C,F_{\Omega}(\rho^{\ast}\Omega_{Y/C}^{1})).$$
Par la formule de Riemann-Roch, la différence $\chi(C,F_{\Omega}(\Omega_{X/C}^{1}))- \chi(C,F_{\Omega}(\rho^{\ast}\Omega_{Y/C}^{1}))$ est égale à $\mathrm{deg}(F_{\Omega}(\Omega^{1}_{X/C}))-\mathrm{deg}(F_{\Omega}(\rho^{\ast}\Omega^{1}_{Y/C}))$.
Il reste à calculer $\chi(C,F_{\Omega}(\rho^*L_{Y/\tgo_D}))=-\chi(C,F_{\Omega}(\rho^* \hc^{-1}(L_{Y/\tgo_D})))$. D'après le lemme \ref{lem:HY}, on a 
$$\rho^*\hc^{-1}(L_{Y/\tgo_D})=\rho^*\pi_Y^* a^*\Omega^1_{\car_D/C}=\pi^* a^*\Omega^1_{\car_D/C}.$$
Par la formule des projections, on a donc
$$F_{\Omega}(\rho^* \hc^{-1}(L_{Y/\tgo_D}))=a^*\Omega^1_{\car_D/C} \otimes_{\oc_C} \pi_*(\Omega^1_{X/k})^W.$$ 
Or  $(\pi_*\Omega^{1}_{X/k})^W=\Omega^{1}_{C/k}$ (cf. lemme \ref{S:annexes} assertion  2) et on peut identifier  $a^*\Omega^1_{\car_D/C}$ à $\bigoplus_{i=1}^n \oc_C(-e_i D)$ (cf. démonstration du \ref{lem:pas-d-obstruction}). Donc, à l'aide de la formule de Riemann-Roch, on a
\begin{eqnarray*}
  \chi(C,F_{\Omega}(\rho^*L_{Y/\tgo_D}))&=&-\sum_{i=1}^n \chi(C,\Omega^1_{C/k}(-e_iD))\\
&=& (e_1+\ldots+e_n)\deg(D)-n(g-1) \\
&=& \dim(\AAA_G).
\end{eqnarray*}
Cela termine la démonstration.
\end{preuve}

\begin{lemme} \label{lem:FOmegaY}On a l'égalité suivante
$$\mathrm{deg}(F_{\Omega}(\rho^{\ast}\Omega^{1}_{Y/C}))=\deg(F_{\Omega}(f^{\ast}\Omega^{1}_{\mathfrak{t}_{D}/C}))-\deg(F_{\Omega}(f^{\ast}\chi^{\ast}\Omega^{1}_{\mathfrak{car}_{D}/C} )).$$

  \end{lemme}

  \begin{preuve}
    Appliquons le foncteur $f^*$ à la suite exacte (\ref{eq:OmegatD}). On obtient  la suite 
$$0 \longrightarrow f^{\ast}\chi^{\ast}\Omega^{1}_{\mathfrak{car}_{D}/C}\longrightarrow f^{\ast}\Omega^{1}_{\mathfrak{t}_{D}/C}\longrightarrow f^* \Omega^{1}_{\tgo_D/\car_D}\longrightarrow 0$$
qui est exacte (l'exactitude à gauche se voit comme celle de la suite \ref{eq:OmegatD2}). Il résulte de (\ref{eq:omegaY}) qu'on a  $f^* \Omega^{1}_{\tgo_D/\car_D}=\rho^{\ast}\Omega^{1}_{Y/C}$. En appliquant le foncteur exact $F_{\Omega}$ à la suite ci-dessus, on obtient la relation voulue.
  \end{preuve}

 \begin{lemme}\label{lem:defaut}
    On a 
$$\deg(F_{\Omega}(f^{\ast}\Omega^{1}_{\mathfrak{t}_{D}/C}))-\deg( (\tgo^\vee\otimes_k\pi_*(\oc_X))^W\otimes_{\oc_C}   \Omega^1_{C/k}(-D) )=\deg((\tgo^\vee\otimes_k\pi_*\Omega^{1}_{X/C})^W).$$
  \end{lemme}

  \begin{preuve}
    Identifions $\Omega^1_{\tgo_D/C}$ à $\chi^*p^* (\tgo^\vee\otimes_k \oc_C(-D))$ et $f^{\ast}\Omega^{1}_{\mathfrak{t}_{D}/C}$ à $\tgo^\vee\otimes_k \pi^*(\oc_C(-D))$. Par la formule des projections, on a
    \begin{equation}
      \label{eq:explicit1}
      F_{\Omega}(f^{\ast}\Omega^{1}_{\mathfrak{t}_{D}/C})=(\tgo^\vee\otimes_k \pi_*(\Omega^1_{X/k}))^W \otimes_{\oc_C}\oc_C(-D).
    \end{equation}
et 
\begin{equation}
      \label{eq:explicit2}
(\tgo^\vee\otimes_k  \pi_*\oc_X)^W \otimes_{\oc_C} \Omega^{1}_{C/k}.=\pi_*(\tgo^\vee\otimes_k\pi^* \Omega^{1}_{C/k})^W.
\end{equation}

Partant de la suite exacte (\ref{eq:hurwitz}), on voit qu'on a une suite exacte
  $$0 \longrightarrow \pi_*(\tgo^\vee\otimes_k\pi^* \Omega^{1}_{C/k})^W(-D) \longrightarrow  \pi_*(\tgo^\vee\otimes_k\Omega^{1}_{X/k})^W(-D) \longrightarrow  \pi_*(\tgo^\vee\otimes_k\Omega^{1}_{X/C})^W(-D) \longrightarrow 0.$$
  qu'on réécrit à l'aide de  (\ref{eq:explicit1}) et (\ref{eq:explicit2})
$$ 0 \longrightarrow  F_{\Omega}(f^{\ast}\Omega^{1}_{\mathfrak{t}_{D}/C}) \longrightarrow (\tgo^\vee\otimes_k  \pi_*\oc_X)^W \otimes_{\oc_C} \Omega^{1}_{C/k}(-D)  \longrightarrow \pi_*(\tgo^\vee\otimes_k\Omega^{1}_{X/C})^W(-D) \longrightarrow 0.$$
L'égalité annoncée s'en déduit puisqu'on a 
$$\deg(\pi_*(\tgo^\vee\otimes_k\Omega^{1}_{X/C})^W(-D))=\deg(\pi_*(\tgo^\vee\otimes_k\Omega^{1}_{X/C})^W$$
vu que  $\pi_*(\tgo^\vee\otimes_k\Omega^{1}_{X/C})^W$ est un $\oc_C$-module de torsion.

\end{preuve}

\begin{lemme}\label{lem:Delta}
  On a 
$$\deg( (\tgo^\vee\otimes_k\pi_*(\oc_X))^W\otimes_{\oc_C}   \Omega^1_{C/k}(-D) )-\deg(F_{\Omega}(f^{\ast}\chi^{\ast}\Omega^{1}_{\mathfrak{car}_{D}/C} ))=\delta(X,\pi,f).$$
\end{lemme}

\begin{preuve} On va montrer l'isomorphisme suivant
  \begin{equation}
    \label{eq:larelation}
    F_{\Omega}(f^{\ast}\chi^{\ast}\Omega^{1}_{\mathfrak{car}_{D}/C})= (\tgo^\vee \otimes_k \pi_{Y,*}\oc_Y )^W\otimes_{\oc_{C}} \Omega^1_{C/k}(-D)
  \end{equation}
  qui donne immédiatement le résultat. 

Partons de la relation (cf. lemme \ref{lem:annexes} assertion 1)
$$\Omega^1_{\car_D/C}=\chi_*(\Omega^1_{\tgo_D/C})^W.$$
Or $\Omega^1_{\tgo_D/C}$ s'identifie à $\tgo^\vee \otimes_k \chi^*p^* \oc_C(-D)$. On a donc
$$\Omega^1_{\car_D/C}=(\tgo^\vee \otimes_k \chi_* \oc_{\tgo_D})^W\otimes_{\oc_{\car_D}} p^*\oc_C(-D).$$
Par conséquent, on a
\begin{eqnarray}\label{eq:OY1}
\nonumber   a^* \Omega^1_{\car_D/C}&=&(\tgo^\vee \otimes_k a^* \chi_* \oc_{\tgo_D})^W\otimes_{\oc_{C}} \oc_C(-D)\\
&=& (\tgo^\vee \otimes_k \pi_{Y,*}\oc_Y )^W\otimes_{\oc_{C}} \oc_C(-D)
\end{eqnarray}
puisque $a^* \chi_* \oc_{\tgo_D}=\pi_{Y,*}\oc_Y$.
En utilisant le lemme \ref{lem:HY}, on a 
$$f^{\ast}\chi^{\ast}\Omega^{1}_{\mathfrak{car}_{D}/C}=\pi^*a^* \Omega^1_{\car_D/C}.$$
Par conséquent, par la formule des projections, on a 
\begin{eqnarray}\label{eq:OY2}
\nonumber   F_{\Omega}(f^{\ast}\chi^{\ast}\Omega^{1}_{\mathfrak{car}_{D}/C} )&= &a^* \Omega^1_{\car_D/C}\otimes_{\oc_C} \pi_*(\Omega^{1}_{X/k})^W\\
&=& a^* \Omega^1_{\car_D/C}\otimes_{\oc_C} \Omega^{1}_{C/k},
\end{eqnarray}
la dernière égalité résultant de $\pi_*(\Omega^{1}_{X/k})^W=\Omega^{1}_{C/k}$ (cf. lemme \ref{S:annexes} assertion 2).
En combinant les lignes (\ref{eq:OY1}) et (\ref{eq:OY2}), on obtient bien (\ref{eq:larelation}).

\end{preuve}

\end{paragr}

\begin{paragr}[Quelques lemmes annexes.] \label{S:annexes} On rassemble dans le lemme suivant quelques résultats qui ont été utilisés dans ce qui précède.
  
\begin{lemme}\label{lem:annexes}
Soit $c\in C(k)$. On a les relations suivantes
\begin{enumerate}
\item $\chi_*(\Omega^1_{\tgo_D/C})^W=\Omega^1_{\car_D/C}$ ;
\item $\pi_*(\Omega^1_{X/k})^W=\Omega^1_{C/k}$ ;
\item $\dim_k((\tgo^\vee\otimes_k \pi_*\Omega^1_{X/C})^W_c)=\dim(\tgo)-\dim(\tgo^{W_x})$ où $W_x\subset W$ est le stabilisateur d'un point $x\in X(k)$ au-dessus de $c$ ;
\item $\dim_k(F_{\Omega}(\Omega^1_{X/C})_c)=\left\lbrace \begin{array}{l} 0 \text{ si } X \text{est étale au-dessus de } c \ ; \\ 1 \text{ sinon}.
    \end{array}\right.$
  \end{enumerate}
  \end{lemme}

  \begin{preuve} Montrons l'assertion 1. On déduit de  (\ref{eq:OmegatD}) la suite exacte
$$0 \longrightarrow (\chi_* \chi^{\ast}\Omega^{1}_{\mathfrak{car}_{D}/C})^W\longrightarrow 
(\chi_* \Omega^{1}_{\mathfrak{t}_{D}/C})^W.$$
Mais, par la formule des projections, on voit que la flèche d'adjonction $\Omega^{1}_{\mathfrak{car}_{D}/C} \to \chi_* \chi^{\ast}\Omega^{1}_{\mathfrak{car}_{D}/C}$ identifie $\Omega^{1}_{\mathfrak{car}_{D}/C}$ à $(\chi_* \chi^{\ast}\Omega^{1}_{\mathfrak{car}_{D}/C})^W$. En choisissant une trivialisation du fibré $\lc_D$ dans un voisinage d'un point de $C(k)$, on est ramené à prouver 
$$\chi_*(\Omega^1_{\tgo/k})^W=\Omega^1_{\car/k},$$
où l'on note $\chi$ le morphisme canonique $\tgo \to \car$. Soit $k(\tgo)$ le corps des fractions de l'anneau des fonctions régulières sur $\tgo$. Comme $\chi$ est génériquement étale, on a
$$  k(\tgo)\otimes_{k[\tgo]}\Omega^1_{\tgo/k}= k(\tgo)\otimes_{k[\tgo]^W}\Omega^1_{\car/k}$$
d'où
\begin{equation}
  \label{eq:omegaWgen}
    (k(\tgo)\otimes_{k[\tgo]}\Omega^1_{\tgo/k})^W= (k(\tgo)\otimes_{k[\tgo]^W}\Omega^1_{\car/k})^W=k(\tgo)^W\otimes_{k[\tgo]^W}\Omega^1_{\car/k})
  \end{equation}
  Soit $t_1,\ldots,t_n$ et $\sigma_1,\ldots,\sigma_n$ des éléments homogènes de $k[\tgo]$ tels que $k[\tgo]\simeq k[t_1,\ldots,t_n]$ et $k[\tgo]^W\simeq k[\sigma_1,\ldots,\sigma_n]$. Soit $d\Sigma\in \chi_*(\Omega^1_{\tgo/k})^W$. Par (\ref{eq:omegaWgen}), il existe $b_1,\ldots, b_n$ dans $k(t)^W$ tels que $d\Sigma= \sum_{i=1}^n b_i d\sigma_i$. Soit $i$ un entier $1\leq i \leq n$. La forme 
$$d\sigma_1 \wedge \ldots \wedge d\Sigma  \wedge \ldots \wedge d\sigma_n \in \bigwedge^n\Omega^1_{\tgo/k},$$
où $d\Sigma$ remplace $d\sigma_i$, est invariante sous $W$. Par ailleurs, cette forme est égale à 
$$b_i d\sigma_1 \wedge \ldots  \wedge d\sigma_n =b_i J dt_1 \wedge \ldots  \wedge dt_n,$$
où $J$ est le déterminant jacobien $\det(\frac{\partial \sigma_i}{\partial t_j})$. On a donc $b_i J \in k[\tgo]$. L'invariance sous $W$ de la forme ci-dessus  implique que $b_j J$ est anti-invariant c'est-à-dire on a pour tout $w\in W$
$$w(b_i J) = \det(w) b_i J$$
où $\det(w)$ est le déterminant de $w$ agissant sur $\tgo$. Or l'ensemble des éléments anti-invariants de $k[\tgo]$ est $k[\tgo]^WJ$ (cf. \cite{Bki} proposition 6 Chap. V). Cela prouve qu'on a $b_i\in k[\tgo]^W$ d'où l'assertion 1.

\medskip

Passons aux autres assertions. Durant toute la démonstration, on utilisera les notations suivantes. On note $c$ un point de $C(k)$. Soit $\oc_{C,c}$ le complété de l'anneau local de $C$ en $c$ et $F_c$ le corps des fractions de $\oc_{C,c}$· La $F_v$-algèbre $\oc_X\times_{\oc_C} F_c$ est un produit fini de corps qui sont des extensions finies et modérément ramifiées de $F_v$. Ils correspondent bijectivement aux points de $X(k)$ au-dessus de $c$. Soit $x$ un tel point et $F_x$ le corps correspondant. La clôture intégrale de $\oc_{C,c}$ dans $F_x$ est le complété $\oc_{X,x}$ de l'anneau local de $X$ en $x$. Le groupe de Weyl $W$ agit sur  la fibre de $X$ au-dessus de $c$. Soit $W_x\subset W$ le stabilisateur de $x$. Puisque la caractéristique de $k$ ne divise pas l'ordre de $W$, ce sous-groupe est cyclique ; il  ne dépend du choix de $x$ qu'à conjugaison près. Le choix d'une uniformisante $\varpi$ de $\oc_{X,x}$ permet d'identifier cet anneau à l'anneau de série formelle $k[[\varpi]]$ et l'action de $W_x$ sur  $\oc_{X,x}$ s'identifie à l'action continue de $W_x$ sur   $k[[\varpi]]$ qui agit  par un caractère primitif $\zeta$ sur l'uniformisante $\varpi$. Soit $m$ l'ordre de $W_x$. Alors le sous-anneau $\oc_{C,c}$ s'identifie au sous-anneau $k[[\varpi^m]]$.

\medskip
Montrons l'assertion 2. Il suffit de la prouver dans le voisinage formel $\Spec(\oc_{C,c})$ du point $c$ auquel cas elle se réduit à l'isomorphisme
 $$(\pi_* \Omega^1_{k[[\varpi]]/k})^{W_x}=\Omega^1_{k[[\varpi^m]]/k}.$$
Or $\pi_*\Omega^1_{k[[\varpi]]/k}$ est le $k[[\varpi^m]]$-module $k[[\varpi]]d\varpi$. Donc son module des invariants sous $W_x$ est 
$$k[[\varpi^m]]\varpi^{m-1} d\varpi=k[[\varpi^m]]d\varpi^m$$ 
qui est bien le module $\Omega^1_{k[[\varpi^m]]/k}$.

\medskip

Montrons l'assertion 3. On a 
$$(\tgo^\vee\otimes_k\pi_*\Omega^{1}_{X/C})^W_c\cong (\tgo^\vee\otimes_k\pi_*\Omega^{1}_{\oc_{X,x}/\oc_{C,c}})^{W_x}$$
et $\pi_*\Omega^{1}_{\oc_{X,x}/\oc_{C,c}}$ s'identifie au $k[[\varpi^m]]$-module $k[[\varpi]]d\varpi/\varpi^{m-1}k[[\varpi]]d\varpi$. 
 Donc, comme représentation de $W_x$, c'est la somme des caractères $\zeta^i$  pour $1\leq i< m$ avec multiplicité $1$. Soit $\tgo_i\subset \tgo$ la composante isotypique de $\tgo$ sur laquelle $W_x$ agit par le caractère $\zeta^i$ pour $0\leq i< m$. On a donc
\begin{eqnarray*}
  \dim_k((\tgo^\vee\otimes_k\pi_*\Omega^{1}_{\oc_{X,x}/\oc_{C,c}})^{W_x})&=&\sum_{i=1}^{m-1}\dim(\tgo_{m-i})\\
&=& \dim(\tgo)-\dim(\tgo^{W_x})
\end{eqnarray*}
d'où l'assertion 3.

\medskip 

Prouvons l'assertion 4. Si $X$ est étale au-dessus de $c$, on a $\Omega^1_{X/C,x}=0$ d'où l'annulation cherchée. Si $X$ n'est pas étale au-dessus de $c$, on a $m\geq 2$. Par définition, on a  $F_{\Omega}(\Omega^1_{X/C})=\pi_*(\Omega^1_{X/C}\otimes_{\oc_X} \Omega^1_{X/k})^W$. On a donc 
\begin{eqnarray*}
  F_{\Omega}(\Omega^1_{X/C})_c&=& \pi_*(\Omega^1_{X/C}\otimes_{\oc_X} \Omega^1_{X/k})^W_c\\
&\simeq & (k[[\varpi]]d\varpi\otimes_k d\varpi /\varpi^{m-1}k[[\varpi]]d\varpi\otimes_k d\varpi)^{W_x}\\
&\simeq & k[[\varpi^m]]/\varpi^{m}k[[\varpi^m]]( \varpi^{m-2}d\varpi\otimes_k d\varpi)
  \end{eqnarray*}
d'où l'assertion sur la dimension.
  
  \end{preuve}

\end{paragr}

\begin{paragr}[La formule de Ngô-Bezrukavnikov.] --- Ngô a relié l'invariant local $\delta_c$ à la dimension d'une fibre de Springer affine (cf. \cite{ngo2} proposition 3.7.5 et corollaire 3.8.3). Par ailleurs,  Bezrukavnikov a démontré (cf. \cite{Bezruk}) une formule énoncée par Kazhdan-Lusztig (cf. \cite{KL}) pour cette dimension. Combinant les deux, Ngô a obtenu le lemme suivant dont on donne ici une preuve directe.

  \begin{lemme}
    \label{lem:Bezruk}
Soit $(X,\pi,f)\in \BBB_G(k)$ et $a\in \AAA_G(k)$ son image par le morphisme défini au lemme \ref{lem:3assertions}. Pour tout $c\in C(k)$, on a la formule suivante
$$\delta_c(X,\pi,f)=1/2(\val_c(D^G(a))-\dim_k(\tgo)+\dim_k(\tgo^{W_x}))$$
où $W_x\subset W$ est le stabilisateur d'un point $x\in X(k)$ au-dessus de $c$.
  \end{lemme}
  
  \begin{preuve}Soit $(X,\pi,f)\in \BBB_G(k)$ et $a\in \AAA_G(k)$ comme ci-dessus. Soit $c\in C(k)$ et $x\in X(k)$ au-dessus de $c$. On reprend les notations de la démonstration du lemme \ref{lem:annexes}. Par définition de $\delta_c(X,\pi,f)$, il suffit de prouver les deux relations suivantes (où $m$ est l'ordre de $W_x$ )
    \begin{equation}
      \label{eq:bez1}
\dim_k\big(\frac{\tgo\otimes_k \oc_{X,x}}{\pi^*(\tgo\otimes_k \pi_{Y,*}\oc_Y)^W_x}\big)=m\val_c(D^G(a))/2      
    \end{equation}
et
\begin{equation}
      \label{eq:bez2}
\dim_k\big(\frac{\tgo\otimes_k \oc_{X,x}}{\pi^*(\tgo\otimes_k \pi_{*}\oc_X)^W_x}\big)=m(\dim_k(\tgo)-\dim_k(\tgo^{W_x}))/2.
    \end{equation}
Comme on peut trouver une forme bilinéaire non dégénérée invariante par $W$, les formules ci-dessus sont équivalentes à celles obtenues lorsqu'on remplace $\tgo$ par son dual $\tgo^\vee$. 
Par les formules (\ref{eq:larelation}) et (\ref{eq:OY2}), on a
$$\pi^*(\tgo^\vee\otimes_k \pi_{Y,*}\oc_Y)^W=\pi^*a^*\Omega^1_{\car/C}=f^*\chi^*\Omega^1_{\car/C}.$$
et l'inclusion $\pi^*(\tgo^\vee\otimes_k \pi_{Y,*}\oc_Y)^W_x\hookrightarrow \tgo^\vee\otimes_k \oc_{X,x}$ s'identifie à l'inclusion de $(f^*\chi^*\Omega^1_{\car/C})_x$ dans $(f^*\Omega^1_{\tgo/C})_x$. La dimension cherchée est donc égale à $m$ fois la valuation du déterminant jacobien $J$ considérée dans la démonstration du lemme \ref{lem:annexes}. Or on sait qu'à un facteur dans $k^\times$ près, le carré de ce jacobien est égal au discriminant $D^G(a)$ (cf. proposition 5 de \cite{Bki}). La formule (\ref{eq:bez1}) s'en déduit. 

Prouvons maintenant la formule (\ref{eq:bez2}). On a 

 $$\pi^*(\tgo\otimes_k \pi_{*}\oc_X)^W_x=k[[\varpi]] \tgo_0 \oplus \bigoplus_{i=1}^{m-1} \varpi^{m-i}k[[\varpi]] \tgo_i.$$
On a donc la formule suivante 
$$\dim_k\big(\frac{\tgo\otimes_k \oc_{X,x}}{\pi^*(\tgo\otimes_k \pi_{*}\oc_X)^W_x}\big)=\sum_{i=1}^{m-1} (m-i)\dim_k(\tgo_i)$$
ce qui donne le résultat voulu vu que $\dim_k(\tgo_i)=\dim_k(\tgo_{m_i})$ par l'existence d'une forme bilinéaire  non dégénérée invariante par $W$.
  \end{preuve}

\end{paragr}

\begin{paragr}[Une majoration de $\kappa$ par le discriminant.] --- Le but de ce paragraphe est d'obtenir une majoration de $\kappa$ par un multiple du discriminant. Auparavant, nous allons donner une autre expression pour $\kappa$.

  \begin{lemme}\label{lem:kappa} Soit $(X,\pi,f)\in \BBB_G(k)$ et $c$ un point fermé de $C$. Alors
$$\kappa_c(X,\pi,f)=\lng(\pi_*(\hc^{-1}(L_{X/Y})\otimes_{\oc_X}\Omega_{X/C}^1)^W_c).$$
      \end{lemme}

 \begin{preuve}On a un triangle distingué 
$$\rho^{\ast}L_{Y/\tgo_D}\rightarrow L_{X/\tgo_D}\rightarrow L_{X/Y} \rightarrow$$
où l'on a introduit le complexe
$$L_{X/Y}=[\rho^{\ast}\Omega_{Y/C}^{1}\rightarrow \Omega_{X/C}^{1}]$$
placé en degrés $-1$ et $0$. On a donc une suite exacte
$$0\rightarrow \hc^{-1}(\rho^*L_{Y/\tgo_D}) \rightarrow \hc^{-1}(L_{X/\tgo_D}) \rightarrow  \hc^{-1}(L_{X/Y}) \rightarrow 0$$
vu que $\hc^{0}(L_{Y/\tgo_D})=0$. Il s'ensuit qu'on a une suite exacte 
$$ \pi_*(\hc^{-1}(\rho^*L_{Y/\tgo_D})\otimes_{\oc_X}\Omega_{X/C}^1)^W \rightarrow  \pi_*(\hc^{-1}(L_{X/\tgo_D}) \otimes_{\oc_X}\Omega_{X/C}^1)^W\rightarrow   \pi_*(\hc^{-1}(L_{X/Y})\otimes_{\oc_X}\Omega_{X/C}^1)^W \rightarrow 0.$$
D'après le lemme \ref{lem:HY}, on a 
$$\hc^{-1}(\rho^*L_{Y/\tgo_D})=\pi^*a^* \Omega^1_{\car_D/C}$$
d'où
 $$\pi_*(\hc^{-1}(\rho^*L_{Y/\tgo_D})\otimes_{\oc_X}\Omega_{X/C}^1)^W=a^* \Omega^1_{\car_D/C}\otimes_{\oc_X} \pi_*(\Omega_{X/C}^1)^W=0$$
car $\pi_*(\Omega_{X/C}^1)^W=0$ (cf. la ligne (\ref{eq:piOmegaW=0})). Il s'ensuit qu'on a 
$$\pi_*(\hc^{-1}(L_{X/\tgo_D}) \otimes_{\oc_X}\Omega_{X/C}^1)^W\simeq   \pi_*(\hc^{-1}(L_{X/Y})\otimes_{\oc_X}\Omega_{X/C}^1)^W.$$
Le lemme se déduit alors de la définition de $\kappa_c$, cf. l.(\ref{eq:kappalocal}).

\end{preuve}

La borne suivante, bien qu'assez grossière, nous suffira par la suite.

  \begin{proposition} \label{prop:ineg-kappa}Soit $(X,\pi,f)\in \BBB_G(k)$ et $a\in \AAA_G(k)$ le point correspondant par le morphisme $\psi$ défini à l'assertion 3 du lemme \ref{lem:3assertions}. Pour tout tout point fermé $c$ de $C$, on a l'inégalité
$$\kappa_c(X,\pi,f)\leq |W|/2 \val_c(D^G(a)).$$
      \end{proposition}

\begin{preuve}Soit $c$ un point fermé de $C$. Par le lemme \ref{lem:kappa}, on a 
\begin{eqnarray*}
  \kappa_c(X,\pi,f)&=&\lng(\pi_*(\hc^{-1}(L_{X/Y})\otimes_{\oc_X}\Omega_{X/C}^1)^W_c)\\
&\leq&  \lng(\pi_*(\hc^{-1}(L_{X/Y})\otimes_{\oc_X}\Omega_{X/k}^1)^W_c),
\end{eqnarray*}
où l'inégalité résulte de la surjection évidente
$$\hc^{-1}(L_{X/Y})\otimes_{\oc_X}\Omega_{X/k}^1 \rightarrow  \hc^{-1}(L_{X/Y})\otimes_{\oc_X}\Omega_{X/C}^1\rightarrow 1.$$
Comme on a également une injection évidente
$$1 \rightarrow\hc^{-1}(L_{X/Y})\otimes_{\oc_X}\Omega_{X/k}^1 \rightarrow \rho^*\Omega^1_{Y/C}\otimes_{\oc_X}\Omega_{X/k}^1$$
et qu'on a (cf. l.(\ref{eq:omegaY}))
$$\rho^*\Omega^1_{Y/C}=f^*\Omega^1_{\tgo_D/\car_D},$$
il vient 
$$\kappa_c(X,\pi,f)\leq \lng(\pi_*(f^*\Omega^1_{\tgo_D/\car_D}\otimes_{\oc_X} \Omega^1_{X/k})^W_c).$$
Soit $x$ un point de $X$ au-dessus de $c$. Soit $W_x\subset W$ le stabilisateur de $x$ dans $W$ et $m$ l'ordre de $W_x$. On peut donc majorer la longueur ci-dessus par la longueur du $\oc_{X,x}$-module  $(f^*\Omega^1_{\tgo_D/\car_D}\otimes_{\oc_X} \Omega^1_{X/k})_x$ (où l'on a enlevé les invariants sous $W_x$). Or celle-ci est la longueur de $(f^*\Omega^1_{\tgo_D/\car_D})_x$ donc elle est égale à la valuation (en $x$) du déterminant jacobien $J$ considéré dans la démonstration du lemme \ref{lem:annexes}, dont le carré est égal au discriminant $D^G(a)$. Cette longueur est donc 
$$m/2 \val_c(D^G(a)).$$
On obtient alors la majoration cherchée.

  \end{preuve}

\end{paragr}

\begin{paragr} Dans ce paragraphe nous allons prouver le lemme suivant.

  \begin{lemme}\label{lem:transverse}
    Soit $a\in \AAA_G(k)$ et $c\in C(k)$ tel que $\val_c(D^G(a))=1$. Alors la courbe camérale $Y_a$ associée à $a$ est lisse au-dessus de $c$.
  \end{lemme}

  \begin{preuve}
    D'après la formule de Ngô-Bezrukavnikov (cf. lemme \ref{lem:Bezruk}), il existe un entier $\delta$ tel que
$$1=\val_c(D^G(a))=2\delta+\dim_k(\tgo)-\dim_k(\tgo^{W_x})$$
où $W_x\subset W$ est le stabilisateur d'un point $x\in X(k)$ au-dessus de $c$. On a donc $\delta=0$ et $\dim_k(\tgo)-\dim_k(\tgo^{W_x})=1$. Soit $M$ le sous-groupe de Levi de $G$ défini comme le centralisateur de $\tgo^{W_x}$. Comme $\tgo^{W_x}$ est de codimension $1$ dans $\tgo$ et qu'il est central dans $M$, le rang du groupe dérivé de $M$ est $1$ (le rang nul correspond à $M=T$ ce qui est exclu : le groupe $W_x$ serait alors trivial). En particulier, le groupe $W_x$ est égal au groupe de Weyl de $M$ d'ordre $2$. Le morphisme $X\to C$ dans un complété formel de $x$ est donné par l'inclusion $k[[\varpi^2]]\hookrightarrow k[[\varpi]]$. Avec les notations de (\ref{eq:memento}), on a, par le choix d'une trivialisation du fibré en droites associé à $D$, un morphisme $f : \Spec(  k[[\varpi]]) \to \tgo\times_k \Spec(k[[\varpi^2]])$. En prenant le quotient par $W^M$, on obtient une section  $\Spec(  k[[\varpi^2]]) \to \car_M \times_k \Spec(k[[\varpi^2]])$. On est donc tout de suite ramené au cas où $G=M$. En utilisant les décompositions $\tgo=\zgo_{M}\oplus \tgo_{M_{\der}}$ et $\car_M=\zgo_{M}\oplus \car_{M_{\der}}$, où l'on note $\zgo_M$ et $\tgo_{M_{\der}}$ les algèbres de Lie du centre de $M$ et du tore $T\cap M_{\der}$, on est ramené au cas où $M$ est semi-simple de rang $1$. Un calcul direct montre que le germe local de la courbe $Y$ est de la forme $\Spec(k[[u,y]]/(y^2=a(u))$ avec $\val_u(a(u))=1$ qui est bien lisse.
  \end{preuve}

\end{paragr}

\begin{paragr}[Démonstration du théorème \ref{thm:Abon}.] --- \label{S:demo-Abon}Le morphisme d'oubli de la donnée $t$ 
$$\Ac_G \to \AAA_G$$
fait de $\Ac_G$ un  revêtement fini, étale et galoisien de groupe de Galois $W$ de $\AAA_G$. En composant par ce morphisme les fonctions $\delta$ et $\kappa$ des définitions \ref{def:kappa}, on obtient des fonctions encore notées  $\delta$ et $\kappa$. La première est semi-continue supérieurement (cf. proposition \ref{prop:semi-continue}) et la seconde est seulement constructible.

On est  enfin en mesure de donner la démonstration  du théorème \ref{thm:Abon}. On en reprend les notations du théorème \ref{thm:Abon} : en particulier, $a=(h_a,t)\in \Ac_G(k)$ vérifie les assertions 1 et 2 du théorème \ref{thm:Abon}. La condition 
$$\kappa \leq \deg(D)-2g +2$$
définit donc un ensemble constructible $\Ac^{\kappa-\mathrm{bon}}$ de $\Ac_G$. On va montrer que cet ensemble contient un voisinage du point $a$. D'après \cite{EGAIII1} chap.0  proposition 9.2.5, il faut et il suffit que pour toute partie  fermée irréductible $\mathcal{Y}$ de  $\Ac_G$ contenant $a$, l'ensemble $\Ac^{\kappa-\mathrm{bon}}\cap \mathcal{Y}$ soit dense dans $\mathcal{Y}$. Pour cela, il suffit de prouver que $\Ac^{\kappa-\mathrm{bon}}\cap \mathcal{Y}$ contient le point générique de $\mathcal{Y}$. Il suffit encore de montrer le lemme suivant.

\begin{lemme} \label{lem:steppe1}
  Soit $K$ une extension algébriquement close de $k$ et $b$ un morphisme de $\Spec(K[[u]])$ dans $\Ac_G$ de fibre spéciale $b_s=a\times_k K$. Soit $b_\eta$ le morphisme déduit de $b$ en fibre générique. Alors $b_\eta$ définit un point de $\Ac^{\kappa-\mathrm{bon}}$.
\end{lemme}

\begin{preuve} Soit $h_b$ la  section du fibré $\car_D\times_k K[[u]]$ au-dessus de la courbe relative $C\times_k K[[u]]$ sous-jacente à la donnée de $b$. Soit $C_\eta$ et $C_s$ les fibres générique et spéciale de cette courbe. Soit $S\subset C$ l'ensemble fini de points fermés de $C$ qui vérifie les assertions $1$ et $2$ du théorème \ref{thm:Abon}. Dans la suite, on identifie $S$ à un ensemble de points fermés de $C_s=C\times_k K$.

  \begin{lemme} Pour tout point fermé $c_\eta$ de $C_\eta$ qui se spécialise en un point de $C_s-S$ on a 
$$\kappa_{c_\eta}(b_\eta)=0$$
où $\kappa_{c_\eta}$ est l'invariant local défini en (\ref{eq:kappalocal}).
      \end{lemme}

      \begin{preuve}
      Si $c_\eta$ se spécialise en un point en lequel $h_a(C_s)$ ne coupe pas  $\mathfrak{R}_M^G$ et coupe transversalement ou ne coupe pas $\mathfrak{D}^M$, alors $h(C_\eta)$ en $c_\eta$ ne coupe pas ou coupe transversalement $\mathfrak{D}^G$. Dans tous les cas, le  lemme \ref{lem:transverse} implique que  la courbe camérale associé à $b_\eta$ est lisse au-dessus de $c_\eta$ d'où la nullité du $\kappa$ en ce point (cf. proposition \ref{prop:kappanul}).

Si $c_\eta$ se spécialise en un point $c_s$ en lequel $h(C_s)$ ne coupe pas $\mathfrak{D}^M$ et coupe transversalement ou ne coupe pas $\mathfrak{R}^G_M$, alors $h(C_\eta)$ en ce point coupe le diviseur $\mathfrak{D}^G$ avec une multiplicité inférieure ou égale à $2$. Si cette multiplicité est nulle ou vaut $1$, on conclut comme précédemment que $\kappa$ en $c_\eta$ est nul.
Supposons donc que cette multiplicité vaut $2$. Soit $Y$ la courbe camérale associée à $h_{b}$. Soit $Y_s$ et $Y_\eta$ ses fibres spéciale et générique. Soit $y_s$ un point au-dessus de $c_s$. Par hypothèse, il existe une unique racine $\pm\al$ de $T$ dans $G$ tel que $Y_s$ coupe en $y_s$ le diviseur défini par $\al$. En particulier, le stabilisateur dans $W$ de $y_s$ est d'ordre $2$. Soit  $Y_\eta$ la courbe camérale associée à $h_{b_\eta}$ et $y_\eta$ un point au-dessus de $c_\eta$. Ce point $y_\eta$ se spécialise en un point $y_s\in Y_s$ au-dessus de $c_s$. Le stabilisateur de $y_\eta$ dans $W$ est inclus dans le stabilisateur de $y_s$ donc son ordre est au plus $2$. Quitte à remplacer $K((u))$ par une extension radicielle et $K[[u]]$ par sa normalisée dans cette extension, on peut supposer que la normalisée $X_\eta$ de $Y_\eta$ est lisse (cf. \cite{EGAIV4} proposition 17.15.14). Soit $x_\eta$ un point de $X_\eta$ au-dessus de $y_\eta$. La stabilisateur de $x_\eta$ dans $W$ est inclus dans le stabilisateur de $y_\eta$ donc est un groupe d'ordre au plus $2$. Par la formule de Ngô-Bezrukavnikov (lemme \ref{lem:Bezruk}), on sait que la valuation du discriminant et la différence $\dim(\tgo)-\dim(\tgo^{W_{x_\eta}})$ ont même parité. Or notre hypothèse est que la valuation du discriminant est $2$. Donc la différence $\dim(\tgo)-\dim(\tgo^{W_{x_\eta}})$ est paire. Si $|W_x|=2$, cette quantité vaut $1$. Donc $W_x=\{1\}$ et $X_\eta$ est étale au-dessus de $c_\eta$. Par conséquent, $\kappa$ s'annule en $c_\eta$ (cf. proposition \ref{prop:kappanul}).
\end{preuve}

Par (\ref{eq:sumkappaloc}) et le lemme précédent, on a
$$\kappa(b_\eta)=\sum_{c_\eta\in S_\eta} \kappa_{c_\eta}(b_\eta)$$
où $S_\eta$ est l'ensemble des points $c_\eta$ de $C_\eta$ qui se spécialisent en des points de $S$. La proposition \ref{prop:ineg-kappa} donne alors la majoration
$$\kappa(b_\eta) \leq |W|/2 \sum_{c_\eta\in S_\eta} \val_{c_\eta}(D^G(b_\eta)).$$
La somme dans le majorant ci-dessus est égale $\sum_{c\in S} \val_{c}(D^G(b_s))$. Or $b_s=a\times_k K$. L'hypothèse 2 du théorème \ref{thm:Abon} entraîne la nouvelle  majoration 
$$\kappa(b_\eta)\leq \deg(D)-2g+2$$
qui prouve bien que $b_\eta$ définit un point de $\Ac^{\kappa-\mathrm{bon}}$ ce qui termine la démonstration du lemme \ref{lem:steppe1}.
\end{preuve}

On a donc prouvé que $\Ac^{\kappa-\mathrm{bon}}$ contient un voisinage ouvert de $a$ qu'on note $U$. Montrons ensuite que $U \subset \Ac^{\mathrm{bon}}_G$. Soit $\delta$ un entier et $\Ac'_{\delta}$ une composante irréductible d'une strate $\Ac_{\delta}$. Il s'agit de voir que si $U\cap \Ac'_{\delta}$ n'est pas vide alors on a 
$$\dim_k(\Ac'_{\delta})\leq \dim(\Ac_G)-\delta.$$
A un revêtement fini et étale près, l'intersection $U\cap \Ac'_{\delta}$ est l'image par le morphisme $\psi$ considéré au lemme \ref{lem:bij-points} d'une partie localement fermée de $\BBB_G$. En tout point $b\in \BBB_G(k)$ de cette partie, on a $\delta(b)=\delta$ et $\kappa(b)<\deg(D)-2g+2$. Les propositions \ref{prop:lissite} et \ref{prop:inegalitedelta} montrent que  $\mathbb{B}_G$ est lisse en un tel point de dimension $\leq \dim_k(\Ac_G)-\delta$. On a donc 
$$\dim_k(\Ac'_{\delta})=\dim_k(U\cap \Ac'_{\delta}) \leq \dim(\Ac_G)-\delta$$
ce qu'il fallait démontrer. Cela achève la démonstration du théorème \ref{thm:Abon}.
\end{paragr}

\section{L'énoncé cohomologique prolongeant celui de Ngô}\label{sec:cohomologie}

\begin{paragr}[Le théorème de décomposition] --- Dans la suite, on note avec un indice $0$ un objet sur $\Fq$ et l'omission de l'indice $0$ indique l'extension des scalaires à $k$. Soit $A_0$ un $\mathbb{F}_{q}$-schéma de type fini intègre. Suivant nos conventions, on a $A=A_0\otimes_{\Fq} k$. Soit $d_A$ la dimension de $A$. Pour tout point $a$ de $A$, soit $d_{a}$ la dimension du fermé irréductible $\overline{\{a\}}$ de $A$. Soit $i_{a}=\overline{\{a\}}\hookrightarrow X$ et $j_{a}=\{a\}\hookrightarrow \overline{\{a\}}$ les inclusions.  On dira qu'un système local de $\Qlb$-espaces vectoriels sur le schéma $\{a\}$ est admissible s'il se prolonge en un système local $\mathcal{F}_{V}$ sur un ouvert dense $V$ de $\overline{\{a\}}$.  Pour un tel système local $\mathcal{F}_{a}$, on notera $j_{a,!\ast}\mathcal{F}_{a}$ le prolongement intermédiaire $j_{V,!\ast}\mathcal{F}_{V}$ du système local $\mathcal{F}_{V}$ à $\overline{\{a\}}$ tout entier pour n'importe quel ouvert dense assez petit $j_{V}:V\hookrightarrow \overline{\{a\}}$. 

\begin{proposition} (\cite{BBD} théorème 5.3.8)\label{prop:BBD}
Soit $K_0$ un $\overline{\mathbb{Q}}_{\ell}$-faisceau pervers pur de poids $w$ sur $A_0$. Le faisceau pervers $K$ sur $A$, obtenu par extension des scalaires de $\Fq$ à $k$, admet une décomposition canonique en somme directe
$$
K=\bigoplus_{a\in A}i_{a,\ast}j_{a,!\ast}\mathcal{F}_{a}[d_{a}]
$$
où pour tous les $a$ sauf un nombre fini, $\mathcal{F}_{a}=(0)$.
\end{proposition}

\begin{preuve}
L'unicité résulte du fait que 
$$\mathrm{Hom}(i_{a,\ast}j_{a,!\ast}\mathcal{F}_{a}[d_{a}], i_{a',\ast}j_{a',!\ast}\mathcal{F}_{a'}[d_{a'}])=(0)$$
quels que soient $a\not=a'$ dans $A$.  En effet la restriction à $\overline{\{a\}}\cap \overline{\{a\}}'$ de $j_{a,!\ast}\mathcal{F}_{a}[d_{a}]$ est dans ${}^{\mathrm{p}}D_{c}^{\leq -1}(\overline{\{a\}}\cap \overline{\{a\}}',\overline{\mathbb{Q}}_{\ell})$ et donc $i_{a}^{\ast}i_{a',\ast}j_{a,!\ast}\mathcal{F}_{a}[d_{a}]$ est dans ${}^{\mathrm{p}}D_{c}^{\leq -1}(\overline{\{a\}},\overline{\mathbb{Q}}_{\ell})$.

Pour l'existence, on procède par récurrence sur la dimension de $A$.  Soit $j:U\hookrightarrow A$ un ouvert dont le réduit est lisse sur $\mathbb{F}_{q}$, qui contient tous les points génériques des composantes irréductibles de $A$ de dimension $d_{A}$ et sur lequel $j^{\ast}{}^{\mathrm{p}}\mathcal{H}^{-d_{A}}(K)$ est un système local et $j^{\ast}K=j^{\ast}{}^{\mathrm{p}}\mathcal{H}^{-d_{A}}(K)[d_{A}]$.  Il existe alors une filtration 
$$(0)=K_{0}\subset K_{1}\subset K_{2}\subset K_{3}=K $$
avec $K_{1}/K_{0}$ et $K_{3}/K_{2}$ à supports dans $A-U$ et $K_{2}/K_{1}=j_{!,\ast}j^{\ast}{}^{\mathrm{p}}\mathcal{H}^{-d_{A}}(K)[d_{A}]$. Or pour tout système local pur $\mathrm{F}$ sur $U$ et tout $\overline{\mathbb{Q}}_{\ell}$-faisceau pervers $L$ pur de m\^{e}me poids sur $A$ supporté par $A-U$, on a
$$\mathrm{Ext}^{1}(L,j_{!,\ast}\mathcal{F})=\mathrm{Ext}^{1}(j_{!,\ast}\mathcal{F},L)=(0)$$
(sur $A$ les $\mathrm{Ext}^{0}$ sont nuls comme on l'a déjà rappelé et les $\mathrm{Ext}^{1}$ ne peuvent pas avoir $1$ comme valeur propre de Frobenius d'après le théorème de Deligne). Il s'ensuit que $j_{!,\ast}j^{\ast}{}^{\mathrm{p}}\mathcal{H}^{-d_{A}}(K)[d_{A}]$ est un facteur direct de $K$ avec un supplémentaire supporté par $A-U$.
\end{preuve}
\end{paragr}

\begin{paragr} \label{S:fMA} Soit $f_0:M_0\rightarrow A_0$  un morphisme de $\mathbb{F}_{q}$-schémas séparés de type fini.  Soit $f:M\rightarrow A$ le morphisme obtenu par changement de base de $\Fq$ à $k$. On suppose que $f_0$ est propre, purement de dimension relative $d_{f}$, et que $M_0$ est lisse, purement de dimension $d_{M}=d_{f}+d_{A}$ sur $\mathbb{F}_{q}$.  On considère le faisceau pervers gradué et pur en chaque degré sur $A_0$
$$K^{\bullet}_0={}^{\mathrm{p}}\mathcal{H}^{\bullet}(Rf_{0,\ast}\overline{\mathbb{Q}}_{\ell , M_0}[d_{M}])=\bigoplus_{n}{}^{\mathrm{p}}\mathcal{H}^{n}(Rf_{0,\ast}\overline{\mathbb{Q}}_{\ell , M_0}[d_{M}]).$$
Soit $K^{\bullet}$ le faisceau pervers gradué sur $A$ déduit de $K^{\bullet}_0$ par changement de base de $\Fq$ à $k$. Par la proposition \ref{prop:BBD}, on a
$$K^{\bullet}=\bigoplus_{a\in A}i_{a,\ast}j_{a,!\ast}\mathcal{F}_{a}^{\bullet}[d_{a}]$$
où chaque $\mathcal{F}_{a}^{\bullet}$ est un système local gradué admissible sur $\{a\}$.  La dualité de Poincaré assure que $\mathcal{F}_{a}^{-n}=(\mathcal{F}_{a}^{n})^{\vee}(d_{a})$.

On définit \emph{le socle} de $Rf_{\ast}\overline{\mathbb{Q}}_{\ell , M}[d_{M}]$ comme l'ensemble fini 
$$\mathrm{Socle}(Rf_{\ast}\overline{\mathbb{Q}}_{\ell , M}[d_{M}])\subset A$$
des $a$ tels que $\mathcal{F}_{a}^{\bullet}\not=(0)$.  Pour tout $a$ dans ce socle, soit $n_{a}$, resp. $n_a^-$  le plus grand entier $n$, le plus petit, tel que $\mathcal{F}_{a}^{n}\not=(0)$. On définit \emph{l'amplitude} de $\mathcal{F}_{a}^{\bullet}$ par
$$\mathrm{Ampl}(\mathcal{F}_{a}^{\bullet})=n_a -n_a^- \leq 0.$$

\begin{lemme} \label{lem:ineg-codim}
  Soit $a\in \mathrm{Socle}(Rf_{\ast}\overline{\mathbb{Q}}_{\ell , M}[d_{M}])$. On a les inégalités
$$\mathrm{Ampl}(\mathcal{F}_{a}^{\bullet})\leq 2(d_{f}-d_{A}+d_{a})$$
et
$$\mathrm{codim}_{A}(a):=d_{A}-d_{a}\leq d_{f}.$$
\end{lemme}

\begin{preuve} Écrivons comme ci-dessus $\mathrm{Ampl}(\mathcal{F}_{a}^{\bullet})=n_a -n_a^-$. Par dualité de Poincaré, on a $n_a=-n_a^-$ de sorte qu'on a 
$$\mathrm{Ampl}(\mathcal{F}_{a}^{\bullet})=2n_a.$$
La restriction à $\{a\}$ de $R^{d_{M}+n_{a}-d_{a}}f_{\ast}\overline{\mathbb{Q}}_{\ell , M}$, qui contient $\mathcal{F}_{a}^{n_{a}}$ comme facteur direct, est non nulle. Par suite, on a $d_{M}+n_{a}-d_{a}\leq 2d_{f}$ ou encore
$$\mathrm{Ampl}(\mathcal{F}_{a}^{\bullet})\leq 2(d_{f}-d_{A}+d_{a}).$$
En particulier, on a nécessairement
$\mathrm{codim}_{A}(a):=d_{A}-d_{a}\leq d_{f}$.
\end{preuve}
\end{paragr}

\begin{paragr}[Retour sur un résultat clef de Ngô.] ---  Soit $g_0:P_0\rightarrow A_0$ un schéma en groupes commutatifs lisse purement de dimension relative $d_{f}$ sur $A_0$, à fibres connexes. On suppose que $P_0$ agit sur  $M_0$ au-dessus de $A_0$.  
Soit 
$$V_{\ell}P_0=R^{2d_{f}-1}g_{0,!}\overline{\mathbb{Q}}_{\ell}(d_{f}).$$
C'est un faisceau sur $A_0$ dont la fibre en un point géométrique $\bar{a}$ de $A_0$ est le $\Qlb$-module de Tate $V_{\ell}P_{0,\bar{a}}$ de la fibre de $P_0$ en $\bar{a}$. On a une action de l'algèbre de Lie commutative $V_{\ell}P_0$ sur  $Rf_{0,*} \Qlb$
$$V_{\ell}P_{0}\otimes Rf_{0,*} \Qlb \to  Rf_{0,*}\Qlb[-1].$$
Pour tout point géométrique $\bar{a}$ de $A_0$, cette action induit d'une part une structure de 
$$\wedge^{\bullet} V_{\ell}P_{0,\bar{a}}=\bigoplus_{i\in \NN}\wedge^{i} V_{\ell}P_{0,\bar{a}}$$
module gradué sur  la cohomologie $H^\bullet(M_{0,\overline{a}},\Qlb)$. Elle induit d'autre part une structure de 
$$\wedge^{\bullet} V_{\ell}P_{0}=\bigoplus_{i\in \NN}\wedge^{i} V_{\ell}P_{0}$$
module sur 
$${}^{\mathrm{p}}\mathcal{H}^{\bullet}(Rf_{0,*} \Qlb)=\bigoplus_{n\in \ZZ}{}^{\mathrm{p}}\mathcal{H}^{n}(Rf_{0,*} \Qlb).$$
En particulier, pour tout $a$ dans le socle de $Rf_* \overline{\mathbb{Q}}_{\ell , M}[d_{M}]$, elle induit une structure de $\wedge^{\bullet} V_{\ell}P_{0,a}$ sur $\fc^\bullet_a$.

\medskip

Pour chaque point géométrique $\overline{a}$ de $A_0$ de corps résiduel noté $k(\bar{a})$, la fibre $P_{0,\overline{a}}$ de $P_0$ admet un dévissage canonique
\begin{equation}
 \label{eq:Chevalley}1\rightarrow P_{0,\overline{a}}^{\mathrm{aff}}\rightarrow P_{0,\overline{a}}\rightarrow P_{0,\overline{a}}^{\mathrm{ab}}\rightarrow 1
\end{equation}
o\`{u} $P_{0,\overline{a}}^{\mathrm{ab}}$ est un $k(\overline{a})$-schéma abélien et $P_{0,\overline{a}}^{\mathrm{aff}}$ est un $k(\overline{a})$-schéma en groupes affine.  La dimension de $P_{0,\overline{a}}^{\mathrm{aff}}$ ne dépend que de l'image $a$ de $\overline{a}$ et sera notée $\delta_{a}^{\mathrm{aff}}$.  La fonction $a\mapsto\delta_{a}^{\mathrm{aff}}$ est semi-continue supérieurement.  En particulier, si $a$ est un point de $A$, la fonction $b\mapsto\delta_{b}^{\mathrm{aff}}$ est constante de valeur $\delta_{a}^{\mathrm{aff}}$ sur un ouvert dense de $\overline{\{a\}}$. Soit $\overline{a}$ un  point géométrique de $A_0$. On a un dévissage analogue au niveau des modules de Tate
\begin{equation}
 \label{eq:Tate}
  0\rightarrow V_{\ell}P_{0,\overline{a}}^{\mathrm{aff}}\rightarrow V_{\ell}P_{0,\overline{a}}\rightarrow V_{\ell}P_{0,\overline{a}}^{\mathrm{ab}}\rightarrow 0
\end{equation}
où $V_{\ell}P_{0,\overline{a}}^{\mathrm{ab}}$ est de dimension $2(d_{f}-\delta_{a}^{\mathrm{aff}})$ et $V_{\ell}P_{0,\overline{a}}^{\mathrm{aff}}$ est de dimension égale à la dimension de la partie torique de $P_{0,\overline{a}}^{\mathrm{aff}}$. 

\medskip

On vient de voir précédemment que la cohomologie $H^\bullet(M_{0,\overline{a}},\Qlb)$ est munie d'une structure de $\wedge^{\bullet} V_{\ell}P_{0,\overline{a}}$-module gradué. Tout scindage de la suite exacte (\ref{eq:Tate}) induit une action de 
$$\wedge^{\bullet} V_{\ell}P_{\overline{a}}^{\mathrm{ab}}$$
sur $H^\bullet(M_{0,\overline{a}},\Qlb)$. Cette action dépend \emph{a priori} du scindage. Si le point géométrique  $\overline{a}$ est localisé en un point fermé $a$ de $A_0$, le corps résiduel $k(a)$ est fini et l'extension de $P_{0,\overline{a}}^{\mathrm{ab}}$ par $P_{0,\overline{a}}^{\mathrm{aff}}$ donnée par (\ref{eq:Chevalley}) est d'ordre fini. Elle est donc, à une isogénie près, scindée et le scindage obtenu est nécessairement unique. En particulier la suite (\ref{eq:Tate}) est canoniquement scindée (cf. \cite{ngo2} proposition 7.5.3).

Soit $a$ un point du  socle de $Rf_* \overline{\mathbb{Q}}_{\ell , M}[d_{M}]$. En fait, le dévissage  (\ref{eq:Chevalley}) existe pour $a$ après un changement de base radiciel. En particulier, on peut définir le dévissage (\ref{eq:Tate}) pour $a$. Par un argument de poids (cf. \cite{ngo2} §7.4.8), l'action de  $\wedge^{\bullet} V_{\ell}P_{0,a}$ sur $\fc^\bullet_a$ se factorise en une action de  $\wedge^{\bullet} V_{\ell}P_{0,a}^{\mathrm{ab}}$ sur $\fc^\bullet_a$. Ainsi $\fc^\bullet_a$ est un module gradué sur  $\wedge^{\bullet} V_{\ell}P_{0,a}^{\mathrm{ab}}$.

Suivant Ngô (cf. \cite{ngo2} § 7.1.4), on dit que $V_{\ell}P_0$ est \emph{polarisable} si, localement pour la topologie étale de $A_0$, il existe une forme bilinéaire alternée
$$
V_{\ell}P_0\times V_{\ell}P_0\rightarrow \overline{\mathbb{Q}}_{\ell}(1)
$$
telle que, pour chaque point géométrique $\overline{a}$ de $A_0$, le noyau de la forme bilinéaire induite sur $V_{\ell}P_{\overline{a}}$ soit $V_{\ell}P_{\overline{a}}^{\mathrm{aff}}$ et que la forme induite sur $V_{\ell}P_{\overline{a}}^{\mathrm{ab}}$ soit non dégénérée. La proposition suivante est une reformulation d'un résultat de Ngô.

\begin{proposition}\label{thm:amplitude}Soit $f_0:M_0\rightarrow A_0$ un morphisme de $\mathbb{F}_{q}$-schémas séparés de type fini. Soit $g_0:P_0\rightarrow A_0$ un schéma en groupes commutatifs lisse purement de dimension relative $d_{f}$ sur $A_0$, à fibres connexes. Soit $P_0\times M_0 \to M_0$ une action de $P_0$ sur $M_0$ au-dessus de $A_0$.

Les hypothèses suivantes 
\begin{enumerate}
\item  
  \begin{enumerate}
\item  $f_0$ est propre, purement de dimension relative $d_{f}$ ;
\item  $M_0$ est lisse, purement de dimension $d_{M}=d_{f}+d_{A}$ sur $\mathbb{F}_{q}$ ;
\end{enumerate}
\item  le module de Tate  $V_{\ell}P_0$ est \emph{polarisable} ;
\item  pour tout point géométrique $\bar{a}$ localisé en un point fermé $a$ de $A_0$, l'action de $\wedge^{\bullet} V_{\ell}P_{\overline{a}}^{\mathrm{ab}}$ sur  $H^\bullet(M_{\overline{a}},\Qlb)$, induite par le scindage canonique de  (\ref{eq:Tate}) décrit ci-dessus, fait de  $H^\bullet(M_{\overline{a}},\Qlb)$ un  $\wedge^{\bullet} V_{\ell}P_{\overline{a}}^{\mathrm{ab}}$-module libre ;

\end{enumerate}
entrainent que pour tout point $a\in \mathrm{Socle}(Rf_{\ast}\overline{\mathbb{Q}}_{\ell , M}[d_{M}])$, la fibre en $a$ de $\fc_a^\bullet$ est un  $\wedge^{\bullet} V_{\ell}P_{0,a}^{\mathrm{ab}}$-module libre.
\end{proposition}

\begin{preuve} Il s'agit essentiellement du contenu de la proposition 7.4.10 de \cite{ngo2}. Les hypothèses sous lesquelles se place Ngô dans cette proposition ne sont pas exactement celles que nous adoptons. Plus précisément, nous ne supposons pas que le morphisme $f_0$ est projectif. En fait, dans la démonstration de Ngô, la projectivité de $f_0$ n'intervient qu'au travers de la proposition 7.5.1 de \cite{ngo2} qui a pour conclusion notre hypothèse 3.
  
\end{preuve}

À la suite de Ngô, on en déduit le corollaire suivant.

\begin{corollaire} \label{cor:amplitude}
Sous les hypothèses de la proposition \ref{thm:amplitude}, pour tout point $a\in A$ qui appartient au socle de $Rf_{\ast}\overline{\mathbb{Q}}_{\ell , M}[d_{M}]$, on a les inégalités
$$
\mathrm{Ampl}(\mathcal{F}_{a}^{\bullet})\geq 2(d_{f}-\delta_{a}^{\mathrm{aff}})
$$
et 
$$
\delta_{a}^{\mathrm{aff}}\geq \mathrm{codim}_{A}(a).
$$
\end{corollaire}

\begin{preuve} La première inégalité est une conséquence évidente de la proposition \ref{thm:amplitude} et la seconde résulte du lemme \ref{lem:ineg-codim}.
  \end{preuve}

Le lemme suivant nous sera utile pour vérifier les hypothèses de la proposition \ref{thm:amplitude}.

\begin{lemme}\label{lem:var-ab}
  Soit $\la : B' \to B$ une isogénie entre variétés abéliennes définies sur un corps algébriquement clos $K$. Soit $X$ un $K$-schéma (ou plus généralement un $K$-champ de Deligne-Mumford) muni d'une action de $B'$ et d'un morphisme propre $f:X\to B$ qui est $B'$-équivariant lorsque $B'$ agit sur $B$ par translation via l'isogénie $\al$. Alors l'action de  
$$H_\bullet(B',\Qlb) = \wedge^{\bullet} V_{\ell}B'$$
sur $H^\bullet(X,\Qlb)$ induite par l'action de $B'$ sur $X$ fait de $H^\bullet(X,\Qlb)$ un $\wedge^{\bullet} V_{\ell}B'$-module libre.
\end{lemme}

\begin{preuve}
  On considère le diagramme cartésien
$$\xymatrix{X_0\times_K B' \ar[rr]^{\al'} \ar[d]_{\mathrm{pr}_{B'}} & & X  \ar[d]^{f} \\ B' \ar[rr]^{\al} & &   B   }$$
où $\mathrm{pr}_{B'}$ est la projection sur le second facteur, $X_0$ est la fibre de $f$ à l'origine de $B$ et $\al'$ est la restriction à $X_0\times_K B'$ de l'action $X\times_K B \to X$.
La cohomologie de $X_0\times_K B$ qui est donnée par le produit tensoriel $H^\bullet(X)\otimes H^\bullet(B')$ est évidemment un module libre sur $\wedge^{\bullet} V_{\ell}B'$.  La cohomologie de $X$ s'identifie aux invariants sous $\ker(\al)$ dans la cohomologie de $X_0\times_K B'$. Or $\ker(\al)$ agit trivialement sur  $H^\bullet(B')$ par l'argument habituel d'homotopie. Ainsi la cohomologie de $X$ qui s'identifie à $H^\bullet(X)^{\ker(\al)}\otimes H^\bullet(B')$ est bien libre comme module sur $\wedge^{\bullet} V_{\ell}B'$.
\end{preuve}

\end{paragr}

\begin{paragr}[Modèle de Néron.] --- On revient aux notations utilisées depuis le début de l'article. Soit  $a\in \Ac_G$ un point fermé. Soit $\nu :J_{a}\hookrightarrow \tilde{J}_{a}$ le modèle de Néron du schéma en groupes $J_{a}$ sur $C$.  La flèche $\nu$ est un isomorphisme au dessus de l'ouvert $U_{a}$ de sorte que  le quotient $\tilde{J}_{a}/J_{a}$ est à support ponctuel. On a la proposition suivante due à Ngô (cf. \cite{ngo2} corollaire 4.8.1).

\begin{proposition}\label{prop:Neron} Soit $a=(h_a,t)\in \Ac_G(k)$. Soit $\rho_{a}:X_{a}\rightarrow Y_{a}$ la normalisée de la courbe $Y_{a}$. 
Il existe un unique isomorphisme
$$
\tilde{J}_{a} \buildrel\sim\over\longrightarrow 
 (X_{\ast}(T)\otimes_{\mathbb{Z}}
\pi_{a,\ast}\rho_{a,\ast}\mathbb{G}_{\mathrm{m},X_{a}})^{W},
$$
où l'exposant $W$ désigne les points fixes sous $W$, qui s'insère dans le carré cartésien
$$
\xymatrix{J_{a}\ar[r]\ar@{^{(}->}[d]_{\nu} &(X_{\ast}(T)\otimes_{\mathbb{Z}}
\pi_{a,\ast}\mathbb{G}_{\mathrm{m},Y_{a}})^{W}\ar@{^{(}->}[d]\cr
\tilde{J}_{a}\ar[r]^-{\sim}&
(X_{\ast}(T)\otimes_{\mathbb{Z}}
\pi_{a,\ast}\rho_{a,\ast}\mathbb{G}_{\mathrm{m},X_{a}})^{W}}
$$
où 
\begin{itemize}
\item la flèche horizontale du haut est obtenue par changement de base par $h_a$  à partir du morphisme  $$J^G_D \longrightarrow (X_*(T)\otimes_\ZZ (\chi^T_G)_*\mathbb{G}_{m,\tgo_D})^W$$
déduit de (\ref{eq:JGversT});
\item la flèche verticale de droite est induite par la flèche d'adjonction 
$$\mathbb{G}_{\mathrm{m},Y_{a}}\hookrightarrow 
\rho_{a,\ast}\mathbb{G}_{\mathrm{m},X_{a}}.$$
\end{itemize}

En particulier, l'algèbre de Lie de $\tilde{J}_{a}$ est l'algèbre de Lie commutative de fibré vectoriel sous-jacent
$$
(X_{\ast}(T)\otimes_{\mathbb{Z}}\pi_{a,\ast}\rho_{a,\ast}\mathcal{O}_{X_{a}})^{W}.
$$
\end{proposition}

\end{paragr}

\begin{paragr}[Application à la fibration de Hitchin.] --- On suppose désormais que le groupe $G$ est \emph{semi-simple}. Soit $(G_0,T_0)$ un couple formé d'un groupe $G_0$ réductif connexe sur $\Fq$ et de $T_0\subset G_0$ un sous-$\Fq$-tore maximal et déployé sur $\Fq$.  Soit $\tau$ l'automorphisme de Frobenius de $k$ donné par $x\mapsto x^q$. On suppose que le couple obtenu à partir de $(G_0,T_0)$ après extension des scalaires de $\Fq$ à $k$ est le couple $(G,T)$. Le couple $(G,T)$ est alors muni de l'action naturelle de $\tau$.

On a fixé au §\ref{S:cbe} une courbe $C$ projective, lisse et connexe sur $k$ ainsi qu'un diviseur $D$ et un point $\infty\in C(k)$. On suppose désormais que $C$ provient par changement de base d'une courbe  $C_0$ sur $\Fq$. On note encore $\tau$ l'automorphisme $1\times \tau$ de $C=C_0\times_{\Fq} k$. On suppose que $D$ et $\infty$ sont fixes sous $\tau$. 

Le schéma $\Ac_G$ est alors muni de l'action naturelle du Frobenius $\tau$. Celle-ci laisse stable les sous-schémas fermés $\Ac_M$ pour $M\in \lc^G(T)$. De même, le champ $\mc_G$ est muni de l'action naturelle de $\tau$.  Soit $\xi\in \ago_T$ \emph{en position générale} (cf. remarque \ref{rq:position}). Soit $f^\xi : \mc^\xi_G \to \Ac_G$ le morphisme de Hitchin tronqué qui est propre de source lisse sur $k$ (cf. §\ref{S:xistab}). Alors le sous-champ ouvert $\mc^\xi_G$ est stable par $\tau$ et le morphisme $f^\xi$ est $\tau$-équivariant. De manière équivalente, on aurait pu directement définir des objets $f_0^\xi : \mc_{G_0}^\xi \to \Ac_{G_0}$ sur $\Fq$ qui redonnent après extension des scalaires de $\Fq$ à $k$ les objets analogues notés sans l'indice $0$. 

 Soit $d_{\mc^{\xi}}$ la dimension relative de $f^\xi$. Rappelons qu'on a défini un ouvert $\Ac_G^{\mathrm{bon}}$ de $\Ac_G$ (cf. définition \ref{def:Abon}). 
  
\begin{theoreme}\label{thm:support}
Le socle de la restriction de  $Rf_{*}^\xi \Qlb [d_{\mc^{\xi}}]$ à $\Ac^{\mathrm{bon}}$ est entièrement contenu dans l'ouvert elliptique de $\Ac^{\el}\cap \Ac^{\mathrm{bon}}$, de sorte que si 
$$j:\mathcal{A}^{\mathrm{ell}}\cap \mathcal{A}^{\mathrm{bon}} \hookrightarrow \mathcal{A}^{\mathrm{bon}}$$
 est l'inclusion de cet ouvert, on a un isomorphisme canonique
$${}^{\mathrm{p}}\mathcal{H}^{i}(Rf_{\ast}^\xi \overline{\mathbb{Q}}_{\ell}[d_{\mathcal{M}^{\xi}}])\cong j_{!\ast}j^{\ast}{}^{\mathrm{p}}\mathcal{H}^{i}(Rf_{\ast}^\xi \overline{\mathbb{Q}}_{\ell}[d_{\mathcal{M}^{\xi}}])$$
sur $\mathcal{A}^{\mathrm{bon}}$ pour tout entier $i$.
\end{theoreme}

\begin{preuve} On va appliquer le corollaire \ref{cor:amplitude} à la fibration de Hitchin $\mc^\xi$ muni de l'action du champ de Picard $\Jc^1$. On va en déduire que tout $a$ qui est dans le socle de $Rf_{*}^\xi \Qlb [d_{\mc^{\xi}}]$ vérifie l'inégalité $\delta_{a}^{\mathrm{aff}}\geq \mathrm{codim}_{\Ac}(a)$. On montrera ensuite qu'on a $\delta_{a}^{\mathrm{aff}}\leq \delta_{a}$ avec égalité si et seulement si $a$ est elliptique. Donc, si $a$ n'est pas elliptique, $a$ ne peut pas être dans $\Ac^{\mathrm{bon}}$.

\medskip

Nous allons tout d'abord montrer que pour tout $a\in \Ac_G(k)$, la fibre $\mc_a^\xi$ vérifie les hypothèses requises par la proposition \ref{thm:amplitude} et  le corollaire \ref{cor:amplitude}. Pour cela, on applique le lemme \ref{lem:var-ab} non pas à $\mc_a$ mais à un champ homéomorphe ce qui suffit pour démontrer la liberté de la cohomologie de $\mc_a$. 

Soit $M\in \lc^G$ tel que $a$ est l'image d'un point $(h_{a_M},t)\in \Ac_M^{\el}(k)^\tau$. Soit $X_a=\eps_M(h_{a_M}) \in \mgo(F)$ où $\eps_M$ est la section de Kostant et $F$ est le corps des fonctions de $C$. Soit $\AAA$ l'anneau des adèles de $F$ et $\oc$ le produit des complétés des anneaux locaux des points fermés de $C$. Soit $S$ l'ensemble fini de points fermés de $C$ complémentaire de l'ouvert $h_{a_M}^{-1}(\car_{M,D}^{\reg})$. Pour tout $v\in S$, soit $\Sc_v$ la fibre de Springer affine associée à $v$ et $X_a$. Soit $F_v$ le corps des fractions du complété de l'anneau local $\oc_v$ en $v$. Soit $J$ le centralisateur de $X_a$ dans $M\times_k F$. Le ind-schéma $J(F_v)/J(\oc_v)$ agit sur $\Sc_v$. Soit $\Xgo$ le produit contracté de $J(\AAA)/J(\oc)\times \prod_{v\in S} \Sc_v$ par l'action diagonale de $\prod_v J(F_v)/J(\oc_v)$. Le quotient champêtre de $\Xgo$ par l'action (sur le premier facteur) de $J(F)$ est homéomorphe à la fibre de Hitchin $\mc_a$. L'homéomorphisme entre $[J(F)\back \Xgo]$ et $\mc_a$ est équivariant sous l'action du champ de Picard  $\Jc_a=[J(F) \back J(\AAA)/ J(\oc) ]$.

Ce dernier admet comme quotient $[J(F) \back J(\AAA)/ J(\oc) \prod_v J(F_v)]$ qui est isogène au champ $\Jc_a^{\mathrm{ab}}$ des  $\tilde{J}_a$-torseurs où  $\tilde{J}_a$ est le modèle de Néron de $J_a$. Le champ $\Jc_a$ agit sur  $[J(F)\back \Xgo]$ et le morphisme évident de $[J(F)\back \Xgo]$  vers $[J(F) \back J(\AAA)/ J(\oc) \prod_v J(F_v)]$ est $J_a$-équivariant. 

 Plus généralement, on peut définir un $\Xgo^\xi$ «tronqué» dont le quotient par $J(F)$ est muni d'une action par $\Jc_a^1$ et homéomorphe de manière $\Jc_a^1$-équivariante à  $\mc_a^\xi$. On a  un morphisme $\Jc_a^1$-équivariant de $[J(F)\back \Xgo^\xi]$ vers   $[J(F) \back J(\AAA)/ J(\oc) \prod_v J(F_v)]$ dont l'image est isogène au champ $\Jc^{1,\mathrm{ab}}_a$ image de $\Jc_a^1$ dans  $\Jc_a^{\mathrm{ab}}$. On obtient une action de  $\Jc^{1,\mathrm{ab}}_a$ sur $[J(F)\back \Xgo^\xi]$ en choisissant, comme il est loisible,  une section à isogénie près du morphisme canonique $\Jc^{1,\mathrm{ab}}_a\to \Jc^{1}_a$. 

Les hypothèses du lemme \ref{lem:var-ab} sont donc vérifiées pour le quotient  $[J(F)\back \Xgo^\xi]$ et on en déduit la liberté voulue de la cohomologie de $[J(F)\back \Xgo^\xi]$ donc de $\mc_a^\xi$.

\medskip

Comparons maintenant pour un point  $a$ de $\Ac_G$ les invariants $\delta_a$ et $\delta_a^{\mathrm{aff}}$. Le champ $\Jc^1_a$ est un champ de Deligne-Mumford dont les groupes d'automorphismes sont tous le groupe fini $H^0(C,J_a)$. Soit $\Jc^0_a$ la composante neutre de $\Jc^1_a$. On rappelle $\tilde{J}_a$ le modèle de Néron de $J_a$. Soit $\Jc^{0,\ab}_a$ le schéma abélien qui représente le foncteur des classes d'isomorphismes des $\tilde{J}_a$-torseurs. Le morphisme $\Jc^0     \longrightarrow \Jc^{0,\ab}$, qui envoie un $J_a$-torseur sur la classe du $\tilde{J}_a$-torseur obtenu lorsqu'on le pousse par le morphisme canonique $J_a\to \tilde{J}_a$, s'inscrit dans la suite exacte
$$1 \longrightarrow \Jc^{0,\aff}  \longrightarrow \Jc^0     \longrightarrow \Jc^{0,\ab}  \longrightarrow 0$$
où le noyau  $\Jc^{0,\aff}$ est le quotient d'un groupe algébrique affine par le groupe fini $H^0(C,J_a)$ agissant trivialement. 
  On va calculer la dimension de ces objets par un calcul d'algèbre de Lie. 
La suite exacte 
$$ 1 \longrightarrow J_a  \longrightarrow \tilde{J}_a  \longrightarrow \tilde{J}_a/J_a \longrightarrow 1 $$
donne la suite exacte longue
$$ 1 \longrightarrow H^0(C,J_a)  \longrightarrow H^0(C,\tilde{J}_a)  \longrightarrow H^0(C,\tilde{J}_a/J_a)  \longrightarrow H^1(C,J_a)  \longrightarrow H^1(C,\tilde{J}_a)  \longrightarrow 1$$
puisque $ \tilde{J}_a/J_a$ étant à support ponctuel a un $H^1$ trivial.
Soit $\delta_a^{\aff}=\dim( \Jc^{0,\aff})$. On a donc
\begin{eqnarray*}
  \delta_a^{\aff}&=&  \dim(\Jc^0) -  \dim(\Jc^{0,\ab})\\
&=&  \dim(H^1(C,J_a))- \dim(H^0(C,J_a)) -  \dim(H^1(C,\tilde{J}_a))\\
&=&  \dim(H^0(C,\tilde{J}_a/J_a))-\dim(H^0(C,\tilde{J}_a ))
\end{eqnarray*}
la dernière égalité provenant de la suite exacte longue ci-dessus.

D'après \ref{prop:Neron}, on a 
\begin{eqnarray*}
  \dim(H^0(C,\tilde{J}_a/J_a))&=& \dim_k(H^0(C,(\mathfrak{t}\otimes_{k}\pi_{Y_a,\ast}(\rho_{a,\ast}\mathcal{O}_{X_a}/ \mathcal{O}_{Y_a}))^{W}))\\
&=& \delta_a
\end{eqnarray*}
où $\delta_a$ est la valeur en $a$ de la fonction $\delta$ considérée à la section \ref{sec:Abon}. 
Il résulte aussi de  \ref{prop:Neron} qu'on a 
$$H^0(C,\tilde{J}_a)=T^{W_a}$$
où $W_a$ est le sous-groupe de $W^G$ défini au §\ref{S:Wa}. 
Combinant les résultats précédents, on a
$$ \delta_a^{\aff}= \delta_a-\dim(T^{W_a}).$$

Supposons maintenant $a\notin \Ac_G^{\el}$. Il s'ensuit que pour un certain $M\in \lc^G$, on a $a\in \Ac_M$ et donc $W_a\subset W^M$ (cf. proposition \ref{prop:AM-Wa}). En particulier, le centre de $M$ est inclus dans $T^{W_a}$. On a donc $\dim(T^{W_a})>0$ et 
$$ \delta_a^{\aff} <\delta_a.$$
Supposons, de plus $a\in \Ac^{\mathrm{bon}}_G$ (cf. définition \ref{def:Abon}). On a donc $\delta_a\leq  \mathrm{codim}_{\Ac}(a)$ ce qui, combiné avec l'inégalité précédente, donne 
\begin{equation}
  \label{eq:<}
  \delta_a^{\aff} < \mathrm{codim}_{\Ac}(a).
\end{equation}
Supposons que  $a\notin \Ac_G^{\el}$ soit dans le  socle de la restriction de  $Rf_{*}^\xi \Qlb [d_{\mc^{\xi}}]$ à $\Ac^{\mathrm{bon}}$. Par le corollaire \ref{cor:amplitude} (dont on a vérifié les hypothèses en début de démonstration)), on a l'inégalité 
$$
\delta_{a}^{\mathrm{aff}}\geq \mathrm{codim}_{\Ac}(a).
$$
Mais celle-ci contredit manifestement (\ref{eq:<}). Donc tout  $a\in \Ac_G$ qui est dans le  socle de la restriction de  $Rf_{*}^\xi \Qlb [d_{\mc^{\xi}}]$ à $\Ac^{\mathrm{bon}}$ appartient à l'ouvert  $\Ac_G^{\el}$.

\end{preuve}

\end{paragr}

\begin{paragr}[Les espaces $\Ac_{G_s,\rho}$.] --- Dans ce paragraphe, on définit l'espace $\Ac_{G_s,\rho}$ associés à une «donnée endoscopique géométrique» de $G$ (cf. définition \ref{def:endogeo} ci-dessous) et on en donne quelques propriétés. On ne suppose pas $G$ semi-simple. Les notations sont celles du paragraphe précédent. Soit $a\in \Ac_G$. On a introduit au §\ref{S:Wa} le sous-groupe $W_a$ de $W^G$ qui est le stabilisateur de la composante irréductible de $Y_a$ passant par le point $\infty_{Y_a}$. Soit $W_a^0$ le sous-groupe normal de $W_a$ engendré par les stabilisateurs dans $W_a$ des points de la composante irréductible de $Y_a$ passant par le point $\infty_{Y_a}$. Si $\bar{a}$ est un point géométrique de $\Ac_G$ localisé en $a$, l'ouvert $V_{\bar{a}}$, pointé par $\infty_{Y_{\bar{a}}}$ (avec les notations du § \ref{S:Wa}) est un revêtement étale de l'ouvert $U_{\bar{a}}$ pointé par $\infty$. La fibre de ce revêtement en $\infty$ est munie d'une action du groupe fondamental  $\pi_1(U_{\bar{a}},\infty)$. En utilisant le point $\infty_{Y_{\bar{a}}}$ de cet fibre, on obtient un morphisme 
$$\pi_1(U_{\bar{a}},\infty) \to W_a.$$
L'image de l'inertie $$\Ker(\pi_1(U_{\bar{a}},\infty)\to \pi_1(C\times_k k(\bar{a}),\infty)$$
est précisément le sous-groupe $W_a^0$.
On en déduit un morphisme
  
  \begin{equation}
    \label{eq:hom-Wa}
        \pi_1(C,\infty) \to W_a/W_a^0.
      \end{equation}
 On a une suite exacte 
$$1 \longrightarrow \pi_1(C,\infty) \longrightarrow \pi_1(C_0,\infty) \longrightarrow \Gal(k/\Fq) \longrightarrow 1$$
qui est canoniquement scindée par le point rationnel $\infty$. On notera que l'homomorphisme 
(\ref{eq:hom-Wa}) est stable par $\Gal(k/\Fq)$.

      \begin{definition}\label{def:GsetWs}
        Pour tout $s\in \Tc$, on introduit les groupes suivants
        \begin{itemize}
        \item $W_s$ est le sous-groupe de $W^G$ défini par
$$W_s=\{w\in W^G \mid w\cdot s=s\} \ ;$$
        \item $W_s^0$ est le sous-groupe de $W^G$ engendré par les réflexions simples associées aux racines $\al\in \Phi^G$ telles que $\al^\vee(s)=1$ ;
        \item $G_s$ est le groupe réductif connexe sur $k$ de dual de Langlands $Z_{\Gc}(s)^0$ le centralisateur connexe de $s$ dans $\Gc$.
        \end{itemize}
\end{definition}

\begin{remarque}
  Le sous-groupe $W_s^0$ est un sous-groupe normal de $W_s$. Il s'identifie d'ailleurs au groupe de Weyl de $G_s$.
\end{remarque}

\begin{definition}
  \label{def:endogeo}
On appelle \emph{donnée endoscopique géométrique} de $G$ la donnée d'un couple $(s,\rho)$ avec 
\begin{enumerate}
\item $s\in \Tc$ ;
\item $\rho :\pi_{1}(C,\infty)\rightarrow W_{s}/W_{s}^{0}$ est un  homomorphisme continu et fixe par  $\Gal(k/\Fq)$.
\end{enumerate}
\end{definition}

Soit  $(s,\rho)$ une donnée endoscopique géométrique de $G$ et soit 
$$\Ac_{s,\rho}$$
la partie de $\Ac$ formée des $a$ pour lesquels les conditions suivantes sont satisfaites
\begin{itemize}
\item $W_a \subset W_s$ ;
\item $W_a^0 \subset W_s^0$
\item l'homomorphisme composé 
$$\pi_1(C,\infty) \to W_a/W_a^0 \to  W_{s}/W_{s}^{0}$$
est égal à $\rho$.
\end{itemize}

\begin{proposition}\label{prop:unferme}
    Pour  une donnée endoscopique géométrique $(s,\rho)$ de $G$, la partie $\Ac_{s,\rho}$ est fermée dans $\Ac_G$ et stable par le Frobenius $\tau$. 
\end{proposition}

\begin{preuve}
  Soit $G_s$ le groupe réductif connexe sur $k$ associé à $s$ (cf. définition \ref{def:GsetWs}). Le groupe $W_s$, qui agit sur $T$, agit aussi sur l'ensemble $\Phi_T^{G_s}$ des racines de $T$ dans $G_s$. Le sous-groupe $W_s^0$ s'identifie au  groupe de Weyl de $(G_s,T)$. Soit $\Delta_{T}^{G_s}$ un ensemble de racines simples de $T$ dans $G_s$. Soit $W_s^{\mathrm{ext}}$ le sous-groupe de $W_s$ qui laisse  $\Delta_{T}^{G_s}$ fixe. On a alors une décomposition en produit semi-direct $W_s=W_s^0 \rtimes W_s^{\mathrm{ext}}$. L'ensemble $\Delta_{T}^{G_s}$ n'est pas unique mais on passe d'un choix à l'autre par l'action d'un unique élément de $W_s^0$ de sorte que changer $\Delta_{T}^{G_s}$ revient à conjuguer $W_s^{\mathrm{ext}}$ par un élément $W_s^0$. On complète la donnée de $\Delta_{T}^{G_s}$ en un épinglage $(B_s,T,\{X_\al\}_{\al\in \Delta_T^{G_s}})$ de $G_s$. Par ce choix, on identifie le groupe $W_s^{\mathrm{ext}}$ à un sous-groupe encore noté $W_s^{\mathrm{ext}}$  du groupe des automorphismes de $G_s$ qui laissent l'épinglage stable. Soit $\mathcal{W}_s^{\mathrm{ext}}$ le $W_s^{\mathrm{ext}}$-torseur au-dessus de $C$ muni d'une trivialisation de sa fibre en $\infty$ déduit de 
$$\rho : \pi_1(C,\infty) \to  W_{s}/W_{s}^{0}\simeq W_s^{\mathrm{ext}} $$
Comme $W_s^{\mathrm{ext}}$ laisse $W_s^0$ invariant, ce groupe agit sur $\car_{G_s}=\tgo//W_s^0$ et on peut introduire le $C$-schéma $\car_{G_{s,\rho}}$ qui est le produit contracté
$$\car_{G_{s,\rho}}=\car_{G_s}\times^{W_s^{\mathrm{ext}}}_k \mathcal{W}_s^{\mathrm{ext}}.$$
On introduit également sa variante 
$$\car_{G_{s,\rho},D}=\lc_D \times_C^{\mathbb{G}_{m,C}}\car_{G_{s,\rho}}.$$
Ce schéma fibré au-dessus de $C$ est muni d'un isomorphisme de sa fibre en $\infty$ sur $\car_{G_s}$ (déduite de la trivialisation canonique de $\lc_D$ et de la trivialisation fixée de $\mathcal{W}_s^{\mathrm{ext}}$). 
En imitant la construction du §\ref{S:esp-car}, on introduit le schéma quasi-projectif $\Ac_{G_{s,\rho}}^{G\text{-}\reg}$ qui classifie les couples $a=(h_a,t)$ où 
\begin{itemize}
\item $h_a\in H^0(C,\car_{G_{s,\rho},D})$ ;
\item $t\in \mathfrak{t}^{G-\textrm{reg}}$ avec $\chi_{G_s}(t)=h_a(\infty)$.
\end{itemize}
La proposition résulte alors de la proposition suivante, la stabilité par Frobenius étant évidente.

\end{preuve}

\begin{proposition}\label{prop:immersionNgo}
 Les notations sont celles de la démonstration de la proposition \ref{prop:unferme}. Le morphisme évident 
$$\chi^{G_{s,\rho}}_G : \Ac_{G_{s,\rho}}^{G\text{-}\reg} \to \Ac_G$$
est une immersion fermée d'image $\Ac_{G_{s,\rho}}$.
\end{proposition}

\begin{preuve}
 On renvoie le lecteur à  \cite{ngo2}, propositions 6.3.2 et 6.3.3.
\end{preuve}

En particulier, pour une donnée endoscopique géométrique $(s,\rho)$ dont l'homomorphisme $\rho$ est  l'homomorphisme trivial noté $1$, le fermé $\Ac_{s,1}$ est l'image du morphisme évident
$$\chi^{G_{s}}_G : \Ac_{G_{s}}^{G\text{-}\reg} \to \Ac_G$$
où  $\Ac_{G_{s}}^{G\text{-}\reg}$ est l'ouvert «$G$-régulier» de l'espace $\Ac_{G_s}$ pour le groupe réductif $G_s$ de la définition \ref{def:GsetWs} (cf. §\ref{S:esp-car}). Comme cet espace $\Ac_{G_s}$ n'interviendra plus  par la suite, on préfère réserver cet notation à l'ouvert $G$-régulier : on pose un peu abusivement
$$\Ac_{G_{s}}=\Ac_{G_{s}}^{G\text{-}\reg}.$$
On introduit  comme dans la définition \ref{def:partieell} un ouvert elliptique $\Ac_{G_s}^{\el}$ de $\Ac_{G_s}$, qui est non vide pour des raisons de dimension. 
 
\begin{proposition} \label{lem:Gs-ss} La partie localement fermée $\chi_G^{G_s}(\Ac_{G_s})\cap\Ac_G^{\el}$ est non vide si et seulement si on a l'égalité
$$Z_{\Gc}^0=Z_{\Gc_s}^0$$
entre les centres connexes des groupes duaux $\Gc$ et $\Gc_s$.

En cas d'égalité, on a, de plus,
$$\chi_G^{G_s}(\Ac_{G_s})\cap\Ac_G^{\el}=\chi_G^{G_s}(\Ac_{G_s}^{\el}).$$
  \end{proposition}

  \begin{remarque}
    Lorsque $G$ est semi-simple, le groupe dual $\Gc$ l'est aussi et la condition sur les centres est équivalente au fait que $\Gc_s$ donc $G_s$ est semi-simple.
  \end{remarque}

  \begin{preuve} Soit $H=G_s$, $a_H\in \Ac_H$ et $a=\chi^H_G(a)$. On a alors $W_a=W_{a_H}$ par une variante du lemme \ref{lem:WaM}. Aussi, dans les notations  on ne distinguera pas $a$ et $a_H$ et on omettra $\chi^H_G$. La démonstration résulte des deux lemmes  ci-dessous.

\begin{lemme} Pour tout $a\in \Ac_G^{\el}\cap \Ac_H$, on a $a\in \Ac_H^{\el}$ et $Z_{\Hc}^0=Z_{\Gc}^0$.
\end{lemme}

\begin{preuve}
Soit $M'$ le plus petit élément de $\lc^H(T)$ tel que $W_a\subset W^{M'}$. Soit $M\in \lc^G(T)$ dont le dual $\Mc$ s'identifie au sous-groupe de Levi $Z_{\Gc}(Z_{\Mc'}^0)$ de $\Gc$ (il s'agit du centralisateur dans $\Gc$ du tore central maximal de $\Mc'$). On a donc $W_a\subset W^{M'}\subset W^M$ et $a\in \Ac_M$ (par la proposition \ref{prop:AM-Wa}). Comme $a\in \Ac_G^{\el}$, on a $M=G$ et $ Z_{\Mc'}^0\subset Z_{\Gc}$ d'où l'égalité $Z_{\Hc}^0=Z_{\Mc'}^0= Z_{\Gc}$ et $a\in \Ac_{H}^{\el}$ (par une application de la proposition \ref{prop:AM-Wa} à $H$ et son sous-groupe de Levi $M'$). 
\end{preuve}

\begin{lemme} Si $Z_{\Hc}^0=Z_{\Gc}^0$, on a $\Ac_H^{\el}\subset \Ac_G^{\el}$.
\end{lemme}

\begin{preuve}
  Soit $a\in \Ac_H^{\el}$ et  $M$ le plus petit élément de $\lc^G(T)$ tel que $W_a\subset W^{M}$. Par la proposition \ref{prop:immersionNgo}, on a $W_a\subset W^H$. On a donc $W_a\subset W^{M'}$ où $M'\in \lc^H(T)$ est défini dualement par $\Mc'=\Mc\cap \Hc$. Comme $a$ est elliptique dans $H$, on a $M'=H$ donc $\Hc\subset \Mc$ donc $Z_{\Mc}^0\subset Z_{\Hc}^0$. La condition $Z_{\Hc}^0=Z_{\Gc}^0$ implique que $M=G$ c'est-à-dire $a\in \Ac_G^{\el}$.
\end{preuve}
  \end{preuve}

  \begin{proposition} \label{prop:finitude} Supposons $G$ semi-simple. L'ensemble des données endoscopiques géométriques $(s,\rho)$ de $G$ telles qu'on ait
\begin{equation}
  \label{eq:non-vide}
  \Ac_{s,\rho}\cap \Ac_G^{\el}\not=\emptyset
\end{equation}
est fini.
  \end{proposition}

  \begin{preuve}
Lorsque $s$ parcourt $\Tc$, il n'y a qu'un nombre fini de possibilités pour les groupes $\Gc_s$, $W_s$ et $W_s^0$. Comme le groupe $\pi_1(C,\infty)$ est topologiquement engendré par $2g$ générateurs, il n'y a qu'un nombre fini de morphismes $\rho$ possibles. Il suffit donc de prouver qu'il existe un sous-ensemble fini de $\Tc$ qui contient tous les $s$ qui appartiennent à une donnée endoscopique $(s,\rho)$ pour laquelle la condition (\ref{eq:non-vide}) est satisfaite. Soit donc un tel $s$ et  $Z_{\Gc_s}^{W_s}$ le groupe diagonalisable des points fixes sous $W_s$ du centre de $\Gc_s$. Comme ce groupe contient évidemment $s$, il suffit de voir qu'il est fini autrement dit que sa composante neutre notée $\hat{S}$ est triviale. Soit $\Mc$ le centralisateur dans $\Gc$ du tore $\hat{S}$. Son dual $M$ est  un sous-groupe de Levi de $G$. Comme $W_s$ centralise $\hat{S}$, on a $W_s\subset W^M$. Par l'hypothèse (\ref{eq:non-vide}), il existe $a\in \Ac_{s,\rho}\cap \Ac_G^{\el}$. On a donc $W_a\subset W_s\subset W^M$. Par la proposition \ref{prop:AM-Wa}, on a $a\in \Ac_M$. Comme $a$ est elliptique dans $G$, il vient $M=G$. Donc $\hat{S}$ est un tore central de $\Gc$ qui est supposé semi-simple. Ainsi $\hat{S}$ est trivial ce qu'il fallait démontrer.

  \end{preuve}

\end{paragr}

\begin{paragr}[Le théorème cohomologique.] --- On continue avec les mêmes notations et hypothèses. On suppose ici que $G$ est \emph{semi-simple}. Soit $s\in \Tc$ et $G_s$ le groupe réductif connexe sur $k$ de la définition \ref{def:GsetWs}.  On munit $G_s$ de l'action du Frobenius $\tau$ qui induit l'action triviale sur le groupe des caractères de $T$.

 Pour $k$-schéma $X$ muni d'une action du Frobenius $\tau$, on note $Perv(X)$ la catégorie des faisceaux pervers $F$ sur $X$ muni d'un isomorphisme $\tau^*F\simeq F$. Voici un théorème dû à Ngô (cf. \cite{ngo2} théorèmes 6.4.1 et 6.4.2).

 \begin{theoreme} \label{thm:Ngo}On suppose $G$ semi-simple. Soit $s\in \Tc$ tel que  $G_s$ est semi-simple.
 Soit $\chi_{G}^{G_s,\el} : \Ac_{G_s}^{\el} \hookrightarrow \Ac_G^{\el}$ la restriction de $\chi_{G}^{G_s}$ à l'ouvert $\Ac_{G_s}^{\el}$.

 On a alors l'égalité suivante dans le groupe de Grothendieck de $Perv(\Ac_{G_s}^{\el})$
   $$(\chi_{G}^{G_s,\el})^* {}^{\mathrm{p}}\mathcal{H}^{\bullet}(Rf_{G,\ast}^{\el}\Qlb)_s= {}^{\mathrm{p}}\mathcal{H}^{\bullet}(Rf_{G_s,\ast}^{\el}\Qlb)_1[-2r_s](-r_s)$$
où
\begin{itemize}
\item $r_s=\mathrm{codim}_{\Ac_G}(\Ac_{G_s})$ ;
\item ${}^{\mathrm{p}}\mathcal{H}^{\bullet}(Rf_{G,\ast}^{\el}\Qlb)_s$ est la $s$-partie de ${}^{\mathrm{p}}\mathcal{H}^{\bullet}(Rf_{G,\ast}^{\el}\Qlb)$ définie au théorème \ref{thm:sdecomposition} ;
\item ${}^{\mathrm{p}}\mathcal{H}^{\bullet}(Rf_{G_s,\ast}^{\el}\Qlb)_1$ est la $1$-partie de ${}^{\mathrm{p}}\mathcal{H}^{\bullet}(Rf_{G_s,\ast}^{\el}\Qlb)$ définie au théorème \ref{thm:sdecomposition} pour le groupe $G_s$ (avec $1$ l'élément neutre de $\Tc$).
\end{itemize}
 \end{theoreme}

 \begin{remarque}
   Comme $G$ et $G_s$ sont semi-simples, il en est de même de leurs duaux $\Gc$ et $\Gc_s=Z_{\Gc}(s)^0$. Par conséquent les centres de ces derniers sont finis. On a donc $Z_{\Gc}^0=Z_{\Gc_s}^0=\{1\}$ et $s$ est un élément de torsion dans $\Tc$.
 \end{remarque}

Soit $M\in \lc^G(T)$. On a défini au §\ref{S:OmegaM} un ouvert $\uc_M$ de $\Ac_G$. Cet ouvert est clairement $\tau$-stable. L'ouvert $\Ac_G^{\mathrm{bon}}$ de la définition \ref{def:Abon} est également $\tau$-stable. Pour tout $s\in \Tc$, le fermé $\bigcup_{\rho} \Ac_{G_{s,\rho}}$, où $\rho$ parcourt l'ensemble fini des morphismes continus et non triviaux 
$$\rho : \pi_1(C,\infty) \to W_s/W_s^0,$$
est $\tau$-stable.
Soit $\bar{s}\in (\Tc/Z_{\Mc}^0)_{\tors}$ et
\begin{equation}
  \label{eq:ucs}
  \uc_{\bar{s}}=\uc_M \cap \Ac_G^{\mathrm{bon}} \cap (\Ac_G-\bigcup_{(s,\rho)} \Ac_{s,\rho})
\end{equation}
où la réunion est prise sur l'ensemble \emph{fini} (cf. proposition \ref{prop:finitude}) des données endoscopiques géométriques $(s,\rho)$ qui vérifient
\begin{itemize}
\item $s\in \bar{s}Z_{\Mc}^0$ ;
\item $\Ac_{G}^{\el}\cap \Ac_{s,\rho}\not=\emptyset$ ;
\item $\rho$ est un homomorphisme  non trivial.
\end{itemize}
Alors $\uc_{\bar{s}}$ est un ouvert $\tau$-stable de $\Ac_G$. Soit 
$$\uc_{\bar{s}}^{\el}=\uc_{\bar{s}}\cap \Ac_G^{\el}.$$
Pour tout $s\in \bar{s}Z_{\Mc}^0$ tel que  $G_s$ est semi-simple, on a un diagramme commutatif, où les flèches notées $\chi$ sont des immersions fermées et les autres flèches sont des immersions ouvertes,
$$\xymatrix{\uc_{\bar{s}}^{\el}\cap \Ac_{G_s}^{\el}  \ar[r]^{i^{\bar{s}}_{G_s}} \ar[d]^{\chi^s_{\bar{s}}} & \Ac_{G_s}^{\el} \ar[d]^{\chi^{G_s,\el}_G}\\ \uc_{\bar{s}}^{\el}  \ar[r]^{i^{\bar{s}}_{G}} \ar[d]^{j_{\bar{s}}} & \Ac_{G}^{\el} \ar[d]^{j} \\  \uc_{\bar{s}}  \ar[r]^{i^{\bar{s}}_{M}} 
 & \uc_M}$$

\begin{theoreme}\label{thm:cohomologique}
  On suppose $G$ semi-simple et $\xi\in \ago_T$ en position générale. Pour tout $\bar{s}\in (\Tc/Z_{\Mc}^0)_{\tors}$, on a l'égalité suivante dans le groupe de Grothendieck de $Perv(\uc_{\bar{s}})$
$$(i_M^{\bar{s}})^*({}^{\mathrm{p}}\mathcal{H}^{\bullet}( Rf^{\xi,M\text{-}\mathrm{par}}_* \Qlb)_{\bar{s}})=\bigoplus_s (j_{\bar{s}})_{!*}(\chi_{\bar{s}}^s)_* (i_{G_s}^{\bar{s}})^*  {}^{\mathrm{p}}\mathcal{H}^{\bullet}(Rf_{G_s,\ast}^{\el}\Qlb)_1[-2r_s](-r_s)  $$
où 
\begin{itemize}
\item la somme porte sur les $s\in \bar{s}Z_{\Mc}^0$ tel que  $G_s$ est semi-simple ;
\item la $\bar{s}$-partie ${}^{\mathrm{p}}\mathcal{H}^{\bullet}( Rf^{\xi,M\text{-}\mathrm{par}}_* \Qlb)_{\bar{s}}$ est celle définie au théorème \ref{thm:sdecomposition} ;
\item les autres notations sont celles utilisées ci-dessus ou dans le théorème \ref{thm:Ngo}.
\end{itemize}
\end{theoreme}

\begin{preuve}
  D'après le théorème \ref{thm:support}, on a un isomorphisme canonique 
$$(i_M^{\bar{s}})^*({}^{\mathrm{p}}\mathcal{H}^{\bullet}( Rf^{\xi,M\text{-}\mathrm{par}}_* \Qlb)_{\bar{s}})\cong  (j_{\bar{s}})_{!*} (j_{\bar{s}})^* (i_M^{\bar{s}})^* ({}^{\mathrm{p}}\mathcal{H}^{\bullet}( Rf^{\xi,M\text{-}\mathrm{par}}_* \Qlb)_{\bar{s}}).$$
D'après le corollaire \ref{cor:dec-ell}, on a 
$$(j_{\bar{s}})^*(i_M^{\bar{s}})^* ({}^{\mathrm{p}}\mathcal{H}^{\bullet}( Rf^{\xi,M\text{-}\mathrm{par}}_* \Qlb)_{\bar{s}})= \bigoplus_s (i_G^{\bar{s}})^*{}^{\mathrm{p}}\mathcal{H}^{\bullet}(Rf_{G,\ast}^{\el}\Qlb)_s$$
où la somme est prise sur les $s\in \bar{s}Z_{\Mc}^0$ d'ordre fini. Pour tout tel $s$, d'après Ngô (cf. \cite{ngo2} §6.4),  ${}^{\mathrm{p}}\mathcal{H}^{\bullet}(Rf_{G,\ast}^{\el}\Qlb)_s$ est à support dans le fermé $(\bigcup_{\rho} \Ac_{G_{s,\rho}})\cap \Ac_G^{\el} \subset \Ac_G^{\el}$ où $\rho$ parcourt les homomorphismes continus de $\pi_1(C,\infty)$ dans $W_s/W_s^0$. Par conséquent, $(i_G^{\bar{s}})^*{}^{\mathrm{p}}\mathcal{H}^{\bullet}(Rf_{G,\ast}^{\el}\Qlb)_s$ est à support dans le fermé $\uc_{\bar{s}}^{\el}\cap \Ac_{G_s}^{\el}$. Si ce fermé n'est pas vide, l'inersection $\Ac_{G_s}\cap\Ac_G^{\el}$ n'est pas vide ce qui implique que $G_s$ est semi-simple (cf. proposition \ref{lem:Gs-ss} et la remarque qui suit). Pour un tel $s$, on a 
$$(i_G^{\bar{s}})^*{}^{\mathrm{p}}\mathcal{H}^{\bullet}(Rf_{G,\ast}^{\el}\Qlb)_s=(\chi_{\bar{s}}^s)_* (i_{G_s}^{\bar{s}})^*(\chi_G^{G_s,\el})^* {}^{\mathrm{p}}\mathcal{H}^{\bullet}(Rf_{G,\ast}^{\el}\Qlb)_s.$$
On utilise le théorème de Ngô (cf. théorème \ref{thm:Ngo}) pour conclure.
\end{preuve}

\end{paragr}

\section{Les $s$-intégrales orbitales pondérées}\label{sec:sIOP}
 
\begin{paragr} \label{S:sIOP-notations}On a fixé au §\ref{S:cbe} une courbe $C$ projective, lisse et connexe sur $k$ ainsi qu'un diviseur $D$ et un point $\infty\in C(k)$. On suppose désormais que $C$ provient par changement de base d'une courbe  $C_0$ sur $\Fq$. Soit $\tau$ l'automorphisme de Frobenius de $k$ donné par $x\mapsto x^q$. On note encore $\tau$ l'automorphisme $1\times \tau$ de $C=C_0\times_{\Fq} k$. On suppose que $D$ et $\infty$ sont fixes sous $\tau$. 

Soit $F$ et $F_0$ les corps des fonctions respectivement de $C$ et $C_0$. On a $F=F_0\otimes_{\Fq} k$ d'où une action de $\tau$ sur $F$ donnée par $1\otimes \tau$ dont $F_0$ est le sous-corps des points fixes. Pour tout point fermé $v_0$ de $C_0$ soit $F_{0,v_0}$ le complété de $F_0$ en  $v_0$ et $\oc_{0,v_0}$ son anneau d'entiers. Des notations analogues valent pour $F$. Le produit tensoriel complété   
$$F\hat{\otimes}_{F_0} F_{0,v_0}=\prod_{v\in V, \ v|v_0}F_v$$ 
se décompose en le produit des complétés  $F_v$ pour les points fermés $v$ de $C$ au-dessus de $v_0$. On munit ce produit tensoriel de l'action continue de $\tau$ donnée par $\tau\otimes 1$. Cette action admet $F_{0,v_0}$ comme sous-anneau de points fixes.  On en déduit une action de $\tau$ sur  l'anneau des adèles $\AAA$ de $F$ dont l'anneau des points fixes est l'anneau $\AAA_0$ des adèles de $F_0$. Cette action respecte le produit
$$\oc=\prod_v \oc_v$$
sur tous les points fermés de $C$. On définit de même $\oc_0$. On a d'ailleurs $\oc_0=\AAA_0 \cap \oc$.

On munit la donnée radicielle du couple $(G,T)$ du §\ref{S:G} de l'action triviale de $\tau$. Le couple $(G,T)$ provient alors par changement de base d'un couple analogue défini sur $\Fq$ où le tore  est déployé sur $\Fq$. Cela munit le couple $(G,T)$ d'une action de $\tau$. Tous les éléments de $\pc^G(T)$ ou de $\lc^G(T)$ c'est-à-dire les sous-groupes paraboliques ou de Levi de $G$ sont stables sous l'action de $\tau$. Pour tout schéma $S$ défini sur $k$ et toute $k$-algèbre $A$ tous deux munis d'une action de $\tau$, on note $S(A)^\tau$ le sous-ensemble des points fixes sous l'action de $\tau$. Si $S$ est, de plus, un schéma en groupes (et $\tau$ agit bien entendu comme un automorphisme de schémas en groupes), $S(A)^\tau$ est un groupe et on note $S(A)_\tau$ l'ensemble des classes de $\tau$-conjugaison (qui est donnée par $g\cdot x= gx\tau(g)^{-1}$) pointé par la classe de l'élément neutre. Si, de plus, $S$ est un schéma en groupes abéliens, $S(A)_\tau$ est un groupe abélien, à savoir le groupe des co-invariants.
\end{paragr}

\begin{paragr} Voici une autre variante d'un résultat déjà évoqué de Kottwitz (cf. \cite{Kot-isocrystal} lemme 2.2).
   
\begin{proposition}\label{prop:alaKot} Soit $U_0$ un ouvert de $C_0$ qui contient le point géométrique $\infty$. Soit $X_*$, $A$ et $B$ les foncteurs de la catégories des schémas en tores sur $U_{0}$  vers celles des groupes abéliens définis respectivement pour tout schéma en tores $\tc$ sur $U_{0}$ par
$$X_*(\tc)=X_*(\tc_{\infty}),$$
où le membre de droite désigne le groupe des cocaractères de la fibre de $\tc$ au point géométrique $\infty$,
$$A(\tc)=X_*(\tc_{\infty})_{\pi_1(U_0,\infty)},$$
où l'indice signifie que l'on prend les co-invariants sous le groupe fondamental de $U_0$ pointé par $\infty$, et
$$B(\tc)=(\tc(F) \back \tc(\AAA))_{\tau},$$
où l'indice $\tau$ désigne le groupe des co-invariants.

Soit $\varpi_\infty$ une uniformisante de $F_\infty$. Pour tout $\tc$ comme ci-dessus, le morphisme qui à $\la\in X_*(\tc_{\infty})$ associe l'image $\varpi_\infty^\la$ dans le groupe $(\tc(F) \back \tc(\AAA))_{\tau}$ ne dépend pas du choix de $\varpi_\infty$ et définit un  morphisme de foncteur
\begin{equation}
  \label{eq:X*B}
  X_* \to B.
\end{equation}
Il existe un unique isomorphisme de foncteur $A\to B$ qui, composé avec le morphisme canonique $X_* \to A$, donne le morphisme  (\ref{eq:X*B}).
\end{proposition}

\begin{preuve}
Elle repose sur les méthodes de Kottwitz (cf. \cite{Kot-isocrystal} preuve du lemme 2.2).
 Le foncteur $B$ a les deux propriétés fondamentales suivantes : il est exact à droite et pour un tore induit on a $B(\tc)=\ZZ$.  L'exactitude à droite se montre comme dans \emph{loc. cit.} \S1.9 : le point clef est que le groupe $H^1(F,\tc)$ est trivial pour tout $F$-tore (et de même si l'on remplace $F$ par un complété $F_v$). La seconde propriété résulte du fait suivant : soit $E_0$ une extension finie de $F$ et $E=E_0\otimes_{F_0} F$ on a $(E^\times \back \AAA_{E}^\times)_\tau \simeq \ZZ$ où  $\AAA_{E}$ est l'anneau des adèles de $E$. Par un lemme de Shapiro, quitte à remplacer $\tau$ par une puissance, on se ramène au cas où $E$ est un corps. Soit $\AAA_{E}^1$  le noyau de l'homomorphisme degré $\AAA_{E}^\times \to \ZZ$. Le quotient $E^\times\back \AAA_{E}^1/ \oc_E^\times$, où $\oc_E^\times$ est le produit des éléments inversibles des complétés des anneaux locaux, est le groupe des $k$-points de la jacobienne de la courbe projective, lisse et connexe sur $k$ dont $E$ est le corps des fonctions. D'après le théorème de l'isogénie de Lang, ces points sont dans l'image de  $1-\tau$. Il en est de même des points de $\oc_{E}^\times$ (cf. \cite{Kot-isocrystal} preuve de la proposition 2.3). Il est alors facile de voir que l'homomorphisme degré induit une bijection $(E^\times \back \AAA_{E}^\times)_\tau \simeq\ZZ$.

Le morphisme $X_*\to B$ ne dépend pas du choix de $\varpi_\infty$ par une variante du  théorème de l'isogénie de Lang. Il résulte de la preuve de \emph{loc. cit.} lemme 2.2 que le morphisme canonique 
$$\Hom(A,B)\to \Hom(X_*,B)$$
est un isomorphisme qui envoie isomorphisme sur isomorphisme. De plus, pour qu'un morphisme $\Hom(X_*,B)$ soit un isomorphisme, il faut et il suffit qu'il envoie un générateur de $X_*(\mathbb{G}_{m,U_{0}})$ sur un générateur de $(F^\times \back \AAA^\times_F)_\tau$. La proposition s'en déduit.
\end{preuve}

\end{paragr}

\begin{paragr}[Caractéristique $a_M$ et schéma en groupes $J_{a_M}$.] --- \label{S:JaM}L'automorphisme de Frobenius $\tau$ agit sur le schéma $\Ac_G$ (cf. §\ref{S:esp-car}) et respecte les sous-schémas fermés $\Ac_M$ pour $M\in \lc^G$. Soit $M\in \lc^G(T)$ un sous-groupe de Levi et  $a_M=(h_{a_M},t)\in \Ac_M(k)^\tau$. L'ouvert 
  $$U_{a_M}=h_{a_M}^{-1}(\car_{M,D}^{G\text{-}\reg})$$
est stable par $\tau$ et se descend donc en un ouvert noté $U_{a_M,0}$ de $C_0$.

Soit $\eps_M$ une section de Kostant du morphisme caractéristique  $\chi_M$ (cf. §\ref{S:Kostant}). Sur $\car_M$, on dispose du schéma en groupes commutatifs lisse (cf.   §\ref{S:centralisateurs-reg})
$$J^M=\eps_M^* I^M.$$
Par la construction du §\ref{S:actionGm}, on en déduit un schéma en groupes $J^M_D$ sur $\car_{M,D}$. Soit
$$J_{a_M}=h_{a_M}^* J^M_D.$$
C'est un schéma en groupes commutatifs lisse sur $C$ : c'est même un schéma en tores sur l'ouvert $U_{a_M}$ (et même sur un ouvert en général plus gros). On peut vérifier qu'il se descend en un schéma   en groupes commutatifs lisse sur $C_0$, qui est un schéma en tores sur $U_{a_M,0}$. 
\end{paragr}

\begin{paragr}[L'invariant $\inv_{a_M}$.] --- La suite exacte courte
$$1  \longrightarrow J_{a_M}(\AAA) \longrightarrow G(\AAA) \longrightarrow  J_{a_M}(\AAA)\back G(\AAA)  \longrightarrow 1$$
donne une suite exacte longue en cohomologie d'ensembles pointés

\begin{equation}
  \label{eq:suite-longue}
  1  \longrightarrow J_{a_M}(\AAA_0) \longrightarrow G(\AAA_0) \longrightarrow  (J_{a_M}(\AAA)\back G(\AAA))^\tau  \longrightarrow \Ker(J_{a_M}(\AAA)_\tau \to G(\AAA)_\tau)  \longrightarrow 1
\end{equation}
où le morphisme de cobord
\begin{equation}
  \label{eq:cobord}
  (J_{a_M}(\AAA)\back G(\AAA))^\tau  \longrightarrow J_{a_M}(\AAA)_\tau
\end{equation}
associe à un élément de $G(\AAA)$ qui relève un élément de  $(J_{a_M}(\AAA)\back G(\AAA))^\tau$ la classe de $\tau$-conjugaison de $g\tau(g)^{-1}$ dans $J_{a_M}(\AAA)$. Notons que l'exactitude pour le troisième terme signifie que  deux éléments de $(J_{a_M}(\AAA)\back G(\AAA))^\tau$ ont même image dans $J_{a_M}(\AAA)_\tau$ si et seulement s'ils diffèrent par une translation à droite par un élément de  $G(\AAA_0)$.

La proposition \ref{prop:alaKot}, appliquée à l'ouvert $U_{a_M,0}$ et au schéma en tores sur $U_{a_M,0}$ déduit de $J_{a_M}$ par descente donne un isomorphisme
$$(J_{a_M}(F)\back J_{a_M}(\AAA))_\tau \simeq X_*(J_{a_M,\infty})_{\pi_1(U_{a_M,0},\infty)}.$$

Les notations et les constructions du §\ref{S:cameral} valent aussi pour le groupe $M$. On a donc un revêtement fini, étale $V_{a_M}\to U_{a_M}$ de groupe de Galois $W^M$ et un point fermé $\infty_{Y_{a_M}}$ de $V_{a_M}$. Tous ces objets sont stables par $\tau$ et se descendent à des objets sur $\Fq$, notées par un indice $0$.

En raisonnant comme dans la démonstration de la proposition \ref{prop:factorisation-faisceau}, on obtient des identifications $X_*(J_{a_M,\infty})\simeq X_*(T)$ et 
$$X_*(J_{a_M,\infty})_{\pi_1(U_{a_M,0},\infty)}  \simeq X_*(T)_{W_{\infty}}$$
où $W_\infty$ est le stabilisateur dans $W^M$ de la composante connexe de $V_{a_M,0}$ qui contient le point  $\infty_{Y_{a_M},0}$. La suite exacte 
$$1  \longrightarrow \pi_1(U_{a_M},\infty)  \longrightarrow \pi_1(U_{a_M,0},\infty)  \longrightarrow \Gal(k/\Fq)  \longrightarrow 1$$
est scindé canoniquement par la donnée du point $\infty_0\in C_0(\Fq)$ déduit de $\infty$. On a donc un produit semi-direct
$$ \pi_1(U_{a_M,0},\infty)= \pi_1(U_{a_M},\infty)\rtimes \Gal(k/\Fq).$$
L'action de $\pi_1(U_{a_M,0},\infty)$ sur la fibre du revêtement  $V_{a_M,0}\to U_{a_M,0}$ au-dessus du point géométrique $\infty$ est triviale sur le facteur $\Gal(k/\Fq)$. On en déduit que le groupe $W_{\infty}$ n'est autre que le groupe $W_{a_M}\subset W^M$ (défini comme au §\ref{S:Wa} mais relativement au groupe $M$). En rassemblant les isomorphismes précédents, on obtient un isomorphisme
\begin{equation}
  \label{eq:iso-alaKot}
  (J_{a_M}(F)\back J_{a_M}(\AAA))_\tau \simeq  X_*(T)_{W_{a_M}}.
\end{equation}

\begin{definition}\label{def:invaM}
  Soit $a_M\in \Ac_M(k)^\tau$. Soit
$$\inv_{a_M} : (J_{a_M}(\AAA)\back G(\AAA))^\tau \to X_*(T)_{W_{a_M}}$$
le morphisme défini par composition par le  morphisme cobord (\ref{eq:cobord}) suivi du morphisme canonique $(J_{a_M}(\AAA))_\tau \to  (J_{a_M}(F)\back J_{a_M}(\AAA))_\tau$ suivi de l'isomorphisme (\ref{eq:iso-alaKot}). 
\end{definition}

Soit $M_{\scnx}$ le revêtement simplement connexe du groupe dérivé de $M$. Soit $T_{M_{\scnx}}$ l'image réciproque de $T$ dans  $M_{\scnx}$ : c'est un sous-tore maximal de $M_{\scnx}$.

\begin{proposition} \label{prop:image} Soit $a_M\in \Ac_M(k)^\tau$.
  Le morphisme $\inv_{a_M}$ est à valeurs dans l'image du morphisme canonique
$$X_*(T_{M_{\scnx}})_{W_{a_M}}\to X_*(T)_{W_{a_M}}.$$
\end{proposition}

Avant de donner la démonstration, énonçons quelques lemmes auxiliaires.

\begin{lemme}\label{lem:ker=} Soit $a_M\in \Ac_M(k)^\tau$.  On a  
\begin{equation}
    \label{eq:ker=}
    \Ker(J_{a_M}(\AAA)_\tau \to G(\AAA)_\tau)=   \Ker(J_{a_M}(\AAA)_\tau \to M(\AAA)_\tau).
  \end{equation}
\end{lemme}

\begin{preuve}   Seule l'inclusion $\subset$ n'est pas triviale. Par le lemme de Shapiro, quitte à remplacer $\tau$ par une de des puissances, on se ramène en un problème local en  un point fermé de $C$ fixe par $\tau$. Soit $g\in G(F_v)$ tel que $g\tau(g)^{-1}\in M(F_v)$. Soit $P\in\pc(M)$ et $B\in \pc(T)$ un sous-groupe de Borel tel que $B\subset P$. Un  couple $(P,M)$ qui contient $(B,T)$ est dit standard. Le couple $\Int(g^{-1})(P,M)$ est $\tau$-stable donc conjugué par un élément de $G(F_v)^\tau$ à un couple standard. Quitte à translater $g$ à droite  par un élément de $G(F_v)^\tau$, on est ramené au cas où $\Int(g^{-1})(P,M)$ est standard donc égal à $(P,M)$. Il s'ensuit qu'on peut supposer $g\in M(F_v)$. Si, de plus, $g\tau(g)^{-1}\in M(\oc_v)$, par une variante du théorème de Lang il existe $m\in M(\oc_v)$ tel que  $g\tau(g)^{-1}=m\tau(m)^{-1}$. Cela prouve le lemme.
\end{preuve}

On définit un schéma en groupes commutatifs sur $\mgo$ par
$$I^{M_{\scnx}}=\{(m,X)\in M_{\scnx}\times_k \mgo \mid \Ad(m)X=X \}.$$
On en déduit un schéma en groupes commutatifs lisse sur $\car_M$
$$J^{M_{\scnx}}=\eps_M^* I^{M_{\scnx}}$$
et une variante sur $\car_{M,D}$ noté $J^{M_{\scnx}}_D$. Soit
$$J_{a_M}^{M_{\scnx}}=h_{a_M}^* J^{M_{\scnx}}_D.$$
Notons que le morphisme évident $I^{M_{\scnx}}\to I^M$ induit des morphismes $J^{M_{\scnx}} \to J^M$ et  $J_{a_M}^{M_{\scnx}}\to J_{a_M}$. 

\begin{lemme} \label{lem:ker<}Soit $a_M\in \Ac_M(k)^\tau$. On a l'inclusion
$$\Ker(J_{a_M}(\AAA)_\tau \to M(\AAA)_\tau) \subset  \mathrm{Im}( J_{a_M}^{M_{\scnx}}(\AAA)_\tau \to  J_{a_M}(\AAA)_\tau  ),$$
où $\mathrm{Im}$ désigne l'image. 
\end{lemme}

\begin{preuve}
Comme dans la preuve du lemme \ref{lem:ker=}, on se ramène en un problème local en un point fermé $v$ de $C$ stable par $\tau$.  Par changement de base, on voit $J_{a_M}$ et $M$ comme des groupes sur $F_v$.  Soit $x\in J_{a_M}(F_v)$ et $m\in M(F_v)$ tel que $x=m\tau(m)^{-1}$.  On a alors 
$$\Coker(M_{\scnx} \to M) \simeq \Coker (J_{a_M}^{M_{\scnx}} \to J_{a_M}).$$ 
Comme   $J_{a_M}^{M_{\scnx}}$ est un tore  au-dessus de $F_v$, on sait que, pour un tel corps $F_v$, le groupe $H^1(F_v,J_{a_M}^{M_{\scnx}})$ est trivial. Il s'ensuit qu'il existe $m'\in M_{\scnx}(F_v)$ et $j\in J_{a_M}(F_v)$ tels que $m=jm'$. Donc, quitte à $\tau$-conjuguer $x$ par $j$, on peut supposer que  $m\in M_{\scnx}(F_v)$ et donc $x\in  J_{a_M}^{M_{\scnx}}(F_v)$. Pour presque tout $v$,  $J_{a_M}^{M_{\scnx}}$ est un tore au-dessus  $\oc_v$ et un raisonnement analogue montre que si  $x\in J_{a_M}(\oc_v)$ et $m\in M(\oc_v)$ vérifient $x=m\tau(m)^{-1}$ alors $x$ est $\tau$-conjugué par un élément de $J_{a_M}(\oc_v)$ à un élément de $J_{a_M}^{M_{\scnx}}(\oc_v)$.
\end{preuve}

\begin{preuve}(de la proposition \ref{prop:image}) Comme ci-dessus, on a une identification 
$$X_*(J_{a_M,\infty}^{M_{\scnx}})_{\pi_1(C_0,\infty)}\simeq X_*(T_{M_{\scnx}})_{W_{a_M}}.$$
Par la proposition \ref{prop:alaKot}, on obtient un diagramme commutatif où les flèches horizontales sont des isomorphismes
$$\xymatrix{(J_{a_M}^{M_{\scnx}}(F)\back J_{a_M}^{M_{\scnx}}(\AAA))_\tau   \ar[r] \ar[d] & X_*(T_{M_{\scnx}})_{W_{a_M}}  \ar[d] \\ (J_{a_M}(F) \back J_{a_M}(\AAA))_\tau \ar[r]  &  X_*(T)_{W_{a_M}} }$$
Soit $g\in (J_{a_M}(\AAA)\back G(\AAA))^\tau$. Par la suite exacte (\ref{eq:suite-longue}), la classe de $g\tau(g)^{-1}$ dans $(J_{a_M}(\AAA))_\tau$ appartient au noyau $\Ker(J_{a_M}(\AAA)_\tau \to G(\AAA)_\tau)$, donc  à l'image de $(J_{a_M}^{M_{\scnx}}(\AAA))_\tau$ par la combinaison des lemmes \ref{lem:ker=} et \ref{lem:ker<}. La proposition \ref{prop:image} est alors claire.
\end{preuve}

\end{paragr}

\begin{paragr}[Mesures de Haar.] --- \label{S:Haar} Le groupe $J_{a_M}(\AAA)_\tau$ est dénombrable : il apparaît en effet dans la suite exacte
$$J_{a_M}(F)_\tau \longrightarrow J_{a_M}(\AAA)_\tau \longrightarrow (J_{a_M}(F)\back J_{a_M}(\AAA))_\tau$$
dont le premier et le troisième terme sont dénombrables (cf. la proposition \ref{prop:alaKot}). On munit alors le noyau $\Ker(J_{a_M}(\AAA)_\tau \to G(\AAA)_\tau)$ de la mesure de comptage. Les groupes localement compacts $G(\AAA_0)$ et $J_{a_M}(\AAA_0)$ sont munis des mesures de Haar qui donnent le volume $1$ respectivement au sous-groupe compact $G(\oc_0)$ et au sous-groupe compact maximal de $J_{a_M}(\AAA_0)$. Par la suite (\ref{eq:suite-longue}), on en déduit une mesure sur   $(J_{a_M}(\AAA)\back G(\AAA))^\tau$ qui, par construction, est $G(\AAA_0)$-invariante pour l'action à droite de ce groupe. Plus précisément, l'ensemble $(J_{a_M}(\AAA)\back G(\AAA))^\tau$ est une réunion disjointe et dénombrable d'orbites sous $G(\AAA_0)$ ; l'orbite de $g\in (J_{a_M}(\AAA)\back G(\AAA))^\tau$ est isomorphe comme ensemble mesurable au quotient  $(g^{-1}J_{a_M}(\AAA_0)g) \back G(\AAA_0)$ muni de la mesure quotient (la conjugaison par $g$ induit une bijection de $g^{-1}J_{a_M}(\AAA_0)g$ sur $J_{a_M}(\AAA_0)$ ce qui détermine, par transport, une mesure sur  $g^{-1}J_{a_M}(\AAA_0)g$).
  \end{paragr}

  \begin{paragr}[Poids d'Arthur.] ---  \label{S:poids-Arthur}Pour tout $P\in \pc^G(M)$, on a la décomposition de Levi $P=MN_P$ et la décomposition d'Iwasawa $G(\AAA)=P(\AAA)G(\oc)$. On définit alors une application $H_P$ de $G(\AAA)$ dans $X_*(M)$ ainsi : pour tout $\la\in X^*(P)=X^*(M)$ et tout $g\in G(\AAA)$
$$\la(H_P(g))=-\deg(\la(p))$$
où $ p\in P(\AAA)$ est tel que $g\in pG(\oc)$. 

\begin{remarque}\label{rq:HM}
  La fonction $H_P$ restreinte à $M(\AAA)$ ne dépend plus de $P$ : on la note $H_M$. On a, de plus, pour tous $m\in M(\AAA)$ et $g\in G(\AAA)$ l'égalité $H_P(mg)=H_M(m)+H_p(g)$.
\end{remarque}

On rappelle (cf. §\ref{S:parab}) que $\ago_T$ est muni d'un produit scalaire $W^G$-invariant et que tous ses sous-espaces sont munis de la mesure euclidienne. La définition suivante dépend bien évidemment du choix de ce produit scalaire.

\begin{definition}
  \label{def:poids}
Pour tout $g\in G(\AAA)$, le poids d'Arthur $v_M^G(g)$ est le volume dans $\ago_M^G$ de l'enveloppe convexe des projections sur $\ago_M^G$ des points $-H_P(g)$ pour $P\in \pc^G(M)$. 
\end{definition}

\begin{remarques} On a $\ago_M=a_M^G\oplus \ago_G$ et la projection considérée ci-dessus est relative à cette décomposition. Le poids est par construction à droite par $G(\oc)$. Il est aussi invariant à gauche par $M(\AAA)$ en vertu de la remarque \ref{rq:HM}. Notons aussi qu'on a $v_G^G(g)=1$ pour tout $g\in G(\AAA)$.
  
\end{remarques}
   
  \end{paragr}

  \begin{paragr}[les $s$-intégrales orbitales pondérées.] --- \label{S:sIOP}On a $D=\sum_{v\in V} d_v v$ où la somme est prise sur les points fermés de $C$ et les entiers $d_v$ sont pairs, positifs et presque tous nuls (par hypothèse, cf. §\ref{S:cbe}, $D$ est effectif et pair). Soit $\varpi^{-D}=(\varpi_v^{-d_v})_v \in \AAA$ où la famille est prise sur tous les points fermés $v$ de $C$. Soit $\mathbf{1}_D$ la fonction caractéristique de $\varpi^{-D}\ggo(\oc)$. Pour tout $a_M\in \Ac_M(k)$, soit $a_{M,\eta}$ la restriction de $a_M$ au point générique de $C$ qu'on voit comme un point de $\car_M(F)$. Soit 
$$X_{a_M}=\eps_M(a_{M,\eta}) \in \mgo(F).$$
Notons que si $a_M$ est fixe par $\tau$, il en est de même de $X_{a_M}$. Notons aussi que la fibre de $J_{a_M}$ au-dessus de $\eta$ est le centralisateur de $X_{a_M}$ dans $M\times_k F$.

    \begin{definition}
      \label{def:sIOP}
Soit $M\in \lc^G(T)$ et $a_M\in \Ac_M(k)^\tau$. Soit $s\in \Tc^{W_{a_M}}$. La $s$-intégrale orbitale pondérée associée à $a_M$ est définie par

$$J_M^{G,s}(a_M)=q^{-\deg(D)|\Phi^G|/2}\int_{(J_{a_M}(\AAA)\back G(\AAA))^\tau} \bg s, \inv_{a_M}(g)\bd \mathbf{1}_D(\Ad(g^{-1})X_{a_M}) v_M^G(g) \, dg.$$
   \end{definition}
   
   \begin{remarques}
     On peut montrer que l'intégrande est à support compact dans $(J_{a_M}(\AAA)\back G(\AAA))^\tau$. L'intégrale est donc finie. On a choisi d'ajouter le facteur $q^{-\deg(D)|\Phi^G|/2}$ afin d'obtenir des formules locales uniformes.
   \end{remarques}

    \begin{proposition} \label{prop:inv-ZM} Avec les notations ci-dessus, l'application
$$s\in  \Tc^{W_{a_M}}\mapsto J_M^{G,s}(a_M)$$
est invariante par translation par $Z_{\Mc}$.
    \end{proposition}

    \begin{preuve}
      On a $W_{a_M}\subset W^M$ et donc $Z_{\Mc}\subset \Tc^{W_{a_M}}$. Le morphisme canonique 
$$\Hom_{\ZZ}(X_*(T)_{W_{a_M}},\CC^\times) \to \Hom_{\ZZ}(X_*(T_{M_{\scnx}})_{W_{a_M}},\CC^\times)$$
n'est autre que le morphisme canonique $\Tc^{W_{a_M}} \to \Tc^{W_{a_M}}/Z_{\Mc}$. La proposition résulte alors de la proposition \ref{prop:image}.
    \end{preuve}

  \end{paragr}

\section{Un premier calcul de trace de Frobenius} \label{sec:premiercalcul}

\begin{paragr}[Énoncé du théorème.] --- \label{S:scomptage} On reprend les notations et les hypothèses de la section \ref{sec:sIOP}. On suppose de plus que le groupe $G$ est \emph{semi-simple}. Soit $\xi\in \ago_T$ \emph{en position générale} au sens de la remarque \ref{rq:position}. Soit $M\in \lc^G(T)$. Sur l'ouvert $\Omega_M$ du §\ref{S:OmegaM} et  pour tout entier $i$ et tout $s\in (\Tc/Z_{\Mc}^0)_{\tors}$, on a introduit au théorème \ref{thm:sdecomposition} la $s$-partie de la cohomologie perverse de $Rf_*^{M\textrm{-}\parab}\Qlb$ notée ${}^{\mathrm{p}}\mathcal{H}^{i}( Rf^{\xi,M\text{-}\mathrm{par}}_* \Qlb)_{s}$. Soit $a\in \Ac_G(k)^\tau$. On suppose que $a$ est l'image par l'immersion fermée $\Ac_M\to \Ac_G$ d'un élément $a_M\in \Ac_M^{\el}(k)^\tau$. On a défini au §\ref{S:JaM} un schéma en groupes $J_{a_M}$ sur $C$. On définira plus tard un  réel positif
$$c(J_{a_M})$$
associé à la fibre générique $J_{a_M}$ (cf. l'égalité (\ref{eq:cJ}) dans la proposition \ref{prop:comptage} ci-dessous).

Soit ${}^{\mathrm{p}}\mathcal{H}^{i}( Rf^{\xi,M\text{-}\mathrm{par}}_* \Qlb)_{a,s}$ la fibre en $a$ de la $s$-partie. Pour tout endomorphisme $\sigma$ gradué (de degré $0$) de  ${}^{\mathrm{p}}\mathcal{H}^{\bullet}( Rf^{\xi,M\text{-}\mathrm{par}}_* \Qlb)_{a,s}$, soit
$$\trace(\sigma,{}^{\mathrm{p}}\mathcal{H}^{\bullet}( Rf^{\xi,M\text{-}\mathrm{par}}_* \Qlb)_{a,s})=\sum_{i\in \ZZ} (-1)^i \trace(\sigma,{}^{\mathrm{p}}\mathcal{H}^{i}( Rf^{\xi,M\text{-}\mathrm{par}}_* \Qlb)_{a,s}).$$

Le but de cette section est de démontrer le théorème suivant.

\begin{theoreme}\label{thm:scomptage}
  Soit $a_M\in \Ac_M^{\el}(k)^\tau$ et $a=\chi_G^M(a_M)\in \Ac_G$. Soit $s\in \Tc^{W_a}$ et $\bar{s}$ son image dans $\Tc/Z_{\Mc}^0$. Alors $\bar{s}$ est un élément de torsion, le groupe  $\Tc^{W_a}/Z_{\Mc}^0$ est fini et on l'égalité suivante
$$\trace(\tau^{-1}, {}^{\mathrm{p}}\mathcal{H}^{\bullet}( Rf^{\xi,M\text{-}\mathrm{par}}_* \Qlb)_{a,\bar{s}})=\frac{c(J_{a_M})}{\vol(\ago_M/X_*(M)) \cdot |\Tc^{W_a}/Z_{\Mc}^0|} \cdot q^{\deg(D)|\Phi^G|/2}   J_M^{G,s}(a_M).$$
\end{theoreme}

\begin{remarque}
  \label{rq:scomptage}
Puisque  $a_M$ est elliptique dans $\Ac_M$, on a l'inclusion $W_a\subset W^M$ (cf. proposition \ref{prop:AM-Wa}) et le groupe $(\Tc/Z_{\Mc}^0)^{W_a}$ est fini puisqu'il est le dual du groupe fini $X_*(T\cap M_{\der})_{W_a}$ (cf. lemme \ref{lem:finitude}). Cela explique pourquoi $\bar{s}$ est de torsion et le groupe  $\Tc^{W_a}/Z_{\Mc}^0$ est fini.
\end{remarque}

La démonstration du théorème \ref{thm:scomptage} sera donnée au §\ref{S:demo-scomptage}. 

\begin{corollaire}
Sous les hypothèses du théorème \ref{thm:scomptage}, la trace
  $$\trace(\tau^{-1}, {}^{\mathrm{p}}\mathcal{H}^{\bullet}( Rf^{\xi,M\text{-}\mathrm{par}}_* \Qlb)_{a,\bar{s}})$$
ne dépend que de l'image de $s$ dans $\Tc/Z_{\Mc}$.
\end{corollaire}

\begin{preuve}
  C'est une conséquence du théorème \ref{thm:scomptage} ci-dessus et de la proposition \ref{prop:inv-ZM}.
\end{preuve}

Dans le paragraphe suivant, on rassemble quelques résultats utiles à la démonstration du théorème \ref{thm:scomptage}.
\end{paragr}

 \begin{paragr} \label{S:comptage-aux}On continue avec les notations du paragraphe \ref{S:scomptage} précédent. On rappelle que $G$ est semi-simple. Soit $M\in \lc^G(T)$ et $X\in \mgo(F_0)$ un élément semi-simple et $G$-régulier. Soit $J$ son centralisateur dans $G$ : c'est un sous-$F$-tore de $M$ stable par $\tau$. On suppose que $J$ vérifie les deux propriétés suivantes :

   \begin{enumerate}
   \item le $F$-tore $J/Z_{M}$ est anisotrope ;
   \item le tore $J \times_F F_\infty$ est déployé sur $F_\infty$.
   \end{enumerate}

Soit
 $$K=G(\oc_0).$$

Soit $j\in J(\AAA)$ et 
 $$\wc=\{(gK,\delta)\in G(\AAA)/K \times J(F) \mid \delta j \tau(g) =g\}.$$
 Le groupe $J(F)$ opère sur $\wc$ par 
 $$\gamma\cdot (gK,\delta)=(\gamma gK,\gamma \delta \tau(\gamma)^{-1}).$$
Soit $\xi\in \ago_T$ (pour le moment, il n'est pas nécessaire de supposer $\xi$ en position générale) et $\wc^\xi$ le sous-ensemble de $\wc$ formé des couples $(gK,\delta)$ qui vérifient les deux conditions suivantes
\begin{equation}
  \label{eq:condition1}
  \Ad(g)^{-1}(X)\in \varpi^{-D}\ggo(\oc)
\end{equation}
et
\begin{equation}
  \label{eq:condition2}
\xi_M\in \mathrm{cvx}(-H_P(g))_{P\in \pc(M)}
\end{equation}
où $\xi_M$ est la projection orthogonale de $\xi$ sur $\ago_M$ et   $\mathrm{cvx}$ désigne l'enveloppe convexe.

Comme $J(F)$ centralise $X$ et que la fonction $H_P$ est invariante à gauche par $M(F)$ et \emph{a fortiori} par $J(F)$, le groupe $J(F)$ laisse le sous-ensemble $\wc^\xi$ globalement stable. Soit $[J(F)\back \wc^\xi]$ le groupoïde quotient. Le but de cette section est de donner une formule pour le cardinal du groupoïde  $[J(F)\back \wc^\xi]$ défini
\begin{equation}
  \label{eq:cardinal-def}
  | [J(F)\back \wc^\xi] |= \sum_{x\in J(F)\back \wc^\xi} |\Aut_{J(F)}(x)|^{-1}
\end{equation}
 où l'on somme sur un système de représentants des orbites de $J(F)$ dans $\wc^\xi$ les inverses des ordres des stabilisateurs (avec la convention que l'inverse de l'infini est $0$).

 Pour tout $g \in G(\AAA)$, soit $\mathbf{1}_{M,g}$ la fonction sur $\ago_M$ caractéristique de l'enveloppe convexe des points $-H_P(g)$ pour $P\in \pc(M)$. Soit 
$$w^\xi_M( g)=|\{\mu \in X_*(M) \ | \ \mathbf{1}_{M,g}(\xi_M+\mu)=1\}| \ ;$$
c'est un entier. 

\begin{lemme}\label{lem:inv-MA}
  La fonction $w^\xi_M$ est invariante à gauche par $M(\AAA)$.
\end{lemme}

\begin{preuve} On renvoie le lecteur à \cite{LFPI} lemme 11.7.1.
\end{preuve}

Soit 
$$(J(\AAA)\back G(\AAA))^{\tau,j}$$
la partie ouverte et fermée de $(J(\AAA)\back G(\AAA))^{\tau}$ formée des $g$ tels que la classe de $g\tau(g)^{-1}$ dans $(J(F)\back J(\AAA))_\tau$ soit égale à celle de $j$. On munit $J(\AAA_0)$ de la mesure de Haar qui donne le volume $1$ au sous-groupe compact maximal et $J(F_0)$ de la mesure de comptage. Rappelons que $G(\AAA_0)$ est muni de la mesure de Haar qui donne le volume $1$ à $K$. Comme au §\ref{S:Haar}, on en déduit une mesure sur $(J(\AAA)\back G(\AAA))^{\tau}$ et donc sur $(J(\AAA)\back G(\AAA))^{\tau,j}$.

\begin{proposition}
  \label{prop:comptage}
Soit   $J(\AAA_0)^1=J(\AAA_0)\cap \Ker(H_M)$ où $H_M$ est le morphisme défini à la remarque \ref{rq:HM} et
\begin{equation}
  \label{eq:cJ}
  c(J)=\vol(J(F_0)\back J(\AAA_0)^1)\cdot | \Ker(J(F)_\tau \to J(\AAA)_\tau) |.
\end{equation}
Alors $c(J)$ est fini et on a l'égalité
$$| [J(F)\back \wc^\xi] |= c(J)\cdot \int_{(J(\AAA)\back G(\AAA))^{\tau,j}} \mathbf{1}_D(\Ad(g^{-1})X)  \, w^\xi_M( g) \, dg$$
où l'intégrale est finie.
\end{proposition}

Avant de donner la démonstration énonçons deux corollaires.

\begin{corollaire} Rappelons que l'ensemble $\wc$ dépend du choix de $j\in J(\AAA)$. Le cardinal du  groupoïde $[J(F)\back \wc^\xi]$ ne dépend que de la classe de $j$ dans $(J(F)\back J(\AAA))_\tau$.
  \end{corollaire}

  \begin{preuve}
    L'ensemble $(J(\AAA)\back G(\AAA))^{\tau,j}$ ne dépend que de  la classe de $j$ dans $(J(F)\back J(\AAA))_\tau$. Le résultat est alors une conséquence immédiate de la proposition \ref{prop:comptage}.
  \end{preuve}

\begin{corollaire} \label{cor:comptage}Supposons, de plus, $\xi$ en position générale (cf. remarque \ref{rq:position}). On a alors
$$| [J(F)\back \wc^\xi] |=  \frac{c(J)}{\vol(\ago_M/X_*(M))}\cdot \int_{(J(\AAA)\back G(\AAA))^{\tau,j}} \mathbf{1}_D(\Ad(g^{-1})X)  \,  v_M^G(g)\, dg$$
\end{corollaire}

\begin{preuve} Pour  $\xi$ en position générale, on a l'égalité entre
$$  \int_{(J(\AAA)\back G(\AAA))^{\tau,j}} \mathbf{1}_D(\Ad(g^{-1})X)  \, w^\xi_M( g) \, dg$$
et 
$$\vol(\ago_M/X_*(M))^{-1} \cdot\int_{(J(\AAA)\back G(\AAA))^{\tau,j}} \mathbf{1}_D(\Ad(g^{-1})X)  \,  v_M^G(g)\, dg.$$
Celle-ci repose sur les formules qui relient les poids $w^\xi_M( g)$ et $v_M^G(g)$ et se démontre exactement comme le théorème 11.14.2 et son corollaire 11.14.3 de \cite{LFPI}. Le résultat se déduit alors de la proposition \ref{prop:comptage}.
\end{preuve}

\begin{preuve} Puisque $J$ est anisotrope modulo $Z_M$ (cf. assertion 1), le quotient  $J(F_0)\back J(\AAA_0)^1$ est compact donc de volume fini. Le groupe $\Ker^1(\Gamma_{F_0},J(F^{\sep}))$ du lemme \ref{lem:H1} est fini  pour tout $F_0$-tore $J$. Le lemme \ref{lem:H1} implique la finitude de $\Ker(J(F)_\tau \to J(\AAA)_\tau)$. Ainsi la constante $c(J)$ est finie.

Soit 
$$\psi : G(\AAA) \to G(\AAA)$$
l'application définie par $\psi(g)=g\tau(g)^{-1}j^{-1}$. L'application $(gK,\delta)\mapsto J(F)gK$ induit une bijection de $J(F)\back \wc$ sur $J(F)\back \psi^{-1}(J(F))/ K$. De plus, le centralisateur de $(gK,\delta)$ dans $J(F)$ est $J(F_0)\cap g K g^{-1}$. L'égalité (\ref{eq:cardinal-def}) se récrit donc
 $$  | [J(F)\back \wc^\xi] |=\sum_{g\in J(F)\back \psi^{-1}(J(F))/ K  } |J(F_0)\cap g K g^{-1}|^{-1}\, \mathbf{1}_D(\Ad(g^{-1})X)\,\mathbf{1}_{M,g}(\xi_M) .$$
Le groupe $G(\AAA_0)$ agit par translation à droite sur $G(\AAA)$ sans point fixe. L'ensemble $\psi^{-1}(F)$ est donc une réunion disjointe et dénombrable d'orbites sous $G(\AAA_0)$, toutes isomorphes à $G(\AAA_0)$. Par transport, on en déduit une mesure sur  $\psi^{-1}(F)$ invariante à droite par $G(\AAA_0)$. On munit $J(F)$ de la mesure de comptage et $J(F)\back \psi^{-1}(J(F))$ de la mesure quotient. On a alors
\begin{equation}
  \label{eq:|wc|}
  | [J(F)\back \wc^\xi] |= \int_{J(F)\back \psi^{-1}(J(F))}  \mathbf{1}_D(\Ad(g^{-1})X)\,\mathbf{1}_{M,g}(\xi_M) \, dg.
\end{equation}
Puisque $J_\infty$ est déployé (cf. assertion 2 ci-dessus), le morphisme $H_M$ se restreint en un morphisme surjectif de $J(\AAA_0)$ sur $X_*(M)$ de noyau $J(\AAA_0)^1$.  Pour tout $g\in G(\AAA)$, on a l'égalité (cf. \cite{LFPI} démonstration de la proposition 11.9.1) 
$$\int_{J(F_0)\back J(\AAA_0)  }  \mathbf{1}_{M,tg}(\xi_M) \, dt= \vol(J(F_0)\back J(\AAA_0)^1) \cdot w^\xi_M( g).$$
La suite exacte 
$$0 \longrightarrow J(F)\longrightarrow J(\AAA) \longrightarrow J(F)\back J(\AAA)\longrightarrow 0$$
donne la suite exacte en cohomologie
$$0\longrightarrow J(F_0)\back J(\AAA_0) \longrightarrow (J(F)\back J(\AAA))^\tau  \longrightarrow \Ker(J(F)_\tau \to J(\AAA)_\tau)    \longrightarrow 0 .$$
On munit alors $(J(F)\back J(\AAA))^\tau$ de la mesure dont le quotient par celle de $J(F_0)\back J(\AAA_0)$ donne la mesure de comptage sur $\Ker(J(F)_\tau \to J(\AAA)_\tau)$. En utilisant la suite exacte ci-dessus et l'invariance de $w^\xi_M$ à droite par $M(\AAA)$, on obtient
\begin{equation}
  \label{eq:wc}
  \int_{(J(F)\back J(\AAA))^\tau  }  \mathbf{1}_{M,tg}(\xi_M) \, dt= c(J) \cdot w^\xi_M( g).
\end{equation}

Le groupe $(J(F)\back J(\AAA))^\tau$ agit sans points fixes sur $J(F)\back \psi^{-1}(J(F))$. L'application qui, à $g\in  \psi^{-1}(J(F))$, associe sa classe dans $J(\AAA)\back G(\AAA)$ passe au quotient en une application surjective $  J(F)\back \psi^{-1}(J(F)) \to (J(\AAA)\back G(\AAA))^{\tau,j}$ dont les fibres sont les orbites sous $(J(F)\back J(\AAA))^\tau$. On a donc 
$$(J(F)\back J(\AAA))^\tau\back J(F)\back \psi^{-1}(J(F))\simeq (J(\AAA)\back G(\AAA))^{\tau,j}$$
et cette bijection est compatible à nos choix de mesures. En combinant (\ref{eq:|wc|}) et (\ref{eq:wc}), on obtient la proposition. L'intégrale est finie car  l'intégrande est à support compact.

\end{preuve}

\begin{lemme} \label{lem:H1}Soit $F^{\sep}$ une clôture séparable de $F$ et $\Gamma_{F_0}$ le groupe de Galois de $F^{\sep}$ sur $F_0$. Soit 
$$\Ker^1(\Gamma_{F_0},J(F^{\sep}))=\Ker(H^1(\Gamma_{F_0},J(F^{\sep}))\to  H^1(\Gamma_{F_0},J(\AAA_{0}\otimes_{F_0} F^{\sep}))$$
On a une identification naturelle
$$\Ker^1(\Gamma_{F_0},J(F^{\sep})) =  \Ker(J(F)_\tau \to J(\AAA)_\tau)$$
  \end{lemme}

  \begin{preuve}
    On a la suite exacte 
$$0 \longrightarrow \Gamma_{F}  \longrightarrow \Gamma_{F_0}  \longrightarrow \Gal(F/F_0)  \longrightarrow 0$$
où  $\Gamma_{F}$  est le groupe de Galois de $F^{\sep}$ sur $F$. Le groupe de Galois $\Gal(F/F_0)$ de $F$ sur $F_0$ s'identifie au groupe de Galois de $k$ sur $\Fq$ donc à $\hat{\ZZ}$. Soit $W_{F_0}$ l'image réciproque de $\ZZ$ dans $\Gamma_{F_0}$. On munit $W_{F_0}$  de la topologie qui fait de $ \Gamma_{F}$ un sous-groupe ouvert de $W_{F_0}$. On a donc une suite exacte
$$0 \longrightarrow \Gamma_{F}  \longrightarrow W_{F_0}  \longrightarrow \ZZ  \longrightarrow 0.$$
Cette suite exacte donne le diagramme commutatif suivant à lignes exactes d'inflation-restriction 

$$\xymatrix{ 0 \ar[r] &  H^1(\ZZ, J(F)) \ar[r] \ar[d] &  H^1(W_{F_0}, J(F^{\sep}))  \ar[r] \ar[d] & H^1(\Gamma_{F}, J(F^{\sep})) \ar[d] \\
 0 \ar[r] &  H^1(\ZZ, J(\AAA)) \ar[r]  &  H^1(W_{F_0}, J(\AAA\otimes_{F}F^{\sep}))  \ar[r]  & H^1(\Gamma_{F}, J(\AAA\otimes_{F}F^{\sep})) }$$
On sait bien que pour un tel corps $F$ les termes de droite sont nuls. Ceux de gauche s'identifient, respectivement de haut en bas, à $J(F)_\tau$ et $J(\AAA)_\tau$. Il s'ensuit qu'on a l'identification
\begin{equation}
  \label{eq:identif1}
\Ker(J(F)_\tau \to J(\AAA)_\tau)= \Ker( H^1(W_{F_0}, J(F^{\sep}))\to  H^1(W_{F_0}, J(\AAA\otimes_{F}F^{\sep})).
\end{equation}
Puisque $J$ est stable par $\tau$, on peut le voir comme un tore sur $F_0$. Soit $E_0$ une extension finie de $F_0$ dans $F^{\sep}$ qui déploie $J$. La suite exacte
 $$0 \longrightarrow W_{E_0}  \longrightarrow W_{F_0}  \longrightarrow \Gal(E_0/F_0) \longrightarrow 0$$
donne comme précédemment le  diagramme commutatif suivant à lignes exactes d'inflation-restriction 
$$\xymatrix{ 0 \ar[r] &  H^1(\Gal(E_0/F_0)   , J(E_0)) \ar[r] \ar[d] &  H^1(W_{F_0}, J(F^{\sep}))  \ar[r] \ar[d] & H^1(W_{E_0}, J(F^{\sep})) \ar[d] \\
 0 \ar[r] &  H^1(\Gal(E_0/F_0), J(\AAA_{0}\otimes_{F_0} E_0)) \ar[r]  &  H^1(W_{F_0}, J(\AAA\otimes_{F}F^{\sep}))  \ar[r]  & H^1(W_{E_0}, J(\AAA\otimes_{F}F^{\sep})) }$$
Montrons que la flèche verticale de droite est injective. Puisque $J$ est déployé sur $E_0$, il suffit de le prouver pour le tore $\mathrm{G}_{m,E_0}$. Par (\ref{eq:identif1}), on est ramené à démontrer que la flèche $F^\times_\tau \to \AAA_\tau$ est injective. Mais tout $x\in F^\times$ d'image triviale dans $\AAA_\tau$  est une unité en toute place,  donc un tel $x$ appartient à $k^\times$ et est de la forme $y\tau(y)^{-1}$ pour un élément $y \in k^\times$.

Cette injectivité étant acquise, on en déduit l'égalité entre 
$$  \Ker(H^1(\Gal(E_0/F_0)   , J(E_0))\to  H^1(\Gal(E_0/F_0), J(\AAA_{0}\otimes_{F_0} E_0))) $$
et 
$$\Ker( H^1(W_{F_0}, J(F^{\sep}))\to  H^1(W_{F_0}, J(\AAA\otimes_{F}F^{\sep}))$$
ce qui, combiné avec (\ref{eq:identif1}), donne le lemme.

  \end{preuve}

\end{paragr}

\begin{paragr}[Démonstration du théorème \ref{thm:scomptage}.] --- \label{S:demo-scomptage}Les notations sont celles du §\ref{S:scomptage}. Le groupe fini $\Gamma(\Omega_M,\pi_0(\Jc^1))$ des sections globales sur $\Omega_M$ du faisceau $\pi_0(\Jc^1)$ agit sur ${}^{\mathrm{p}}\mathcal{H}^{\bullet}( Rf^{\xi,M\text{-}\mathrm{par}}_* \Qlb)$ (cf. \ref{S:sdecomposition}). Comme ce groupe fini reçoit le groupe $X_*(T\cap M_{\der})$, on a également une action de ce dernier sur la cohomologie perverse. Soit $a_M\in \Ac_M^{\el}(k)^\tau$ et $a$ son image dans $\Ac_G(k)^\tau$. La décomposition (\ref{eq:s-decomposition}) du théorème \ref{thm:sdecomposition} donne l'égalité suivante de traces pour tout $\la\in X_*(T\cap M_{\der})$
$$\trace((\la\tau)^{-1}, {}^{\mathrm{p}}\mathcal{H}^{\bullet}( Rf^{\xi,M\text{-}\mathrm{par}}_* \Qlb)_a) =\sum_{s\in (\Tc/Z_{\Mc}^0)_{\tors}} \bg s, \lambda \bd^{-1} \trace(\tau^{-1}, {}^{\mathrm{p}}\mathcal{H}^{\bullet}( Rf^{\xi,M\text{-}\mathrm{par}}_* \Qlb)_{s,a}) $$
  où l'accouplement $\bg \cdot, \cdot \bd$ est l'accouplement canonique entre $\Tc/Z_{\Mc}^0$ et $X_*(T\cap M_{\der})$.

L'assertion 1 du théorème \ref{thm:sdecomposition} montre que, dans la somme ci-dessus, le terme associé à $s$ est nul sauf si $s$ appartient au groupe fini $(\Tc/Z_{\Mc}^0)^{W_a}$ (pour un commentaire sur ce groupe, cf. la remarque \ref{rq:scomptage}). En conséquence, la trace de $(\la\tau)^{-1}$ sur la cohomologie perverse ne dépend que de la classe de $\la$ dans le groupe des co-invariants $X_*(T\cap M_{\der})_{W_a}$. Par inversion de Fourier sur le groupe fini et commutatif $(\Tc/Z_{\Mc}^0)^{W_a}$ de dual $X_*(T\cap M_{\der})_{W_a}$, on obtient la formule suivante pour tout $s \in (\Tc/Z_{\Mc}^0)^{W_a}$
\begin{eqnarray}\label{eq:tracetau}
  \trace(\tau^{-1}, {}^{\mathrm{p}}\mathcal{H}^{\bullet}( Rf^{\xi,M\text{-}\mathrm{par}}_* \Qlb)_{s,a})= \\
| (\Tc/Z_{\Mc}^0)^{W_a} |^{-1} \sum_{\la \in X_*(T\cap M_{\der})_{W_a}} \bg s, \lambda \bd \, \trace((\la\tau)^{-1}, {}^{\mathrm{p}}\mathcal{H}^{\bullet}( Rf^{\xi,M\text{-}\mathrm{par}}_* \Qlb)_a).\nonumber
\end{eqnarray}

Soit $\la\in  X_*(T\cap M_{\der})$. On a les égalités suivantes 
\begin{eqnarray}\label{eq:lefschetz}
  \trace((\la\tau)^{-1}, {}^{\mathrm{p}}\mathcal{H}^{\bullet}( Rf^{\xi,M\text{-}\mathrm{par}}_* \Qlb)_a)\nonumber &=&  \trace((\la\tau)^{-1}, H^\bullet(\mc_a^\xi,\Qlb))\\  &=& |\mc_a^\xi(k)^{\lambda \tau}|.
\end{eqnarray}
La première égalité relie la trace de $(\la\tau)^{-1}$ sur la cohomologie perverse à la trace (alternée) de $(\la\tau)^{-1}$ sur la cohomologie de la fibre en $a$ de $\mc^\xi$ (notée $\mc^\xi_a$). Elle se démontre comme dans \cite{laumon-Ngo} lemme 3.10.1. La seconde égalité résulte d'une version champêtre de la formule des traces de Grothendieck-Lefschetz. Elle relie cette trace au cardinal, noté entre deux barres, du groupoïde des $k$-points de $\mc_a^\xi$ qui sont fixes sous $\lambda\tau$. Nous allons maintenant décrire ce groupoïde. On reprend les notations de la section  \ref{sec:sIOP}. On y a défini un schéma en groupes $J_{a_M}$ sur $C$, stable par $\tau$ (cf. §\ref{S:JaM}), ainsi qu'un élément $X_{a_M}\in \mgo(F_0)$ (cf. §\ref{S:sIOP}). On a le résultat suivant.

\begin{lemme}\label{lem:groupoide1}
  Le groupoïde $\mc_a^\xi(k)$ des $k$-points de la fibre en $a$ de $\mc^\xi$ est équivalent au groupoïde quotient de l'ensemble des $g\in G(\AAA)\back G(\oc)$ tels que
  \begin{equation}
    \label{eq:cond1}
     \Ad(g)^{-1}(X)\in \varpi^{-D}\ggo(\oc)
  \end{equation}
et
\begin{equation}
  \label{eq:cond2}
\xi_M\in \mathrm{cvx}(-H_P(g))_{P\in \pc(M)}
\end{equation}
par l'action par translation à gauche du groupe $J_{a_M}(F)$.
\end{lemme}

\begin{preuve}
  C'est une conséquence facile de la proposition 7.6.3 de \cite{LFPI}.
\end{preuve}

Le groupe $X_*(T\cap M_{\der})$ agit sur la fibre $\mc_a^\xi$ via le morphisme $X_*(T\cap M_{\der}) \to \Jc^1_a$ considéré au §\ref{S:X1} (cf. aussi §\ref{S:X-J}) et l'action de $\Jc^1_a$ sur $\mc_a^\xi$ (cf. \ref{S:actionJ1}). Soit $\varpi_\infty$ une uniformisante de $\oc_\infty$. Rappelons que pour construire ce morphisme, on identifie canoniquement $J_a$ à $T$ sur $\oc_\infty$ (cf. proposition \ref{prop:Ja=T-isocan}) et l'élément $\varpi_\infty^\la\in T(\oc_\infty)\simeq J_a(\oc_\infty)$ s'interprète comme une donnée de recollement qui définit un $J_a$-torseur sur $C$.  En fait, on a aussi une identification canonique de $J_{a_M}$ à $T$ sur $\oc_\infty$ et on note $j_\lambda\in J_{a_M}(\oc_\infty)$ l'élément correspondant à $\varpi_\infty^\la$ dans cette identification. 

\begin{lemme}\label{lem:groupoide2}
   Le groupoïde $\mc_a^\xi(k)^{\lambda \tau}$ est équivalent au groupoïde quotient de l'ensemble des couples $(g,\delta)\in G(\AAA)\back G(\oc)\times J_{a_M}(F)$ qui vérifient les conditions (\ref{eq:cond1}) et  (\ref{eq:cond2}) du lemme \ref{lem:groupoide1} pour $g$ et la condition
$$j_\lambda \delta\tau(g)=g$$
par l'action de $J_{a_M}(F)$ donnée par la translation à gauche sur $g$ et la $\tau$-conjugaison sur $\delta$. 
\end{lemme}

\begin{preuve}
Elle résulte de la description du lemme \ref{lem:groupoide1}, de la définition des points fixes d'un groupoïde sous un automorphisme (cf. par exemple \cite{LFPI} §11.2) et de la description adélique de l'action de $\Jc_a^1$ (laissée au lecteur).
\end{preuve}

Voici une variante du lemme \ref{lem:groupoide2} qui va nous permettre d'utiliser les résultats du §\ref{S:comptage-aux}.

\begin{lemme}\label{lem:groupoide3}
  Le groupoïde $\mc_a^\xi(k)^{\lambda \tau}$ est équivalent au groupoïde quotient $[J_{a_M}(F)\back \wc^\xi]$ défini au §\ref{S:comptage-aux} relatif aux données $X_{a_M}$ et $j_\lambda$.
\end{lemme}

\begin{preuve}
La seule difficulté est de voir que le morphisme évident de  $[J_{a_M}(F)\back \wc^\xi]$ vers le groupoïde quotient considéré au lemme \ref{lem:groupoide2} est essentiellement surjectif. On utilise pour cela la variante du théorème de Lang qui montre que $1-\tau$ est une application surjective de $G(\oc)$ dans $G(\oc)$.
\end{preuve}

On en déduit le lemme suivant.

\begin{lemme} \label{lem:tracetaula} Pour tout $\la\in  X_*(T\cap M_{\der})$, on a 
  
$$ \trace((\la\tau)^{-1}, {}^{\mathrm{p}}\mathcal{H}^{\bullet}( Rf^{\xi,M\text{-}\mathrm{par}}_* \Qlb)_a)=     \frac{c(J_{a_M})}{\vol(\ago_M/X_*(M))}\cdot \int_{(J_{a_M}(\AAA)\back G(\AAA))^{\tau,j_\lambda}} \mathbf{1}_D(\Ad(g^{-1})X)  \,  v_M^G(g)\, dg.$$

\end{lemme}

\begin{preuve}
  L'égalité cherchée est une simple combinaison de l'égalité (\ref{eq:lefschetz}) et du corollaire \ref{cor:comptage} (rappelons qu'on suppose que $\xi$ est en position générale). 
\end{preuve}

\begin{lemme}  \label{lem:tracetau}Pour tout $s \in (\Tc/Z_{\Mc}^0)^{W_a}$, on a l'égalité entre  
$$ \trace(\tau^{-1}, {}^{\mathrm{p}}\mathcal{H}^{\bullet}( Rf^{\xi,M\text{-}\mathrm{par}}_* \Qlb)_{s,a})$$
et
$$\frac{c(J_{a_M})\cdot |\Ker(  X_*(T\cap M_{\der})_{W_a}\to X_*(T)_{W_a})|}{\vol(\ago_M/X_*(M)) \cdot | (\Tc/Z_{\Mc}^0)^{W_a} |}  \cdot q^{\deg(D)|\Phi^G|/2}  \cdot  J_M^{G,s}(a_M)$$
  \end{lemme}

  \begin{preuve} Soit  $s \in (\Tc/Z_{\Mc}^0)^{W_a}$. L'égalité (\ref{eq:tracetau}) combinée au lemme \ref{lem:tracetaula} donne la formule suivante pour  $\trace(\tau^{-1}, {}^{\mathrm{p}}\mathcal{H}^{\bullet}( Rf^{\xi,M\text{-}\mathrm{par}}_* \Qlb)_{s,a})$ 

$$  \frac{c(J_{a_M})}{\vol(\ago_M/X_*(M)) \cdot | (\Tc/Z_{\Mc}^0)^{W_a} |} \sum_{\la \in X_*(T\cap M_{\der})_{W_a}} \bg s, \lambda \bd \, \int_{(J_{a_M}(\AAA)\back G(\AAA))^{\tau,j_\lambda}} \mathbf{1}_D(\Ad(g^{-1})X)  \,  v_M^G(g)\, dg.  $$
Le quotient $(J_{a_M}(\AAA)\back G(\AAA))^{\tau,j_\lambda}$ n'est autre que l'ensemble des $g\in (J_{a_M}(\AAA)\back G(\AAA))^{\tau}$ tels que $\inv_{a_M}(g)$ est égal à l'image de $\la$ dans $X_*(T)_{W_a}$. Par la proposition \ref{prop:image}, l'image de l'invariant $\inv_{a_M}$ est incluse dans l'image du morphisme $X_*(T\cap M_{\der})\to X_*(T)_{W_a}$. On en déduit que la somme 
$$\sum_{\la \in X_*(T\cap M_{\der})_{W_a}} \bg s, \lambda \bd \, \int_{(J_{a_M}(\AAA)\back G(\AAA))^{\tau,j_\lambda}} \mathbf{1}_D(\Ad(g^{-1})X)  \,  v_M^G(g)\, dg$$
est égale à
$$ |\Ker(  X_*(T\cap M_{\der})_{W_a}\to X_*(T)_{W_a})| \cdot  q^{\deg(D)|\Phi^G|/2}\cdot   J_M^{G,s}(a_M)$$
ce qui donne le résultat.
  \end{preuve}

Le lemme \ref{lem:tracetau} donne le théorème en vertu du lemme \ref{lem:ordres} suivant.

\begin{lemme}\label{lem:ordres} On a l'égalité

$$  \frac{|\Ker(  X_*(T\cap M_{\der})_{W_a}\to X_*(T)_{W_a})|}{| (\Tc/Z_{\Mc}^0)^{W_a} |}=| \Tc^{W_a}/Z_{\Mc}^0 |^{-1}$$

\end{lemme}

\begin{preuve}
La suite exacte
$$0 \longrightarrow  X_*(T\cap M_{\der}) \longrightarrow    X_*(T) \longrightarrow   X_*(M) \longrightarrow   0$$ 
donne la suite exacte de co-invariants (rappelons que $W_a\subset W^M$ agit trivialement sur $X_*(M)$)
$$0 \longrightarrow I   \longrightarrow X_*(T\cap M_{\der})_{W_a} \longrightarrow    X_*(T)_{W_a} \longrightarrow   X_*(M) \longrightarrow   0$$ 
avec $I=\Ker(  X_*(T\cap M_{\der})_{W_a}\to X_*(T)_{W_a})$. En appliquant le foncteur exact $\Hom_\ZZ(\cdot,\CC^\times)$ à cette suite, on obtient la suite exacte de groupes 
$$1 \longrightarrow Z_{\Mc}^0  \longrightarrow \Tc^{W_a} \longrightarrow   (\Tc/Z_{\Mc}^0)^{W_a}  \longrightarrow  \Hom_{\ZZ}(I,\CC^\times) \longrightarrow   0.$$
On en déduit que l'ordre de $\Hom_{\ZZ}(I,\CC^\times)$ est égal au rapport entre les ordres des groupes finis  $(\Tc/Z_{\Mc}^0)^{W_a}$ et $\Tc^{W_a}/ Z_{\Mc}^0$ ce qui donne   l'égalité cherchée puisque cet ordre est aussi celui de $I$.
\end{preuve}
\end{paragr}

\section{Avatars stables des intégrales orbitales pondérées} \label{sec:avatarstable}

\begin{paragr}\label{S:avatar1} Dans cette section, on reprend dans notre contexte le formalisme introduit par Arthur pour stabiliser les intégrales orbitales pondérées. Comme dans tout l'article, $G$ est un groupe réductif, connexe défini sur $k$ et $T$ un sous-tore maximal. On a une bijection naturelle $M\mapsto \Mc$ (cf. \ref{S:dual-Langlands}) entre $\lc^G(T)$ et les sous-groupes de Levi de $\Gc$ qui contiennent $\Tc$.  Pour tout $s\in \Tc$ et tout sous-groupe $\Hc$ de $\Gc$, on note $\Hc_s$ la composante neutre du centralisateur de $s$ dans $\Hc$ et $H_s$ son dual sur $k$.
\end{paragr}

\begin{paragr} \label{S:avatar2} Soit $s_M \in \Tc$ et $\Mc'=\Mc_{s_M}$. On a $\Mc'=\Mc$ si et seulement si $s_M$ est central dans $\Mc$. Soit $\ec^G_M(s_M)$ l'ensemble (fini) des sous-groupes $\Hc$ de $\Gc$ pour lesquels il existe un entier $n\geq 1$ et une famille finie $\underline{s}=(s_1,\ldots,s_n)$ d'éléments de $\Tc$ qui vérifie 
\begin{itemize}
\item $s_1\in s_M Z_{\Mc}$ et $s_i\in Z_{\Mc'}$ pour $i\geq 2$ ;
\item $\Hc$ est la composante neutre du groupe
$$\Gc_{s_1}\cap \ldots \cap \Gc_{s_n}.$$
\end{itemize}

\begin{remarques}\label{rq:scentral} Tous les groupes de $\ec^G_M(s_M)$ sont réductifs, connexes et de même rang, à savoir la dimension de $\Tc$.
 
  Pour tout $\Hc\in \ec^G_M(s_M)$, on a $\Mc\cap \Hc=\Mc_{s_M}$. En particulier, $\Gc\in   \ec^G_M(s_M)$ si et seulement si $s_M\in Z_{\Mc}$.
\end{remarques}

On définit de même un ensemble $\ec_M^L(s_M)$ pour tout $L\in \lc^L(M)$. On a une inclusion évidente $\ec_M^L(s_M)\subset \ec_M^G(s_M)$. Plus généralement, pour tout $s\in s_MZ_{\Mc}$, on définit $\ec_{M'}^{G_s}(s_M)$. On a d'ailleurs $\ec_{M'}^{G_s}(s_M)=\ec_{M'}^{G_s}(1)$ où $1$ est l'élément neutre de $\Tc$.

\begin{definition}
  Soit  $\Lc_1\subset \Lc_2$ deux sous-groupes réductifs et connexes  de $\Gc$ qui contiennent $\Tc$. On dit que $\Lc_1$ est \emph{elliptique} dans $\Lc_2$ si l'une des deux conditions équivalentes suivantes est satisfaite
  \begin{enumerate}
  \item  les centres connexes de $\Lc_1$ et $\Lc_2$ sont égaux ;
  \item le quotient $Z_{\Lc_1}/Z_{\Lc_2}$ est fini.
  \end{enumerate}
\end{definition}

Pour tout $L\in \lc^G(M)$, soit $\ec^L_{M,\el}(s_M)\subset \ec^L_M(s_M)$ le sous-ensemble des sous-groupes $\Hc\in   \ec^L_M(s_M)$ qui sont elliptiques dans $\Lc$.

\begin{lemme}\label{lem:reunion-EML}
  On a l'égalité suivante où la réunion, dans le membre de droite, est disjointe
$$\ec_M^G(s_M)=\bigcup_{L\in \lc(M)} \ec_{M,\el}^L(s_M)$$
\end{lemme}

\begin{preuve}
 Soit $\Hc\in \ec_M^{G}(s_M)$. Le centralisateur du centre connexe de $\Hc$ est un sous-groupe de Levi $\Lc$ de $\Gc$ qui contient $\Hc$. Par suite, on a $Z_{\Hc}^0\subset Z_{\Lc}\subset Z_{\Hc}$. On a donc $  Z_{\Lc}^0= Z_{\Hc}^0$ et $\Hc$ est elliptique dans $\Lc$. Tout autre sous-groupe de Levi $\Lc'$ qui contient $\Hc$ comme sous-groupe elliptique vérifie $Z_{\Lc'}^0= Z_{\Hc}^0$ et donc $\Lc'=\Lc$.
\end{preuve}

\begin{proposition}\label{prop:SH}
Soit $s_M\in \Tc$ et $\Mc'=\Mc_{s_M}$.
 Pour toute application 
$$J \ :\ \ec^G_{M,\el}(s_M) \to \CC,$$
il existe une unique application
$$S \ :\ \ec^G_{M,\el}(s_M) \to \CC$$
telle que la relation suivante 
 \begin{equation}
   \label{eq:SH}
   S(\Hc)= J(\Hc)-\sum_{s\in Z_{\Mc'}/Z_{\Hc}, \ s\not=1 } |Z_{\Hc_s}/Z_{\Hc}|^{-1} S(\Hc_s)
 \end{equation}
soit satisfaite pour tout $\Hc\in  \ec^G_{M,\el}(s_M)$.
\end{proposition}

\begin{remarque}\label{rq:finitude}
  Soit $\Hc\in \ec_{\el}^G(s_M)$. Pour tout $s\in Z_{\Mc'}$, le groupe  $\Hc_s$ appartient à $\ec^G_M(s_M)$ et ne dépend que de la classe de $s$ dans  $Z_{\Mc'}/Z_{\Hc}$. Si $\Hc_s$  n'est pas elliptique dans $\Hc$, le quotient $Z_{\Hc_s}/Z_{\Hc}$ n'est pas fini et, par convention on pose  $|Z_{\Hc_s}/Z_{\Hc}|^{-1}=0$. Ainsi la somme sur $s$ ne fait intervenir que des $S(\Hc_s)$ pour lesquels $\Hc_s$ appartient à  $\ec^G_{M,\el}(s_M)$. Comme il n'y a qu'un nombre fini de tels $\Hc_s$, le support de la somme sur $s$ est fini, puisqu'inclus dans une réunion finie de quotients finis $Z_{\Hc_s}/Z_{\Hc}$.

\end{remarque}
\medskip

\begin{preuve} Soit $\Hc\in \ec_{\el}^G(s_M)$. Remarquons que si $s\in Z_{\Mc'}$ et $s\notin Z_{\Hc}$ alors $\dim(\Hc_s)< \dim(\Hc)$. En particulier, si $\Hc$ est de dimension minimale, 
 la relation (\ref{eq:SH}) se réduit à 
$$S(\Hc)= J(\Hc).$$
Une récurrence sur la dimension des groupes $\Hc$  donne l'existence et l'unicité de $S$.
\end{preuve}

\begin{corollaire}\label{lem:SH2}
Soit $s_M\in \Tc$.
 Pour toute application 
$$J \ :\ \ec^G_{M}(s_M) \to \CC,$$
il existe une unique application
$$S \ :\ \ec^G_{M}(s_M) \to \CC$$
telle que pour tout $L\in \lc(M)$, la restriction de $S$ à $\ec_{M,\el}^L(s_M)$ soit l'application déduite de la  restriction de $J$ à $\ec_M^L(s_M)$ par le lemme \ref{prop:SH}.
\end{corollaire}

\begin{preuve}
  Elle découle du lemme \ref{lem:reunion-EML} et de la proposition \ref{prop:SH}.
  \end{preuve}
\end{paragr}

\begin{paragr}[Avatar stable de $J_M^{G,1}(a_M)$.] --- On reprend les notations de la section \ref{sec:sIOP}. On introduit la définition suivante. On note $1$ l'élément neutre de $\Tc$.

  \begin{definition}\label{def:IOPst}
    Soit $M\in \lc^G(T)$ et $a_M\in \Ac_M(k)^\tau$. L'application
$$H\in \ec_M^G(1) \mapsto S_M^H(a_M)$$
est l'application déduite de l'application de la définition \ref{def:sIOP} 
$$ H\in \ec_M^G(1) \mapsto J_M^{H,1}(a_M)$$
par le corollaire \ref{lem:SH2}.
  \end{definition}

\end{paragr}

\section{Un second calcul de trace de Frobenius}

\begin{paragr}
  On reprend les hypothèses et les notations des sections \ref{sec:cohomologie}, \ref{sec:sIOP}, \ref{sec:premiercalcul}. Le résultat principal de cette section est le théorème suivant qui est un calcul d'une trace de Frobenius sur un prolongement intermédiaire en terme de la variante stable des intégrales orbitales pondérées $J_M^{G,1}(a_M)$ de la définition \ref{def:sIOP}.

  \begin{theoreme}\label{thm:traceFrob} Soit $G$ un groupe semi-simple, $M\in \lc^G$ et $a_M\in \Ac_M^{\el}(k)^\tau$. Soit $a=\chi^M_G(a_M)$. On suppose que $a\in \Ac_G^{\mathrm{bon}}$. Soit $j$ l'immersion ouverte
$$j : \Ac_G^{\el}\to \Ac_G.$$
On a alors l'égalité entre
$$\trace(\tau^{-1}, (j_{!*}({}^{\mathrm{p}}\mathcal{H}^{\bullet}( Rf^{\el}_{G,*} \Qlb)_{1}))_a)$$
et 
$$\frac{c(J_{a_M})  \cdot q^{\deg(D)|\Phi^G|/2}}{\vol(\ago_M/X_*(M)) \cdot |\Tc^{W_a}/Z_{\Mc}^0|\cdot |Z_{\Mc}^0\cap Z_{\Gc}|}\cdot    S_M^{G}(a_M).$$
où
\begin{enumerate}
\item ${}^{\mathrm{p}}\mathcal{H}^{\bullet}( Rf^{\el}_{G,*} \Qlb)_{1}$ est la $1$-partie de ${}^{\mathrm{p}}\mathcal{H}^{\bullet}( Rf^{\el}_{G,*} \Qlb)$ (cf. théorème \ref{thm:sdecomposition}) ;
\item son prolongement intermédiaire à $\Ac_G$ est notée $j_{!*}$ ;
\item les notations dans le membre de droite, à l'exception de $S_M^G(a_M)$ sont celles du théorème \ref{thm:scomptage} ;
\item $S_M^G(a_M)$ est la variante de l'intégrale $J_M^{G,1}(a_M)$ introduite à la  définition \ref{def:IOPst}.
\end{enumerate}
  \end{theoreme}

Ce théorème est non-trivial. S'il repose en partie sur des arguments de comptage (c'est-à-dire sur le   théorème \ref{thm:scomptage}), sa démonstration requiert aussi le théorème cohomologique \ref{thm:cohomologique}. Ce dernier point explique la condition $a\in \Ac_G^{\mathrm{bon}}$. La démonstration se trouve au paragraphe suivant.
\end{paragr}

\begin{paragr}[Démonstration du théorème \ref{thm:traceFrob}.] --- \label{S:preuve-traceFrob} On reprend les hypothèses du  théorème \ref{thm:traceFrob}. Pour tout $H\in \ec_{M,\el}^G(1)$, on définit $\tilde{S}_M^H(a_M)$ par l'égalité suivante
$$    \trace(\tau^{-1}, (j_{!*}({}^{\mathrm{p}}\mathcal{H}^{\bullet}( Rf^{\el}_{H,*} \Qlb)_{1}))_a)= $$
$$\frac{c(J_{a_M})  \cdot q^{\deg(D)|\Phi^G|/2}}{\vol(\ago_M/X_*(M)) \cdot |\Tc^{W_a}/Z_{\Mc}^0|\cdot |Z_{\Mc}^0\cap Z_{\Gc}|}\cdot \tilde{S}_M^{H}(a_M).$$

Il s'agit de voir qu'on a  

\begin{equation}
  \label{eq:S=Stilde}
  \tilde{S}_M^{H}(a_M)= S_M^{H}(a_M).
\end{equation}

Comme $\Ac_G^{\mathrm{bon}}\cap \Ac_H \subset \Ac_H^{\mathrm{bon}}$ (cf. lemme \ref{lem:AHbon}), on peut supposer, en raisonnant par récurrence sur la dimension des groupes, que l'égalité (\ref{eq:S=Stilde}) vaut pour tout $H\in \ec_{M,\el}^G(1)$ tel que $H\not=G$. Montrons que le point $a=\chi^M_G(a_M)$ appartient à l'ouvert (notation de (\ref{eq:ucs}) pour $\bar{s}=1$)
$$\uc_1=\uc_M\cap \Ac_G^{\mathrm{bon}}\cap (\Ac_G-\bigcup_{(s,\rho)} \Ac_{s,\rho})$$
où la réunion est prise sur l'ensemble fini (cf. proposition \ref{prop:finitude} des données endoscopiques géométriques $(s,\rho)$ de $G$ qui vérifient
\begin{itemize}
\item $s\in Z_{\Mc}^0$ ;
\item $\Ac_{G}^{\el}\cap \Ac_{s,\rho}\not=\emptyset$ ;
\item $\rho$ est un homomorphisme  non trivial.
\end{itemize}

Comme $\Ac_M^{\el}\subset \uc_M$ (cf. la définition de $\uc_M$ donnée au §\ref{S:OmegaM}), il reste à voir qu'on a $a\notin  \Ac_{s,\rho}$ pour tout couple $(s,\rho)$ qui vérifie les trois conditions ci-dessus. Soit $(s,\rho)$ une donnée endoscopique géométrique de $G$ tel que $s\in Z_{\Mc}^0$ et $a\in \Ac_{s,\rho}$.  Puisque $a=\chi^M_G(a_M)$, on a $W_a\subset W^M$  (cf. proposition \ref{prop:AM-Wa}). Or  la condition $s\in Z_{\Mc}^0$ entraîne $W^M\subset W^0_s$. Donc $\rho$ est nécessairement trivial.

On peut donc prendre la trace de $\tau^{-1}$ sur la fibre en $a$ de l'égalité du théorème \ref{thm:cohomologique}. Après simplification par le facteur $\frac{c(J_{a_M})  \cdot q^{\deg(D)|\Phi^G|/2}}{\vol(\ago_M/X_*(M)) \cdot |\Tc^{W_a}/Z_{\Mc}^0|}$ , on trouve l'égalité
$$ J_M^{G,1}(a_M)=\sum_{s} |Z_{\Mc}^0\cap Z_{\Gc_s}|^{-1} \cdot  \tilde{S}_M^{G_s}(a_M).$$
où la somme est prise sur les $s\in Z_{\Mc}^0$ tels que $\Gc_s$ soit semi-simple. Comme les termes dans la somme de droite sont invariants par translation par $Z_{\Mc}^0\cap Z_{\Gc}$, on a aussi
$$ J_M^{G,1}(a_M)=\sum_{s} |Z_{\Mc}^0\cap Z_{\Gc}|\cdot |Z_{\Mc}^0\cap Z_{\Gc_s}|^{-1} \cdot  \tilde{S}_M^{G_s}(a_M),$$
où, cette fois-ci, la somme est prise sur les $s\in Z_{\Mc}^0/(Z_{\Mc}^0\cap Z_{\Gc})$ tels que $\Gc_s$ soit semi-simple. Puisque $\Mc$ est un sous-groupe de Levi de $\Gc$ et de $\Gc_s$, on a $Z_{\Mc}=Z_{\Gc}Z_{\Mc}^0$ et  $Z_{\Mc}=Z_{\Gc_s}Z_{\Mc}^0$. En particulier, on a $Z_{\Gc_s}\subset Z_{\Gc}( Z_{\Mc}^0 \cap Z_{\Gc_s})$. On en déduit les isomorphismes 
$$ Z_{\Mc}^0/(Z_{\Mc}^0\cap Z_{\Gc})\cong Z_{\Mc}/ Z_{\Gc}$$
et
$$(Z_{\Mc}^0\cap Z_{\Gc_s})/(Z_{\Mc}^0\cap Z_{\Gc})\cong Z_{\Gc_s}/ Z_{\Gc},$$
on obtient
$$ J_M^{G,1}(a_M)=\sum_{s\in Z_{\Mc}/ Z_{\Gc} } |Z_{\Gc_s}/ Z_{\Gc}|^{-1}\tilde{S}_M^{G_s}(a_M),$$
où, par la  convention habituelle, l'inverse du cardinal d'un ensemble infini est nul. 
Par hypothèse de récurrence, on a $\tilde{S}_M^{G_s}(a_M)=S_M^{G_s}(a_M)$ pour $s\not=1$. Il vient alors
\begin{eqnarray*}
  \tilde{S}_M^{G}(a_M)&=&J_M^{G,1}(a_M)-\sum_{s\in Z_{\Mc}/ Z_{\Gc}, s\not=1 } |Z_{\Gc_s}/ Z_{\Gc}|^{-1}\tilde{S}_M^{G_s}(a_M)\\
&=&J_M^{G,1}(a_M)-\sum_{s\in Z_{\Mc}/ Z_{\Gc, s\not=1} } |Z_{\Gc_s}/ Z_{\Gc}|^{-1}S_M^{G_s}(a_M).
\end{eqnarray*}
Or cette dernière expression est précisément $S_M^G(a_M)$ ce qui termine la démonstration.

\end{paragr}

\begin{paragr} Le lemme suivant a servi dans la démonstration du théorème \ref{thm:traceFrob}.

  \begin{lemme} \label{lem:AHbon} Soit $H$ un groupe dual d'un sous-groupe réductif $\Hc$ de $\Gc$ qui contient $\Tc$. Soit $\Ac_H  \hookrightarrow  \Ac_G$ l'immersion fermée canonique (induite par l'inclusion $W^H \subset W^G$). On a alors
$$\Ac_G^{\mathrm{bon}}\cap \Ac_H \subset \Ac_H^{\mathrm{bon}}.$$
\end{lemme}

\begin{preuve}
D'après la section \ref{sec:Abon},  on dispose sur $\Ac_H$ d'une fonction  $\delta^H$ et sur $\Ac_G$ d'une fonction $\delta^G$. Soit $a\in \Ac_G^{\mathrm{bon}}$. On a donc $\mathrm{codim}_{\Ac_G}(a)\geq \delta^G_a$. Supposons de plus $a\in \Ac_H$. On vérifie la formule suivante
$$\mathrm{codim}_{\Ac_G}(\Ac_H)=(\dim(G)-\dim(H))\deg(D)/2=\delta^G_a-\delta_a^H.$$
La première égalité provient du calcul de la dimension de $\Ac_G$ par la formule de Riemann-Roch et la seconde résulte de la formule de Ngô-Bezrukavnikov, cf. lemme \ref{lem:Bezruk} (cf. aussi \cite{ngo2} §4.17). On en déduit aussitôt l'inégalité $\mathrm{codim}_{\Ac_H}(a)\geq \delta^H_a$ qu'il fallait démontrer.
\end{preuve}

\end{paragr}

\section{Une forme globale du lemme fondamental pondéré}\label{sec:LFPglobal}

\begin{paragr} Le théorème suivant, qui est le résultat principal de cette section, est une forme globale du lemme fondamental pondéré d'Arthur. 

  \begin{theoreme}
    \label{thm:lfpglobal}
Soit $M\in \lc^G$ et $s_M\in \Tc$. Soit $M'$ de dual $\Mc'=\Mc_{s_M}$. Soit $a_M\in \Ac_{M'}^{\el}(k)^\tau$. Par abus, on note encore $a_M$ l'élément $\chi^{M'}_M(a_M)$. On fait les hypothèses suivantes :
\begin{itemize}
\item $\Mc'\subsetneq\Mc$ (c'est-à-dire $s_M\notin Z_{\Mc}$) ;
\item $\Mc'$ est elliptique dans $M$ (ce qui équivaut à  $a_M\in \Ac_M^{\el}(k)$) ;
\item l'élément $\chi^{M'}_G(a_M)$ appartient à $\Ac_G^{\mathrm{bon}}$.
\end{itemize}
On a alors l'égalité suivante
$$J_{M}^{G,s_M}(a_M)=|Z_{\Mc'}/Z_{\Mc}| \sum_{s\in s_M Z_{\Mc}/Z_{\Gc}} |Z_{\Gc_s}/Z_{\Gc}|^{-1} S_{M'}^{G_s}(a_M).$$
  \end{theoreme}

  \begin{paragr}[Démonstration du théorème \ref{thm:lfpglobal}.] --- Commençons par le cas d'un groupe semi-simple.  Montrons le lemme suivant.

    \begin{lemme}
      Soit $a=\chi^{M'}_G(a_M)$. Alors $a$ appartient à l'ouvert $\uc_{s_M}$ défini en (\ref{eq:ucs}).
    \end{lemme}

    \begin{preuve}
      Par hypothèse, on sait que $a$ appartient à $\Ac_G^{\mathrm{bon}}$. L'hypothèse $\Mc'$ est elliptique dans $M$ implique  que $a\in \Ac_M^{\el}$ (cf. proposition \ref{lem:Gs-ss}) et donc $a$ appartient à l'ouvert $\uc_M$ du §\ref{S:OmegaM}. Soit  $(s,\rho)$ une donnée endoscopique géométrique (cf. définition \ref{def:endogeo}) telle que les trois conditions suivantes soient satisfaites :
\begin{itemize}
\item $s\in s_M Z_{\Mc}^0$ ;
\item $\Ac_{G}^{\el}\cap \Ac_{s,\rho}\not=\emptyset$ ;
\item $a\in  \Ac_{s,\rho}$.
\end{itemize}
 Puisque $a_M\in \Ac_{M'}$, on a $W_a\subset W^{M'}$ (par une variante de la proposition \ref{prop:AM-Wa}). Comme on a $W^{M'}\subset W_s^0$, on a aussi $W_a\subset W_s^0$ et $\rho$ est nécessairement trivial. Cela termine la démonstration.
\end{preuve}

Prenons $\xi$ en position générale. Le théorème \ref{thm:cohomologique} spécialisé au point $a\in \uc_{s_M}$ entraîne une identité de traces de Frobenius qui s'explicite ainsi par les théorèmes \ref{thm:scomptage} et \ref{thm:traceFrob} (cf. aussi lemme \ref{lem:AHbon}) :
    $$\frac{c(J_{a_M})\cdot q^{\deg(D)|\Phi^G|/2}}{\vol(\ago_M/X_*(M)) \cdot |\Tc^{W_a}/Z_{\Mc}^0|}    J_M^{G,s_M}(a_M)=$$
$$\sum_{s} \frac{c(J_{a_M})  \cdot q^{\deg(D)|\Phi^G|/2}}{\vol(\ago_{M'}/X_*(M')) \cdot |\Tc^{W_a}/Z_{\Mc'}^0|\cdot |Z_{\Mc'}^0\cap Z_{\Gc_s}|}\cdot    S_{M'}^{G_s}(a_M),$$
où la somme est prise sur les $s\in s_M Z_{\Mc}^0$ tels que $G_s$ soit semi-simple. L'hypothèse $M'$ elliptique entraîne les égalités $Z_{\Mc}^0=Z_{\Mc'}^0$, $\ago_M=\ago_{M'}$ (en tant qu'espace euclidien) et $X_*(M)=X_*(M')$. Après simplification, on obtient l'égalité 
$$ J_M^{G,s_M}(a_M)=\sum_{s} |Z_{\Mc}^0\cap Z_{\Gc_s}|^{-1}    S_{M'}^{G_s}(a_M)$$
pour le même ensemble de sommation, qu'on peut réécrire (comme dans la démonstration du théorème \ref{thm:traceFrob}, cf. §\ref{S:preuve-traceFrob})
$$ J_M^{G,s_M}(a_M)=\sum_{s} \frac{[Z_{\Mc}^0\cap Z_{\Gc}|}{|Z_{\Mc}^0\cap Z_{\Gc_s}|}    S_{M'}^{G_s}(a_M)$$
où cette fois-ci la somme est sur les $s\in s_MZ_{\Mc}^0/ Z_{\Mc}^0\cap Z_{\Gc}$ tels que $G_s$ est elliptique. La suite exacte 
$$1 \longrightarrow  ( Z_{\Mc}\cap Z_{\Gc_s})/Z_{\Gc} \longrightarrow  Z_{\Gc_s}/Z_{\Gc} \longrightarrow  Z_{\Mc'}/Z_{\Mc}\longrightarrow 1$$
(où la surjectivité provient des égalités  $Z_{\Mc'}=Z_{\Gc_s} Z_{\Mc'}^0$ et $Z_{\Mc'}^0=Z_{\Mc}^0$) et l'isomorphisme
$$(Z_{\Gc_s}\cap Z_{\Mc})/Z_{\Gc}\cong (Z_{\Mc}^0\cap Z_{\Gc_s})/(Z_{\Mc}^0\cap Z_{\Gc})$$ 
entraînent la formule
$$\frac{[Z_{\Mc}^0\cap Z_{\Gc}|}{|Z_{\Mc}^0\cap Z_{\Gc_s}|}=\frac{|Z_{\Mc'}/Z_{\Mc}|}{|Z_{\Gc_s}/Z_{\Gc}|}.$$
Tous les ordres dans la formule ci-dessus sont finis si $G_s$ est elliptique. Dans le second membre dette formule, le dénominateur est infini si $G_s$ n'est pas elliptique auquel cas son inverse vaut par convention $0$. En utilisant l'isomorphisme $Z_{\Mc}^0/ Z_{\Mc}^0\cap Z_{\Gc}\cong Z_{\Mc}/Z_{\Gc}$, on obtient l'égalité 
$$ J_M^{G,s_M}(a_M)=|Z_{\Mc'}/Z_{\Mc}|\sum_{s\in s_M Z_{\Mc}/Z_{\Gc}} |Z_{\Gc_s}/Z_{\Gc}|^{-1}  S_{M'}^{G_s}(a_M)$$
qu'il fallait démontrer.

Le passage d'un groupe semi-simple à un groupe réductif quelconque est laissé au lecteur.
  \end{paragr}

\end{paragr}

\section{Formalisme d'Arthur de la stabilisation des intégrales orbitales pondérées}\label{sec:formalisme}

\begin{paragr} Cette section reprend les notations et les définitions des §§ \ref{S:avatar1} et \ref{S:avatar2}. Elle est indépendante des autres parties de l'article. On développe dans un contexte abstrait le formalisme d'Arthur de la stabilisation des intégrales orbitales pondérées  et de leurs formules de descente et de scindage (cf.  \cite{Ar-transfer} par exemple). On donne plusieurs propositions qui nous permettront ensuite de déduire le lemme fondamental pondéré de la variante globale donnée par le théorème \ref{thm:lfpglobal}. Tous les résultats reposent sur les méthodes d'Arthur.
\end{paragr}

\begin{paragr} Soit $M\in \lc^G$. Rappelons que pour tout $s_M\in \Tc$, on a défini au §\ref{S:avatar2} des ensembles finis $\ec^G_{M,\el}(s_M) \subset \ec^G_{M}(s_M)$ de sous-groupes réductifs de $\Gc$ qui contiennent $\Tc$. Pour tout sous-groupe de Levi de $\Gc$ contenant $\Tc$, soit $\lc^{\Gc}(\Mc)$ l'ensemble fini des sous-groupes de Levi de $\Gc$ qui contiennent $\Mc$. Cet ensemble est en bijection  naturelle avec $\lc^G(M)$. On pose $\lc^{\Gc}=\lc^{\Gc}(\Tc)$. Cette définition vaut plus généralement pour un sous-groupe réductif de $\Gc$ et un de ses sous-groupes de Levi. 

\begin{definition} \label{def:JE}Soit $s_M\in \Tc$ et $\Mc'=\Mc_{s_M}$.
 Pour toute application 
$$J \ :\ \ec^G_{M,\el}(s_M) \to \CC$$
Soit
$$J^{\ec}(\Gc)=  |Z_{\Mc'}/Z_{\Mc}| \sum_{s\in s_M Z_{\Mc}/Z_{\Gc}}  |Z_{\Gc_s}/Z_{\Gc}|^{-1} S(\Gc_s)$$
où $S$ est l'application déduite de $J$ par la proposition \ref{prop:SH}. 
\end{definition}

\begin{remarque} \label{rq:JE1} Si $s_M\in Z_{\Mc}$, on a $\Mc=\Mc'$ et $\Gc\in \ec^G_{M,\el}(s_M)$. Dans ce cas, on a $J^\ec(\Gc)=J(\Gc)$. En effet la formule ci-dessus se réécrit
$$J^{\ec}(\Gc)=   S(\Gc) +\sum_{s\in s_M Z_{\Mc}/Z_{\Gc},s\not=1}  |Z_{\Gc_s}/Z_{\Gc}|^{-1} S(\Gc_s)$$
et la somme sur $s$ ci-dessus est égale à $J(\Gc)-S(\Gc)$ par la formule (\ref{eq:SH}). 
 \end{remarque}

 \begin{remarque}\label{rq:JE2} Par le lemme \ref{lem:reunion-EML}, on a $\ec_M^G(s_M)=\bigcup_{L\in \lc(M)} \ec_{M,\el}^L(s_M)$ où la réunion est disjointe. Pour toute application $J \ :\ \ec^G_{M}(s_M) \to \CC$, on obtient donc une application
$$J^\ec : \lc^{\Gc}(\Mc) \to \CC$$
donnée par la définition \ref{def:JE}.
    \end{remarque}

\begin{definition} \label{def:stable}Soit $M\in \lc^G$ et $s_M\in \Tc$. On dit qu'une application 
$$J \ :\ \ec^G_{M,\el}(s_M)\cup\{\Gc\} \to \CC$$
est \emph{stabilisante} si elle vérifie l'égalité $J(\Gc)=J^\ec(\Gc)$. 
 On dit qu'une application 
$$J \ :\ \ec^G_{M}(s_M)\cup\lc^{\Gc}(\Mc) \to \CC$$
est \emph{stabilisante} si pour tout $L\in \lc^G(M)$ l'application restreinte à $\ec_{M,\el}^L(s_M)\cup\{\Lc\}$ est stabilisante.
\end{definition}

Voici un premier exemple trivial.

\begin{exemple}
  Si $s_M\in Z_{\Mc}$, toute application $J \ :\ \ec^G_{M,\el}(s_M)\cup\{\Gc\} \to \CC$ est stabilisante en vertu de la remarque \ref{rq:JE1}.
\end{exemple}

En voici un second qui ne l'est pas.

\begin{exemple}\label{ex:2}
  Soit $s_M\in \Tc$ et $\Mc'=\Mc_{s_M}$. On suppose $s_M\notin Z_{\Mc}$. Pour $a_M\in \Ac_{M'}^{\el}(k)^\tau$ qui vérifie les hypothèses du théorème \ref{thm:lfpglobal}, l'application $J$ définie par $J( \Gc)=J_{M}^{G,s_M}(a_M)$ et pour $\Hc\in \ec^G_{M,\el}(s_M)$ par $J(\Hc)=J_{M'}^{H,1}(a_M)$ est stabilisante : ce n'est qu'une reformulation du théorème \ref{thm:lfpglobal}.
\end{exemple}
  
\end{paragr}

\begin{paragr}[Scindage et descente.] --- Dans tout le reste de cette section, on fixe $M\in \lc^G$ et $s_M\in \Tc$. Les intégrales orbitales pondérées admettent des formules de descente et de scindage mises en évidence par Arthur (cf. \cite{Ar-JAMS1}). Le coefficient $d$ de la définition \ref{def:d} ci-dessous joue un rôle clef dans ces formules. Auparavant, pour tout $\Hc\in \ec_M^G(s_M)$, on pose
$$\ago_{\Hc}=X^*(\Hc)\otimes_\ZZ \RR=X^*(Z_{\Hc}^0)\otimes_\ZZ \RR.$$
On a d'ailleurs $\ago_H=\ago_{\Hc}$ si $H$ est le dual de $\Hc$. On a alors des décompositions $\ago_{\Tc}=\ago_{\Tc}^{\Hc}\oplus \ago_{\Hc}$ comme au §\ref{S:parab}.

\begin{definition}
  \label{def:d} Soit $\Hc\in \ec_M^G(s_M)$ et $\Rc\in \lc^{\Hc}(\Tc)$. Pour tous sous-groupes de Levi $\Lc_1$ et $\Lc_2$ dans $\lc^{\Hc}(\Rc)$, on pose 
$$d_{\Rc}^{\Hc}(\Lc_1,\Lc_2)=0$$
sauf si on a une somme directe $\ago_{\Rc}^{\Lc_1}\oplus \ago_{\Rc}^{\Lc_2}=\ago_{\Rc}^{\Hc}$ auquel cas $d$ est la valeur absolue du déterminant de l'isomorphisme canonique $\ago_{\Rc}^{\Lc_1}\times  \ago_{\Rc}^{\Lc_2}\to\ago_{\Rc}^{\Hc}$ écrit dans des bases orthonormales.
\end{definition}

\begin{remarque}
  On a $d_{\Rc}^{\Hc}(\Lc_1,\Lc_2)=d_{\Rc}^{\Hc}(\Lc_2,\Lc_1)$ et ce coefficient est non nul si et seulement si $\ago^{\Gc}_{\Lc_1}\oplus \ago^{\Gc}_{\Lc_2}=\ago_{\Rc}^{\Hc}$. On a aussi $d_{\Rc}^{\Hc}(\Rc,\Hc)=1$.
\end{remarque}

On rassemble dans le lemme suivant quelques faits élémentaires bien connus dont la démonstration est laissée au lecteur.

  \begin{lemme}\label{lem:centres}
  \begin{enumerate}
  \item On a  $a_{\Lc_1}^{\Hc}\cap a_{\Lc_2}^{\Hc}=0$ 
si et seulement si le quotient
$$(Z_{\Lc_1}\cap Z_{\Lc_2})/Z_{\Hc}$$
est fini ;
\item on a $a_{\Rc}=a_{\Lc_1}+a_{\Lc_2}$ si et seulement si $Z_{\Rc}=Z_{\Lc_1}Z_{\Lc_2}$.
\item Soit $s\in \Tc$ tel que $\Gc_s$ est elliptique dans $\Gc$ et $\Mc_s$ est elliptique dans $\Mc$. L'application $\Lc\to \Lc_s$ induit une bijection de l'ensemble des $\Lc\in \lc(\Mc)$ tel que $\Lc_s$ est elliptique dans $\Lc$ sur $\lc^{\Hc_s}(\Mc_s)$.
  \end{enumerate}
\end{lemme}

\begin{definition}
  \label{def:e} Soit $\Hc\in \ec_M^G(s_M)$ et $\Rc\in \lc^{\Hc}(\Tc)$. Pour tous sous-groupes de Levi $\Lc_1$ et $\Lc_2$ dans $\lc^{\Hc}(\Rc)$, on pose 
$$e_{\Rc}^{\Hc}(\Lc_1,\Lc_2)= | (Z_{\Lc_1}\cap Z_{\Lc_2})/Z_{\Hc}|^{-1} d_{\Rc}^{\Hc}(\Lc_1,\Lc_2),$$
avec la convention habituelle que l'inverse du cardinal d'un ensemble infini est nul.
\end{definition}

\begin{remarque}
Par l'assertion 1 du lemme \ref{lem:centres}, on a $e_{\Rc}^{\Hc}(\Lc_1,\Lc_2)=0$ si et seulement si $d_{\Rc}^{\Hc}(\Lc_1,\Lc_2)=0$.
\end{remarque}

Nous pouvons maintenant énoncer les deux principaux résultats de cette section.

\begin{proposition} (Scindage)
    \label{prop:scindage-stab}
Soit $s_M\in \Tc$ tel que $s_M\notin Z_{\Mc}$ et $\Mc'=\Mc_{s_M}$ est un sous-groupe elliptique de $\Mc$. 

Soit $J_1$ et  $J_2$ deux applications de $\ec_M^G(s_M)\cup \lc^{\Gc}(\Mc)$ dans $\CC$. Soit $J$ l'application de $\ec_M^G(s_M) \cup \lc^{\Gc}(\Mc)$ dans $\CC$ définie pour tout $\Hc\in \ec_M^G(s_M)$ par
$$J(\Hc)=\sum_{\Lc_1,\Lc_2\in \lc^{\Hc}(\Mc')} d_{\Mc'}^{\Hc}(\Lc_1,\Lc_2)\cdot J_1(\Lc_1)\cdot J_2(\Lc_2)$$
et tout $\Lc\in \lc^{\Gc}(\Mc)$ par
$$J(\Lc)=\sum_{\Lc_1,\Lc_2\in \lc^{\Lc}(\Mc)} d_{\Mc}^{\Lc}(\Lc_1,\Lc_2)\cdot J_1(\Lc_1)\cdot J_2(\Lc_2).$$
On a
\begin{enumerate}
\item Si $J_1$ et $J_2$ sont stabilisantes alors $J$ est stabilisante ;
\item Si  $J$ et $J_2$ sont stabilisantes et $J_2(\Mc)\not=0$ alors  $J_1$ est stabilisante.
\end{enumerate}
  \end{proposition}

 \begin{proposition} (Descente)
    \label{prop:descente-stab}Soit $s_M\in \Tc$ tel que $s_M\notin Z_{\Mc}$ et $\Mc'=\Mc_{s_M}$ est un sous-groupe elliptique de $\Mc$. Soit $\Rc\in \lc^{\Mc}$ tel que $\Rc_{s_M}$ est elliptique dans $\Rc$. Soit $\Rc'=\Rc_{s_M}$. 

Soit $J_R$ une application  de $\ec_R^G(s_M) \cup \lc^{\Gc}(\Rc)$ dans $\CC$. Soit $J_M$ l'application de $\ec_M^G(s_M) \cup \lc^{\Gc}(\Mc)$ dans $\CC$ définie pour tout $\Hc\in \ec_M^G(s_M)$ par
$$J_M(\Hc)=\sum_{\Lc\in \lc^{\Hc}(\Rc')}d_{\Rc'}^{\Hc}(\Mc',\Lc)J_R(\Lc)$$
et pour tout $\Lc\in \lc^{\Gc}(\Mc)$ par
$$J_M(\Lc)=\sum_{\Lc_1\in \lc^{\Lc}(\Rc)}d_{\Rc}^{\Lc}(\Mc,\Lc_1)J_R(\Lc_1).$$
Si  $J_R$ est stabilisante, alors $J_M$ est stabilisante.
  \end{proposition}
 
\begin{remarque}
    En utilisant  l'assertion 3 du lemme \ref{lem:centres}, on voit que pour tout  que  pour $\Hc\in \ec_M^G(s_M)$ et tout $\Lc\in \lc^{\Hc}(\Rc)$, on a $\Lc\in \ec_R^G(s_M)$.
  \end{remarque}

Le reste de la section  (en particulier les §§ \ref{S:scindage-stab} et \ref{S:descente-stab}) est consacrée aux démonstrations des propositions \ref{prop:scindage-stab} et \ref{prop:descente-stab}.
\end{paragr}

\begin{paragr}[Démonstration de la proposition \ref{prop:scindage-stab}.] --- \label{S:scindage-stab} Elle repose sur la proposition suivante dont le lecteur trouvera pour sa commodité une démonstration au §\ref{S:scindage}.
 
 \begin{proposition}\label{prop:scindage} (Formules de scindage)
On reprend les hypothèses de la proposition \ref{prop:scindage-stab}
  
\begin{enumerate}
\item Soit $S$,   resp $S_i$, l'application déduite de $J$, resp $J_i$, par le corollaire \ref{lem:SH2}. Pour tout $\Hc\in \ec^G_{M}(s_M)$, on a
$$S(\Hc)=\sum_{\Lc_1,\Lc_2 \in \lc^{\Hc}(\Mc')} e_{\Mc'}^{\Hc}(\Lc_1,\Lc_2)  \cdot S_1(\Lc_1)\cdot S_2(\Lc_2) \ ;$$
\item Soit $J^\ec$, resp $J^\ec_i$, l'application déduite de $J$ (cf. définition \ref{def:JE} et remarque \ref{rq:JE2}). Pour tout $\Lc\in \lc^{\Gc}(\Mc)$, on a

$$J^\ec(\Lc)=\sum_{\Lc_1,\Lc_2 \in \lc^{\Lc}(\Mc)} d_{\Mc}^{\Lc}(\Lc_1,\Lc_2)  \cdot J^\ec_1(\Lc_1)\cdot J^\ec_2(\Lc_2)$$
\end{enumerate}
\end{proposition}

On peut maintenant donner la démonstration de la proposition \ref{prop:scindage-stab}. L'assertion 1 est une conséquence évidente de l'assertion 2 de la proposition \ref{prop:scindage} ci-dessus. Prouvons l'assertion 2 de la proposition \ref{prop:scindage-stab}. Par hypothèse, $J$ et $J_2$ sont stabilisantes. On a donc 
$$J(\Gc)=J^\ec(\Gc).$$
On peut réécrire cette égalité en utilisant pour le membre de gauche la définition de $J$ en termes de $J_1$ et $J_2$ et pour le membre de droite l'assertion 2 de la proposition \ref{prop:scindage}. Il vient
$$\sum_{\Lc_1,\Lc_2\in \lc^{\Gc}(\Mc)} d_{\Mc}^{\Gc}(\Lc_1,\Lc_2)\cdot J_1(\Lc_1)\cdot J_2(\Lc_2)=\sum_{\Lc_1,\Lc_2\in \lc^{\Gc}(\Mc)} d_{\Mc}^{\Gc}(\Lc_1,\Lc_2)\cdot J_1^\ec(\Lc_1)\cdot J_2^\ec(\Lc_2).$$
Puisque $J_2$ est stabilisante, on  a $J_2^\ec=J_2$. En raisonnant par récurrence, on peut bien supposer qu'on a $J_1^{\ec}(\Lc_1)=J_1(\Lc_1)$ pour tout $\Lc_1\in \lc^{\Gc}(\Mc)-\{\Gc\}$. En retranchant membre à membre les termes égaux, on voit que l'égalité précédente est équivalente à 
$$J_1(\Gc) J_2(\Mc)= J_1^\ec(\Gc) J_2(\Mc)$$
d'où $J_1(\Gc) = J_1^\ec(\Gc)$ puisqu'on a supposé $J_2(\Mc)\not=0$. Donc  $J_1$ est bien stabilisante.
\end{paragr}

\begin{paragr}[Démonstration de la proposition \ref{prop:scindage-stab}.] --- \label{S:scindage} Commençons par le lemme suivant.

  \begin{lemme}  \label{lem:unebijection}  Soit $s_M\in \Tc$ tel que $\Mc'=\Mc_{s_M}$ soit elliptique dans $\Mc$. Soit $\tilde{\ec}$ l'ensemble des triplets 
$$(s,\Lc_1',\Lc_2')$$
formés de $s\in s_M Z_{\Mc}/Z_{\Gc}$ et $\Lc_1',\Lc_2' \in \lc^{\Gc_s}(\Mc')$ tels que 
\begin{itemize}
\item $\Gc_s$ est elliptique dans $\Gc$ ;
\item  $d_{\Mc'}^{\Gc_s}(\Lc_1',\Lc_2')\not=0$.
\end{itemize}
Soit  $\tilde{\lc}$  l'ensemble des quadruplets 
$$(\Lc_1,\Lc_2,s_1,s_2)$$
formés de $\Lc_1,\Lc_2\in \lc^{\Gc}(\Mc)$ et $s_i\in Z_{\Mc}/Z_{\Lc_i}$ tels que
\begin{itemize}
\item $\Lc_{i,s_i}$ est elliptique dans $\Lc_i$ ;
\item  $d_{\Mc}^{\Gc}(\Lc_1,\Lc_2)\not=0$.
\end{itemize}

Soit $\tilde{\ec}\to \tilde{\lc}$ l'application qui à $(s,\Lc_1',\Lc_2')$ associe $(\Lc_1,\Lc_2,s_1,s_2)$ défini par
\begin{itemize}
\item $\Lc_i$ est l'unique élément de $\lc^{\Gc}(\Mc)$ tel que $\Lc_i'$ est elliptique dans $\Lc_i$ et $\Lc_i'=\Lc_{i,s}$ (cf. lemme \ref{lem:centres} assertion 3);
\item $s_i$ est l'image de $s$ par la projection $s_M Z_{\Mc}/Z_{\Gc} \to s_M Z_{\Mc}/Z_{\Lc_i}$.
\end{itemize}
Alors cette application est surjective et ses fibres ont toutes le même cardinal à savoir  
$$|(Z_{\Lc_1}\cap Z_{\Lc_2})/Z_{\Gc}|.$$
De plus, si  $(s,\Lc_1',\Lc_2')$ a pour image $(\Lc_1,\Lc_2,s_1,s_2)$, on a l'égalité entre 
 $$ |(Z_{\Lc_1}\cap Z_{\Lc_2})/Z_{\Gc}|\cdot|Z_{\Mc'}/Z_{\Mc}| \cdot |Z_{\Gc_s}/Z_{\Gc}|^{-1} \cdot e_{\Mc'}^{\Gc_s}(\Lc_1',\Lc_2') $$
et
$$|Z_{\Mc'}/Z_{\Mc}|^2 \cdot |Z_{\Lc_{1,s}}/Z_{\Lc_1}|^{-1}  \cdot |Z_{\Lc_{2,s}}/Z_{\Lc_2}|^{-1}  \cdot d_{\Mc}^{\Gc}(\Lc_1,\Lc_2).$$

\end{lemme}

\begin{preuve} Remarquons tout d'abord que si  $(s,\Lc_1',\Lc_2')$ a pour image $(\Lc_1,\Lc_2,s_1,s_2)$ alors $d_{\Mc'}^{\Gc_s}(\Lc_1',\Lc_2')=d_{\Mc}^{\Gc}(\Lc_1,\Lc_2)$ (car si $\Lc'$ est elliptique dans $\Lc$, on a $\ago_{\Lc}=\ago_{\Lc'}$.)\\

L'application $\tilde{\ec}\to \tilde{\lc}$ est \emph{surjective}. En effet, si  $(\Lc_1,\Lc_2,s_1,s_2)\in \tilde{\lc}$, on a  $d_{\Mc}^{\Gc}(\Lc_1,\Lc_2)\not=0$ et par le lemme \ref{lem:centres} on a $Z_{\Mc}=Z_{\Lc_1}Z_{\Lc_2}$. Donc $(s_1,s_2)$ se relève en $s\in Z_{\Mc}/Z_{\Gc}$. Le triplet $(s,\Lc_{1,s},\Lc_{2,s})$ appartient à $\tilde{\ec}$  : la seule condition qui n'est pas évidente est $\Gc_s$ elliptique dans $\Gc$. Par hypothèse, $\Lc_{i,s}=\Lc_{i,s_i}$ est elliptique dans $\Lc_i$. On a donc
$$\ago_{\Mc}^{\Gc}=\ago_{\Mc}^{\Lc_1}\oplus \ago_{\Mc}^{\Lc_2}=\ago_{\Mc}^{\Lc_{1,s}}\oplus \ago_{\Mc}^{\Lc_{2,s}}$$
et, passant aux orthogonaux, on a 
$$\ago_{\Gc}=\ago_{\Lc_{1,s}}\cap \ago_{\Lc_{2,s}}$$
Comme l'intersection ci-dessus contient $\ago_{\Gc_s}$, on a $\ago_{\Gc}=\ago_{\Gc_s}$ et $\Gc_s$ est elliptique dans $\Gc$. Il est clair que l'image de  $(s,\Lc_{1,s},\Lc_{2,s})$ est le quadruplet $(\Lc_1,\Lc_2,s_1,s_2)\in \tilde{\lc}$ de départ.

L'assertion sur le cardinal des fibres est évidente. 

Soit $(s,\Lc_1',\Lc_2')$ d'image $(\Lc_1,\Lc_2,s_1,s_2)$. Montrons la dernière assertion. Après simplification par $d_{\Mc'}^{\Gc_s}(\Lc_1',\Lc_2')=d_{\Mc}^{\Gc}(\Lc_1,\Lc_2)\not= 0$, on tombe sur l'égalité entre
$$|(Z_{\Lc_1}\cap Z_{\Lc_2})/Z_{\Gc}|\cdot |Z_{\Mc'}/Z_{\Mc}| \cdot |Z_{\Gc_s}/Z_{\Gc}|^{-1} |(Z_{\Lc_1'}\cap Z_{\Lc_2'})/Z_{\Gc_s}|^{-1}$$
et 
$$|(Z_{\Lc_1}\cap Z_{\Lc_2})/Z_{\Gc}|^{-1}\cdot |Z_{\Mc'}/Z_{\Mc}|^2 \cdot |Z_{\Lc_{1,s}}/Z_{\Lc_1}|^{-1}  \cdot |Z_{\Lc_{2,s}}/Z_{\Lc_2}|^{-1}.$$
Celle-ci résulte du diagramme commutatif à lignes exactes
$$\xymatrix{1  \ar[r]  &  (Z_{\Lc_1}\cap Z_{\Lc_2})/Z_{\Gc} \ar[r] \ar[d]^{\varphi}  & Z_{\Mc}/Z_{\Gc}  \ar[r] \ar[d] &  Z_{\Mc}/Z_{\Lc_1} \times Z_{\Mc}/Z_{\Lc_2}  \ar[r] \ar[d]&  1\\
1  \ar[r]  &  (Z_{\Lc_{1,s}}\cap Z_{\Lc_{2,s}})/Z_{\Gc_s} \ar[r]  & Z_{\Mc'}/Z_{\Gc_s}  \ar[r]  &  Z_{\Mc'}/Z_{\Lc_{1,s}} \times Z_{\Mc'}/Z_{\Lc_{2,s}}  \ar[r] &  1 }$$
et du lemme suivant.
\end{preuve}

\begin{lemme}\label{lem:morphsurj}
Soit $s$ tel que $\Mc=\Mc_s$ est elliptique dans $\Mc$ et $\Gc_s$ est elliptique dans $\Gc$.
Pour tout $\Lc\in \lc^{\Gc}(\Mc)$ tel que $\Lc_s$ est elliptique dans $\Lc$, le morphisme
$$  Z_{\Mc}/Z_{\Lc} \to Z_{\Mc'}/Z_{\Lc_{s}}$$
est surjectif et son noyau a pour cardinal 
 $$|Z_{\Lc_s}/ Z_{\Lc}|\cdot |Z_{\Mc'}/ Z_{\Mc}|^{-1}.$$
\end{lemme}

\begin{preuve}
  La surjectivité est une conséquence des égalités
$$Z_{\Mc'}=Z_{\Gc_s}(Z_{\Mc'})^0=Z_{\Gc_s}(Z_{\Mc})^0.$$
Comme $Z_{\Lc_s}=Z_{\Gc_s} (Z_{\Lc_s}^0)=Z_{\Gc_s} (Z_{\Lc}^0)$, le noyau s'écrit
$$(Z_{\Mc}\cap Z_{\Lc_s}) /Z_{\Lc} =   Z_{\Lc} (Z_{\Gc_s}\cap Z_{\Mc})/ Z_{\Lc}   \simeq  (Z_{\Gc_s}\cap Z_{\Mc})/  (Z_{\Gc_s}\cap Z_{\Lc}).$$
On conclut en remarquant que $(Z_{\Gc_s}\cap Z_{\Lc})/Z_{\Gc} $ est le noyau du morphisme surjectif 
 $$Z_{\Gc_s}/Z_{\Gc} \to Z_{\Lc_s}/Z_{\Lc}.$$

\end{preuve}

On peut maintenant donner la démonstration de la proposition \ref{prop:scindage}. Prouvons l'assertion 1. Soit    $\Hc\in \ec^G_{M}(s_M)$. Soit $\tilde{\ec}$ et $\tilde{\lc}$ les ensembles définis au lemme \ref{lem:unebijection} relativement au groupe $\Hc$, à son sous-groupe de Levi $\Mc'$ et l'élément $s_M=1$. Par définition des applications $S$ et $S_i$, on a 
\begin{eqnarray*}
  J(\Hc) &=& \sum_{\Lc_1,\Lc_2 \in \lc^{\Hc}(\Mc')} d_{\Mc'}^{\Hc}(\Lc_1,\Lc_2) \cdot J_1(\Lc_1)\cdot J_2(\Lc_2)\\
&=& \sum_{\Lc_1,\Lc_2 \in \lc^{\Hc}(\Mc')} d_{\Mc'}^{\Hc}(\Lc_1,\Lc_2)\big[\sum_{s_1\in Z_{\Mc'}/Z_{\Lc_1}} |Z_{\Lc_{1,s_1}}/Z_{\Lc_1}|^{-1} S_1(\Lc_{1,s_1})  \big]\times \\
& & \big[\sum_{s_2\in Z_{\Mc'}/Z_{\Lc_2}} |Z_{\Lc_{2,s_2}}/Z_{\Lc_2}|^{-1} S_2(\Lc_{2,s_2})  \big]\\
&=& \sum_{ (\Lc_1,\Lc_2,s_1,s_2)\in \tilde{\lc} } d_{\Mc'}^{\Hc}(\Lc_1,\Lc_2)\cdot |Z_{\Lc_{1,s_1}}/Z_{\Lc_1}|^{-1}        \cdot  |Z_{\Lc_{2,s_2}}/Z_{\Lc_2}|^{-1} \cdot S_1(\Lc_{1,s_1})\cdot  S_2(\Lc_{2,s_2})
\end{eqnarray*}
En utilisant le  lemme \ref{lem:unebijection}, on obtient
\begin{eqnarray*}
  J(\Hc) &=& \sum_{(s,\Lc_1',\Lc_2')\in \tilde{\ec}}   |Z_{\Hc_{s}}/Z_{\Hc}|^{-1}  e_{\Mc'}^{\Hc_s}(\Lc_1',\Lc_2') \cdot  S_1(\Lc_{1}) \cdot S_2(\Lc_{2}) \\
&=& \sum_{s\in Z_{\Mc'}/Z_{\Hc}}   |Z_{\Hc_{s}}/Z_{\Hc}|^{-1} \big(\sum_{\Lc_1,\Lc_2 \in \lc^{\Hc_s}(\Mc')}  e_{\Mc'}^{\Hc_s}(\Lc_1,\Lc_2) \cdot  S_1(\Lc_{1}) \cdot S_2(\Lc_{2})\big).
\end{eqnarray*}
Par l'unicité dans le corollaire \ref{lem:SH2}, on obtient l'expression voulue pour  $S(\Hc)$.

Montrons ensuite l'assertion 2.  Soit $\tilde{\ec}$ et $\tilde{\lc}$ les ensembles définis au lemme \ref{lem:unebijection} relativement au groupe $\Gc$, à son sous-groupe de Levi $\Mc$ et l'élément $s_M$. Par définition de $J^{\ec}$ et l'assertion 1 qu'on vient de démontrer, on a

\begin{eqnarray*}
  J^\ec(\Gc)&=& |Z_{\Mc'}/Z_{\Mc}| \sum_{s\in s_M Z_{\Mc}/Z_{\Gc}}  |Z_{\Gc_s}/Z_{\Gc}|^{-1} S(\Gc_s)\\
&=& |Z_{\Mc'}/Z_{\Mc}| \sum_{s\in s_M Z_{\Mc}/Z_{\Gc}}  |Z_{\Gc_s}/Z_{\Gc}|^{-1} \sum_{\Lc_1',\Lc_2' \in \lc^{\Gc_s}(\Mc')} e_{\Mc'}^{\Gc_s}(\Lc_1',\Lc_2')  \cdot S_1(\Lc_1')\cdot S_2(\Lc_2')\\ 
&=&   |Z_{\Mc'}/Z_{\Mc}|     \sum_{(s,\Lc_1',\Lc_2')\in \tilde{\ec}}  |Z_{\Gc_s}/Z_{\Gc}|^{-1} \cdot  e_{\Mc'}^{\Gc_s}(\Lc_1',\Lc_2')  \cdot S_1(\Lc_1')\cdot S_2(\Lc_2')\\
\end{eqnarray*}
Le  lemme \ref{lem:unebijection} implique qu'on a 
\begin{eqnarray*}
   J^\ec(\Gc)&=& \sum_{ (\Lc_1,\Lc_2,s_1,s_2)\in \tilde{\lc} }|Z_{\Mc'}/Z_{\Mc}|^2 \cdot |Z_{\Lc_{1,s}}/Z_{\Lc_1}|^{-1}  \cdot |Z_{\Lc_{2,s}}/Z_{\Lc_2}|^{-1}  \cdot d_{\Mc}^{\Gc}(\Lc_1,\Lc_2) \cdot S_1(\Lc_1')\cdot S_2(\Lc_2')\\
&=& \sum_{\Lc_1,\Lc_2 \in \lc^{\Lc}(\Mc)} d_{\Mc}^{\Lc}(\Lc_1,\Lc_2)   \cdot J^\ec_1(\Lc_1)\cdot J^\ec_2(\Lc_2).
\end{eqnarray*}
La dernière égalité résulte de la définition de  $J^\ec_1$ et $J^\ec_2$. Cela prouve donc l'assertion 2.

\end{paragr}

\begin{paragr}[Démonstration de la proposition \ref{prop:descente-stab}.] --- \label{S:descente-stab}
  Énonçons la proposition suivante qui sera démontrée au §\ref{S:descente}.

 \begin{proposition}\label{prop:descente}
    On se place sous les hypothèses de la proposition \ref{prop:descente-stab}.
\begin{enumerate}
\item Soit $S_M$ et $S_R$ les applications  déduites de $J_M$ et  $J_R$ par le corollaire \ref{lem:SH2}. Pour tout $\Hc\in \ec_M^G(s_M)$, on a
$$S_M(\Hc)= \sum_{\Lc\in \lc^{\Hc}(\Rc')}e_{\Rc'}^{\Hc}(\Mc',\Lc) S_R(\Lc) \ ;$$
\item Soit $J^\ec_M$ et $J^\ec_R$ les applications  déduites de $J_M$ et  $J_R$ (cf. définition \ref{def:JE} et remarque \ref{rq:JE2}). Pour tout $\Lc\in \lc^{\Gc}(\Mc)$, on a
$$J^{\ec}_M(\Lc)=\sum_{\Lc_1 \in \lc^{\Lc}(\Mc)} d_{\Rc}^{\Lc}(\Mc,\Lc_1) J^\ec_R(\Lc_1).$$
\end{enumerate}
\end{proposition}
On peut alors donner la preuve de la  proposition \ref{prop:descente-stab}. En utilisant l'assertion 2 de la  proposition \ref{prop:descente} et le fait que $J_R$ est stabilisante, on a 
    \begin{eqnarray*}
      J^{\ec}_M(\Lc)&=&\sum_{\Lc_1 \in \lc^{\Lc}(\Rc)} d_{\Rc}^{\Lc}(\Mc,\Lc_1) J^\ec_R(\Lc_1)\\
&=& \sum_{\Lc_1 \in \lc^{\Lc}(\Rc)} d_{\Rc}^{\Lc}(\Mc,\Lc_1) J_R(\Lc_1)\\
&=& J_M(\Lc).
\end{eqnarray*}
Donc $J_M$ est stabilisante.

\end{paragr}

\begin{paragr}[Démonstration de la proposition \ref{prop:descente}.] ---\label{S:descente}  Elle se déduit du lemme \ref{lem:bij-descente} ci-dessous  (cf. la démonstration de la proposition \ref{prop:scindage}).

\begin{lemme}
    \label{lem:bij-descente}
Soit $\Rc\subset  \Mc\subset \Gc$ deux sous-groupes de Levi de $\Gc$. Soit $s_M\in \Tc$, $\Mc'=\Mc_{s_M}$ et $\Rc'=\Rc_{s_M}$. On suppose que $\Mc'$ et $\Rc'$ sont elliptiques respectivement dans $\Mc$ et $\Rc$.

 Soit $\tilde{\ec}$ l'ensemble des couples $(s',\Lc')$
formés de $s'\in s_M Z_{\Mc}/Z_{\Gc}$ et $\Lc'\in \lc^{\Gc_s}(\Rc')$ tels que 
\begin{itemize}
\item $\Gc_{s'}$ est elliptique dans $\Gc$ ; 
\item  $d_{\Rc'}^{\Gc_{s'}}(\Mc',\Lc')\not=0$.
\end{itemize}
Soit  $\tilde{\lc}$  l'ensemble des couples $(\Lc,s)$ formés de $\Lc\in \lc^{\Gc}(\Rc)$ et $s\in Z_{\Mc}/Z_{\Lc}$ tels que
\begin{itemize}
\item $\Lc_{s}$ est elliptique dans $\Lc$ ;
\item  $d_{\Rc}^{\Gc}(\Mc,\Lc)\not=0$.
\end{itemize}

Soit $\tilde{\ec}\to \tilde{\lc}$ l'application qui à $(s',\Lc')$ associe $(\Lc,s)$ défini par
\begin{itemize}
\item $\Lc$ est l'unique élément de $\lc^{\Gc}(\Rc)$ tel que $\Lc'$ est elliptique dans $\Lc$ et $\Lc'=\Lc_{s'}$ (cf. lemme \ref{lem:centres} assertion 3);
\item $s$ est l'image de $s'$ par la projection $s_M Z_{\Mc}/Z_{\Gc} \to s_M Z_{\Rc}/Z_{\Lc}$.
\end{itemize}
Cette application est surjective et ses fibres ont toutes le même cardinal à savoir  
$$|(Z_{\Mc}\cap Z_{\Rc})/Z_{\Gc}|.$$
De plus, si  $(s',\Lc')$ a pour image $(\Lc,s)$, on a l'égalité entre 
 $$ |(Z_{\Mc}\cap Z_{\Lc})/Z_{\Gc}|  \cdot |Z_{\Mc'}/Z_{\Mc}| \cdot |Z_{\Gc_{s'}}/Z_{\Gc}|^{-1} \cdot e_{\Rc'}^{\Gc_{s'}}(\Mc',\Lc') $$
et
$$|Z_{\Rc'}/Z_{\Rc}| \cdot |Z_{\Lc_{s}}/Z_{\Lc}|^{-1}  \cdot d_{\Rc}^{\Gc}(\Mc,\Lc).$$
  \end{lemme}

  \begin{preuve}
    Elle est analogue à celle du lemme \ref{lem:unebijection}. Les détails sont laissés au lecteur.

  \end{preuve}
\end{paragr}

\section{Du global au local}

\begin{paragr} Le résultat principal est le théorème \ref{thm:LFPsemilocal} qui est une variante de  nature  «locale» du  théorème \ref{thm:lfpglobal}. Les résultats de scindage et de descente de la section \ref{sec:formalisme} jouent un rôle clef ainsi que certains calculs locaux d'intégrales orbitales pondérées. Les notations sont celles utilisées dans tout l'article, en particulier celles des sections \ref{sec:sIOP}, \ref{sec:avatarstable} et  \ref{sec:formalisme}.
  \end{paragr}

  \begin{paragr}[Notations.] --- Soit $V$ un ensemble fini ou infini mais non vide  de points fermés de la courbe $C$ qui est stable par le Frobenius $\tau$. Soit $\AAA_V$ le sous-groupe des adèles de $\AAA$ de composante nulle hors $V$. Le groupe $\AAA_V$ a une structure évidente d'anneau. Si $V$ est fini (resp. infini), $\AAA_V$ est le produit (resp. produit restreint) des complétés du corps des fonctions $F$ en les points de $V$. Soit $\oc_V=\prod_{v\in V}\oc_v$ où $\oc_v$ est le complété de l'anneau local de $C$ en $v$.

Soit $M\in \lc^G$,  $a_M\in \Ac_M(k)^\tau$ et $X_{a_M}=\eps_M(a_{M,\eta})\in \mgo(F)$ (cf. §\ref{S:sIOP}). Soit $J_{a_M}$ le schéma en groupes sur $C$ défini au §\ref{S:JaM} dont la fibre au point générique de $C$ est le centraliseur de $X_{a_M}$ dans $M\times_k F$ (c'est un sous-$F$-tore maximal de  $M\times_k F$). 

On complète les choix de mesures de Haar faits au §\ref{S:Haar}. On munit $J_{a_M}(\AAA_V)^\tau$ et $G(\AAA_V)^\tau$ des mesures de Haar qui donnent le volume $1$ respectivement au sous-groupe compact maximal de $J_{a_M}(\AAA_V)^\tau$ et à $G(\oc_V)^\tau$. On en déduit comme au §\ref{S:Haar} une mesure sur $
(J_{a_M}(\AAA_V)\back G(\AAA_V))^\tau$ invariante à droite par  $G(\AAA_V)^\tau$. Si $V'$ est le complémentaire de $V$ dans l'ensemble des points fermés de $C$ on a une bijection canonique
$$ (J_{a_M}(\AAA_V)\back G(\AAA_V))^\tau \times (J_{a_M}(\AAA_{V'})\back G(\AAA_{V'}))^\tau\to (J_{a_M}(\AAA)\back G(\AAA))^\tau$$
qui préserve les mesures.

Pour tout $s\in \Tc^{W_{a_M}}$, on définit la $s$-intégrale orbitale pondérée sur $V$ associée à $a_M$  par
\begin{eqnarray}\label{eq:JMGV}
  J_{M,V}^{G,s}(a_M)= q^{-\deg(D_V)|\Phi^G|/2} |D^G(X_{a_M})|^{1/2}_V \times \\ \nonumber \int_{(J_{a_M}(\AAA_V)\back G(\AAA_V))^\tau} \bg s, \inv_{a_M}(g)\bd \mathbf{1}_D(\Ad(g^{-1})X_{a_M}) v_M^G(g) \, dg,
\end{eqnarray}
où $D_V$ est la restriction de $D$ à $V$ et où  $|D^G(X_a)|_V$ est le produit sur $V$ des valeurs absolues usuelles du discriminant $D^G$ de $X_a$ (cf. §\ref{S:discr}).

\begin{remarque}
  Lorsque $V$ est l'ensemble de tous les points fermés l'intégrale orbitale pondérée définie par (\ref{eq:JMGV}) n'est autre que l'intégrale de la définition \ref{def:sIOP}.
\end{remarque}
\end{paragr}

\begin{paragr}[L'énoncé principal.] --- On peut maintenant énoncer le principal résultat de cette section. On y utilise les diviseurs $\mathfrak{R}_{M}^{G}$ et $\mathfrak{D}^{M}$ de $\car_{M,D}$ définis au §\ref{S:enonceAbon}.

\begin{theoreme}
  \label{thm:LFPsemilocal}
On se place sous les hypothèses du théorème \ref{thm:lfpglobal}. Soit 
$$h_{a_M} : C \to \car_{M',D}$$
la section associée à $a_M\in \Ac_{M'}(k)$. Soit $V$ un ensemble $\tau$-stable non vide de points fermés de $C$ tel que tout point fermé $v\notin V$ vérifie l'une des deux conditions suivantes :
\begin{enumerate}
\item $h_{a_M}(v)$  est un point lisse de $\mathfrak{R}_{M'}^{G}$ et n'appartient pas à  $\mathfrak{D}^{M'}$ ; de plus,  l'intersection de  $h_a(C)$ et  $\mathfrak{R}_{M'}^{G}$ est transverse en $h_a(v)$ ;
\item $h_{a_M}(v)$ est un point lisse de $\mathfrak{D}^{M'}$ et n'appartient pas à $\mathfrak{R}_{M'}^{G}$ ; de plus,  l'intersection de  $h_a(C)$ et  $\mathfrak{D}^{M'}$ est transverse en $h_a(v)$.
\end{enumerate}

Alors l'application $J_V :\ec_M^G(s_M)\cup \lc^{\Gc}(\Mc) \to \CC$ définie pour tout $\Hc\in  \ec_M^G(s_M)$ par $J_V(\Hc)=J_{M',V}^{H,1}(a_M)$ et pour tout $\Lc\in \lc^{\Gc}(\Mc)$ par  $J_V(\Lc)=J_{M,V}^{L,s_M}(a_M)$ est stabilisante au sens de la définition \ref{def:stable}.
\end{theoreme}

\begin{remarque}\label{rq:semilocal}
  Si $V$ est l'ensemble de tous les points fermés de $C$, le théorème n'est qu'une paraphrase du théorème  \ref{thm:lfpglobal} (cf. exemple \ref{ex:2}). L'ensemble des $v$ qui ne vérifient ni l'assertion 1 ni l'assertion 2 est fini. Autrement dit il existe toujours un ensemble fini $V$ qui vérifie les hypothèses du théorème. En ce sens, le théorème \ref{thm:LFPsemilocal} est une variante locale du théorème  \ref{thm:lfpglobal}.
\end{remarque}

\begin{preuve} Soit $V$ l'ensemble des points fermés de $C$. Soit $V_1$ un ensemble $\tau$-stable de points fermés et $V_2$ son complémentaire dans $V$. On suppose que tout $v\in V_2$ vérifie l'une des assertions 1 et 2 du théorème \ref{thm:LFPsemilocal}. Il s'agit de voir que $J_{V_1}$ est stabilisante. Or $J_V$, $J_{V_1}$ et $J_{V_2}$ vérifient les relations de scindage de la proposition \ref{prop:scindage-stab} comme le montre la proposition \ref{prop:scindage-IOP} ci-dessous. Il suffit donc, par l'assertion 2 de la  proposition \ref{prop:scindage-stab}, de montrer que $J_V$ et $J_{V_2}$ sont stabilisantes et que $J_2(\Mc)\not=0$. On a déjà dit dans la remarque \ref{rq:semilocal} ci-dessus que $J_V$ est stable. Il existe une partition (finie) de $V_2$ en parties $V_i'$ qui vérifient les hypothèses du lemme \ref{lem:final}. On a donc $J_2(\Mc)=\prod_i J_{V'_i}(\Mc)\not=0$ (avec les notations du lemme \ref{lem:final}) et par une application répétée de l'assertion 1  de la  proposition \ref{prop:scindage-stab}, on voit que $J_2$ est stabilisante et qu'on  $J_2(\Mc)\not=0$. 
\end{preuve}

\end{paragr}

\begin{paragr}[Formules de scindage des intégrales orbitales pondérées.] Pour tout $L_1,L_2\in \lc^G(M)$, on pose
$$d_M^G(L_1,L_2)=d_{\Mc}^{\Gc}(\Lc_1,\Lc_2)$$
où le second est celui de la définition \ref{def:d}. Voici les formules de scindage d'Arthur.

\begin{proposition}
\label{prop:scindage-IOP}
  Soit $V_1$ et $V_2$ deux ensembles disjoints et non vides de points fermés de $C$. Soit $V=V_1\cup V_2$. Pour tout $a\in \Ac_M(k)^\tau$ et tout $s\in \Tc^{W_a}$, on a l'égalité
$$J_{M,V}^{G,s}(a_M)=\sum_{L_1,L_2\in  \lc^G(M)}d_M^G(L_1,L_2)\cdot J_{M,V}^{L_1,s}(a_M)\cdot J_{M,V}^{L_2,s}(a_M).$$
\end{proposition}

\begin{preuve}
 Arthur a défini une section $(L_1,L_2)\mapsto (Q_{L_1},Q_{L_2})$ de l'application $\fc(M)\times\fc(M)\to \lc(M)\times\lc(M)$ donné par $(Q_1,Q_2)\mapsto (M_{Q_1},M_{Q_2})$ définie sur l'ensemble des $(L_1,L_2)$ tels que $d_M^G(L_1,L_2)\not=0$. Cette section a la propriété suivante : pour tout $g\in (J_{a_M}(\AAA_{V})\back G(\AAA_{V}))^\tau$ d'image inverse $(g_1,g_2)$ par la bijection  canonique
 \begin{equation}
   \label{eq:bijcan}
    (J_{a_M}(\AAA_{V_1})\back G(\AAA_{V_1}))^\tau \times (J_{a_M}(\AAA_{V_2})\back G(\AAA_{V_2}))^\tau \to (J_{a_M}(\AAA_{V})\back G(\AAA_{V}))^\tau,
  \end{equation}
  on a
  \begin{equation}
    \label{eq:scindage-poids}
    v_M^G(g)=\sum_{L_1,L_2\in  \lc^G(M)}d_M^G(L_1,L_2) \cdot v_M^{Q_{L_1}}(g_1) \cdot v_M^{Q_{L_2}}(g_2),
  \end{equation}
où pour tout $Q\in \fc(M)$ et tout $h\in G(\AAA)$ on pose
$$v_M^Q(h)=v_M^{M_Q}(l)$$
où $l$ est l'unique élément de $M(\AAA)\back M(\oc)$ tel que $g\in l N(\AAA) G(\oc)$ (décomposition d'Iwasawa). Introduisons l'intégrale $J_{M,V}^{Q,s}(a_M)$ donnée par l'égalité (\ref{eq:JMGV}) où l'on a remplacé le poids $v_M^G$ par le poids $v_M^Q$. Comme la bijection (\ref{eq:bijcan}) préserve les mesures, la formule (\ref{eq:scindage-poids}) implique qu'on a 
$$J_{M,V}^{G,s}(a_M)=\sum_{L_1,L_2\in  \lc^G(M)}d_M^G(L_1,L_2)\cdot J_{M,V}^{Q_{L_1},s}(a_M)\cdot J_{M,V}^{Q_{L_2},s}(a_M).$$
Pour conclure, il suffit de montrer que pour tout $L\in \lc(M)$ et tout $Q\in \pc(L)$, on a 
$$J_{M,V}^{Q,s}(a_M)=J_{M,V}^{L,s}(a_M).$$
On montre à l'aide du théorème de l'isogénie de Lang qu'on a une décomposition
\begin{equation}
  \label{eq:decomposition}
  (J_{a_M}(\AAA_{V})\back G(\AAA_{V}))^\tau= (J_{a_M}(\AAA_{V})\back L(\AAA_{V}))^\tau N_Q(\AAA_V)^\tau G(\oc_V)^\tau
\end{equation}
qui est compatible aux mesures pourvu que $N_Q(\AAA_V)^\tau$ soit muni de la mesure de Haar qui donne le volume $1$ au sous-groupe  $(N_Q(\oc_V))^\tau$. Par le changement de variable $g=lnk$ avec des notations qu'on espère évidentes, on voit que $J_{M,V}^{Q,s}(a_M)$ est égale à 
\begin{eqnarray*}
   q^{-\deg(D_V)|\Phi^G|/2} |D^G(X_a)|^{1/2}_V \times \\ \nonumber \int_{(J_{a_M}(\AAA_V)\back L(\AAA_V))^\tau} \int_{N_Q(\AAA_V)^\tau} \int_{ G(\oc_V)^\tau } \bg s, \inv_{a_M}(lnk)\bd \mathbf{1}_D(\Ad((lnk)^{-1})X_{a_M}) v_M^L(l) \, dl\,dn \, dk.
\end{eqnarray*}
Par la définition \ref{def:invaM}, l'invariant  $\inv_{a_M}(lnk)$ ne dépend que de $lnk\tau(lnk)^{-1}=l\tau(l)^{-1}$. On a donc  $\inv_{a_M}(lnk)= \inv_{a_M}(l)$. L'intégrande est alors  $G(\oc_V)^\tau $-invariante et l'intégrale sur $G(\oc_V)^\tau$ qui vaut $1$ peut être omise. On utilise ensuite le changement de variable $\Ad((ln)^{-1})X_{a_M}=\Ad(l^{-1})X_{a_M}+U$ (où $U$ appartient à $\ngo_Q(\AAA_V)^\tau$ qui est muni de la mesure de Haar qui donne le volume $1$ à $\ngo_Q(\oc_V)^\tau$) de jacobien  $|D^G(X_a)|^{-1/2}_V |D^L(X_a)|^{1/2}_V$ pour obtenir que $J_{M,V}^{Q,s}(a_M)$ est égale à
\begin{eqnarray*}
   q^{-\deg(D_V)|\Phi^G|/2} |D^L(X_a)|^{1/2}_V \times \\ \nonumber \int_{(J_{a_M}(\AAA_V)\back L(\AAA_V))^\tau} \int_{\ngo_Q(\AAA_V)^\tau}  \bg s, \inv_{a_M}(l) \bd \mathbf{1}_D(\Ad(l^{-1})X_{a_M}+U)\,  v_M^L(l) \, dl\,dU.
\end{eqnarray*}
Comme on a 
\begin{eqnarray*}
  \int_{\ngo_Q(\AAA_V)^\tau}\mathbf{1}_D(\Ad(l^{-1})X_{a_M}+U)\,dU&=& \mathbf{1}_D(\Ad(l^{-1})X_{a_M})  \int_{\ngo_Q(\AAA_V)^\tau}\mathbf{1}_D(\Ad(l^{-1})X_{a_M}+U)\, dU\\
&=& q^{\deg(D_V)\dim(\ngo_Q)} \mathbf{1}_D(\Ad(l^{-1})X_{a_M})
\end{eqnarray*}
et $ q^{-\deg(D_V)|\Phi^G|/2}q^{\deg(D_V)\dim(\ngo_Q)}=q^{-\deg(D_V)|\Phi^L|/2}$, on obtient bien $J_{M,V}^{Q,s}(a_M)=J_{M,V}^{L,s}(a_M)$ ce qui conclut la démonstration.

\end{preuve}
\end{paragr}

\begin{paragr} \label{S:triv-IOP}Soit $a_M=(h_{a_M},t)\in \Ac_M(k)^\tau$. Soit $h_{a_{M},\eta}$ la restriction de $h_a$ au point générique de la courbe $C$. On voit alors $h_{a_{M},\eta}$ comme un point dans $\car_M(F)$ où $F$ est le corps des fonctions de $F$. Soit $V$ un ensemble $\tau$-stable de points fermés de $C$. Soit
$\varpi^{D_V}=\prod_{v\in v} \varpi_v^{d_v}\in \oc_V$. On a $\varpi^{D_V}h_{a_M,\eta}\in \car_M(\oc_V)$. Soit $X_{a_M,V}=\eps_M(\varpi^{D_V}h_{a_M,\eta})\in \mgo(\oc_V)$ l'image par la section de Kostant $\eps_M$. Rappelons qu'on a $X_{a_M}=\eps_M(h_{a_M,\eta})\in \mgo(F)$. Hors du support de $D$, ces deux éléments ne sont en général pas égaux. Néanmoins, ils satisfont toujours la relation suivante 
  $$X_{a_M,V}= \varpi^{D_V} \Ad(\rho_M(\varpi^{D_V/2})) X_{a_M},$$
où $\rho_M$ est la somme des coracines dans $M$ positives pour le choix de l'épinglage qui détermine la section de Kostant $\eps_M$. On rappelle que le diviseur $D$ est pair. L'élément $\rho_M(\varpi^{D_V/2})$ est $\tau$-stable. Il conjugue les centralisateurs $J_{a_M}$ et $J_{X_{a_M,V}}$. La translation par $\rho_M(\varpi^{D_V/2})$ envoie $(J_{a_M}(\AAA_V)\back G(\AAA_V))^\tau$ bijectivement sur $(J_{a_M,V}(\AAA_V)\back G(\AAA_V))^\tau$ qu'on munit de la mesure obtenue par transport. De même, on note $\inv_{a_M,V}$ l'application définie  sur $(J_{a_M,V}(\AAA_V)\back G(\AAA_V))^\tau$ et qui, composée avec la bijection précédente, redonne $\inv_{a_M}$. Soit $\mathbf{1}_{\ggo(\oc)}$ la fonction caractéristique de $\ggo(\oc)$. Par un changement de variables, on obtient le lemme suivant.

\begin{lemme} \label{lem:triv-IOP}Les notations sont celles utilisées ci-dessus. Pour tout $s\in \Tc^{W_{a_M}}$, on a l'égalité 

$$  J_{M,V}^{G,s}(a_M)= |D^G(X_{a_M,V})|^{1/2}_V  \int_{(J_{a_M,V}(\AAA_V)\back G(\AAA_V))^\tau} \bg s, \inv_{a_M,V}(g)\bd \mathbf{1}_{\ggo(\oc)}(\Ad(g^{-1})X_{a_M,v}) v_M^G(g) \, dg$$
\end{lemme}

\end{paragr}

\begin{paragr}[Calculs locaux d'intégrales orbitales pondérées.] --- Dans tout ce paragraphe, on suppose que $M\in \lc^G$ vérifie $M\not=G$. Soit $a_M=(h_{a_M},t)\in \Ac_M(k)$ et $s\in \Tc^{W_{a_M}}$. On va calculer dans quelques cas simples des intégrales orbitales pondérées locales.

  \begin{lemme}\label{lem:IOPreg}Soit $V$ un ensemble $\tau$-stable de points fermés de $C$.  On suppose que  pour tout $v\in V$, la section $h_{a_M}$ est telle que $h_{a_M}(v)$ n'appartient pas à $\mathfrak{R}_M^G$. Alors
$$  J_{M,V}^{G,s}(a_M)=0.$$
      \end{lemme}

      \begin{preuve}
On utilise les notations du §\ref{S:triv-IOP}. D'après le lemme \ref{lem:triv-IOP}, on a 
$$  J_{M,V}^{G,s}(a_M)= |D^G(X_{a_M,V})|^{1/2}_V  \int_{(J_{a_M,V}(\AAA_V)\back G(\AAA_V))^\tau} \bg s, \inv_{a_M,V}(g)\bd \mathbf{1}_{\ggo(\oc)}(\Ad(g^{-1})X_{a_M,V}) v_M^G(g) \, dg.$$
Soit $g\in G(\AAA_V)$ tel que l'intégrande soit non nulle. On a donc $\Ad(g^{-1})X_{a_M,V}\in \ggo(\oc_V)$. Soit $P\in \pc(M)$. On écrit $g=mnk$ avec $m\in M(\AAA_V)$, $n\in N_P(\AAA_V)$ et $k\in G(\oc_V)$ (décomposition d'Iwasawa). Il vient 
$$\Ad(m^{-1})X_{a_M,V}\in \mgo(\oc_V)$$
et
\begin{equation}
  \label{eq:integrite222}
  \Ad(n^{-1})\Ad(m^{-1})X_{a_M,V} -\Ad(m^{-1})X_{a_M,V}\in \ngo_P(\oc_V).
\end{equation}
L'hypothèse sur $h_{a_M}$ entraîne que l'action adjointe de $\Ad(m^{-1})X_{a_M,V}$ sur $\ngo_P(\oc_V)$ est un $\oc_V$-isomorphisme. Par un argument standard de dévissage, on en déduit que la condition (\ref{eq:integrite222}) implique  $n\in N_P(\oc)$. Mais alors on a $v_M^G(g)=v_M^G(n)=v_M^G(1)=0$ (car $M\not=G$). Donc l'intégrande est toujours nulle et l'intégrale aussi \emph{a fortiori}.
      \end{preuve}
  
      \begin{lemme}\label{lem:IOP1}
Soit $v$ un point fermé de $C$ et $n\geq 1$ le plus petit entier tel que $\tau^n(v)=v$. Soit $V$ l'orbite de $v$ sous l'action de $\tau$. On suppose que $v$ vérifie les conditions suivantes :
\begin{itemize}
\item $h_{a_M}(c)$ est un point lisse de $\mathfrak{R}_M^G$;
\item $h_{a_M}(c)$ n'appartient pas à $\mathfrak{D}^M$ ;
\item l'intersection de  $h_{a_M}(C)$ avec $\mathfrak{R}_M^G$ est transverse en  $h_{a_M}(c)$.
\end{itemize}
On a alors les conclusions suivantes.
\begin{enumerate}
\item Il existe  $Y_v\in \tgo(\oc_v)$ tel que $\chi^T_M(Y_v)=\chi_M(X_{a_M,v})$. Pour tout tel $Y_v$, il existe un unique couple $\{\al,-\al\}\subset \Phi_T^G-\Phi_T^M$ tel que $\val_v(\pm \al(Y_v))=1$.  Soit $w$ l'unique élément de $W^M$ tel que $w\tau(Y_v)=Y_v$. Alors $w(\al)=\al$.

\item Avec les notations ci-dessus, Soit $L\in \lc^G$ défini par $\Phi^L_T=\{\al,-\al\}$ et $(B_L,T,X_\al)$ un épinglage de $L$. Soit $w$ l'unique automorphisme de $L$ qui fixe l'épinglage et qui est de la forme $\Ad(m)$ où $m\in \Norm_{M(\oc_v)}(T)$ est un relèvement de $w$.  
On a 
$$  J_{M,V}^{G,s}(a_M)=d_T^G(M,L)|D^L(Y_v)|^{|V|/2} |V|^{\dim(a_M^G)}\int_{T(F_v)^{w\tau^n}\back L(F_v)^{w\tau^n}} \mathbf{1}_{\lgo(\oc_v)}(\Ad(l^{-1})Y_v) v_T^L(l) \,dl,$$
où les groupes $T(F_v)^{w\tau^n}$ et $ L(F_v)^{w\tau^n}$ sont munis des mesures de Haar qui donnent le volume $1$ respectivement aux groupes  $T(\oc_v)^{w\tau^n}$ et $ L(\oc_v)^{w\tau^n}$.

\end{enumerate}

      \end{lemme}

      \begin{preuve} Puisque $h_{a_M}(c)$ n'appartient pas à $\mathfrak{D}^M$, on a $\chi_M(X_{a_M,v})\in \car_{M}^{M-\reg}(\oc_v)$. Sa réduction mod $\varpi_v$ définit donc un élément de $\car_{M}^{M-\reg}(k)$ qui se relève évidemment en un élément $y$ de $\tgo^{M-\reg}(k)$. Or sur l'ouvert  $\tgo^{M-\reg}$, on sait que $\chi^T_M$ est étale. Il s'ensuit que $y$ se relève en un unique élément $Y_v\in  \tgo^{M-\reg}(\oc_v)$ tel que $\chi^T_M(Y_v)=\chi_M(X_{a_M,v})$. 

Soit $P\in \pc(M)$. Les hypothèses sur l'intersection de $h_{a_M}(C)$ avec $\mathfrak{R}_M^G$ entraînent qu'il existe une unique racine $\al$ de $T$ dans $N_P$ tel que $\val_v(\al(Y_v))=1$, les autres racines de $T$ dans $N_P$ étant de valuation nulle en $Y_v$. En raisonnant avec le sous-groupe parabolique opposé à $P$, on voit que $\al$ et $-\al$ sont les seules racines  hors $M$ telles que $\val_v(\pm \al(Y_v))=1$. Toutes les autres racines y compris celles dans $M$ (vu qu'on a  $h_{a_M}(c)\notin\mathfrak{D}^M$) ont pour valuation $0$ en $Y_v$. Soit $Y_v$ et $\al\in \Phi^G_T$ qui vérifient l'assertion 1. Il existe $w\in W^M$ tel que $w\tau(Y_v)=Y_v$ puisque $ \chi^T_M(Y_v)\in \car_M(\oc_v)^\tau$. On a $\val_v(\al(Y_v))= \val_v(\tau(\al(Y_v)))= \val_v(\tau(\al)(Y_v))$. On a donc nécessairement $w(\al)=\pm\al$ mais comme $\al$ est une racine dans le radical unipotent d'un sous-groupe parabolique $P\in \pc(M)$, il en est de même de $w(\al)$ et donc $w(\al)=\al$. Cela démontre donc l'assertion 1.

Prouvons ensuite l'assertion 2. Partons de l'égalité (cf. §\ref{S:triv-IOP} et lemme \ref{lem:triv-IOP})

$$  J_{M,V}^{G,s}(a_M)= |D^G(X_{a_M,V})|^{1/2}_V  \int_{(J_{a_M,V}(\AAA_V)\back G(\AAA_V))^\tau} \bg s, \inv_{a_M,V}(g)\bd \mathbf{1}_{\ggo(\oc)}(\Ad(g^{-1})X_{a_M,V}) v_M^G(g) \, dg.$$

Soit  $g\in G(\AAA_V)$ et $j\in J_{a_M,V}(\AAA_V)$ tel que $\tau(g)=jg$. On va d'abord prouver que si $\mathbf{1}_{\ggo(\oc_V)}(\Ad(g^{-1})X_{a_M,V})\not=0$  alors $\inv_{a_M,V}(g)=0$. Pour cela, il suffit de prouver que l'image de $j$ dans les coinvariants $J(\AAA_V)_{\tau}$ est triviale. Soit $P\in \pc(M)$ et $g=mnk$ la décomposition d'Iwasawa de $g$ associée à $P$. On a donc $\tau(m)\in jmM(\oc_V)$. Par une application du théorème de l'isogénie de Lang, on voit que, quitte à translater $m$ par un élément de $M(\oc_V)$, on a $\tau(m)=jm$ ce que l'on suppose désormais. Comme dans la preuve du lemme  \ref{lem:triv-IOP}, la condition $\Ad(g^{-1})X_{a_M,V}\in \ggo(\oc_V)$ entraîne
$$\Ad(m^{-1})X_{a_M,V}\in \mgo(\oc_V).$$
Comme  $h_{a_M}(c)$ n'appartient pas à $\mathfrak{D}^M$, on a $\chi_M(X_{a_M,V})\in \car_M^{M-\reg}(\oc_V)$. Par conséquent, on a $m\in J_{a_M,V}(\AAA_V)M(\oc_V)$. Comme $\inv_{a_M,V}$ est invariante par translation à gauche par  $J_{a_M,V}(\AAA_V)$, on peut et on va supposer que $m\in M(\oc_V)$. Mais alors on a $j\in M(\oc_V)\cap J_{a_M,V}(\AAA_V)$. Soit $J$ le $\oc_V$-schéma en groupes défini comme le centralisateur de $X_{a_M,V}$ dans $M\times_k \AAA_V$. Comme $\chi_M(X_{a_M,V})\in \car_M^{M-\reg}(\oc_V)$, c'est en fait un schéma en tores. On a $M(\oc_V)\cap J_{a_M,V}(\AAA_V)=J(\oc_V)$ et une nouvelle application du théorème de Lang montre que $j$ est trivial dans les coinvariants.  Il s'ensuit qu'on peut restreindre l'intégrale aux $g\in (J_{a_M,V}(\AAA_V)\back G(\AAA_V))^\tau$ d'invariant $\inv_{a_M,V}(g)$ trivial c'est-à-dire qu'on peut restreindre l'intégrale à $J_{a_M,V}(\AAA_V)^\tau \back G(\AAA_V)^\tau$. On a donc

\begin{eqnarray*}
            J_{M,V}^{G,s}(a_M)&=& |D^G(X_{a_M,V})|^{1/2}_V  \int_{J_{a_M,V}(\AAA_V)^\tau\back G(\AAA_V)^\tau} \mathbf{1}_{\ggo(\oc_V)}(\Ad(g^{-1})X_{a_M,V}) v_M^G(g) \, dg\\
            &=& |D^G(X_{a_M,v})|^{|V|/2}_v  |V|^{\dim(a_M^G)}\int_{J_{a_M,v}(F_v)^{\tau^n}\back G(F_v)^{\tau^n}} \mathbf{1}_{\ggo(\oc_v)}(\Ad(g^{-1})X_{a_M,v}) v_M^G(g) \, dg,
\end{eqnarray*}
où l'on utilise la bijection évidente $J_{a_M,V}(\AAA_V)^\tau\back G(\AAA_V)^\tau \to J_{a_M,v}(\AAA_v)^{\tau^n}\back G(\AAA_v)^{\tau^n}$. Le facteur   $|V|^{\dim(a_M^G)}$ provient de la comparaison des poids. 

Il est clair aussi qu'on a 
$$|D^G(X_{a_M,v})|_v=|D^G(Y_{v})|_v=|D^L(Y_v)|_v.$$

\emph{Désormais, pour alléger les notations on note $\tau$ au lieu de $\tau^n$.}  Puisque $Y_v$ et $X_{a_M,v}$ ont même image par $\chi_M$ et que cette image appartient à l'ouvert $\car_M^{M-\reg}(\oc_v)$, il existe $m\in M(\oc_v)$ tel que $X_{a_M,v}=\Ad(m^{-1})Y_v$. Il s'ensuit que $m \tau(m)^{-1}$ est un élément du normalisateur de $T$ dans $M(\oc_v)$ et que son image dans $W^M$ est $w$. On en déduit que $m\tau(m)^{-1}$ normalise $L$. L'automorphisme $\Int(m\tau(m)^{-1})$ préserve le couple $(B_L,T)$ de l'épinglage $(N_L,T,X_\al)$ et agit sur $X_\al$ par homothétie par un élément de $\oc_v^\times$. Si l'on change $m$ en $tm$ avec $t\in T(\oc_v)$, le rapport de cette homothétie est multipliée par $\al(t\tau(t)^{-1})$. On voit donc qu'en prenant un élément $t$ judicieux, on peut supposer que  $\Int(m\tau(m)^{-1})$ préserve l'épinglage. 

Il est clair qu'on a 
$$|D^G(X_{a_M,v})|_v=|D^G(Y_{v})|_v=|D^L(Y_v)|_v.$$
Par abus de notations, on note encore $w$ l'automorphisme intérieur de $G$ donné par $\Int(m\tau(m)^{-1})$. L'application $g\mapsto mgm^{-1}$ induit une bijection de $J_{a_M,v}(F_v)^\tau \back G(F_v)^{\tau}$ sur $T(F_v)^{w\tau}\back G(F_v)^{w\tau}$ et l'on a
$$\int_{J_{a_M,v}(F_v)^{\tau}\back G(F_v)^{\tau}} \mathbf{1}_{\ggo(\oc_v)}(\Ad(g^{-1})X_{a_M,v}) v_M^G(g) \, dg=\int_{T(F_v)^{w\tau}\back G(F_v)^{w\tau}} \mathbf{1}_{\ggo(\oc_v)}(\Ad(g^{-1})Y_v) v_M^G(g) \, dg$$
où la mesure sur $T(F_v)^{w\tau}\back G(F_v)^{w\tau}$ se déduit de celle sur $J_{a_M,v}(F_v)^\tau\back G(F_v)^{\tau}$ par transport par la bijection précédente. En utilisant comme précédemment la décomposition d'Iwasawa pour $Q\in \pc(L)$, on voit que si $g\in G(F_v)$ vérifie $\Ad(g^{-1})Y_v\in \ggo(\oc_v)$ alors $g=ls$ avec $l\in L(F_v)$ et $s\in  G(\oc_v)$. Supposons de plus qu'on a  $w\tau(g)=g$. Ainsi $l^{-1}w\tau(l)\in L(\oc_v)$ et par une nouvelle application du théorème de Lang, on voit qu'on peut supposer $w\tau(l)=l$ et  $s\in G(\oc)^{w\tau}$. On a donc
$$\int_{T(F_v)^{w\tau}\back G(F_v)^{w\tau}} \mathbf{1}_{\ggo(\oc_v)}(\Ad(g^{-1})Y_v) v_M^G(g) \, dg=\int_{T(F_v)^{w\tau}\back L(F_v)^{w\tau}} \mathbf{1}_{\lgo(\oc_v)}(\Ad(l^{-1})Y_v) v_M^G(l) \, dl.$$
On utilise ensuite la formule de descente d'Arthur
$$v_M^G(l) =\sum_{R\in \lc^G(T)} d_T^G(M,R) v_T^{Q_R}(l)$$
où $R\in \lc^G(T) \mapsto Q_R\in \pc(R)$ est une certaine section (définie sur les $R$ tels que   $d_T^G(M,R)\not=0$) de l'application $\fc^G(T) \to \lc(T)$ définie par $Q\mapsto M_Q$. On se limite dans ce qui suit à des $R$ tels que  $d_T^G(M,R)\not=0$. Cela entraîne donc $R\not=T$. Le poids  $v_T^{Q_R}$ a été défini dans la démonstration de la proposition \ref{prop:scindage-IOP}. Voici comme il se calcule ici. Le groupe $Q_R\cap L$ est un sous-groupe parabolique de $L$. Comme $L$ est de rang semi-simple $1$, de deux choses l'une
\begin{itemize}
\item soit $Q_R\cap L$ est un sous-groupe de Borel de $L$ ; on a $l\in t N_{Q_R\cap L}(F_v) L(\oc_v)$ avec $t\in T(F_v)$ et donc $v_T^{Q_R}(l)=v_T^R(t)=0$ (puisque $T\not=R$) ;
\item soit $L\subset R$ et alors $v_T^{Q_R}(l)=v_T^R(l)$.
\end{itemize}
On se place désormais dans le second cas. Rappelons que le poids $v_T^R(l)$ est le volume de la projection sur $\ago_T^R$ de l'enveloppe convexe des points $-H_B(l)$ pour $B\in \pc^R(T)$. Or cette enveloppe convexe est ici contenue dans une droite donc le volume est donc nul sauf si $\dim(\ago_T^R)=1$ auquel cas $R=L$. La formule de descente d'Arthur se simplifie donc en 
$$v_M^G(l) =d_T^G(M,L) v_T^{L}(l)$$
ce qui donne le résultat voulu.
      \end{preuve}
\end{paragr}

\begin{paragr}[Calculs locaux d'intégrales orbitales ordinaires.] --- On donne quelques calculs locaux d'intégrales orbitales ordinaires. Soit $v\in V$ un point fermé de $C$ et $V$ l'orbite sous le Frobenius $\tau$ de $v$. Soit  $s\in \Tc$ et  $\Mc'=\Mc_s$. Soit $a_M=(h_{a_M},t)\in \Ac_{M'}(k)^\tau$.

  \begin{lemme}\label{lem:IOreg}
     On suppose que $h_{a_M}(c)$ n'appartient pas à $\mathfrak{D}^{M}$. Alors 
$$J_{M,V}^{M,s}(a_M)=1.$$
  \end{lemme}

  \begin{preuve}
    Partons de l'égalité
$$  J_{M,V}^{M,s}(a_M)= |D^M(X_{a_M,V})|^{1/2}_V  \int_{(J_{a_M,V}(\AAA_V)\back M(\AAA_V))^\tau} \bg s, \inv_{a_M,V}(m)\bd \mathbf{1}_{\mgo(\oc_V)}(\Ad(m^{-1})X_{a_M,V}) \, dm.$$
Soit $m \in M(\AAA_V)$ tel que $\Ad(m^{-1})X_{a_M,V}\in \mgo(\oc_V)$. Comme $X_{a_M,V}$ est encore semi-simple régulier en réduction (vu que $h_{a_M}(c)$ n'appartient pas à $\mathfrak{D}^{M}$), on voit que $m\in J_{a_M,V}(\AAA_V) M(\oc_V)$. Comme dans la démonstration du lemme, on en déduit que  $\inv_{a_M,V}(m)=0$. Comme par ailleurs  $|D^M(X_{a_M,V})|^{1/2}_V=1$, on a
\begin{eqnarray*}
    J_{M,V}^{M,s}(a_M)&= &\int_{J_{a_M,V}(\oc_V)^\tau\back M(\oc_V)^\tau} \mathbf{1}_{\mgo(\oc_V)}(\Ad(m^{-1})X_{a_M,V}) \, dm.\\
&=& \vol(J_{a_M,V}(\oc_V)^\tau)^{-1} \vol(M(\oc_V)^\tau).
\end{eqnarray*}
Or $\vol(M(\oc_V)^\tau)=1$ et $ \vol(J_{a_M,V}(\oc_V)^\tau)=1$ car $J_{a_M,V}(\oc_V)^\tau$ est le sous-groupe compact maximal de $J_{a_M,V}(\AAA_V)^\tau$.
\end{preuve}

  \begin{lemme}\label{lem:IO1}
     On suppose que 
\begin{itemize}
\item $h_{a_M}(c)$ est un point lisse de $\mathfrak{R}_{M'}^M$;
\item $h_{a_M}(c)$ n'appartient pas à $\mathfrak{D}^{M'}$ ;
\item l'intersection de  $h_{a_M}(C)$ avec $\mathfrak{R}_{M'}^M$ est transverse en  $h_{a_M}(c)$.
\end{itemize}
Alors
$$J_{M,V}^{M,s}(a_M)=J_{M',V}^{M',1}=1.$$
  \end{lemme}

  \begin{preuve}
La première égalité n'est qu'un cas particulier du lemme fondamental de Langlands-Shelstad. La second résulte du lemme \ref{lem:IOreg} ci-dessus. Le lecteur qui rechignerait à utiliser le lemme fondamental, peut vérifier directement l'égalité $J_{M,V}^{M,s}(a_M)=1$ : on se ramène à un calcul sur un groupe de rang semi-simple $1$ qu'on fait à la main (cf. aussi \cite{ngo2} §8.5.9).
  \end{preuve}

\begin{lemme}\label{lem:IO2}
     On suppose que 
\begin{itemize}
\item $h_{a_M}(c)$ n'appartient pas à $\mathfrak{R}_{M'}^M$;
\item $h_{a_M}(c)$ est un point lisse de $\mathfrak{D}^{M'}$ ;
\item l'intersection de  $h_{a_M}(C)$ avec $\mathfrak{D}^{M'}$ est transverse en  $h_{a_M}(c)$.
\end{itemize}
Alors
$$J_{M,V}^{M,s}(a_M)=J_{M',V}^{M',1}(a_M)\not=0.$$
  \end{lemme}

  \begin{preuve}
Là encore, la première égalité est un cas particulier du lemme fondamental de Langlands-Shelstad. L'intégrande dans  
$$J_{M',V}^{M',1}= |D^{M'}(X_{a_M,V})|^{1/2}_V  \int_{(J_{a_M,V}(\AAA_V)\back M'(\AAA_V))^\tau} \mathbf{1}_{\mgo'(\oc_V)}(\Ad(m^{-1})X_{a_M,V}) \, dm.$$
est positive et l'intégrale est clairement minoré par $\vol((J_{a_M,V}(\oc_V)\back M'(\oc_V))^\tau)>0$.
Bien sûr, il est aussi possible de mener à bien directement les calculs de $J_{M,V}^{M,s}(a_M)$ et $J_{M',V}^{M',1}$ (on se ramène à un calcul sur un groupe de rang semi-simple $1$, cf. aussi \cite{ngo2} § 8.5.8).
  \end{preuve}

\end{paragr}

\begin{paragr} Soit $V$ un ensemble $\tau$-stable non vide de points fermés de $C$.  Soit  $s\in \Tc$ tel que   $\Mc'=\Mc_s\not=\Mc$. Soit $a_M=(h_{a_M},t)\in \Ac_{M'}(k)^\tau$. Soit $J_V :\ec_M^G(s)\cup \lc^{\Gc}(\Mc) \to \CC$ l'application définie pour tout $\Hc\in  \ec_M^G(s)$ par $J_V(\Hc)=J_{M',V}^{H,1}(a_M)$ et pour tout $\Lc\in \lc^{\Gc}(\Mc)$ par  $J_V(\Lc)=J_{M,V}^{L,s}(a_M)$

  \begin{lemme}\label{lem:final}
On suppose que $V$ satisfait l'une des trois conditions suivantes :
\begin{enumerate}
\item Pour tout $v\in V$, le point   $h_{a_M}(v)$ n'appartient pas à $\mathfrak{D}^{G}$ ;
\item $V$ est l'orbite sous $\tau$ d'un point $v$ qui vérifie l'assertion 1 du théorème \ref{thm:LFPsemilocal} ;
\item $V$ est l'orbite sous $\tau$ d'un point $v$ qui vérifie l'assertion 2 du théorème \ref{thm:LFPsemilocal} .
\end{enumerate}
Alors $J_V$ est stabilisante et $J_V(\Mc)\not=0$.
  \end{lemme}

  \begin{preuve}
Supposons l'assertion 1 vérifiée. D'après le lemme \ref{lem:IOPreg}, on a  $J_V(\Hc)=0$ sauf si $\Hc=\Mc$ ou $\Hc=\Mc'$. De plus,  $J_V(\Mc)=J_V(\Mc')=1$ d'après le lemme \ref{lem:IOreg}. Donc $J_V$ est stabilisante.

Supposons l'assertion 2 vérifiée. Dans ce cas, $h_{a_M}(v)$ est un point lisse de  $\mathfrak{R}_{M'}^G$ et n'appartient  pas à $\mathfrak{D}^{M'}$ et l'intersection de $h_{a_M}(v)$  avec $\mathfrak{R}_{M'}^G$ est transverse en $h_{a_M}(v)$.  Supposons, de plus, pour commencer que $h_{a_M}(v)$ n'appartient pas à   $\mathfrak{R}_{M}^G$. Par conséquent, l'intersection de $h_a(C)$ avec $\mathfrak{R}_{M'}^M$ est transverse
en $h_{a_M}(v)$. On a donc $J_V(\Mc)=J_V(\Mc')=1$ d'après le lemme \ref{lem:IO1}. D'après le lemme \ref{lem:IOPreg}, on a $J_V(\Lc)=0$ pour tout $\Lc\in \lc^{\Gc}(\Mc)$ et $\Lc\not=\Mc$. Pour tout  $\Hc\in  \ec_M^G(s)$, on a $\Hc\cap \Mc=\Mc'$ ; par conséquent, $h_{a_M}(v)$ n'appartient pas à $\mathfrak{R}_{M'}^H$ et on a $J_V(\Hc)=0$ par le  lemme \ref{lem:IOPreg}. Ainsi  $J_V$ est stabilisante.

Supposons toujours l'assertion 2 vérifiée mais cette fois on suppose que $h_{a_M}(v)$ appartient à   $\mathfrak{R}_{M}^G$ : sous nos hypothèses,  $h_{a_M}(v)$ est un point d'intersection transverse de $h_{a_M}(C)$ et $\mathfrak{R}_{M}^G$ et il n'appartient pas à  $\mathfrak{D}^{M}$. 

Soit $Y_v\in \tgo(\oc_v)$, $L\in \lc^G$, $n$ un entier et $w$ un automorphisme de $L$ qui préserve $T$ tels que ces éléments vérifient les conclusions du lemme \ref{lem:IOP1}. Introduisons l'application $J_T : \ec_T^G(1) \to \CC$ définie par $J_T(\Tc)=1$,
$$J_T(\Lc)=|D^L(Y_v)|^{|V|/2} |V|^{\dim(a_M^G)}\int_{T(F_v)^{w\tau^n}\back L(F_v)^{w\tau^n}} \mathbf{1}_{\lgo(\oc_v)}(\Ad(l^{-1})Y_v) v_T^L(l) \,dl,$$
et $J_T(\Hc)=0$ si $\Hc\notin \{\Tc,\Lc\}$. Pour tout $\Hc \in  \lc^{\Gc}(\Mc)$, l'expression 
\begin{equation}
  \label{eq:lexpression}
  \sum_{\Rc \in \lc^{\Hc}(\Tc)} d_{\Tc}^{\Hc}(\Mc,\Rc) J_T(\Rc)
\end{equation}
vaut
\begin{itemize}
\item $J_{T}(\Tc)=1=J_V(\Mc)$ si $\Hc=\Mc$ (cf. lemme \ref{lem:IOreg});
\item $d_{\Tc}^{\Hc}(\Mc,\Lc) J_T(\Lc)=J_V(\Hc)$ si $\Hc\not=\Mc$  (cf. lemme \ref{lem:IOP1}).
\end{itemize}
Donc l'expression (\ref{eq:lexpression}) est égale à $J_V(\Hc)$. Pour tout $\Hc \in  \ec_M^G(s)$, on a $h_{a_M}(v)\in \mathfrak{R}_{M'}^{H}$ si et seulement si $\Lc\subset \Hc$ et  l'expression 
\begin{equation}
  \label{eq:lexpression2}
 \sum_{\Rc \in \lc^{\Hc}(\Tc)} d_{\Tc}^{\Hc}(\Mc',\Rc) J_T(\Rc)=J_V(\Hc)
\end{equation}
vaut
\begin{itemize}
\item $J_{T}(\Tc)=1=J_V(\Mc')$ si $\Hc=\Mc'$ (cf. lemme \ref{lem:IOreg});
\item $d_{\Tc}^{\Hc}(\Mc,\Lc) J_T(\Lc)=J_V(\Hc)$ si $\Lc\subset \Hc$ (cf. lemme \ref{lem:IOP1}) ;
\item $0$ si $\Lc\not\subset \Hc$ ( lemme \ref{lem:IOPreg}).
\end{itemize}
En tout cas l'expression (\ref{eq:lexpression2}) est égale à $J_V$. Comme $J_T$ est (trivialement) stabilisante, $J_V$ est stabilisante (cf. proposition \ref{prop:descente-stab}). Le lemme \ref{lem:IOreg} donne $J_V(\Mc)=J_V(\Mc')=$.

Supposons maintenant l'assertion 3 vérifiée. Dans ce cas, le lemme \ref{lem:IOPreg} indique que $J_V(\Hc)=0$ sauf si $\Hc=\Mc$ ou $\Mc'$. On a $J_V(\Mc)=J_V(\Mc')\not=0$ par le lemme \ref{lem:IO2}. Donc $J_V$ est stabilisante.

  \end{preuve}
\end{paragr}

\section{Fibres de Springer affines tronquées}\label{sec:FSAT}

\begin{paragr}Dans cette section, les notations sont celles utilisées dans tout l'article. On ajoute les précisions suivantes. Soit $F=k((\eps))$ et $\oc=k[[\eps]]$. Soit $\val$ la valuation usuelle de $F$. Soit $F_0=\Fq((\eps))$. On rappelle qu'on note $\tau$ l'automorphisme de Frobenius de $k$ donné par l'élévation à la puissance $q$. On note encore $\tau$ l'automorphisme de Frobenius de $F$ de sorte qu'on a $F^\tau=F_0$. Soit $\overline{F}$ une clôture séparable de $F$ et $I=\Gal(\overline{F}/F)$. Soit $\overline{F}_0$ la clôture algébrique de $F_0$ dans  $\overline{F}$. Soit $\Gamma=\Gal( \overline{F}_0/F_0)$.

On suppose que le couple $(G,T)$ est défini sur $\Fq$. Pour tout $P\in \pc(M)$, on définit comme au §\ref{S:poids-Arthur} une application
$$H_P : G(F) \to X_*(M)$$
définie pour tout $\la\in X^*(P)$ par $\la(H_P(g))=-\val(\la(p))$ où $p\in P(F)$ satisfait $g\in pG(\oc)$ (par la décomposition d'Iwasawa).

On rappelle que $\ago_T$ est muni d'un produit scalaire $W^G$-invariant et que tous ses sous-espaces sont munis de la mesure euclidienne qui s'en déduit (cf. §\ref{S:parab}). Soit 
$$H_P^G : G(F) \to \ago_M^G$$
le composé de l'application $H_P$ avec la projection orthogonale $\ago_M\to \ago_M^G$. Les restrictions à $M(F)$ de ces applications ne dépendent pas du choix de $P\in \pc(M)$ : on les note respectivement $H_M$ et $H_M^G$.

\end{paragr}

\begin{paragr}[Grassmannienne affine tronquée.] --- Soit $\Xgo^G$ la grassmannienne affine : c'est un ind-$k$-schéma projectif dont l'ensemble des $k$-points est le quotient $G(F)/G(\oc)$. Soit $M\in \lc$.  Pour tout $P\in \pc(M)$, l'application $H_P^G$ sur $G(F)$ induit une application $H_P^G :\Xgo^G \to \ago_M^G$. Soit $\xi\in \ago_M^G$. On définit la grassmannienne affine tronquée 
$$\Xgo^{G,\xi}_M$$ 
par la condition suivante : pour tout $x\in \Xgo^G$, on a $x\in \Xgo^{G,\xi}_M$ si et seulement si $\xi$ appartient à l'enveloppe convexe des projections sur $\ago_M^G$ des points $-H_P^G(x)$ pour $P\in \pc(M)$. 

\begin{proposition}
  La grassmannienne affine tronquée $\Xgo^{G,\xi}_M$ est un sous-ind-schéma ouvert de $\Xgo^G$.
\end{proposition}

\begin{preuve}
  On sait bien que l'enveloppe convexe des points $-H_P^G(x)$ est l'intersection des cônes $-H_P^G(x) -\overline{^+\ago_P}$, où $\overline{^+\ago_P}$ est la chambre de Weyl obtuse positive dans $\ago_M$ (c'est un résultat d'Arthur, cf. \cite{dis_series}). Il suffit donc de montrer que la condition  $H_P^G(x)\in -\xi  -\overline{^+\ago_P}$ définit un ouvert de  $\Xgo^G$. Pour cela, on procède comme dans \cite{CL}, démonstration de la proposition 6.4.
\end{preuve}

\end{paragr}

\begin{paragr}[Fibre de Springer affine tronquée.] --- Soit $Y\in \mgo(\oc)$ un élément semi-simple et $G$-régulier. À la suite de Kazhdan-Lusztig (cf. \cite{KL}), on introduit la fibre de Springer affine 
$$\Xgo^G_Y=\{x\in \Xgo^G \mid \Ad(x^{-1})Y\in \ggo(\oc)\}.$$
C'est un sous-ind-schéma fermé de $\Xgo^G$. En fait, $\Xgo^G_Y$ est un $k$-schéma localement de type fini et de dimension finie (cf. \cite{KL}). La fibre de Springer affine tronquée $\Xgo^{G,\xi}_{M,Y}$ est l'ouvert de $\Xgo^G_Y$ défini par
$$\Xgo^{G,\xi}_{M,Y}=\Xgo^G_Y\cap \Xgo^{G,\xi}_M.$$
Soit $J$ le centralisateur de $Y$ dans $G\times_k F$. C'est un sous-$F$-tore maximal de $G\times_k F$. Le groupe $J(F)$ agit sur l'ensemble des $k$-points de  $\Xgo^G_Y$. Il ne respecte pas le tronqué $\Xgo^{G,\xi}_{M,Y}$. Soit $J'$ le sous-$F$-tore de $J$ défini par
$$X_*(J')= \Ker ( \la\in X_*(J) \mapsto H_M^G(\eps^\la)).$$
Le groupe $J'(F)$ agit sur  $\Xgo^G_Y$ et respecte le tronqué $\Xgo^{G,\xi}_{M,Y}$. Soit 
$$\Lambda = X_*(J')^I$$
le groupe des cocaractères du sous-tore déployé maximal de $J'$. Via le morphisme $\la \mapsto \eps^\la$, on identifie le groupe discret $\Lambda$ à un sous-groupe de $J'(F)$.

\begin{proposition}
  Le quotient $\Lambda\back \Xgo^{G,\xi}_M(Y)$ est une variété quasi-projective définie sur $k$.
\end{proposition}

\begin{preuve}
Cela se démontre comme dans \cite{KL}.
  
\end{preuve}

\end{paragr}

\begin{paragr}[Une $s$-intégrale orbitale pondérée locale.] --- On continue avec les notations précédentes en supposant désormais qu'on a $\tau(Y)=Y$. Dans ce cas, $J$ est muni d'une action de $\tau$ et on a une suite exacte en cohomologie d'ensemble pointé 

\begin{equation}
  \label{eq:suite-longue-locale}
  1  \longrightarrow J(F)^\tau \longrightarrow G(F)^\tau \longrightarrow  (J(F)\back G(F))^\tau  \longrightarrow \Ker(J(F)_\tau \to G(F)_\tau)  \longrightarrow 1,
\end{equation}  
où comme au §\ref{S:sIOP-notations} l'indice $\tau$ désigne l'ensemble des classes de $\tau$-conjugaison.
On munit les groupes  $J(F)^\tau$ et $G(F)^\tau$ des mesures de Haar qui donnent la mesure $1$ respectivement au sous-groupe compact maximal de $J(F)^\tau$ et au sous-groupe  $G(\oc)^\tau$ de $G(F)^\tau$. D'après Kottwitz (??????), on a une suite exacte canonique
\begin{equation}
  \label{eq:suite-Kott}
  0 \longrightarrow J(F)_1 \longrightarrow J(F) \longrightarrow X_*(J)_{I}\longrightarrow 0
\end{equation}
qui induit au niveau des co-invariants sous $\tau$ un isomorphisme
\begin{equation}
  \label{eq:iso-Kott}
  J(F)_\tau \simeq X_*(J)_{\Gamma}.
\end{equation}
En particulier, le groupe  $J(F)_\tau$ est discret. On le munit de la mesure de comptage. On en déduit une mesure invariante à droite par $G(F)^\tau$ sur le quotient  $(J(F)\back G(F))^\tau$.

Soit 
$$\inv_Y : (J(F)\back G(F))^\tau \to X_*(J)_{\Gamma}$$
le composé de l'application cobord de (\ref{eq:suite-longue-locale}) avec l'isomorphisme (\ref{eq:iso-Kott}) de Kottwitz.

\begin{definition}
  Pour $Y$ et $M$ comme ci-dessus et $s\in \hat{J}^\Gamma$, la $s$-intégrale orbitale pondérée $J_M^{G,s}(Y)$ est définie par
  \begin{equation}
    \label{eq:def-sIOPlocale}
    J_M^{G,s}(Y)=|D^G(Y)|^{1/2}_F \int_{(J(F)\back G(F))^\tau} \bg s, \inv_Y(g) \bd\, \mathbf{1}_{\ggo(\oc)}(\Ad^{-1}(Y)) \, v_M^G(g) \,dg
  \end{equation}
\end{definition}
où
\begin{itemize}
\item la valeur absolue est $|\cdot|_F=q^{-\val(\cdot)}$ ;
\item l'accouplement est celui entre $\hat{J}^\Gamma$ et $ X_*(J)_{\Gamma}$ ;
\item $\mathbf{1}_{\ggo(\oc)}$ est la fonction caractéristique de $\ggo(\oc)$ ;
\item le poids d'Arthur $v_M^G(g)$ est le volume dans $\ago_M^G$ de l'enveloppe convexe des points $-H_P^G(g)$ pour $P\in \pc(M)$.
\end{itemize}

\end{paragr}

\begin{paragr}[Interprétation géométrique des intégrales orbitales pondérées.] --- Puisque $Y$ est $\tau$-stable, le quotient $\Lambda \back  \Xgo^{G,\xi}_{M,Y}$ est défini sur $\Fq$. Soit $\eta\in \Hom_{\ZZ}(\Lambda,\Qb^\times)$ un caractère d'ordre fini. On suppose que ce caractère est fixe par $\tau$ ou, ce qui revient au même, fixe par $\Gamma$. Par le choix d'un plongement de $\Qb$ dans $\Qlb$, on interprète $\eta$ comme un caractère $\ell$-adique continu de $\Lambda\rtimes \Gal(k/\Fq)$ trivial sur le second facteur. On en déduit un système local $\lc_\eta$ sur  $\Lambda \back  \Xgo^{G,\xi}_{M,Y}$ défini sur $\Fq$ qui devient trivial lorsqu'on le tire par le revêtement $\Xgo^{G,\xi}_{M,Y} \to  \Lambda \back  \Xgo^{G,\xi}_{M,Y}$. L'action de $J'(F)$ sur  $\Xgo^{G,\xi}_{M,Y}$ commute à celle de $\Lambda$. On en déduit une action de $J'(F)$ sur la cohomologie $H^\bullet(\Xgo^{G,\xi}_{M,Y},\lc_\eta)$. Par un lemme d'homotopie, celle-ci se factorise par son groupe de composantes connexes. Ce dernier s'identifie au groupe $X_*(J')_I$ par la suite (\ref{eq:suite-Kott}) appliqué au $F$-tore $J'$. Or le sous-groupe discret $\Lambda$ de $J'(F)$ agit par le caractère $\eta^{-1}$ sur la cohomologie $H^\bullet(\Xgo^{G,\xi}_{M,Y},\lc_\eta)$. Comme l'image de $\Lambda$ dans $X_*(J')_I$ est d'indice fini, on en déduit la décomposition suivante
  \begin{equation}
    \label{eq:decomposition-cohom-FSAT}
    H^\bullet(\Xgo^{G,\xi}_{M,Y},\lc_\eta)=\bigoplus_{s'} H^\bullet(\Xgo^{G,\xi}_{M,Y},\lc_\eta)_{s'}
  \end{equation}
  où la somme est prise sur l'ensemble fini des $s'\in (\hat{J}')^{I}=\Hom(X_*(J')_I,\CC^\times)$ qui induisent le caractère $\eta^{-1}$ sur $\Lambda$ et où  $H^\bullet(\Xgo^{G,\xi}_{M,Y},\lc_\eta)_{s'}$ désigne le facteur direct de $H^\bullet(\Xgo^{G,\xi}_{M,Y},\lc_\eta)$ sur lequel $J'(F)$ agit via le caractère $s'$ de $X_*(J')_I$.  Soit $p_M^G$ la projection orthogonale de $\ago_M$ sur $\ago_M^G$. Le réseau $H_M^G(J(F)^\tau)$ est d'indice fini dans  $p_M^G(X_*(M))$. Soit $\mu_1,\ldots,\mu_n$ un système de représentants du quotient. Soit
$$\psi : \Lambda_\tau \to X_*(J')_\Gamma \simeq J'(F)_\tau$$ 
et
$$\varphi : X_*(J')_\Gamma \simeq J'(F)_\tau  \to  X_*(J)_\Gamma \simeq J(F)_\tau$$
les applications canoniques, les isomorphismes étant ceux donnés par (\ref{eq:iso-Kott}). On notera que les noyaux ainsi que le conoyau de $\psi$ sont finis. Soit 
$$ J(F)^{0}=J(F)\cap \Ker(H_M^G).$$
 Le quotient $J(F)^{0,\tau}\back J(F)^\tau$ est discret et il est muni de la mesure de comptage. On en déduire une mesure de Haar sur  $J(F)^{0,\tau}$.

\begin{theoreme}
  \label{thm:IOP-FSAT}
Soit $s\in \hat{J}^\Gamma$ un élément d'ordre fini. Soit $s'\in  (\hat{J}')^\Gamma$ l'image de $s$ par le morphisme canonique $\hat{J}^\Gamma\to(\hat{J}')^\Gamma$. Soit $\eta\in \Hom_{\ZZ}(\Lambda,\Qb^\times)$ le caractère induit par l'inverse de $s'$.

Pour tout $\xi$ \emph{en position générale}, on a l'égalité 
$$J_M^{G,s}(Y)=c_M^G(Y) \sum_{i=1}^n  \trace(\tau^{-1},H^\bullet(\Xgo^{G,\xi+\mu_i}_{M,Y},\lc_\eta)_{s'})$$
où
$$c_M^G(Y)= \frac{|\Ker(\psi)|}{|\Coker(\psi)|} \cdot |\Ker(\varphi)|\cdot \frac{\vol(\Lambda^\tau \back J(F)^{0,\tau})}{\vol\big(\ago_M^G/ p_M^G(X_*(M))\big)} \cdot |D^G(Y)|_F^{1/2}.$$
\end{theoreme}

\end{paragr}

\begin{paragr}[Démonstration du théorème \ref{thm:IOP-FSAT}.] --- On continue avec les notations précédentes. Pour $g\in G(F)$, soit $\mathbf{1}_{M,g}$ la fonction sur $\ago_M^G$ caractéristique de l'enveloppe convexe des points $-H_P(g)$, $P\in \pc(M)$. Pour tout $\xi\in \ago_M^G$, soit
$$w^{G,\xi}_M(g)=|\{\mu\in p_M^G(X_*(M)) \mid      \mathbf{1}_{M,g}(\xi+\mu)=1 \}|.$$
Le poids $w^{G,\xi}_M$ est invariant à droite par $M(F)$ donc par $J(F)$. On définit alors l'intégrale orbitale pondérée $J_{M}^{G,\xi,s}$ par la formule  (\ref{eq:def-sIOPlocale}) dans laquelle on remplace le poids $v_M^G$ par $w^{G,\xi}_M$.
  
\begin{proposition}
  Sous les hypothèses du théorème \ref{thm:IOP-FSAT}, pour \emph{tout} $\xi\in \ago_M^G$, on a 
$$J_M^{G,\xi,s}(Y)=c_M^G(Y)\cdot\vol\big(\ago_M^G/ p_M^G(X_*(M))\big)\cdot \sum_{i=1}^n  \trace(\tau^{-1},H^\bullet(\Xgo^{G,\xi+\mu_i}_{M,Y},\lc_\eta)_{s'}).$$
\end{proposition}

\begin{preuve}
  Cette proposition généralise au cas tronqué le théorème 15.8 de l'article \cite{GKM} de Goresky-Kottwitz-MacPherson. Notre démonstration s'inspire de la leur. 
Soit $j\in J'(F)$. On confond dans les notations $j$ et son image dans $J'(F)_\tau$. Soit $\xi\in \ago_M^G$. La formule des traces de Grothendieck-Lefschetz donne l'égalité
\begin{eqnarray*}
   \bg s',j\bd \trace((j\tau)^{-1},H^\bullet(\Xgo^{G,\xi}_{M,Y},\lc_\eta))&=&\sum_{x \in  (\Lambda\back \Xgo^{G,\xi}_{M,Y})^{j\tau}} \bg s',\psi(\la_x) j \bd\\
\end{eqnarray*}
où $\la_x$ est la classe dans $\Lambda_\tau$ d'un élément $\la\in \Lambda$ qui vérifie $\la j\tau(x)=x$ et où l'accouplement est l'accouplement canonique entre . L'expression ci-dessus ne dépend que l'image de $j$ dans le conoyau $\Coker(\psi)$ de $\psi$. On a donc
$$\sum_{j\in\Coker(\psi) }\bg s',j\bd \trace((j\tau)^{-1},H^\bullet(\Xgo^{G,\xi}_{M,Y},\lc_\eta))=|\Ker(\psi)| \sum_{j\in J'(F)_\tau}  \bg s',j \bd |\Lambda^\tau\back (\Xgo^{G,\xi}_{M,Y})^{j\tau} |.$$
Par  la décomposition (\ref{eq:decomposition-cohom-FSAT}), le membre de gauche dans l'expression ci-dessus vaut
\begin{equation}
  \label{eq:latrace}
  |\Coker(\psi)| \trace(\tau^{-1},H^\bullet(\Xgo^{G,\xi}_{M,Y},\lc_\eta)_{s'}).
\end{equation}
À l'aide du théorème de l'isogénie de Lang, on voit qu'on a 
\begin{eqnarray*}
  \sum_{j\in J'(F)_\tau}  \bg s',j \bd |\Lambda^\tau\back (\Xgo^{G,\xi}_{M,Y})^{j\tau}| =  \sum_{j\in J'(F)_\tau}  \bg s',j \bd \int_{\Lambda^\tau\back G(F)^{j \tau}} \mathbf{1}_{\ggo(\oc)}(\Ad(g^{-1})Y)  \mathbf{1}_{M,g}(\xi) \, dg\\
=   \vol(\Lambda^\tau \back J(F)^{0,\tau}) \sum_{j\in J'(F)_\tau} \bg s',j \bd \int_{ J(F)^{\tau}  \back G(F)^{j \tau}} \mathbf{1}_{\ggo(\oc)}(\Ad(g^{-1})Y) w_J^\xi(g) \, dg
\end{eqnarray*}
où 
$$ w_J^\xi(g)=\sum_{h\in  J(F)^{0,\tau}\back J(F)^\tau}  \mathbf{1}_{M,g}(\xi)$$
est le nombre de points du réseau $\xi+H_M^G(J(F)^\tau)$ dans l'enveloppe convexe des $-H_P(g)$, $P\in \pc(M)$.
On a donc 
$$\sum_{i=1}^n  w_J^{\xi+\mu_i}(g)=w_M^{G,\xi}(g).$$
Par conséquent lorsqu'on somme sur $\xi+\mu_1,\ldots, \xi+\mu_n$ l'expression (\ref{eq:latrace}), on trouve
\begin{equation}
  \label{eq:expr-intermediaire}
  |\Ker(\psi)| \vol(\Lambda^\tau \back J(F)^{0,\tau}) \sum_{j\in J'(F)_\tau} \bg s',j \bd \int_{ J(F)^{\tau}  \back G(F)^{j \tau}} \mathbf{1}_{\ggo(\oc)}(\Ad(g^{-1})Y) w_M^{G,\xi}(g) \, dg.
\end{equation}
On a $\bg s',j \bd=\bg s,\varphi(j) \bd$ et $J(F)^{\tau}  \back G(F)^{j \tau}$ s'identifie à l'ensemble des  $g\in (J(F)\back G(F))^\tau$ tels que $\inv_Y(g)=\varphi(j)$. En fait, il résulte du lemme \ref{lem:ker<} que l'application $\inv_Y$ est à valeurs dans l'image de $\varphi$. L'expression (\ref{eq:expr-intermediaire}) est donc égale à 
$$ |\Ker(\psi)|\cdot |\Ker(\varphi)| \cdot \vol(\Lambda^\tau \back J(F)^{0,\tau}) |D^G(Y)|^{-1/2} J_M^{G,\xi,s}(Y).$$
La proposition s'en déduit.
\end{preuve}

\begin{proposition}
  Avec les notations du théorème \ref{thm:IOP-FSAT}, pour tout $\xi$ \emph{en position générale}, on a 
$$J_M^{G,\xi,s}(Y)=\vol\big(\ago_M^G/ p_M^G(X_*(M))\big) J_M^{G,s}(Y).$$
\end{proposition}

\begin{preuve}
  Elle est identique à la démonstration du théorème 11.14.2 de \cite{LFPI}.
\end{preuve}
\end{paragr}

\begin{paragr}[Constance locale des fibres de Springer affines tronquées.] --- Soit $M\in \lc$ et $a\in \car_M(\oc)$ un élément génériquement $G$-régulier. On fixe une section de Kostant dans $\mgo$ notée $\eps^M$. On pose $X_a=\eps^M(a)$. On note $J_a$ le centralisateur de $X_a$ dans $M\times_k \oc$ : c'est un $\oc$-schéma en groupes lisse. Soit $\xi\in \ago_M^G$. On pose
$$\Xgo_M^{G,\xi}(a)=\Xgo_{M,X_a}^{G,\xi}.$$
Le but de ce paragraphe est de démontrer la proposition suivante qui généralise au cas tronqué la proposition 3.5.1 de \cite{ngo2}.

\begin{proposition}\label{prop:approx}
  Soit $a\in \car_M(\oc)\cap \car_M^{\reg}(F)$. Il existe un entier $n$ tel que pour tout $a'\in  \car_M(\oc)$ tel que 
$$a'=a\mod \eps^n$$
il existe $h\in M(\oc)$ tel que
\begin{enumerate}
\item $\Ad(h)$ induit un isomorphisme de schéma en groupes de $J_a$ sur $J_{a'}$ ;
\item la translation à gauche par $h$ induit un isomorphisme de $\Xgo_M^{G,\xi}(a)$ sur $\Xgo_M^{G,\xi}(a')$.
\end{enumerate}
De plus, si $a$ et $a'$ sont fixes par une puissance de $\tau$, on peut supposer de plus que $h$ est aussi fixe par la même  puissance de $\tau$.
\end{proposition}

\begin{preuve} Soit $r$ et $n_0$ des entiers pour lesquels les conclusions du lemme \ref{lem:approximation} valent pour l'élément $X_a$ et le groupe $M$. Soit $v=\val(\det(\ad(X_a)_{|\ngo})$ : c'est un entier. Soit $\zeta :  \ngo \to N$  un isomorphisme de variétés affines et $i$  un entier qui satisfont les conclusions du lemme \ref{lem:iso-nilp}. Soit $c \geq 0$ un entier tel que 
  \begin{equation}
    \label{eq:condition-c}
    \forall Y\in \eps^{c}\mgo(\oc)  \textrm{ et }  \forall n\in \zeta(\eps^{-iv}\oc) \ \textrm{ on   a } \Ad(n^{-1})Y -Y \in \ngo(\oc).
  \end{equation}
  Soit $c_1 \geq \max(1,n_0,c+r)$. On suppose de plus que $c_1$ est assez grand pour que pour tout $a'\in \car_M(\oc)$ congru à $a$ modulo $\eps^{c_1}$ on ait $a'\in \car^{\reg}(F)$ et $\val(\det(\ad(X_a)_{|\ngo})=v$. D'après le lemme d'approximation \ref{lem:approximation}, pour tout $a'\in \car_M(\oc)$ tel que $a'=a\mod \eps^{c_1}$, il existe $h\in M(\oc)$ tel que
  \begin{equation}
    \label{eq:cond-h1}
    h=1 \mod \eps^{c}
  \end{equation}
et
\begin{equation}
  \label{eq:cond-h2}
  [\Ad(h)X_a,X_{a'}]=0.
\end{equation}

La conjugaison par $h$ induit un morphisme de $J_a$ sur le centralisateur $I_{\Ad(h)X_a}$ de $\Ad(h)X_a$ qui n'est autre que $J_{a'}$ puisque $X_a$ et $X_{a'}$ sont génériquement semi-simples réguliers. D'après la proposition \ref{prop:IP}, ce morphisme se prolonge. En raisonnant avec le morphisme inverse $\Ad(h^{-1})$, on voit que $\Ad(h)$ induit un isomorphisme de $J_a$ sur $J_{a'}$. Soit $gG(\oc)\in \Xgo_{M}^{G,\xi}(a)$. Montrons que $hgG(\oc)\in \Xgo_{M}^{G,\xi}(a)$. Remarquons tout de suite qu'on a $H_M(h)=0$ puisque $h=1 \mod \eps$. En particulier, la translation par $h$ respecte la troncature. Il suffit donc de prouver qu'on a 
$$\Ad(hg)^{-1}X_{a'}\in \ggo(\oc).$$
Soit $P\in \pc(M)$. Par la décomposition d'Iwasawa, on peut bien supposer qu'on a  $g=mn$ avec $m\in M(F)$ et $n\in N(F)$. Puisque $gG(\oc)\in  \Xgo_{M}^{G,\xi}(a)$, on  a $\Ad(mn)^{-1}X_a\in \ggo(\oc)$. On écrit 
$$\Ad(mn)^{-1}X_a= [\Ad(n^{-1})(\Ad(m)^{-1}X_a)-\Ad(m)^{-1}X_a] +\Ad(m^{-1})X_a$$
où le crochet appartient à $\ngo(F)$ et où $\Ad(m^{-1})X_a\in \mgo(F)$. On a donc les deux conditions
\begin{equation}
  \label{eq:integrite1}
\Ad(m^{-1})X_a\in \mgo(\oc)
\end{equation}
et 
\begin{equation}
 \label{eq:integrite2}
\Ad(n^{-1})(\Ad(m)^{-1}X_a)-\Ad(m)^{-1}X_a\in \ngo(\oc).
\end{equation}
Soit $Y=\Ad(m^{-1})X_a$. La condition (\ref{eq:integrite1}) entraîne que $Y\in \mgo(\oc)$. Soit $I_Y$ le centralisateur de $Y$ dans $M$. L'automorphisme intérieur $\Ad(m^{-1})$ induit un isomorphisme en fibre générique de $J_a$ sur $I_Y$. Par la proposition \ref{prop:IP}, on sait que cet isomorphisme se prolonge en un morphisme de schémas en groupes sur $\oc$ de $J_a$ dans $I_Y$. En particulier, on a
\begin{equation}
  \label{eq:integrite3}
  \Ad(m^{-1}) (\Lie(J_a)(\oc)) \subset I_Y(\oc)\subset \mgo(\oc).
\end{equation}
Comme $\Ad(h^{-1})(X_{a'})$ commute à $X_a$, on a   $\Ad(h^{-1})(X_{a'})\in \Lie(J_a)(\oc)$. Vu (\ref{eq:integrite3}), on obtient
\begin{equation}
  \label{eq:integrite4}
 \Ad(m^{-1})(\Ad(h^{-1})(X_{a'}))\in \mgo(\oc).
\end{equation}
Les conditions (\ref{eq:cond-h1}) et  (\ref{eq:cond-h2}) entraînent qu'on a $X_a-\Ad(h^{-1})X_{a'}\in \eps^{c} \Lie(J_a)(\oc)$.
 Par (\ref{eq:integrite3}), on a 
\begin{equation}
  \label{eq:integrite5}
  \Ad(m^{-1})(X_a-\Ad(h^{-1})X_{a'})\in \eps^{c} \mgo(\oc)
\end{equation}
D'après le lemme \ref{lem:iso-nilp}, les conditions (\ref{eq:integrite1}) et (\ref{eq:integrite2}) et le fait que $v=\val(\det(\ad(X_a)_{\ngo}))$ entraînent que $n\in \zeta(\eps^{-iv}\ngo(\oc))$. Par (\ref{eq:integrite5}) et  (\ref{eq:condition-c}), on a
\begin{equation}
  \label{eq:integrite6}
   \Ad(n^{-1})(\Ad(m^{-1})(X_a-\Ad(h^{-1})X_{a'}))-\Ad(m^{-1})(X_a-\Ad(h^{-1})X_{a'})\in \ngo(\oc).
\end{equation}
Donc la condition (\ref{eq:integrite2}) entraîne l'assertion
\begin{equation}
  \label{eq:integrite7}
   \Ad(n^{-1})(\Ad(m^{-1})(\Ad(h^{-1})X_{a'}))-\Ad(m^{-1})(\Ad(h^{-1})X_{a'})\in \ngo(\oc)
\end{equation}
qui, combinée avec (\ref{eq:integrite4}), donne 
$$\Ad(hmn)^{-1}(X_{a'})\in \ggo(\oc)$$
ce qu'on cherchait à démontrer. Le même raisonnement montre l'implication réciproque. La proposition s'en déduit, la dernière assertion résultant de la dernière assertion de \ref{lem:approximation}.

\end{preuve}

Les deux lemme suivants ont été utilisés dans la démonstration de la proposition \ref{prop:approx}.
 
  \begin{lemme}\label{lem:iso-nilp}
Soit $M\in \lc$ et $P\in \pc(M)$. Il existe un isomorphisme de variétés affines $\zeta : \ngo \to N$, un morphisme $\phi : \ngo \times \mgo \to \ngo$ et un entier $i$ de sorte que pour tout $Y\in \mgo$ tel que $\det(\ad(Y)_{|\ngo})\not=0$ le morphisme $N \to \ngo$ induit par
$$n \mapsto \Ad(n^{-1})Y-Y$$
est un isomorphisme d'inverse $U\in \ngo \mapsto \zeta( (\det(\ad(Y)_{|\ngo}))^{-i} \phi(U,Y))$.
     \end{lemme}

     \begin{preuve}
       Cet énoncé est probablement bien connu mais faute de référence on esquisse une démonstration. Soit $B\subset P$ un sous-groupe de Borel qui contient $T$. Soit $\Delta_B$ l'ensemble des racines simples de $T$ dans $B$ et $\Delta_M\subset \Delta_B$ le sous-ensemble des racines hors $M$. La restriction des racines dans $\Delta_M$ à la composante neutre notée $A_M$ du centre de $M$ est injective. On identifie les éléments de $\Delta_M$ à leur restriction à $A_M$. Soit $\Phi^N$ l'ensemble des racines de $T$ dans $N$. La restriction d'une racine dans $\Phi^N$ à $A_M$ est une combinaison linéaire  à coefficients entiers positifs non tous nuls d'éléments de $\Delta_M$. Pour tout entier $l$, soit $\Phi^N_l$ le sous-ensemble des racines dont la somme des coefficients sur $\Delta_M$ est égale à $l$. Pour $l$ assez grand  l'ensemble $\Phi^N_l$ est vide. Pour tout $\al \in \Phi^N$, on fixe un sous-groupe à $1$ paramètre additif $\zeta_\al : \mathbb{G}_a \to N$ de sorte que $t\zeta_{\al}(u)t^{-1}=\zeta_{\al}(\al(t)u)$ pour $t\in T$ et $u\in \mathbb{G}_a$. 

Pour tout $\al\in \Phi^N$, soit $\ngo_\al\subset \ngo$ le sous-espace de poids $\al$. Posons pour tout entier $l$
$$\ngo_l=\bigoplus_{\al\in \Phi^N_l} \ngo_\al$$
et 
$$\ngo^l=\bigoplus_{i>l}\ngo_i.$$
Par le choix d'une base, on identifie $\ngo_l$ à $\mathbb{G}_a^{\Phi^N_l}$. Le morphisme $\zeta_l=\prod_{\al\in \Phi^N_l} \zeta_\al$  (le produit est le produit dans $N$ selon un ordre fixé une fois pour toutes) envoie $\ngo_l$ sur une sous-variété notée $N_l$ de $N$. Soit $N^l=N_{l+1}\cdot N_{l+2}\cdot \ldots$ (le $\cdot$ désigne le produit dans $N$). On a alors les propriétés suivantes pour tout $l\geq 0$
\begin{enumerate}
\item $\ngo_l$ est stable par l'action adjointe de $M$ ;
\item $[\ngo,\ngo_l]\subset \ngo^l$ ;
\item $N^l$ est un sous-groupe de $N$ stable par conjugaison par $M$ et $[N,N_l]\subset N^l$ ;
\item le quotient $N^{l}/ N^{l+1}$ est commutatif et le morphisme $\zeta_{l+1}$ induit un isomorphisme de groupe de $\ngo_{l+1}$ sur $N^{l}/ N^{l+1}$. On peut de plus supposer que la différentielle de ce morphisme est le morphisme canonique $\ngo_{l+1} \to \ngo^{l}/ \ngo^{l+1}$.
\end{enumerate}
Soit $l\geq 0$, $n\in N^{l}$ et $y\in M$. On a donc $n^{-1}yny^{-1}\in N^{l}$. Il existe $U\in \ngo_{l+1}$ tel que $n\in \zeta_{l+1}(U)N^{l+1}$. En utilisant les propriétés ci-dessus, on voit qu'il existe un morphisme $\psi: \ngo_{l+1}\times M \to \ngo_{l+1}$ tel que   $n^{-1}yny^{-1}\in \zeta_{l+1}(\psi(U,y)) N^{l+1}$. On remarquera que le morphisme $\psi$ est linéaire en la première variable. En prenant une différentielle, on obtient la relation suivante pour tout $Y\in \mgo$
\begin{equation}
  \label{eq:congruence1}
  \Ad(n^{-1})Y-Y=\psi'(U,Y) \mod(\ngo^{l+1})
\end{equation}
pour une certaine application bilinéaire $\psi'$. En différentiant encore, on trouve 
\begin{equation}
  \label{eq:congruence2}
  \psi'(U,Y)=[Y,U]\mod(\ngo^{l+1}).
\end{equation}

Cela dit, nous sommes en mesure de démontrer le lemme. Comme le morphisme $n \mapsto \Ad(n^{-1})Y-Y$ envoie $N^l$ dans $\ngo^l$, on peut raffiner l'énoncé pour les sous-groupes $N^l$. Par récurrence, on suppose l'énoncé raffiné vrai pour $N^1$. Il s'agit alors de le démontrer pour $N^0=N$. Partons de $V\in \ngo$. Soit $V_1$ le projeté de $V$ sur $\ngo_1$ suivant $\ngo=\ngo_1\oplus\ngo^1$. Soit $n=\zeta_{1}(U)$ avec $U\in \ngo_1$. De (\ref{eq:congruence1}) et (\ref{eq:congruence2}), on déduit
$$ \Ad(n^{-1})Y-Y=[Y,U]\mod(\ngo^{1}).$$

Supposons de plus que $\det(\ad(Y)_{|\ngo})\not=0$. Il existe alors un unique $U\in \ngo_1$ tel que $[Y,U]=V_1$. De plus $U$ est un polynôme $\phi_1$ en $V_1$ et $Y$ multiplié par l'inverse de $\det(\ad(Y)_{|\ngo_1})\not=0$. Pour un tel $U$, on a donc
$$   \Ad(n^{-1})Y-Y -V \in \ngo^{1}.$$
En appliquant l'hypothèse de récurrence, un élément de $\ngo_1$ est de la forme $-(\Ad(n_1^{-1})Y-Y)$ avec $n_1\in N^1$ qui est donné par un certain morphisme du type décrit dans l'énoncé. On a alors
\begin{eqnarray*}
  \Ad((nn_1)^{-1})Y-Y -V &=&\Ad(n_1^{-1})( \Ad(n^{-1})Y-Y)+ \Ad(n^{-1})Y-Y -V\\
&=&  \Ad(n^{-1})Y-Y+ \Ad(n^{-1})Y-Y -V\mod(\ngo^2)\\
&=& 0 \mod(\ngo^2).
\end{eqnarray*}
On obtient le lemme par récurrence.
     \end{preuve}

     \begin{lemme}\label{lem:approximation}
       Soit $X\in \ggo(\oc)$ un élément génériquement semi-simple régulier. Alors il existe des entiers $r$ et $n_0$ tels que pour tout $n\geq n_0$ et tout $Y\in \ggo(\oc)$ tel que 
$$Y=X\mod \eps^n$$
il existe $g\in G(\oc)$ tel que
\begin{enumerate}
\item $g=1 \mod \eps^{n-r}$ ;
\item $[\Ad(g)X,Y]=0$.
\end{enumerate}

Si, de plus, $X$ et $Y$ sont fixes par une puissance $\tau$, alors on peut exiger de $g$ d'être fixe par la même puissance de $\tau$.
\end{lemme}

\begin{preuve}
Soit $Y\in \ggo(\oc)$. Le foncteur qui à toute  $\oc$-algèbre $A$ associe l'ensemble
$$\{g\in G(A) \mid [\Ad(g)X,Y]\}$$
est représentable par un schéma sur $\oc$. Supposons  $Y=X \mod \eps^n$ pour $n$ assez grand. La fibre générique de ce $\oc$-schéma est lisse, en particulier si $Y=X$. Il n'en est pas de même en général de sa fibre spéciale. La section unité fournit un point dans $\oc/ \eps^n\oc$. Par une variante du lemme 1 de \cite{Elkik}, en utilisant des approximations successives, on peut relever cette section en un $\oc$-point qui vérifie les conditions voulues. 
\end{preuve}

\end{paragr}

\section{Fin de la démonstration du lemme fondamental pondéré}\label{sec:fin}

\begin{paragr}
  Soit $K_1=\Fq((\eps))$ et $K=k((\eps)$. Soit $\oc_1=\Fq[[\eps]]$. Soit $(G,T)$ un couple formé d'un groupe réductif connexe défini sur $\Fq$ et d'un sous-tore $T$ maximal et déployé sur $\Fq$. Soit $M\in \lc^G(T)$. Soit $s\in \Tc$ et $\Mc'=\Mc_s$. On suppose qu'on a $\Mc\not=\Mc'$ et $Z_{\Mc}^0=Z_{\Mc'}^0$. Soit $M'$ le groupe déployé sur $\Fq$ dual de $\Mc'$. Alors $M'$ est un groupe endoscopique elliptique «non ramifié» de $M$. Pour tout $Y\in \mgo'(K_1)$ et $X\in \mgo(K_1)$ semi-simples et $G$-réguliers, Langlands et Shelstad (cf. \cite{LS}) ont défini un facteur de transfert $\Delta_{M',M}(Y,X)$ ; plus exactement, il s'agit ici d'une variante de leur définition adaptée aux algèbres de Lie et «privée du facteur $\Delta_{\mathrm{IV}}$», que Waldspurger introduit dans  \cite{lem_fond}.  Pour tout $Y\in \mgo'(K_1)$ semi-simple et $G$-régulier, soit
$$J_{M',M}^{G}(Y)=\sum_{X} \Delta_{M',M}(Y,X)|D^G(X)|^{1/2} \int_{J_X(K_1)\back G(K_1)} \mathbf{1}_{\ggo(\oc_1)}(\Ad(g^{-1})X) v_M^G(g) dg.$$
où la somme porte sur l'ensemble des éléments semi-simples $G$-réguliers de $\mgo(K_1)$ pris modulo $M(K_1)$-conjugaison. Le poids $v_M^G$ est le poids d'Arthur considéré précédemment. Les conventions sur les mesures de Haar sont aussi celles considérées précédemment. Notons que $\Delta_{M'}(Y,X)$ est nul sauf sur un sous-ensemble fini de l'ensemble de sommation.

\begin{theoreme} (Lemme fondamental pondéré pour les algèbres de Lie sur les corps locaux de caractéristiques égales)\label{thm:LFP}

On reprend les notations et les hypothèses ci-dessus. Soit $Y\in \mgo'(K_1)$ semi-simple et $G$-régulier. L'application $J :\ec_{M,\el}^G(s) \cup {\Gc} \to \CC$ définie  par
$$J(\Gc)= J_{M',M}^{G}(Y)$$
et pour tout $\Hc\in \ec_{M,\el}^G(s)$ par
$$J(\Hc)= J_{M',M'}^{H}(Y)$$
est stabilisante.
 \end{theoreme}

 \begin{preuve}  Elle est tout-à-fait analogue à la démonstration par Ngô du lemme fondamental ordinaire (cf. \cite{ngo2} §8.6). Soit $Y\in \mgo'(K_1)$ semi-simple et $G$-régulier. On suppose que $\chi_{M'}(Y)\in \car_{M'}(\oc_1)$ sans quoi le théorème est trivialement vrai ($J$ est alors identiquement nulle). Les valeurs de $J$ ne dépendent que de $\chi_{M'}(Y)\in  \car_{M'}(\oc_1)$. Soit $n\geq 1$ un entier, $K_n=\mathbb{F}_{q^n}((\eps))$ et $\oc_n=\mathbb{F}_{q^n}[[\eps]]$. Pour tout entier $n\geq 1$ et tout $b\in \car_{M'}(\oc_n)$, on définit $J_{n,b}$ comme l'application $J$ définie comme ci-dessus mais associée au corps $K_n$ et à un relèvement arbitraire $Y_n\in \mgo'(K_n)$ de $b$ (c'est-à-dire $\chi_{M'}(Y_n)=b$).

Soit $b\in \car_{M'}(\oc_1)$. Il résulte du théorème \ref{thm:IOP-FSAT} qu'il existe des familles presque nulles de complexes $(c_i)_{i\in \NN}$ et $(\la_i)_{i\in \NN}$ qui dépendent de $b$ et $\Gc$ telles que pour tout entier $n\in \NN$ on ait
$$J_{n,b}(\Gc)=\sum_{i\in \NN} c_i \la_i^n.$$
Il en est de même pour tout $\Hc\in  \ec_{M,\el}^G(s)$. Par conséquent, on déduit d'un argument standard (utilisation d'un déterminant de Van der Monde), que l'application $J_{1,b}$ est stabilisante (ce qui est notre but) dès qu'il existe un entier $m$ tel que pour tout entier $n\geq m$, l'application $J_{n,b}$ est stabilisante.

Soit $N$ un entier tel que  pour tout $n\geq 1$ et tout $b'\in \car_{M'}(\oc_n)$ qui vérifie
\begin{equation}
  \label{eq:congruence}
  b'\equiv b \mod \eps^N,
\end{equation}
les fibres de Springer affines $\Xgo_M^{G,\xi}(b)$ et $\Xgo_M^{G,\xi}(b')$ d'une part et  $\Xgo_{M'}^{H,\xi}(b)$ et $\Xgo_{M'}^{H,\xi}(b')$ pour tout $\Hc\in  \ec_{M,\el}^G(s)$ d'autre part vérifient les conclusions de la proposition \ref{prop:approx}.  On suppose également que, sous la condition (\ref{eq:congruence}), on a 
\begin{equation}
  \label{eq:cong-disc}
  D^G(b')\in D^G(b)\oc_n^\times.
\end{equation}
D'après la proposition \ref{prop:approx}, un tel $N$ existe. Pour tout $b'$ qui vérifie (\ref{eq:congruence}), la proposition \ref{prop:approx} et le théorème \ref{thm:IOP-FSAT} impliquent que   pour tout $n\geq 1$, on a 
\begin{equation}
  \label{eq:egalitedesJ}
  J_{n,b'}=J_{n,b}.
\end{equation}

On fixe une courbe projective $C_0$ sur $\Fq$ géométrique connexe et deux points distincts $v$ et $\infty$ dans $C_0(\Fq)$. Soit $g$ son genre.  On fixe un isomorphisme entre $K_1$ et le complété $F_{0,v}$ du corps des fonctions de $C_0$ au point $v$. Ainsi $b$ s'identifie à un élément de $\car_{M'}(\oc_{0,v})$ où $\oc_{0,v}$ est l'anneau des entiers de  $F_{0,v}$. Soit $D$ un diviseur de $C_0$ dont le support ne contient ni $\infty$ ni $v$ et dont le degré vérifie les inégalités suivantes
\begin{equation}
  \label{eq:ineg-pour-kappa}
  \val_v(D^G(b))\leq 2|W^G|^{-1}(\deg(D)-2g+2)
\end{equation}
et
\begin{equation}
  \label{eq:ineg-pour-deg}
N\rang(G)+2g < \deg(D).
\end{equation}
Cette dernière inégalité implique que l'application de restriction
\begin{equation}
  \label{eq:rest-sections}
  H^0(C_0,\car_{M',D})\to \car_{M'}(\oc_{0,v}/ \varpi_v^N)
\end{equation}
est surjective.
Le sous-schéma fermé $\Ac_{M'}(b,N)$ de $\Ac_{M'}$ des $(h_{a},t)$ tels que $h_{a}$ ait même image que $b$ par l'application (\ref{eq:rest-sections}) est de codimension $N\rang(G) $ dans $\Ac_{M'}$. 

Soit  $\Ac_{M'}(b,N)^0$ l'ouvert de $\Ac_{M'}(b,N)$ formé des $a=(h_a,t)$ tels que, pour tout point $c\not=v$ de $C_0$,  si  $h_a(c)$ appartient à $\mathfrak{R}_{M'}^G\cup \mathfrak{D}^{M'}$, ce point  est lisse et l'intersection de $h_a(C)$ avec $\mathfrak{R}_{M'}^G\cup \mathfrak{D}^{M'}$ y est transverse. Par une variante d'un théorème de Bertini (cf. \cite{ngo2} §8.6.6), $\Ac_{M'}(b,N)^0$ est non vide. Comme la codimension de $\Ac_{M'}-\Ac_{M'}^{\el}$ dans $\Ac_{M'}$ est $\geq \deg(D)$ et que la codimension de $\Ac_{M'}(b,N)^0$ dans $\Ac_{M'}$ est $N\rang(G) <\deg(D)$ (cf. (\ref{eq:ineg-pour-deg})), on a $\Ac_{M'}(b,N)^0\cap \Ac_{M'}^{\el}\not=\emptyset$. Cette intersection est même de dimension $>0$. D'après le théorème de Deligne, il existe un entier $m$ tel que pour tout $n\geq m$, il existe $a_M\in \Ac_{M'}(k)^{\tau^n}$ qui appartient à cette intersection. Or un tel point vérifie les hypothèses du théorème \ref{thm:Abon} : l'image de $a_M$ dans $\Ac_G$ appartient donc à l'ouvert $\Ac_G^{\mathrm{bon}}$. Le théorème \ref{thm:LFPsemilocal} implique alors que l'application $J_{n,b'}$, où $b'$ est l'image de $a_M$ dans $\car_{M'}(\oc_n)$, est stabilisante. Or $J_{n,b'}=J_{n,b}$ par l'égalité (\ref{eq:egalitedesJ}). L'application $J_{n,b}$ est donc stabilisante pour tout $n\geq m$. Comme on l'a expliqué en début de démonstration cela entraîne que $J_{b,1}$ est stabilisante ce qu'il fallait démontrer.
 \end{preuve}

\end{paragr}

\bibliographystyle{plain}
\bibliography{biblio}

\bigskip

\begin{flushleft}
Pierre-Henri Chaudouard \\
CNRS et Universit\'{e} Paris-Sud \\
 UMR 8628 \\
 Math\'{e}matique, B\^{a}timent 425 \\
F-91405 Orsay Cedex \\
France \\
\bigskip

Gérard.Laumon \\
CNRS et Universit\'{e} Paris-Sud \\
 UMR 8628 \\
 Math\'{e}matique, B\^{a}timent 425 \\
F-91405 Orsay Cedex \\
France \\
\bigskip

Adresses électroniques :\\
Pierre-Henri.Chaudouard@math.u-psud.fr \\
Gerard.Laumon@math.u-psud.fr\\
\end{flushleft}

\end{document}